%% file: TRRfactor.tex
\documentclass[10pt]{article}

\usepackage[none,bottom,dark]{draftcopy}
\usepackage{amssymb,amsmath,amsthm,mathrsfs}
\usepackage[dvips]{graphics,epsfig}
\usepackage{verbatim}
\usepackage{url}
\usepackage{colordvi,color}
\usepackage{psfrag}
\usepackage{colordvi,color,pspicture}

\graphicspath{{Figures/}}

\usepackage[authoryear,sort&compress]{natbib}

\newcommand{\I}{\mathbf{I}}

\newcommand{\E}{\mathbf{E}}

\newcommand{\U}{\mathcal{U}}

\newcommand{\G}{\Gamma}
\newcommand{\V}{\mathcal{V}}
\newcommand{\A}{\mathcal{A}}
\newcommand{\Y}{\mathcal{Y}}
\newcommand{\y}{\mathsf{y}}
\newcommand{\X}{\mathcal{X}}

\newcommand{\mI}{\mathcal{I}}
\newcommand{\NY}{N_{\Y}}

\newcommand{\NPE}{N_{PE}}

\newcommand{\N}{\mathcal{N}}
\newcommand{\R}{\mathbb{R}}

\newcommand{\msP}{\mathscr{P}}
\newcommand{\RS}{\mathscr{R}_S}

\newcommand{\Var}{\mathbf{Var}}
\newcommand{\Cov}{\mathbf{Cov}}

\DeclareMathOperator{\argsup}{argsup}
\DeclareMathOperator{\PAE}{PAE}
\DeclareMathOperator{\HLAE}{HLAE}

\begin{document}

\title{Technical Report \# 644:\\
Relative Density of the Random $r$-Factor Proximity Catch Digraph
for Testing Spatial Patterns of Segregation and Association}
\author{\it Elvan Ceyhan, Carey E. Priebe, \& John C. Wierman\\
Department of Applied Mathematics and Statistics,\\
The Johns Hopkins University, Baltimore, MD, 21218}

\date{July 29, 2004}

\maketitle
\begin{abstract}
Statistical pattern classification methods based on data-random
 graphs were introduced recently.  In this approach, a
random directed graph is constructed from the data using the
relative positions of the data points from various classes.
Different random graphs result from different definitions of the
proximity region associated with each data point and different
graph statistics can be employed for data reduction.  The approach used in this article
is based on a parameterized family of proximity maps determining an associated family of data-random digraphs.
The relative arc density of the digraph is used
as the summary statistic, providing an alternative to the
domination number employed previously.  An important advantage of
the relative arc density is that, properly re-scaled, it is a
$U$-statistic, facilitating analytic study of its asymptotic
distribution using standard $U$-statistic central limit theory.  The
approach is illustrated with an application to the testing of
spatial patterns of segregation and association. Knowledge of the
asymptotic distribution allows evaluation of the Pitman and
Hodges-Lehmann asymptotic efficacy, and selection of the proximity map
parameter to optimize efficacy.  Notice that the approach presented here also
has the advantage of validity for data in any dimension.
\end{abstract}

\section{Introduction}
Classification and clustering have received considerable attention
in the statistical literature. In recent years, a new classification
approach has been developed which is based on the relative positions
of the data points from various classes.  Priebe et al. introduced
the class cover catch digraphs (CCCD) in $\R$ and gave the exact and
the asymptotic distribution of the domination number of the CCCD (\cite{priebe:2001}).
\cite{devinney:2002a}, Marchette
and Priebe \cite{marchette:2003}, \cite{priebe:2003a}, \cite{priebe:2003b}
applied the concept in higher dimensions and
demonstrated relatively good performance of CCCD in classification.
The methods employed involve data reduction (condensing) by using
approximate minimum dominating sets as prototype sets (since finding
the exact minimum dominating set is an NP-hard problem ---in
particular for CCCD).  Furthermore the exact and the asymptotic
distribution of the domination number of the CCCD are not
analytically tractable in multiple dimensions.

Ceyhan and Priebe introduced the central similarity proximity map
and $r$-factor proximity maps and the associated random digraphs in
\cite{ceyhan:CS-JSM-2003} and \cite{ceyhan:TR-dom-num-NPE-spatial},
respectively. In both cases, the space is partitioned by the
Delaunay tessellation which is the Delaunay triangulation in $\R^2$.
In each triangle, a family of data-random proximity catch digraphs
is constructed based on the proximity of the points to each other.
The advantages of the $r$-factor proximity catch digraphs are that
an exact minimum dominating set can be found in polynomial time and
the asymptotic distribution of the domination number is analytically
tractable. The latter is then used to test segregation and
association of points of different classes in
\cite{ceyhan:TR-dom-num-NPE-spatial}.

In this article, we employ a different statistic, namely the
relative (arc) density, that is the proportion of all possible arcs
(directed edges) which are present  in the data random digraph. This
test statistic has the advantage that, properly rescaled, it is a
$U$-statistic.  Two simple classes of alternative hypotheses, for
segregation and association, are defined in Section
\ref{sec:null-and-alt}. The asymptotic distributions under both the
null and the alternative hypotheses are determined in Section
\ref{sec:asy-norm} by using standard $U$-statistic central limit
theory. Pitman and Hodges-Lehmann asymptotic efficacy are analyzed
in Sections \ref{sec:Pitman} and \ref{sec:Hodges-Lehmann},
respectively.  This test is related to the available tests of
segregation and association in the literature, such as Pielou's test
and Ripley's test.  See discussion in Section \ref{sec:discussion}
for more detail. Our approach is valid for data in any dimension,
but for simplicity of expression and visualization, will be
described for two-dimensional data.

\section{Preliminaries}
\subsection{Proximity Maps}
Let $(\Omega,\mathcal{M})$ be a measurable space and consider a
function $N:\Omega \times \wp(\Omega) \rightarrow \wp(\Omega)$,
where $\wp(\cdot)$ represents the power set functional. Then given
$\Y \subseteq \Omega$, the {\em proximity map} $\NY(\cdot) =
N(\cdot,\Y): \Omega \rightarrow \wp(\Omega)$ associates with each
point $x \in \Omega$ a {\em proximity region} $\NY(x) \subset
\Omega$.  Typically, $N$ is chosen to satisfy $x \in \NY(x)$ for all
$x \in \Omega$. The use of the adjective \emph{proximity} comes form
thinking of the region $\NY(x)$ as representing a neighborhood of
points ``close" to $x$ (\cite{toussaint:1980,jaromczyk:1992}).

\subsection{$r$-Factor Proximity Maps}
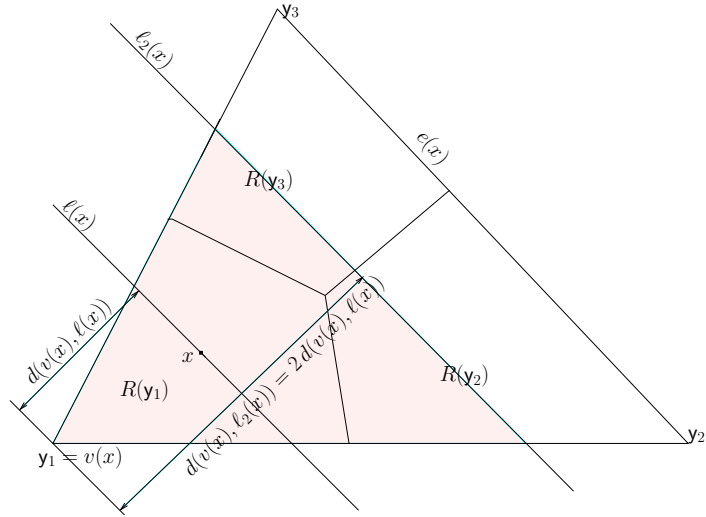
\begin{figure} [ht]
    \centering
   \scalebox{.4}{\input{Nofnu.pstex_t}}
    \caption{Construction of $r$-factor proximity region, $\NY^2(x)$ (shaded region). }
\label{fig:ProxMapDef}
    \end{figure}

We now briefly define $r$-factor proximity maps
(see \cite{ceyhan:TR-dom-num-NPE-spatial} for more details).
Let $\Omega = \R^2$
and let $\Y = \{\y_1,\y_2,\y_3\} \subset \R^2$
be three non-collinear points.
Denote by $T(\Y)$ the triangle ---including the interior--- formed by the three points.
For $r \in [1,\infty]$,
define $\NY^r$ to be the {\em r-factor} proximity map as follows;
see also Figure \ref{fig:ProxMapDef}.
Using line segments from the center of mass of $T(\Y)$ to the midpoints of its edges, we partition $T(\Y)$ into ``vertex regions" $R(\y_1)$, $R(\y_2)$, and $R(\y_3)$.
For $x \in T(\Y) \setminus \Y$, let $v(x) \in \Y$ be the
vertex in whose region $x$ falls, so $x \in R(v(x))$.
If $x$ falls on the boundary of two vertex regions,
 we assign $v(x)$ arbitrarily to one of the adjacent regions.
Let $e(x)$ be the edge of $T(\Y)$ opposite $v(x)$.
Let $\ell(x)$ be the line parallel to $e(x)$ through $x$.
Let $d(v(x),\ell(x))$ be the Euclidean (perpendicular) distance from $v(x)$ to $\ell(x)$.
For $r \in [1,\infty)$, let $\ell_r(x)$ be the line parallel to $e(x)$
such that $d(v(x),\ell_r(x)) = rd(v(x),\ell(x))$ and $d(\ell(x),\ell_r(x)) < d(v(x),\ell_r(x))$.
Let $T_r(x)$ be
the triangle similar to
and with the same orientation as $T(\Y)$
having $v(x)$ as a vertex
and $\ell_r(x)$ as the opposite edge.
Then the {\em r-factor} proximity region
$\NY^r(x)$ is defined to be $T_r(x) \cap T(\Y)$.
Notice that $r \ge 1$ implies $x \in \NY^r(x)$.
Note also that
$\lim_{r \rightarrow \infty} \NY^r(x) = T(\Y)$
for all $x \in T(\Y) \setminus \Y$,
so we define $\NY^{\infty}(x) = T(\Y)$ for all such $x$.
For $x \in \Y$, we define $\NY^r(x) = \{x\}$ for all $r \in [1,\infty]$.

\subsection{Data-Random Proximity Catch Digraphs}
\label{sec:data-random-PCD}
If $\X_n:=\{X_1,X_2,\cdots,X_n\}$ is a set of $\Omega$-valued random variables,
then the $\NY(X_i),\; i=1,\cdots,n$, are random sets.
If the $X_i$ are independent and identically distributed,
then so are the random sets $\NY(X_i)$.

In the case of an $r$-factor proximity map, notice that if $X_i \stackrel{iid}{\sim} F$
and $F$ has a non-degenerate two-dimensional
probability density function $f$ with support$(f) \subseteq T(\Y)$,
 then the special case in the construction
of $\NY^r$ ---
$X$ falls on the boundary of two vertex regions ---
occurs with probability zero.

The proximities of the data points to each other are used to construct a digraph. A digraph is a directed graph; i.e. a graph with directed edges from one vertex to another based on a binary relation.
Define the data-random proximity catch digraph $D$
with vertex set $\V=\{X_1,\cdots,X_n\}$
and arc set $\A$ by
$(X_i,X_j) \in \A \iff X_j \in \NY(X_i)$.  Since this relationship is not symmetric, a digraph is needed rather than a graph.
The random digraph $D$ depends on
the (joint) distribution of the $X_i$ and on
the map $\NY$.

\subsection{Relative Density}
The \emph{relative arc density} of a digraph $D=(\V,\A)$
of order $|\V| = n$,
denoted $\rho(D)$,
is defined as
$$
\rho(D) = \frac{|\A|}{n(n-1)}
$$
where $|\cdot|$ denotes the set cardinality functional (\cite{janson:2000}).

Thus $\rho(D)$ represents the ratio of the number of arcs
in the digraph $D$ to the number of arcs in the complete symmetric
digraph of order $n$, which is $n(n-1)$. For brevity of notation we use \emph{relative density} rather than relative arc density henceforth.

If $X_1,\cdots,X_n \stackrel{iid}{\sim} F$
the relative density
of the associated data-random proximity catch digraph $D$,
denoted $\rho(\X_n;h,\NY)$, is a $U$-statistic,
\begin{eqnarray}
\rho(\X_n;h,\NY) =
  \frac{1}{n(n-1)}
    \sum\hspace*{-0.1 in}\sum_{i < j \hspace*{0.25 in}}   \hspace*{-0.1 in}
      \hspace*{-0.1 in}\;\;h(X_i,X_j;\NY)
\end{eqnarray}
where
\begin{eqnarray}
h(X_i,X_j;\NY)&=& \I\{(X_i,X_j) \in \A\}+ \I\{(X_j,X_i) \in \A\} \nonumber \\
       &=& \I\{X_j \in \NY(X_i)\}+ \I\{X_i \in \NY(X_j)\}.
\end{eqnarray}
We denote $h(X_i,X_j;\NY)$ as $h_{ij}$ for brevity of notation.
Although the digraph is asymmetric, $h_{ij}$ is defined as
the number of arcs in $D$ between vertices $X_i$ and $X_j$,
in order to produce a symmetric kernel with finite variance (\cite{lehmann:1988}).

The random variable $\rho_n := \rho(\X_n;h,\NY)$ depends on $n$ and $\NY$ explicitly
and on $F$ implicitly.
The expectation $\E[\rho_n]$, however, is independent of $n$
and depends on only $F$ and $\NY$:
\begin{eqnarray}
0 \leq \E[\rho_n] = \frac{1}{2}\E[h_{12}] \leq 1 \text{ for all $n\ge 2$}.
\end{eqnarray}
The variance $\Var[\rho_n]$ simplifies to
\begin{eqnarray}
\label{eq:var-rho}
0 \leq
  \Var[\rho_n] =
     \frac{1}{2n(n-1)} \Var[h_{12}] +
     \frac{n-2}{n(n-1)} \Cov[h_{12},h_{13}]
  \leq 1/4.
\end{eqnarray}
A central limit theorem for $U$-statistics
(\cite{lehmann:1988})
yields
\begin{eqnarray}
\sqrt{n}\bigl(\rho_n-\E[\rho_n]\bigr) \stackrel{\mathcal{L}}{\longrightarrow} \N\bigl(0,\Cov[h_{12},h_{13}]\bigr)
\end{eqnarray}
provided $\Cov[h_{12},h_{13}] > 0$.
The asymptotic variance of $\rho_n$, $\Cov[h_{12},h_{13}]$,
depends on only $F$ and $\NY$.
Thus, we need determine only
$\E[h_{12}]$
and
$\Cov[h_{12},h_{13}]$
in order to obtain the normal approximation
\begin{eqnarray}
\rho_n \stackrel{\text{approx}}{\sim}
\N\bigl(\E[\rho_n],\Var[\rho_n]\bigr) =
\N\left(\frac{\E[h_{12}]}{2},\frac{\Cov[h_{12},h_{13}]}{n}\right) \text{ for large $n$}.
\end{eqnarray}

\subsection{Null and Alternative Hypotheses}
\label{sec:null-and-alt}
The phenomenon known as {\em segregation}
involves observations from different classes
having a tendency to repel each other
--- in our case, this means the $X_i$ tend to be located away from all elements of $\Y$.
{\em Association} involves observations from different classes
having a tendency to attract one another,
so that the $X_i$ tend to be located near an element of $\Y$.
See, for instance, \cite{dixon:1994}, \cite{coomes:1999}. For statistical testing for segregation and association, the null hypothesis is generally some form of
{\em complete spatial randomness};
thus we consider
$$H_0: X_i \stackrel{iid}{\sim} \U(T(\Y)).$$
If it is desired to have the sample size be a random variable,
we may consider a spatial Poisson point process on $T(\Y)$
as our null hypothesis.

We define two simple classes of alternatives,
$H^S_{\epsilon}$ and $H^A_{\epsilon}$
with $\epsilon \in \bigl( 0,\sqrt{3}/3 \bigr)$,
for segregation and association, respectively.
For $\y \in \Y$,
let $e(\y)$ denote the edge of $T(\Y)$ opposite vertex $\y$,
and for $x \in T(\Y)$
let $\ell_\y(x)$ denote the line parallel to $e(\y)$ through $x$.
Then define
$T(\y,\epsilon) = \bigl\{x \in T(\Y): d(\y,\ell_\y(x)) \le \epsilon\bigr\}$.
Let $H^S_{\epsilon}$ be the model under which
$X_i \stackrel{iid}{\sim} \U \bigl(T(\Y) \setminus \cup_{\y \in \Y} T(\y,\epsilon)\bigr)$
and $H^A_{\epsilon}$ be the model under which
$X_i \stackrel{iid}{\sim} \U\bigl(\cup_{\y \in \Y} T(\y,\sqrt{3}/3 - \epsilon)\bigr)$.
Thus the segregation model excludes the possibility of
any $X_i$ occurring near a $\y_j$,
and the association model requires
that all $X_i$ occur near a $\y_j$.
The $\sqrt{3}/3 - \epsilon$ in the definition of the
association alternative is so that $\epsilon=0$
yields $H_0$ under both classes of alternatives.

{\bf Remark:}
These definitions of the alternatives
are given for the standard equilateral triangle.
The geometry invariance result of Theorem 1 from Section 3 still holds
under the alternatives, in the following sense.
If, in an arbitrary triangle,
a small percentage $\delta \cdot 100\%$ where $\delta \in (0,4/9)$ of the area is carved
away as forbidden from each vertex using line segments parallel
to the opposite edge, then
under the transformation to the standard equilateral triangle
this will result in the alternative $H^S_{\sqrt{3 \delta / 4}}$.
This argument is for
segregation with $\delta < 1/4$;
a similar construction is available for the other cases.

\section{Asymptotic Normality Under the Null and Alternative Hypotheses}
\label{sec:asy-norm}
First we present a ``geometry invariance" result which allows us to assume $T(\Y)$ is the standard equilateral triangle, $T\bigl((0,0),(1,0),\bigl( 1/2,\sqrt{3}/2 \bigr)\bigr)$, thereby simplifying our subsequent analysis.

{\bf Theorem 1:}
Let $\Y = \{\y_1,\y_2,\y_3\} \subset \R^2$
be three non-collinear points.
For $i=1,\cdots,n$
let $X_i \stackrel{iid}{\sim} F = \U(T(\Y))$,
the uniform distribution on the triangle $T(\Y)$.
Then for any $r \in [1,\infty]$
the distribution of $\rho(\X_n;h,\NY^r)$
is independent of $\Y$,
 hence the geometry of $T(\Y)$.

{\bf Proof:}
A composition of translation, rotation, reflections, and scaling
will transform any given triangle $T_o = T\bigl(\y_1,\y_2,\y_3\bigr)$
into the ``basic'' triangle $T_b = T\bigl((0,0),(1,0),(c_1,c_2) \bigr)$
with $0 < c_1 \le 1/2$, $c_2 > 0$ and $(1-c_1)^2+c_2^2 \le 1$,
preserving uniformity.
The transformation $\phi_e: \R^2 \rightarrow \R^2$
given by $\phi_e(u,v) = \left(u+\frac{1-2\,c_1}{\sqrt{3}}\,v,\frac{\sqrt{3}}{2\,c_2}\,v \right)$
takes $T_b$ to
the equilateral triangle
$T_e = T\bigl((0,0),(1,0),\bigl( 1/2,\sqrt{3}/2 \bigr)\bigr)$.
Investigation of the Jacobian shows that $\phi_e$
also preserves uniformity.
Furthermore, the composition of $\phi_e$ with the rigid motion transformations
maps
     the boundary of the original triangle $T_o$
  to the boundary of the equilateral triangle $T_e$,
     the median lines of $T_o$
  to the median lines of $T_e$,
and  lines parallel to the edges of $T_o$
  to lines parallel to the edges of $T_e$.
Since the joint distribution of any collection of the $h_{ij}$
involves only probability content of unions and intersections
of regions bounded by precisely such lines,
and the probability content of such regions is preserved since uniformity is preserved,
the desired result follows.
$\blacksquare$

Based on Theorem 1 and our uniform null hypothesis,
we may assume that
$T(\Y)$ is the standard equilateral triangle
with $\Y = \bigl\{(0,0),(1,0),\bigl( 1/2,\sqrt{3}/2 \bigr)\bigr\}$
henceforth.

For our $r$-factor proximity map and uniform null hypothesis,
the asymptotic null distribution of $\rho_n(r) = \rho(\X_n;h,\NY^r)$ can be derived as a function of $r$.  Let $\mu(r):=\E[\rho_n(r)]$ and $\nu(r):=\Cov[h_{12},h_{13}]$. Notice that $\mu(r)=\E[h_{12}]/2=P(X_2 \in \NY^r(X_1))$ is the probability of an arc occurring between any pair of vertices.

\subsection{Asymptotic Normality under the Null Hypothesis}
By detailed geometric  probability calculations, provided in Appendix 1, the mean and the asymptotic variance of the relative density of the $r$-factor proximity catch digraph can explicitly be computed.  The central limit theorem for $U$-statistics then establishes the asymptotic normality under the uniform null hypothesis. These results are summarized in the following theorem.

{\bf Theorem 2:}
For $r \in [1,\infty)$,
\begin{eqnarray}
 \frac{\sqrt{n}\,\bigl(\rho_n(r)-\mu(r)\bigr)}{\sqrt{\nu(r)}}
 \stackrel{\mathcal{L}}{\longrightarrow}
 \N(0,1)
\end{eqnarray}
where
\begin{eqnarray}
\label{eq:Asymean}
\mu(r) =
 \begin{cases}
  \frac{37}{216}r^2                                 &\text{for} \quad r \in [1,3/2) \\
  -\frac{1}{8}r^2 + 4 - 8r^{-1} + \frac{9}{2}r^{-2}  &\text{for} \quad r \in [3/2,2) \\
  1 - \frac{3}{2}r^{-2}                             &\text{for} \quad r \in [2,\infty) \\
 \end{cases}
\end{eqnarray}
and
\begin{equation}
\label{eq:Asyvar}
\nu(r) =\nu_1(r) \,\I(r \in [1,4/3)) + \nu_2(r) \,\I(r \in [4/3,3/2))+ \nu_3(r) \,\I(r \in [3/2,2)) + \nu_4(r) \,\I( r \in [2,\infty])
\end{equation}
with
{\small
\begin{align*}
  \nu_1(r) &=\frac{3007\,r^{10}-13824\,r^9+898\,r^8+77760\,r^7-117953\,r^6+48888\,r^5-24246\,r^4+60480\,r^3-38880\,r^2+3888}{58320\,r^4},\\
  \nu_2(r) &=\frac{5467\,r^{10}-37800\,r^9+61912\,r^8+46588\,r^6-191520\,r^5+13608\,r^4+241920\,r^3-155520\,r^2+15552}{233280\,r^4},  \\
  \nu_3(r) &=-[7\,r^{12}-72\,r^{11}+312\,r^{10}-5332\,r^8+15072\,r^7+13704\,r^6-139264\,r^5+273600\,r^4-242176\,r^3\\
& +103232\,r^2-27648\,r+8640]/[960\,r^6],\\
  \nu_4(r) &=\frac{15\,r^4-11\,r^2-48\,r+25}{15\,r^6}.
\end{align*}
}
For $r=\infty$, $\rho_n(r)$ is degenerate.

See Appendix 1 for the proof.

\begin{figure}[ht]
\centering
\psfrag{mu(r)}{\scriptsize{$\mu(r)$}}
\psfrag{r}{\scriptsize{$r$}}
\epsfig{figure=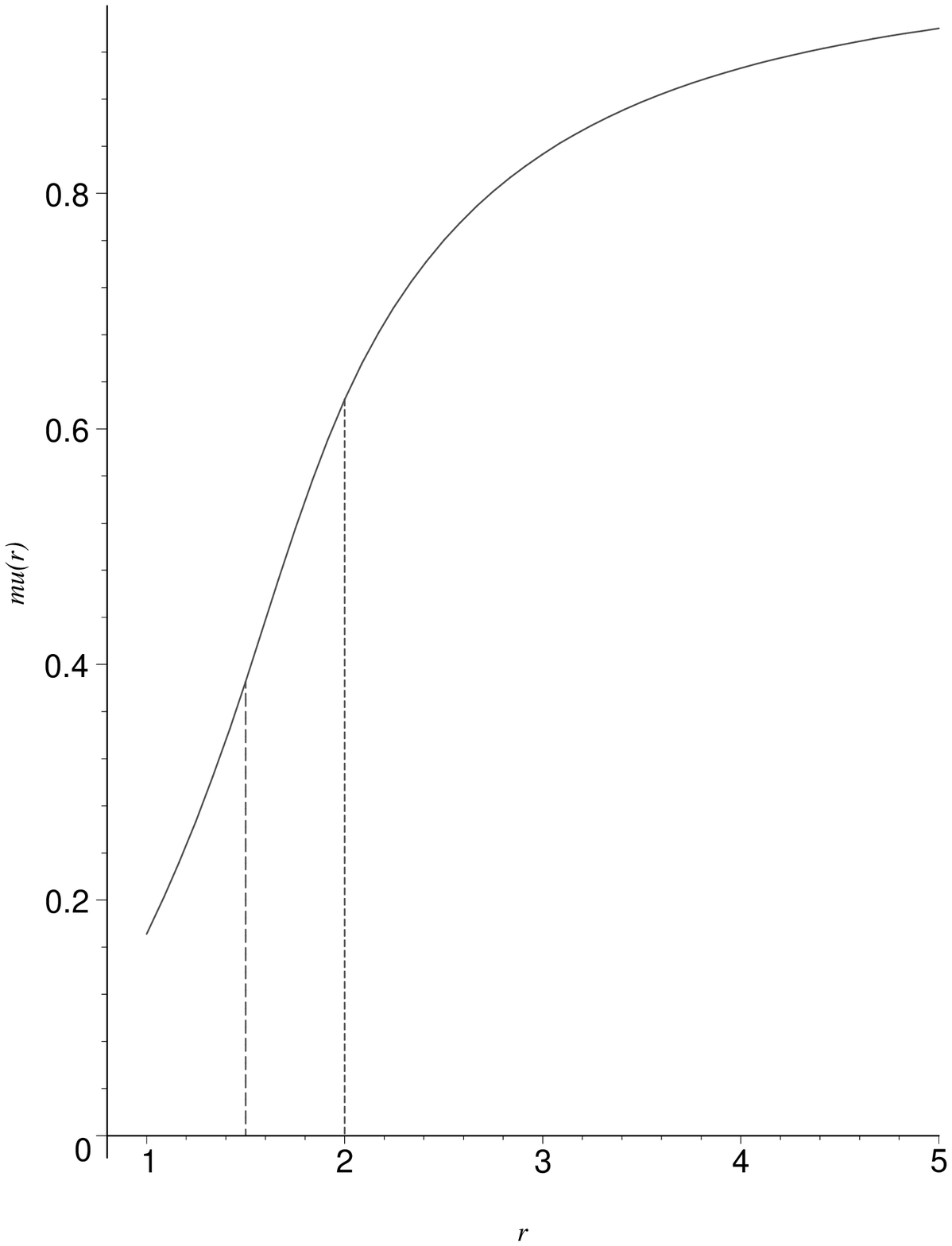, height=180pt, width=180pt}
\psfrag{nu(r)}{\small{$\nu(r)$}}
\epsfig{figure=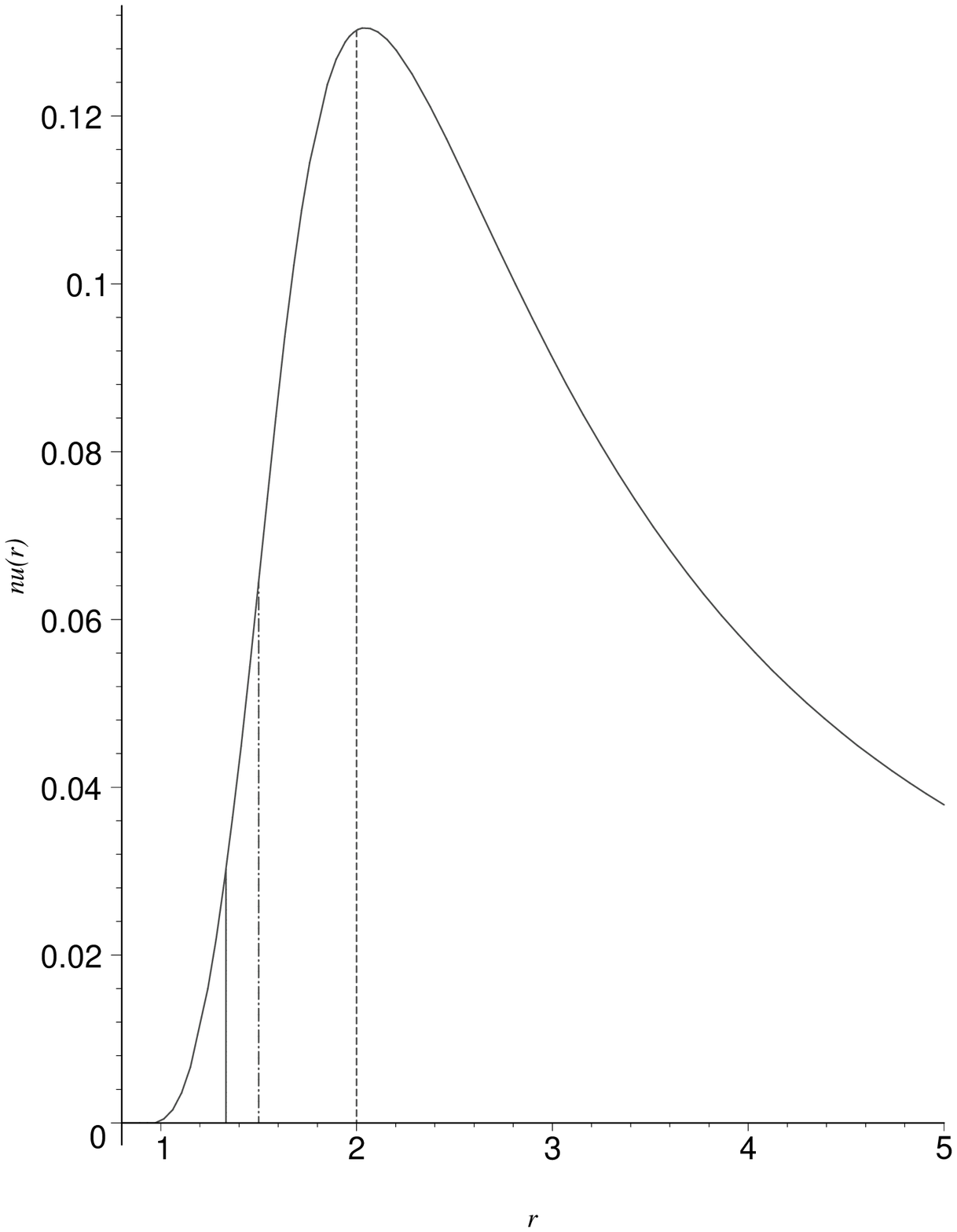, height=180pt, width=180pt}
\caption{
\label{fig:AsyNormCurves}
Asymptotic null mean $\mu(r)$ (left) and variance $\nu(r)$ (right) from Theorem 2. The vertical lines indicate the endpoints of the intervals in the piecewise definition of the functions.  Notice that the vertical axes are differently scaled.
}
\end{figure}

Consider the forms of the mean and asymptotic variance functions, which are depicted in Figure \ref{fig:AsyNormCurves}. Note that $\mu(r)$ is monotonically increasing in $r$, since $\NY^r(x)$ increases with $r$ for all $x \in R_{CM}(\y_j)\setminus \RS(\NY^r,M_C)$, where $\RS(\NY^r,M_C):=\{x \in T(\Y):\;\NY^r(x)=T(\Y)\}$ .  In addition, $\mu(r) \rightarrow 1$ as $r \rightarrow \infty$ (at rate $O\left( r^{-2} \right)$), since the digraph becomes complete asymptotically, which explains why $\rho_n(r)$ becomes degenerate, i.e. $\nu(r=\infty)=0$. Note also that $\mu(r)$ is continuous, with the value at $r=1$, $\mu(1) = 37/216 \approx .1713$.

Regarding the asymptotic variance, note that $\nu(r)$ is also continuous in $r$ with $\lim_{r \rightarrow \infty} \nu(r) = 0$ and $\nu(1) = 34/58320\approx .000583$ and observe that $\sup_{r \ge 1} \nu(r) \approx .1305$ at $\argsup_{r \ge 1} \nu(r) \approx 2.045$.

To illustrate the limiting distribution, $r=2$ yields
$$
\frac{\sqrt{n} \bigl( \rho_n(2) - \mu(2)\bigr)}{\sqrt{\nu(2)}}
=
\sqrt{\frac{192n}{25}} \left(\rho_n(2) - \frac{5}{8}\right)\stackrel{\mathcal{L}}{\longrightarrow} \N(0,1)$$
or equivalently,
$$
\rho_n(2) \stackrel{\text{\scriptsize approx}}{\sim} \N\left(\frac{5}{8},\frac{25}{192n}\right)
.$$

The finite sample variance and skewness may be derived analytically
in much the same way as was $\Cov[h_{12},h_{13}]$
for the asymptotic variance.
In particular, the variance of $h_{12}$ is
\begin{multline*}
\omega(r)=\Var[h_{12}]= \omega_{1,1}(r)\,\I(r \in [1,4/3)) +\\
 \omega_{1,2}(r)\,\I(r \in [4/3,3/2))+\omega_{1,3}(r)\,\I(r \in [3/2,2))+\omega_{1,4}(r)\,\I(r \in [2,\infty))
\end{multline*}
where
{\small
\begin{align*}
\omega_{1,1}(r)&=\frac{-(1369\,r^8+4107\,r^7+902\,r^6-78084\,\,r^5+161784\,\,r^4-182736\,r^3-23328\,r^2+155520\,r-55296)}{11664\,\,(r+2)(r+1)r^2},\\
\omega_{1,2}(r)&=-\frac{1369\,r^7+4107\,r^6+9650\,r^5-98496\,r^4+132624\,\,r^3-79056\,r^2-57888\,r+72576}{11664\,\,(r+2)(r+1)r},\\
\omega_{1,3}(r)&=-\frac{r^{10}+3\,r^9-62\,r^8+968\,r^6-1704\,\,r^5-1824\,\,r^4+5424\,\,r^3-1168\,r^2-3856\,r+2208}{16\,(r+2)(r+1)r^4},\\
\omega_{1,4}(r)&=\frac{3\,r^3+3\,r^2+3\,r-13}{r^4(r+1)}.
\end{align*}
}
\begin{figure}[ht]
\centering
\psfrag{var}{\small{$\omega(r)$}}
\psfrag{r}{\small{$r$}}
\epsfig{figure=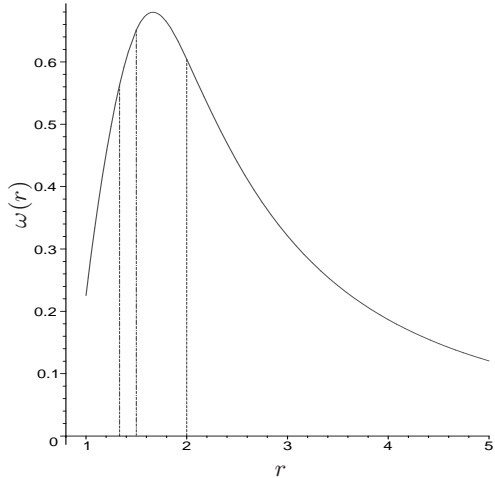, height=180pt, width=180pt}
\caption{
\label{fig:Var(r)-NYr}
$\Var[h_{12}]=\omega(r)$ as a function of $r$ for $r \in [1,5] $.}
\end{figure}

In Figure \ref{fig:Var(r)-NYr} is the graph of $\omega(r)$ for $r \in [1,5]$. Note that $\omega(r=1)=2627/11664\approx .2252$ and $\lim_{r \rightarrow \infty}\omega(r)=0$ (at rate $O\left( r^{-2} \right)$), $\argsup_{r \in [1,\infty)} \omega(r) \approx 1.66$ with $\sup_{r \in [1,\infty)} \omega(r) \approx .6796$.

In fact,
the exact distribution of $\rho_n(r)$
is, in principle, available
by successively conditioning on the values of $X_i$.
Alas,
while the joint distribution of $h_{12},h_{13}$ is available,
the joint distribution of $\{h_{ij}\}_{1 \leq i < j \leq n}$,
and hence the calculation for the exact distribution of $\rho_n(r)$,
is extraordinarily tedious and lengthy for even small values of $n$.

\begin{figure}[ht]
\centering
\psfrag{Density}{ \Huge{\bfseries{Density}}}
\rotatebox{-90}{ \resizebox{1.84 in}{!}{ \includegraphics{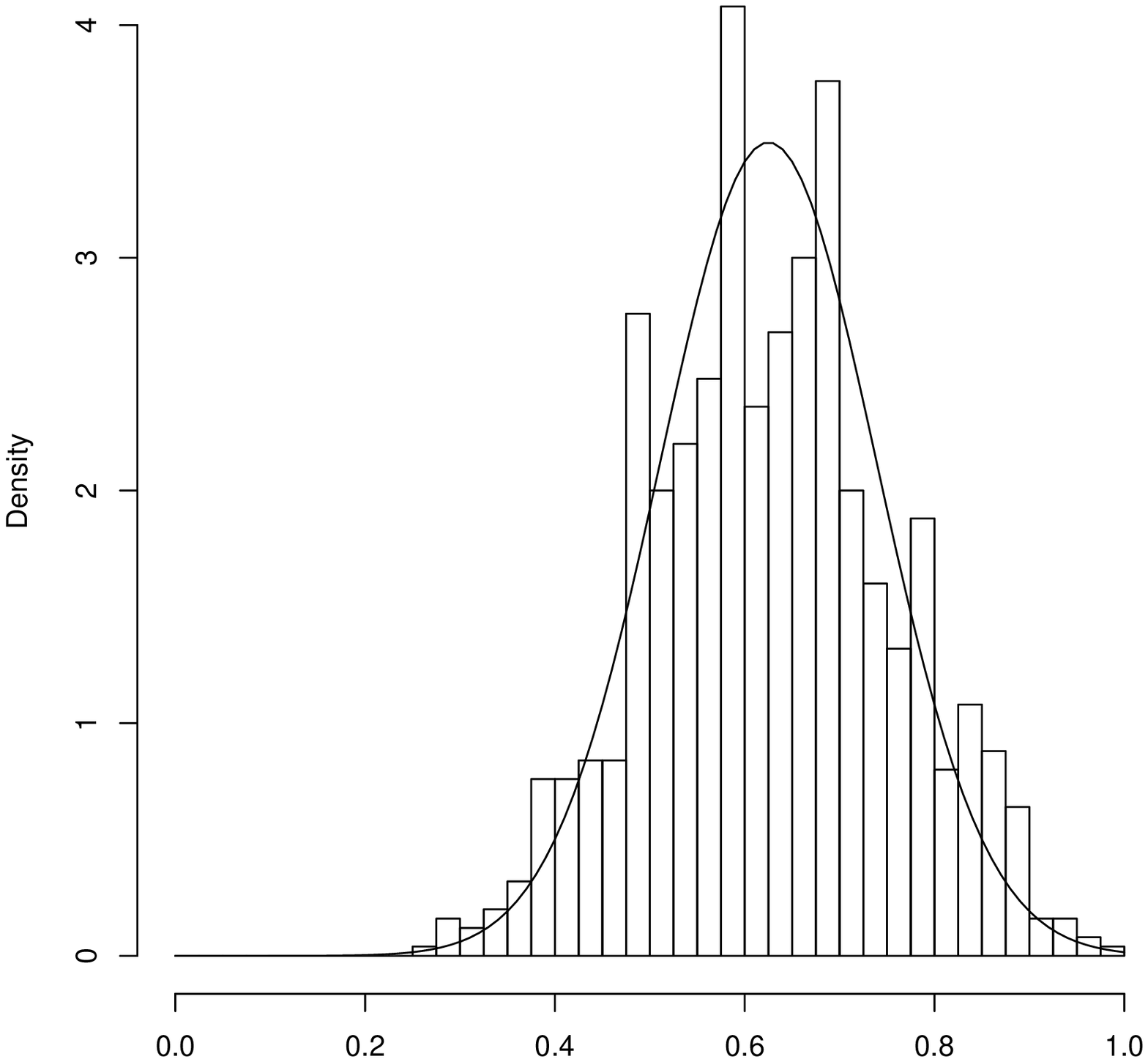}}}
\rotatebox{-90}{ \resizebox{1.84 in}{!}{ \includegraphics{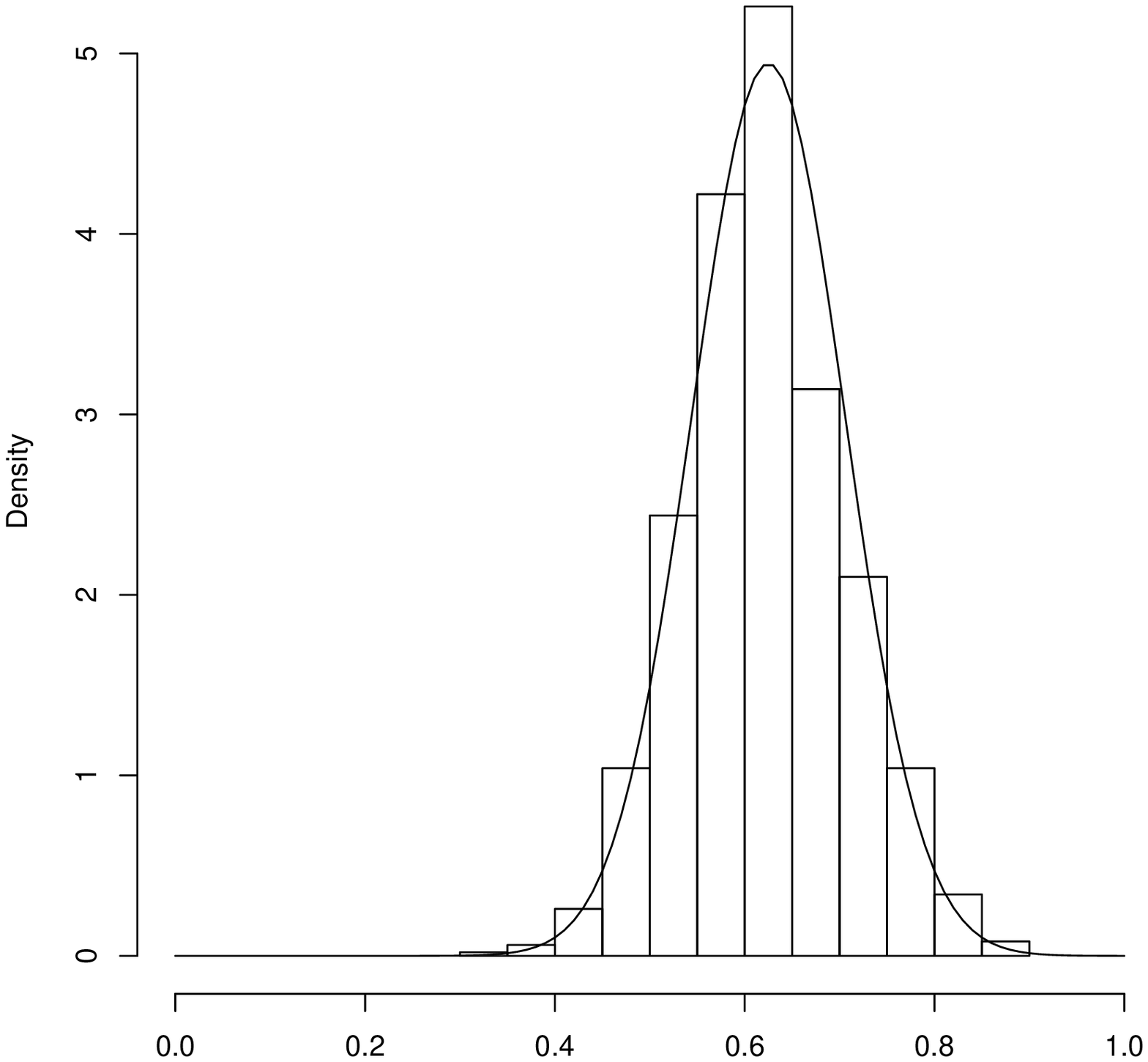}}}
\rotatebox{-90}{ \resizebox{1.84 in}{!}{ \includegraphics{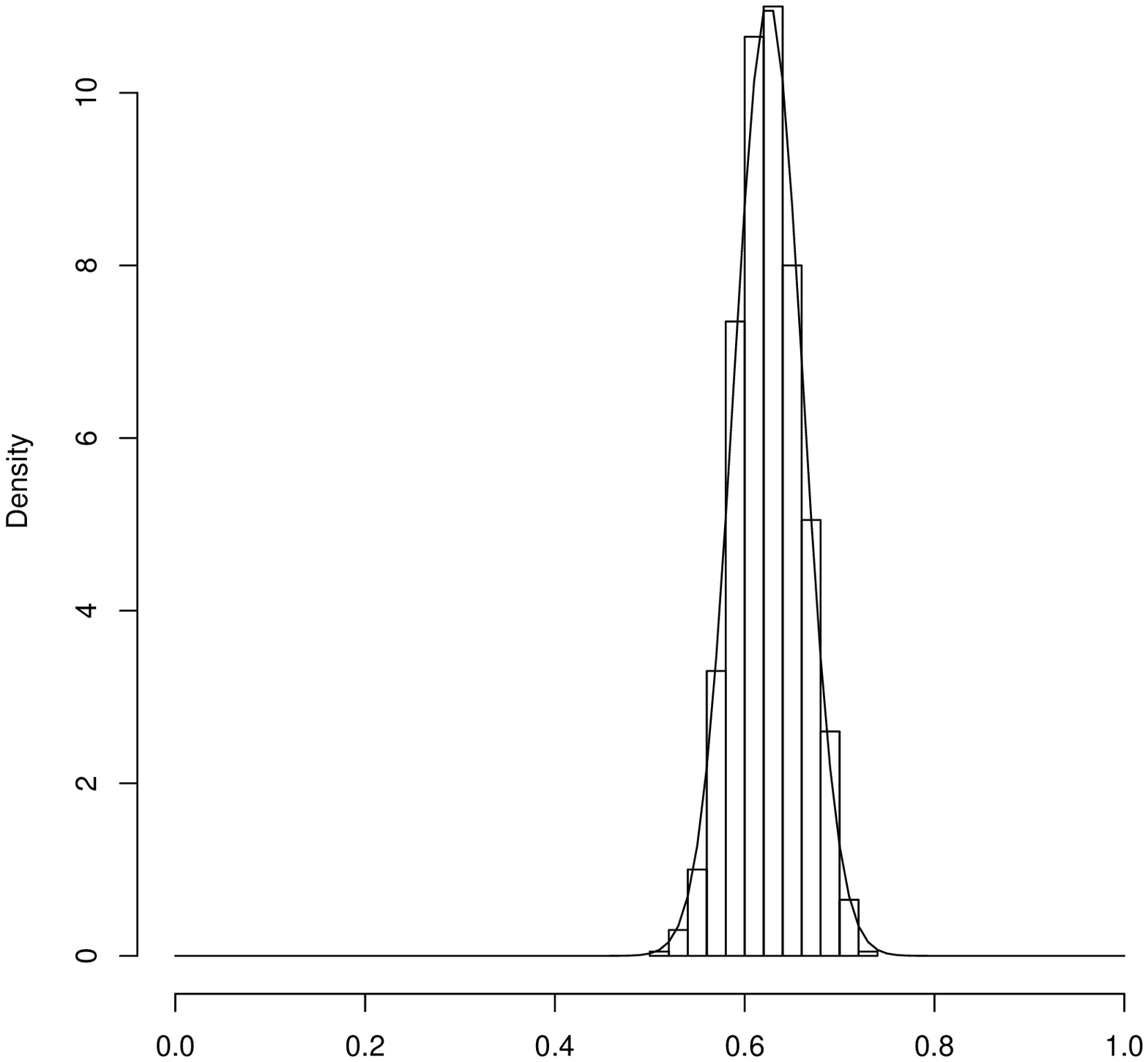}}}
\caption{
\label{fig:NormSkew}
Depicted are the distributions of
$\rho_n(2) \stackrel{\text{\scriptsize approx}}{\sim} \N\left(\frac{5}{8},\frac{25}{192n}\right)$
for $10,20,100$ (left to right).
Histograms are based on 1000 Monte Carlo replicates. Solid curves represent the approximating normal densities given in Theorem 2. Note that the vertical axes are differently scaled.
}
\end{figure}

\begin{figure}[ht]
\centering
\psfrag{Density}{ \Huge{\bfseries{Density}}}
\rotatebox{-90}{ \resizebox{2.1 in}{!}{ \includegraphics{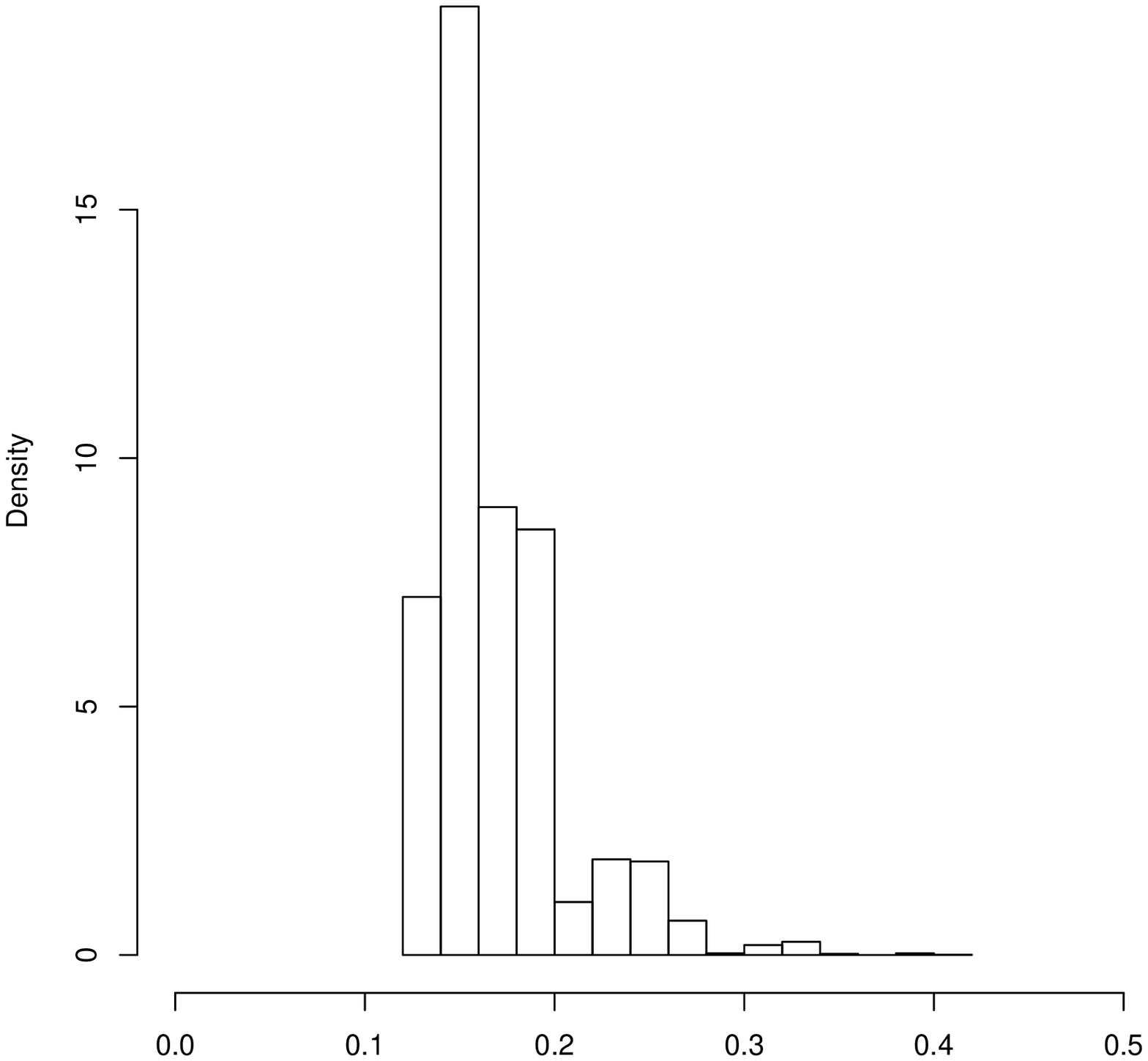}}}
\rotatebox{-90}{ \resizebox{2.1 in}{!}{ \includegraphics{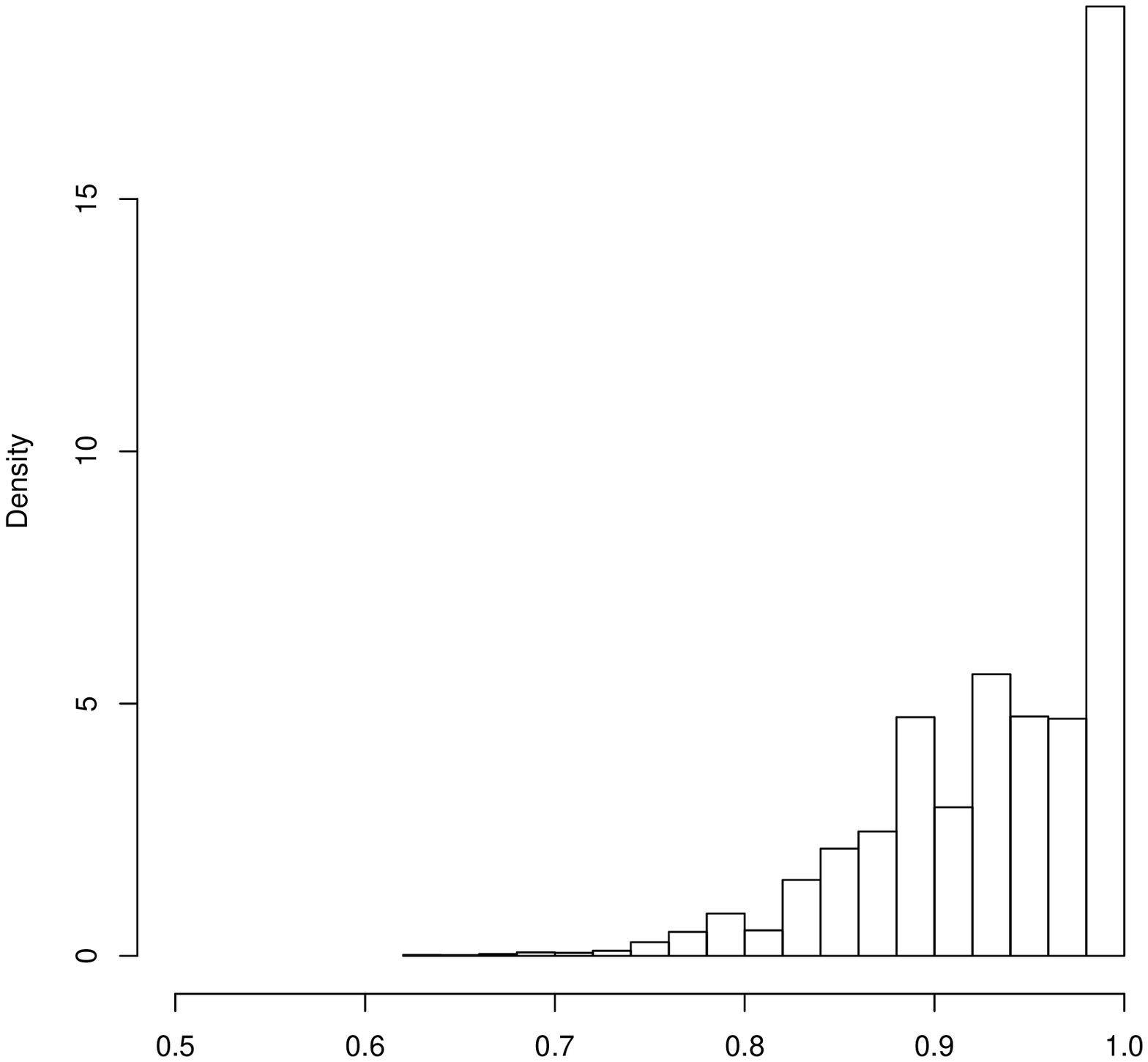}}}
\caption{
\label{fig:NormSkew1}
Depicted are the histograms for 10,000 Monte Carlo replicates of $\rho_{10}(1)$ (left) and $\rho_{10}(5)$ (right) indicating severe small sample skewness for extreme values of $r$.}
\end{figure}

Figure \ref{fig:NormSkew}
indicates that, for $r=2$,
the normal approximation is accurate even for small $n$
(although kurtosis may be indicated for $n=10$).
Figure \ref{fig:NormSkew1} demonstrates,
however, that severe skewness obtains for small values of $n$ and extreme values of $r$.

\vspace*{0.5 in}

Letting $H_n(r) = \sum_{i=1}^n h\bigl( X_i,X_{n+1} \bigr)$,
the exact distribution of $\rho_n(r)$
can be written as the recurrence
$$(n+1)\,n\,\rho_{n+1}(r) \stackrel{d}{=} n\,(n-1)\,\rho_n(r) + H_n(r)$$
by noting that the conditional random variable
$H_n(r)|X_{n+1}$
is the sum of $n$ independent and identically distributed random variables.
Alas, this calculation is also tedious for large $n$.

\subsection{Asymptotic Normality Under the Alternatives}
Asymptotic normality of relative density of the proximity catch digraphs under the alternative hypotheses of segregation and association can be established by the same method as under the null hypothesis. Let $\E^S_{\epsilon}[\cdot]$ ( $\E^A_{\epsilon}[\cdot]$) be the expectation with respect to the uniform distribution under the segregation ( association ) alternatives with $\epsilon \in \left( 0,\sqrt{3}/3 \right)$.

{\bf Theorem 3}
Let $\mu_S(r,\epsilon)$ be the mean $\E^S_{\epsilon}[h_{12}]$ and $\nu_S(r,\epsilon)$ be the covariance, $\Cov^S_{\epsilon}[h_{12},h_{13}]$ for $r \in [1,\infty]$ and $\epsilon \in \bigl[ 0,\sqrt{3}/3 \bigr)$ under $H^S_{\epsilon}$, $\sqrt{n}\bigl( \rho_n(r)-\mu_S(r,\epsilon) \bigr) \stackrel {\mathcal L}{\longrightarrow} \N(0,\nu_S(r,\epsilon))$ for the values of $(r,\epsilon)$ for which $\nu_S(r,\epsilon)>0$.  Likewise for $H^A_{\epsilon}$.

{\bfseries Sketch of Proof:}
Under the alternatives, i.e. $\epsilon>0$ ,
$\rho_n(r)$ is a $U$-statistic
with the same symmetric kernel $h_{ij}$ as in the null case. Under $H^S_{\epsilon}$,
the mean $\mu_S(r,\epsilon)=\E^S_{\epsilon}[\rho_n(r)] = \E^S_{\epsilon}[h_{12}]/2$,
now a function of both $r$ and $\epsilon$, is again in $[0,1]$.
The asymptotic variance $\nu_S(r,\epsilon)=\Cov^S_{\epsilon}[h_{12},h_{13}]$, also a function of both $r$ and $\epsilon$,
is bounded above by $1/4$, as before. Thus asymptotic normality obtains provided $\nu_S(r,\epsilon) > 0$;
otherwise $\rho_n(r)$ is degenerate. Likewise for $H^A_{\epsilon}$.

The explicit forms of $\mu_S(r,\epsilon)$ and $\mu_A(r,\epsilon)$ are given, defined piecewise, in Section \ref{sec:pi_a-r-epsilon}. Sample values of $\mu_S(r,\epsilon)$, $\nu_S(r,\epsilon)$, and $\mu_A(r,\epsilon)$, $\nu_A(r,\epsilon)$ are given in Section \ref{sec:Hodges-Lehmann-Seg} under segregation with $\epsilon=\sqrt{3}/8,\sqrt{3}/4,\,2\,\sqrt{3}/7$ and in Section \ref{sec:Hodges-Lehmann-Agg} under association with $\epsilon=5\,\sqrt{3}/24,\,\sqrt{3}/12,\,\sqrt{3}/21$.
Note that under $H^S_{\epsilon}$,
$$\nu_S(r,\epsilon)>0 \text{ for } (r,\epsilon) \in \left[ 1,\sqrt{3}/(2 \epsilon) \right) \times \left( 0,\sqrt{3}/4 \right] \cup \left[ 1,\sqrt{3}/\epsilon-2 \right) \times \left( \sqrt{3}/4,\sqrt{3}/3 \right),$$
and under $H^A_{\epsilon}$,
$$\nu_A(r,\epsilon)>0\text{ for }(r,\epsilon)\in (1,\infty) \times \left( 0,\sqrt{3}/3 \right) \cup \{1\} \times \left( 0,\sqrt{3}/12 \right). \;\;\blacksquare$$

Notice that under the association alternatives
any $r \in (1,\infty)$ yields asymptotic normality
for all $\epsilon \in \left( 0,\sqrt{3}/3 \right)$,
while under the segregation alternatives
only $r=1$ yields this universal asymptotic normality.

\subsection{The Test and Analysis}
The relative density of the proximity catch digraph
is a test statistic for the segregation/association alternative;
rejecting for extreme values of $\rho_n(r)$ is appropriate
since under segregation we expect $\rho_n(r)$ to be large,
while under association we expect $\rho_n(r)$ to be small.
Using the test statistic
\begin{equation}
R = \frac{\sqrt{n} \bigl( \rho_n(r) - \mu(r) \bigr)}{\sqrt{\nu(r)}},
\end{equation}
the asymptotic critical value
for the one-sided level $\alpha$ test against segregation
is given by
\begin{equation}
z_{\alpha} = \Phi^{-1}(1-\alpha)
\end{equation}
where $\Phi(\cdot)$ is the standard normal distribution function.
Against segregation, the test rejects for $R>z_{\alpha}$ and against association,
the test rejects for $R<z_{1-\alpha}$.

\subsection{Consistency}
{\bf Theorem}
The test against $H^S_{\epsilon}$ which rejects for $R>z_{1-\alpha}$
and
the test against $H^A_{\epsilon}$ which rejects for $R<z_{\alpha}$
are consistent for $r \in [1,\infty)$ and $\epsilon \in \left( 0,\sqrt{3}/3 \right)$.

{\bfseries Proof:}
Since the variance of the asymptotically normal test statistic,
under both the null and the alternatives,
converges to 0 as $n \rightarrow \infty$
(or is degenerate),
it remains to show that the mean under the null, $\mu(r)=\E[\rho_n(r)]$, is less than (greater than) the mean under the alternative, $\mu_S(r,\epsilon)=\E^S_{\epsilon}[\rho_n(r)]$ against segregation ($\mu_A(r,\epsilon)=\E^A_{\epsilon}[\rho_n(r)]$ against association) for $\epsilon > 0$.
Whence it will follow that power converges to 1 as $n \rightarrow \infty$.

 Detailed analysis of $\mu_S(r,\epsilon)$ in Appendix 2 indicates that under segregation $\mu_S(r,\epsilon)>\mu(r)$ for all $\epsilon >0$ and $r \in [1,\infty)$. Likewise, detailed analysis of $\mu_A(r,\epsilon)$ in Appendix 2 indicates that under association $\mu_A(r,\epsilon)<\mu(r)$ for all $\epsilon >0$ and $r \in [1,\infty)$.  Hence the desired result follows for both alternatives. $\blacksquare$

{\bf Remark:}
In fact, the analysis of $\mu_S(r,\epsilon)$ and $\mu_A(r,\epsilon)$ under the alternatives reveals more than what is required for consistency.  Under segregation, the analysis indicates that $\mu_S(r,\epsilon_1) < \mu_S(r,\epsilon_2)$ for $\epsilon_1<\epsilon_2$.  Likewise, under association, the analysis indicates that $\mu_A(r,\epsilon_1) > \mu_A(r,\epsilon_2)$ for $\epsilon_1<\epsilon_2$.  $\square$

\subsection{Monte Carlo Power Analysis Under Segregation}
\label{sec:monte-carlo-seg}
In segregation alternatives with $\epsilon>0$, we implement the above described Monte Carlo experiment for various values of $r$ (for which $\rho_n(r)$ is non-degenerate).
Recall that $\rho_n(r)$ is degenerate for large $r$ at each $\epsilon>0$.  In particular, $\rho_n\left( r,\sqrt{3}/8 \right)$ is degenarate for $r \ge 4 $, $\rho_n\left( r,\sqrt{3}/4 \right)$ is degenarate for $r \ge 2 $, and $\rho_n\left( r,2\,\sqrt{3}/7 \right)$ is degenarate for $r \ge 3/2$.

Let $\rho_k(n)$ be the empirical relative density for experiment $k$ and $\rho_{(j)}(n)$ be the $j^{th}$ (ordered) empirical relative density for $j=1,\ldots, N$.  Then for each $r$ value, we estimate the empirical critical value $\widehat{C}^S_n:=\rho_{(\lceil (1-\alpha)\,N \rceil)}(n)$ and the empirical significance level $\widehat{\alpha}^S_{mc}(n):=\frac{1}{N}\sum_{j=1}^{N}\I\left( \rho_{j}(n) > \widehat{C}^S_n \right)$ under $H_0$ and the empirical power $\widehat{\beta}^S_{mc}(n,\epsilon):=\frac{1}{N}\sum_{j=1}^{N}\I\left( \rho_{j} > \widehat{C}^S_n \right)$ under $H^S_{\epsilon}$ with $\epsilon=\sqrt{3}/8,\,\sqrt{3}/4,\,2\,\sqrt{3}/7$.

\begin{figure}[]
\centering
\psfrag{kernel density estimate}{ \Huge{\bfseries{kernel density estimate}}}
\psfrag{relative density}{ \Huge{\bfseries{relative density}}}
\rotatebox{-90}{ \resizebox{1.7 in}{!}{ \includegraphics{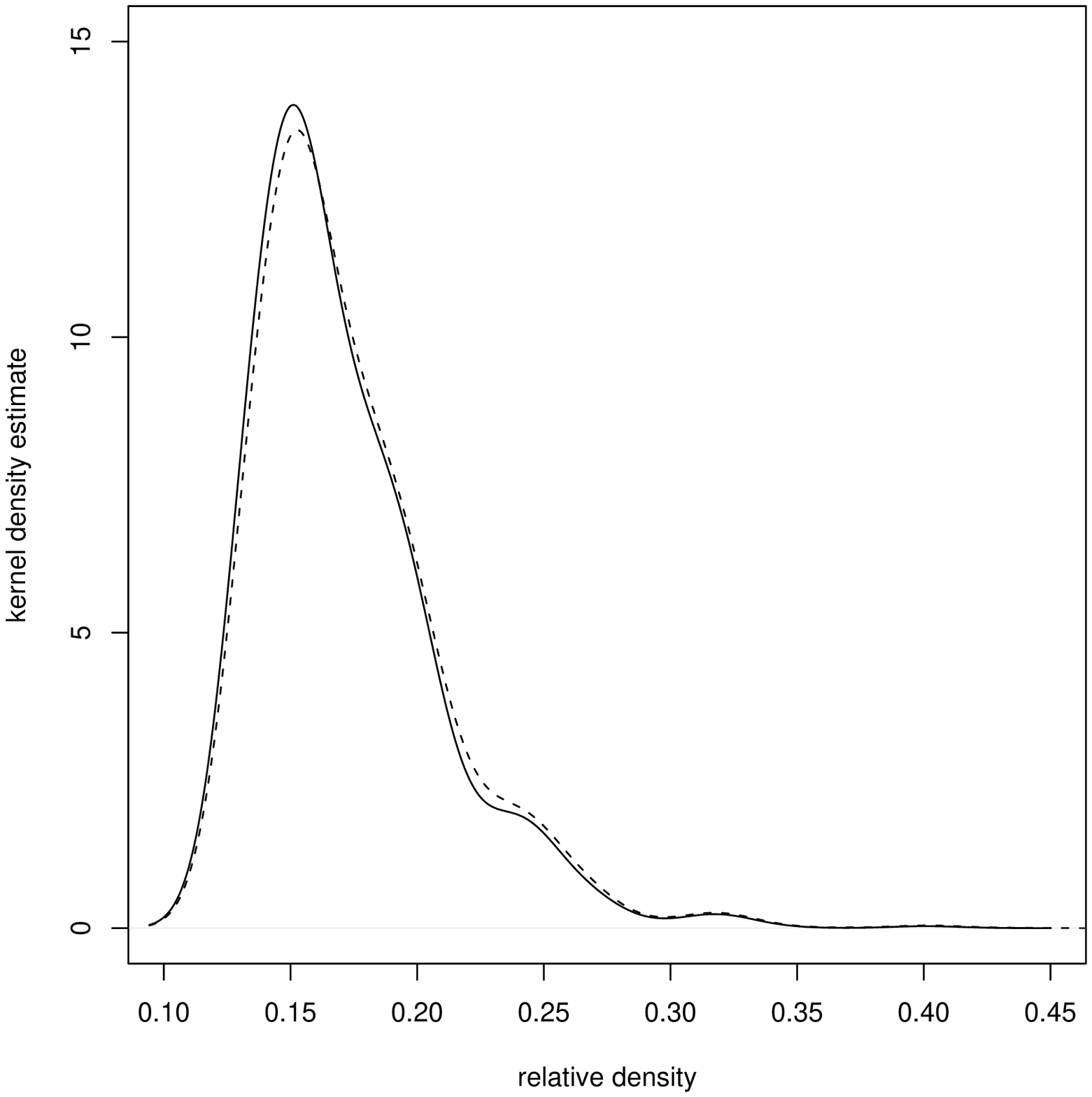}}}
\rotatebox{-90}{ \resizebox{1.7 in}{!}{ \includegraphics{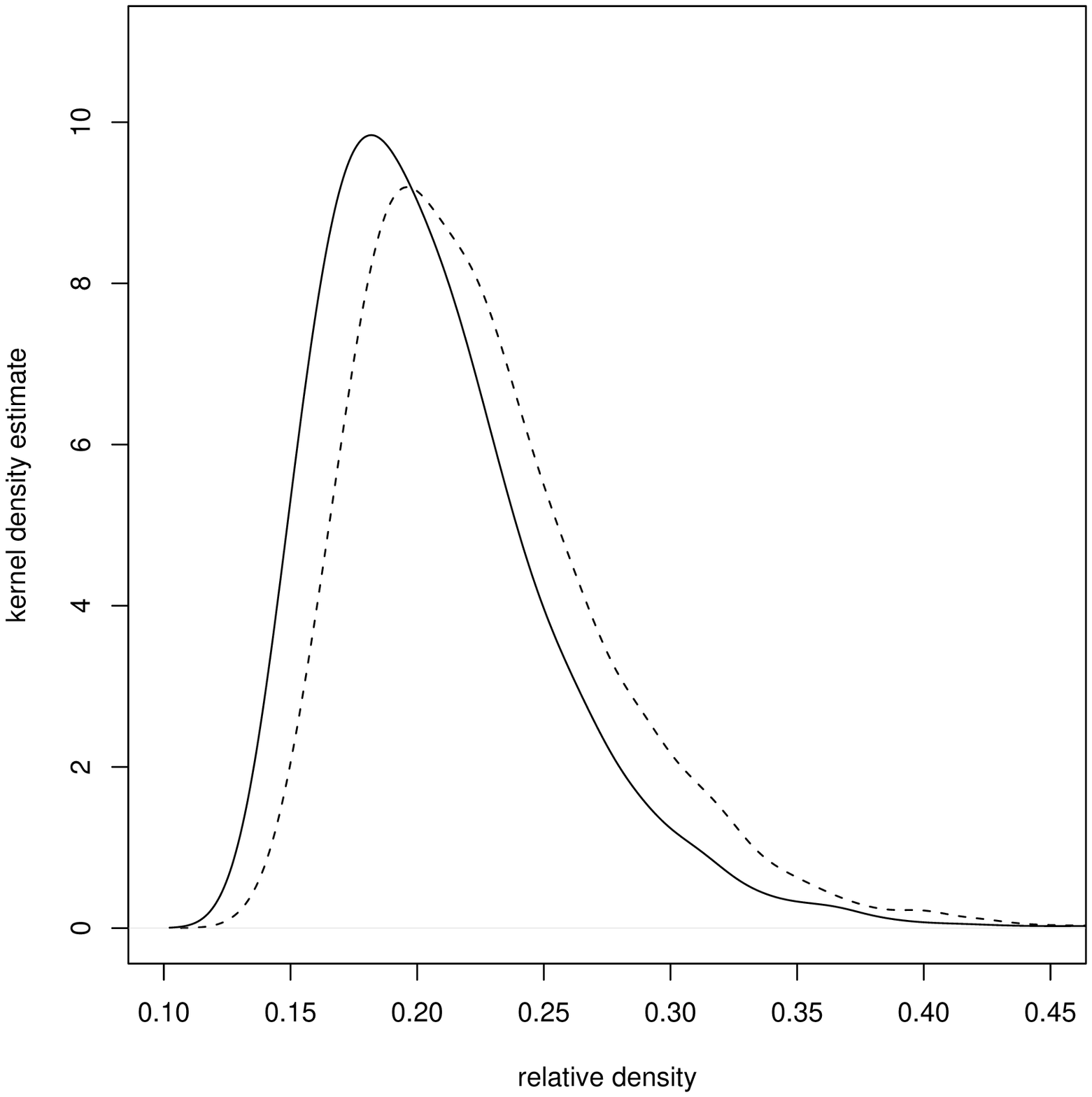}}}
\rotatebox{-90}{ \resizebox{1.7 in}{!}{ \includegraphics{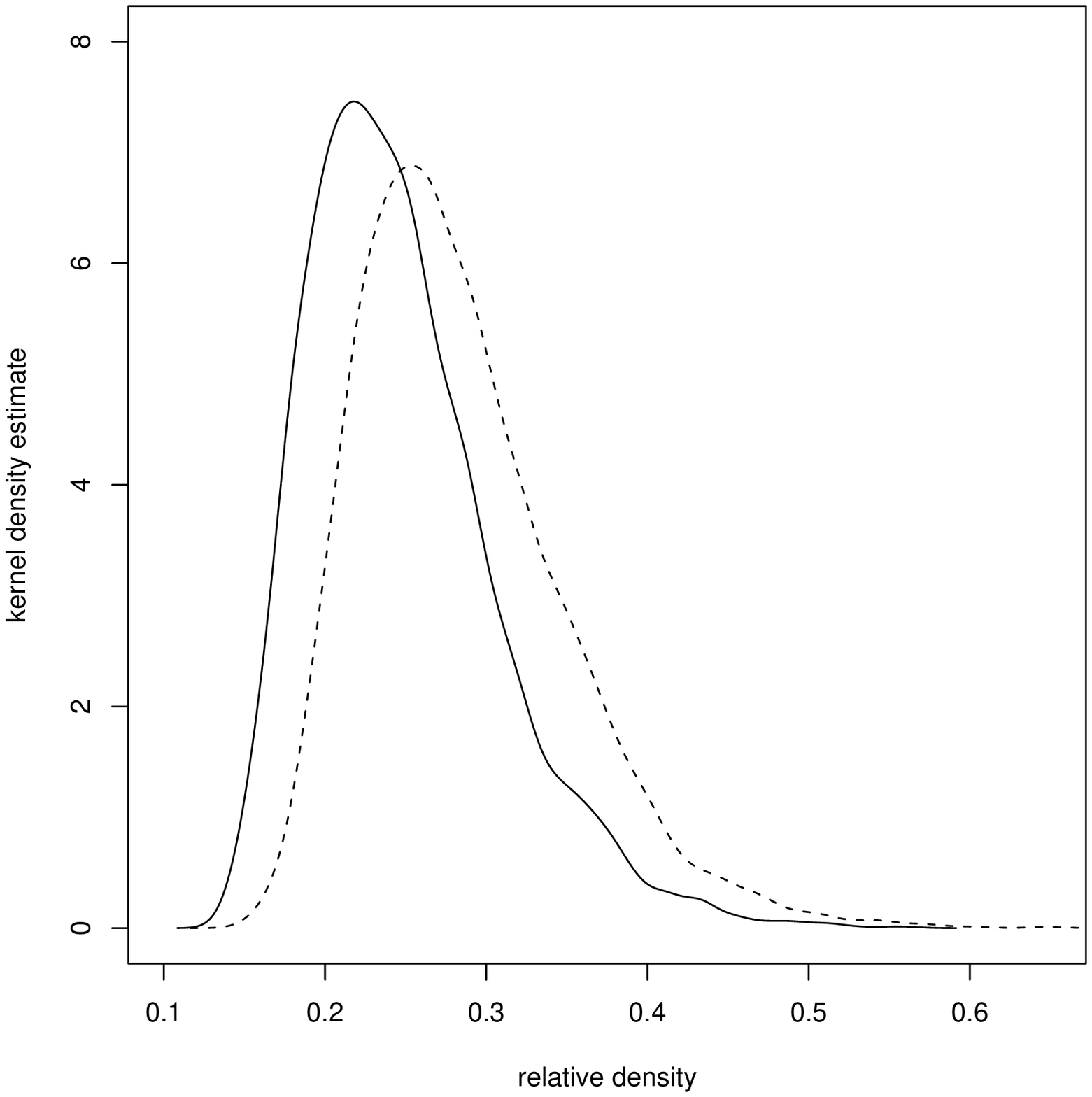}}}
\rotatebox{-90}{ \resizebox{1.7 in}{!}{ \includegraphics{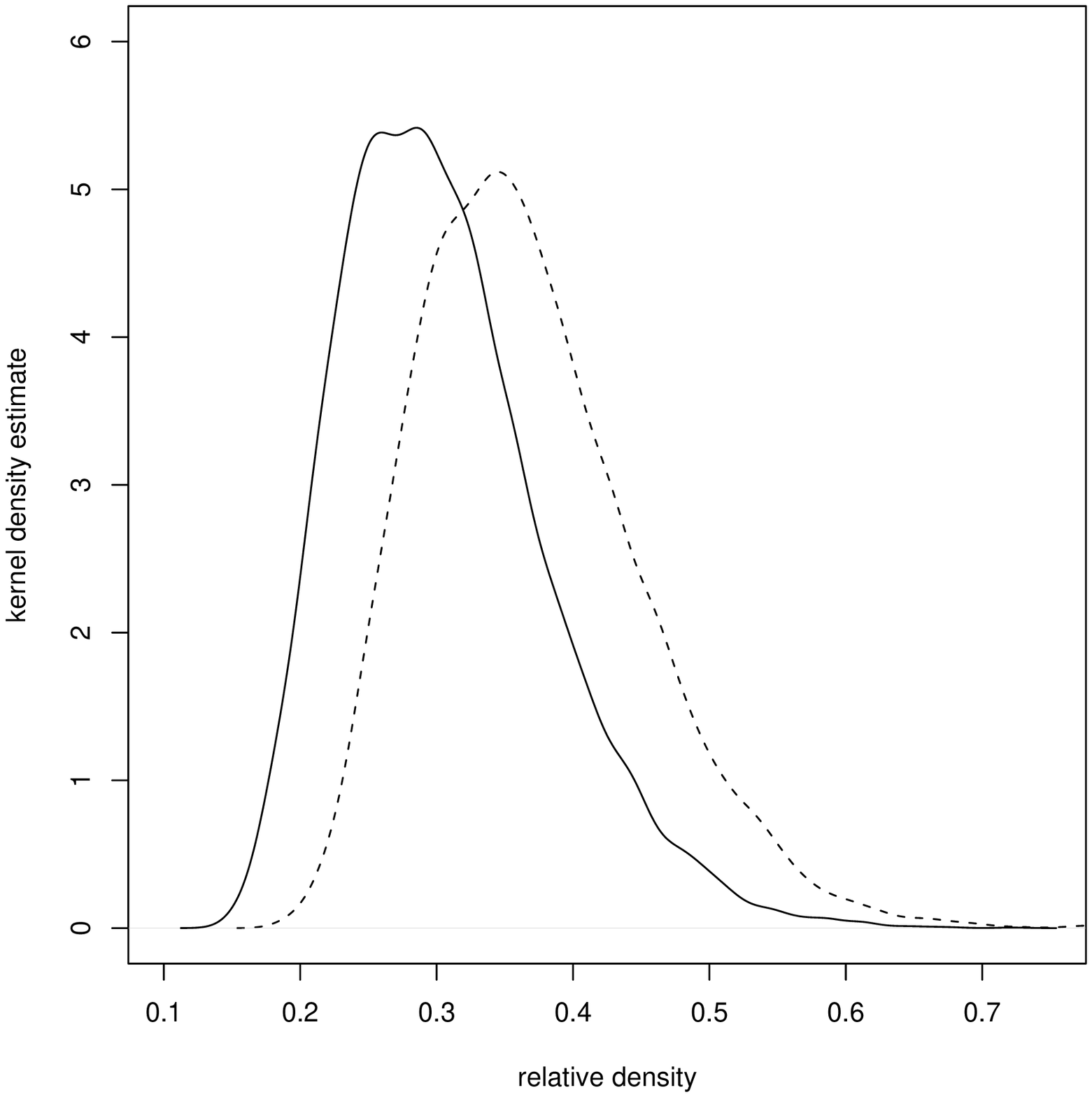}}}
\rotatebox{-90}{ \resizebox{1.7 in}{!}{ \includegraphics{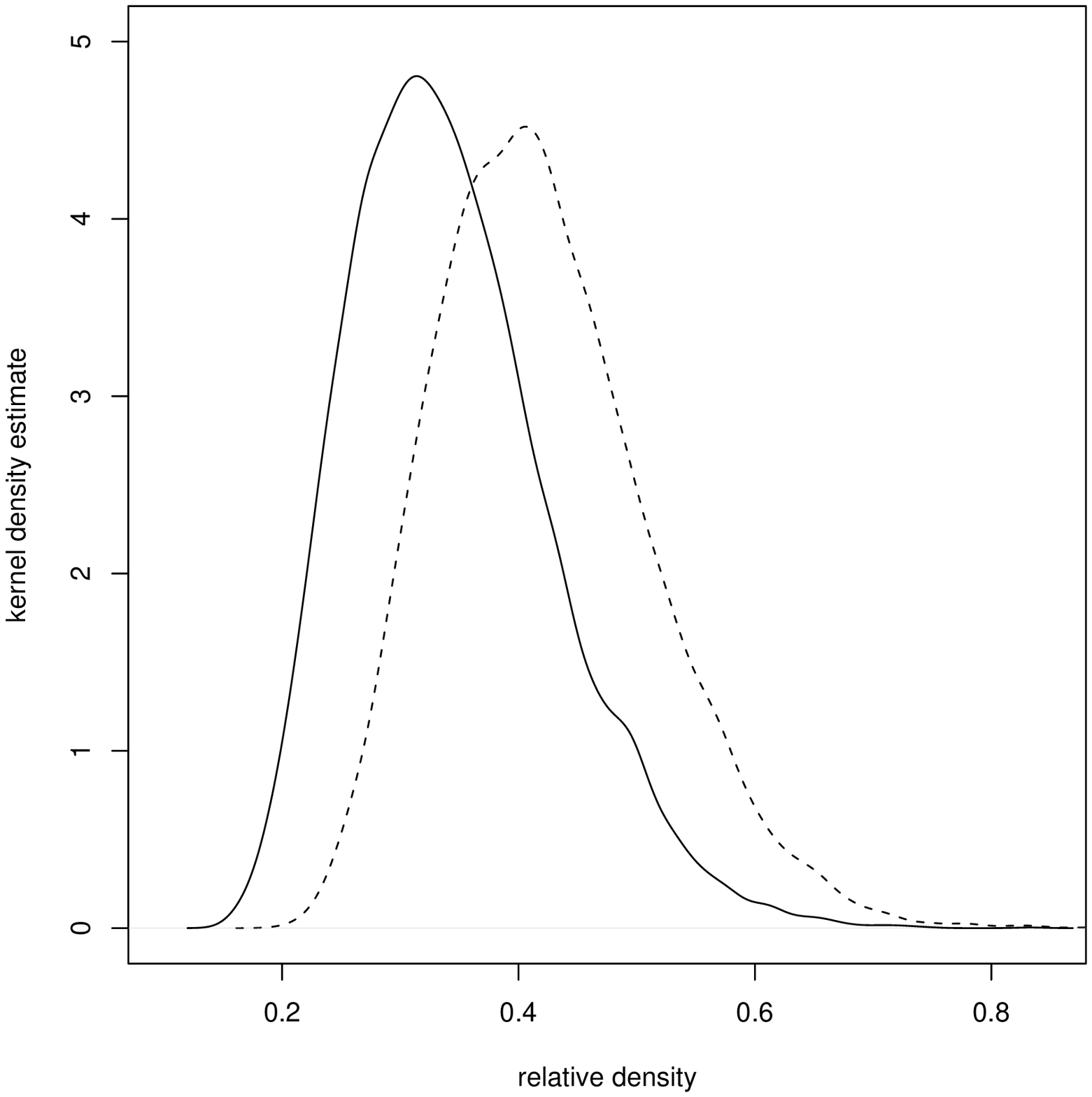}}}
\rotatebox{-90}{ \resizebox{1.7 in}{!}{ \includegraphics{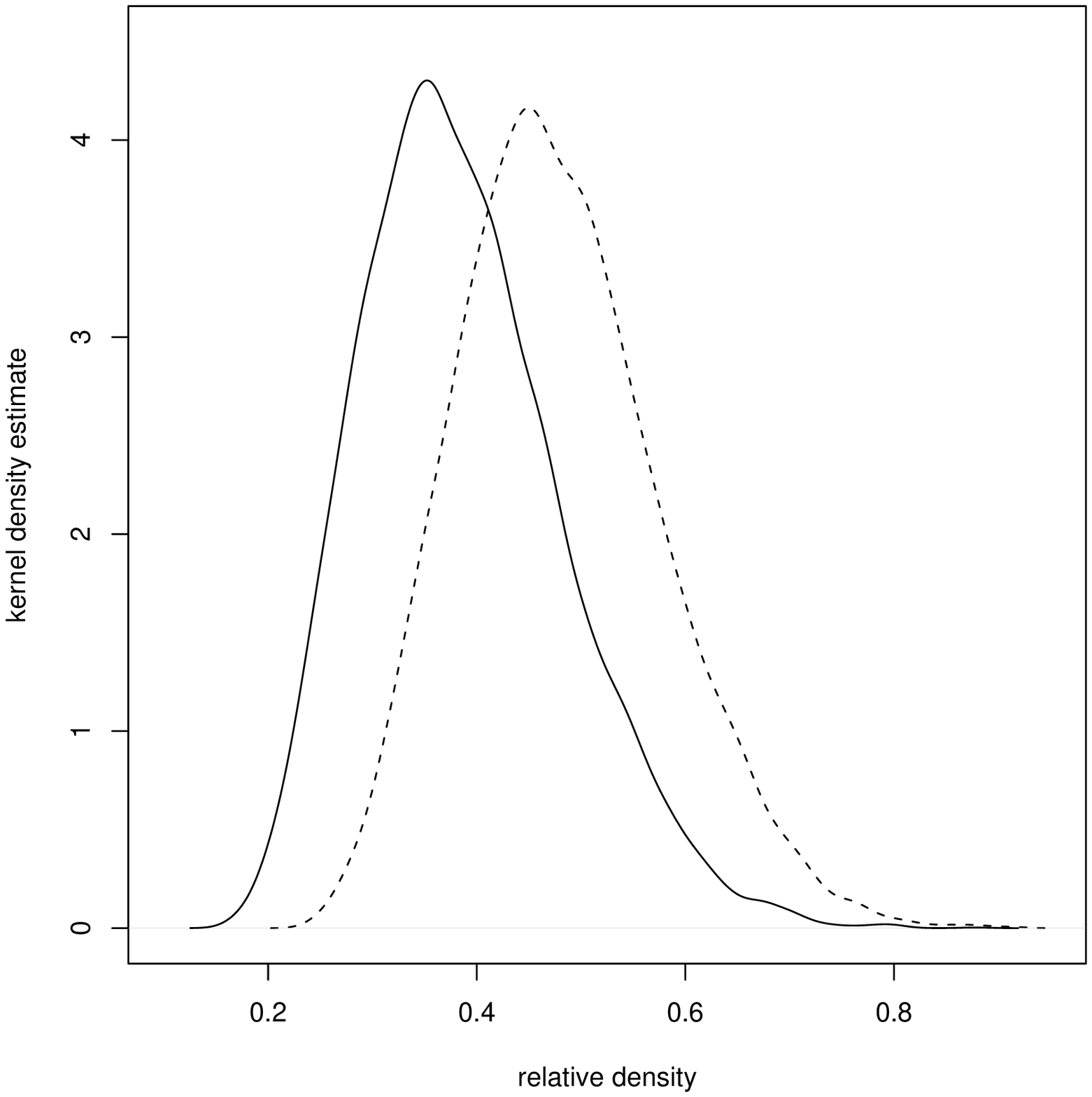}}}
\rotatebox{-90}{ \resizebox{1.7 in}{!}{ \includegraphics{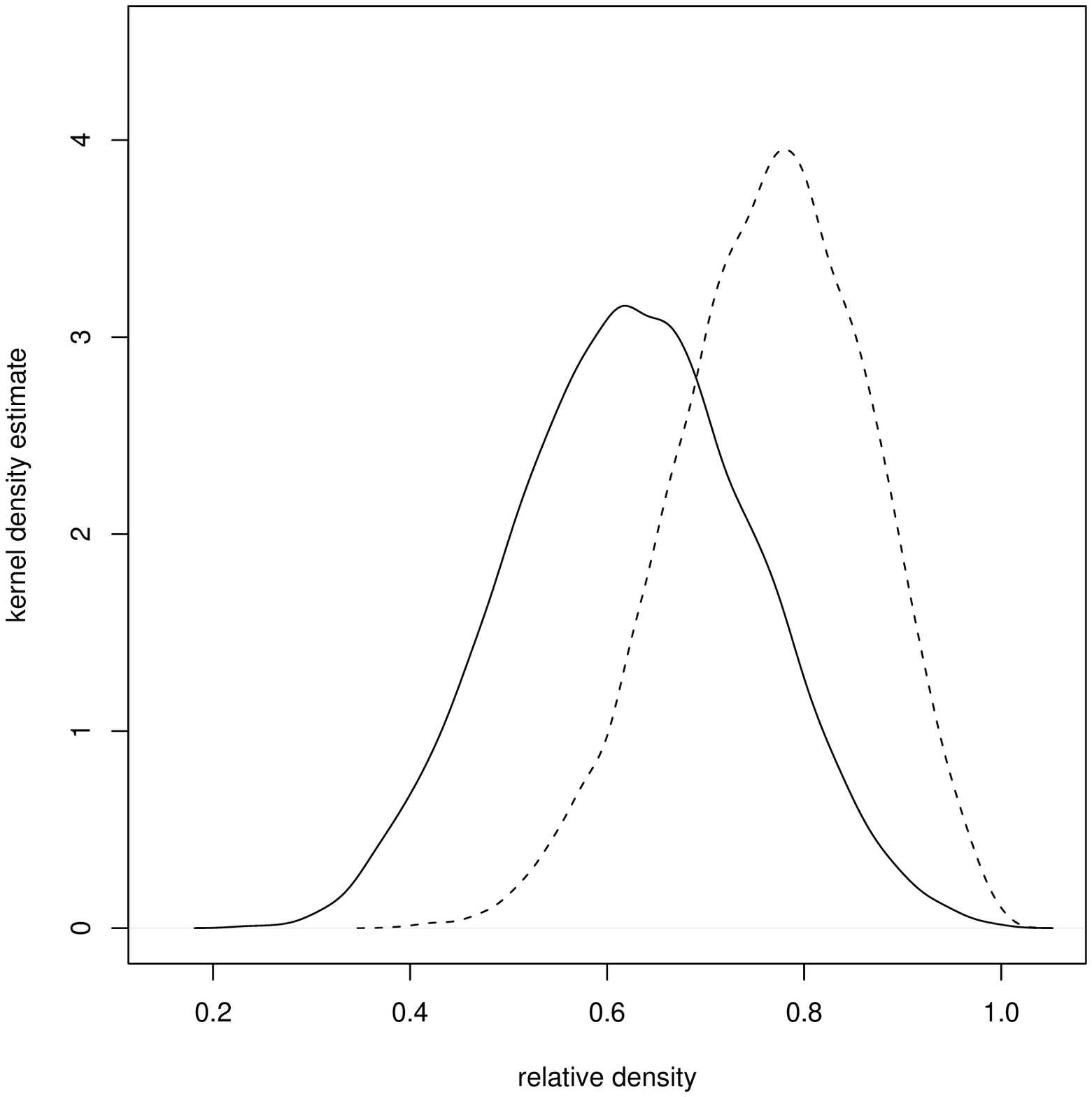}}}
\rotatebox{-90}{ \resizebox{1.7 in}{!}{ \includegraphics{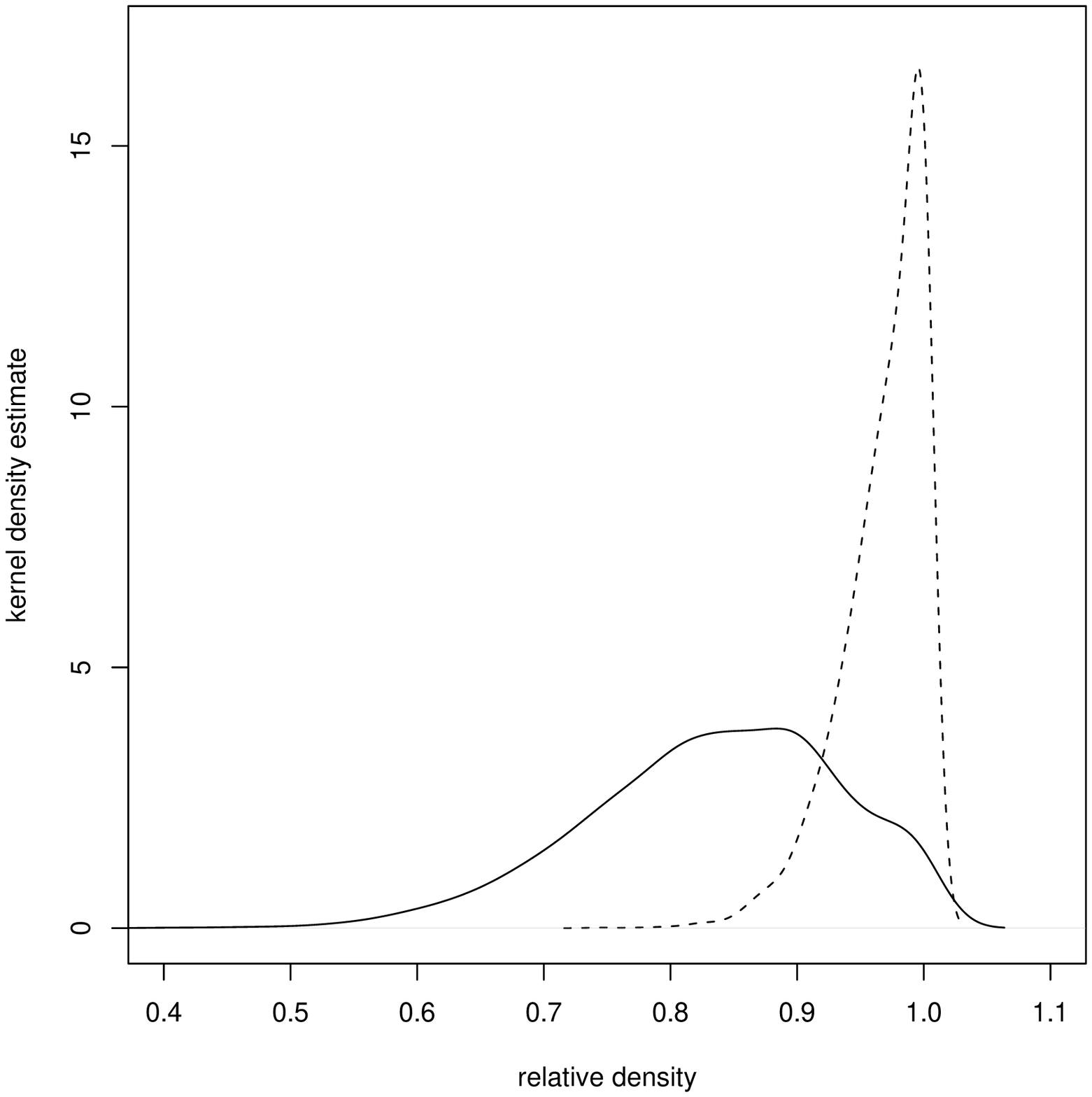}}}
\caption{\label{fig:segsim1-8}
Kernel density estimates for the null (solid) and the segregation alternative $H^S_{\sqrt{3}/8}$ (dashed) for $r=1,\,11/10,\,6/5,\,4/3,\,\sqrt{2},\,3/2,\,2,$ and $3$ (left-to-right).
}
\end{figure}

For segregation with $\epsilon=\sqrt{3}/8 \approx .2165$, we run the Monte Carlo experiments for eight $r$ values: $1,\,11/10,\,6/5,\,4/3,\,\sqrt{2},\,3/2,\,2,$ and $3$. In Figure \ref{fig:segsim1-8} are the kernel density estimates for the null case and the segregation alternative with $\epsilon=\sqrt{3}/8$, for the eight $r$ values with $n=10$ and $N=10,000$.  Observe that under both $H_0$ and $H^S_{\sqrt{3}/8}$, kernel density estimates are skewed right for $r=1,\,11/10 $, (with skewness increasing as $r$ gets smaller) and kernel density estimates are almost symmetric for $r=6/5,\,4/3,\,\sqrt{2}, 3/2,\, 2 $, with most symmetry occurring at $r=3/2 $, kernel density estimate is skewed left for $r=3$ (with skewness increasing as $r$ gets larger).

The empirical critical values, empirical significance levels, and empirical power estimates under $H^S_{\sqrt{3}/8}$ are presented in Table \ref{tab:emp-val-S1}.

\begin{table}[]
\centering
\begin{tabular}{|c|c|c|c|c|c|c|c|c|}
\hline
$r$  & 1 & 11/10 & 6/5 & 4/3 & $\sqrt{2}$ & 3/2 & 2 & 3 \\
\hline
$\widehat{C}^S_n$ & $0.2\bar{4}$ & .3 & $.3\bar{5}$ & $.\bar 4$ & .5 & $.\bar 5$ & $.8\bar{2}$ & $0.9\bar{8}$\\
\hline
$\widehat{\alpha}^S_{mc}(10)$      & .0324 & .0403 & .0484 & .0442 &  .0446 & .0492 & .049 & .0389\\
\hline
$\widehat{\beta}^S_{mc}\left(10,\,\sqrt{3}/8 \right)$ & .0381 &  .0787 & .122 & .1571  &  .1719 & .1955 & .2791 & .2901\\
\hline
\end{tabular}
\caption{
\label{tab:emp-val-S1}
The empirical critical values, empirical significance levels, and empirical power estimates under $H^S_{\sqrt{3}/8}$, $N=10,000$, and $n=10$ at $\alpha=.05$.}
\end{table}

\begin{figure}[]
\centering
\psfrag{kernel density estimate}{ \Huge{\bfseries{kernel density estimate}}}
\psfrag{relative density}{ \Huge{\bfseries{relative density}}}
\rotatebox{-90}{ \resizebox{2.2 in}{!}{ \includegraphics{segsim2.ps}}}
\rotatebox{-90}{ \resizebox{2.2 in}{!}{ \includegraphics{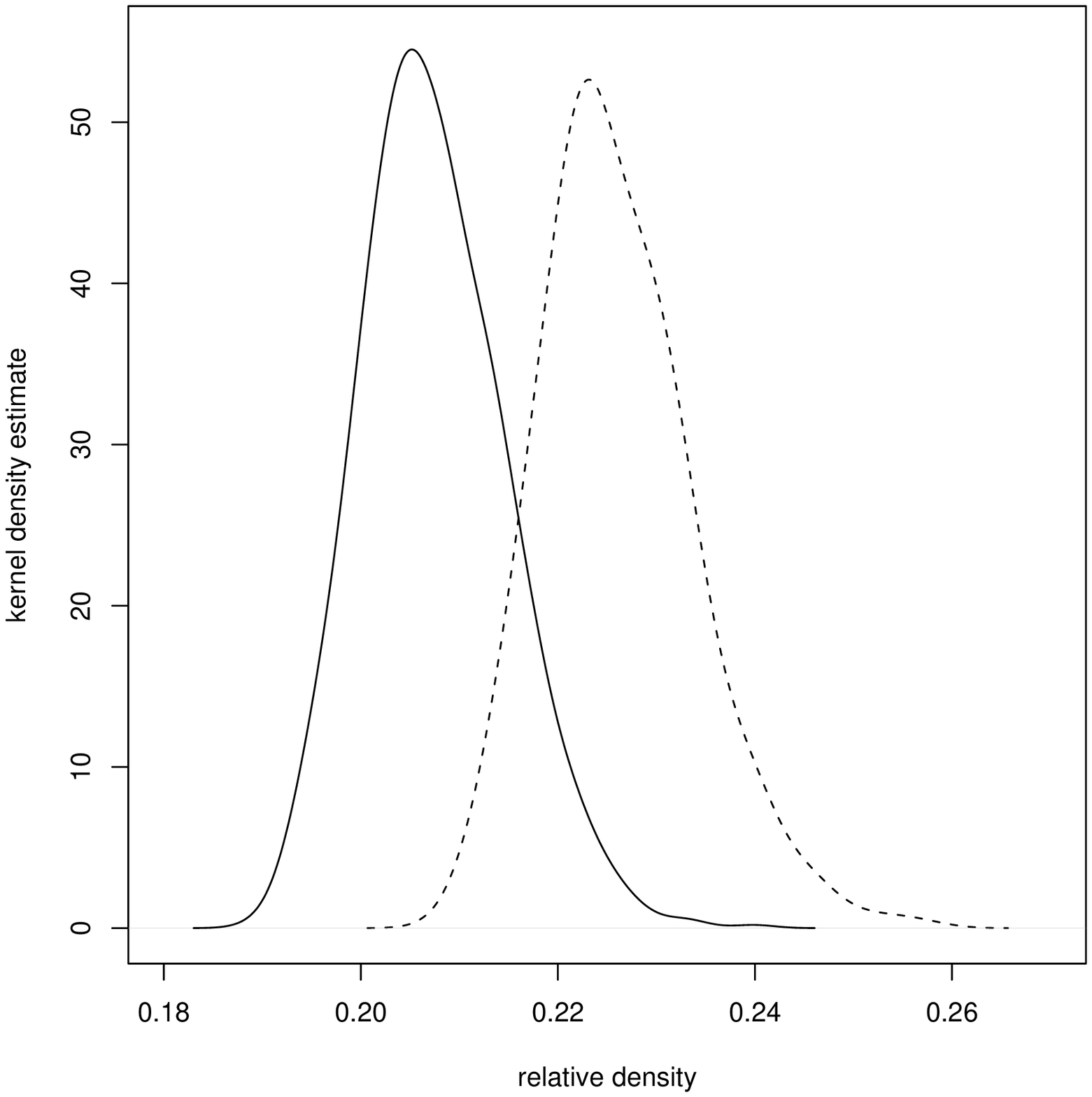}}}
\caption{ \label{fig:SegSimPowerPlots}
Two Monte Carlo experiments against the segregation alternative $H^S_{\sqrt{3}/8}$.
Depicted are kernel density estimates for $\rho_n(11/10)$ for
$n=10$ (left) and $n=100$ (right) under the null (solid) and alternative (dashed).
}
\end{figure}

In Figure \ref{fig:SegSimPowerPlots}, we present a Monte Carlo investigation
against the segregation alternative $H^S_{\sqrt{3}/8}$ for $r=11/10$, and $n=10$, $N=10,000$ (left), $n=100$, $N=1000$ (right).
With $n=10$, the null and alternative probability density functions
for $\rho_{10}(11/10)$ are very similar, implying small power
(10,000 Monte Carlo replicates yield $\widehat{\beta}^S_{mc}\left( 10,\,\sqrt{3}/8 \right)= 0.0787$, $\widehat{C}^S_n = 0.1\bar{5}$, and $\widehat{\alpha}^S_{mc}(10)= 0.0484$).
With $n=100$,
there is more separation
between null and alternative probability density functions;
for this case, 1000 Monte Carlo replicates yield $\widehat{\beta}^S_{mc}\left(100,\, \sqrt{3}/8 \right) = 0.77$, $\widehat{C}^S_n = 0.2203$, and $\widehat{\alpha}^S_{mc}(100)= 0.05$.
Notice also that the probability density functions are more skewed for $n=10$,
while approximate normality holds for $n=100$.

\begin{figure}[]
\centering
\psfrag{kernel density estimate}{ \Huge{\bfseries{kernel density estimate}}}
\psfrag{relative density}{ \Huge{\bfseries{relative density}}}
\rotatebox{-90}{ \resizebox{1.7 in}{!}{ \includegraphics{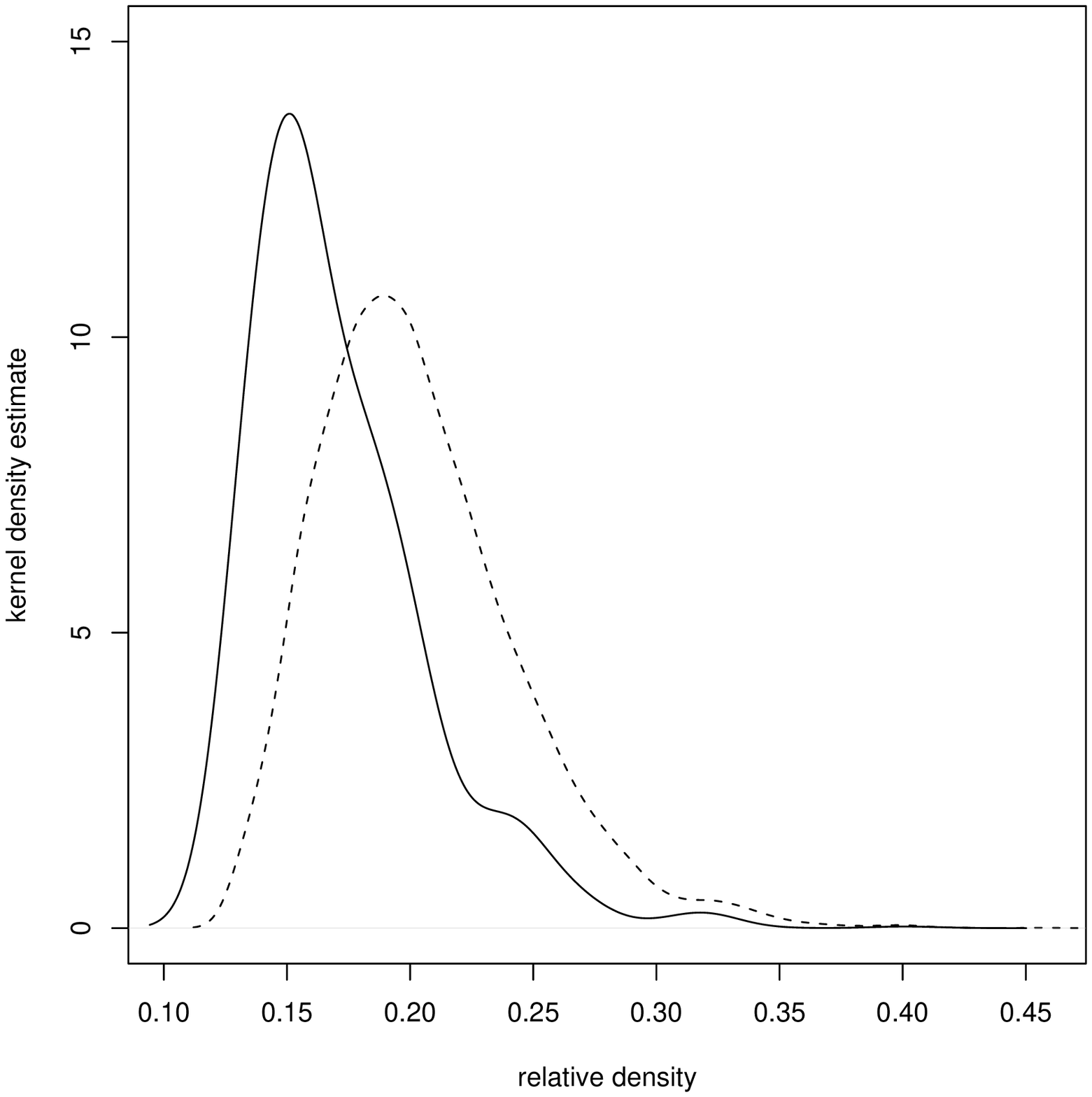}}}
\rotatebox{-90}{ \resizebox{1.7 in}{!}{ \includegraphics{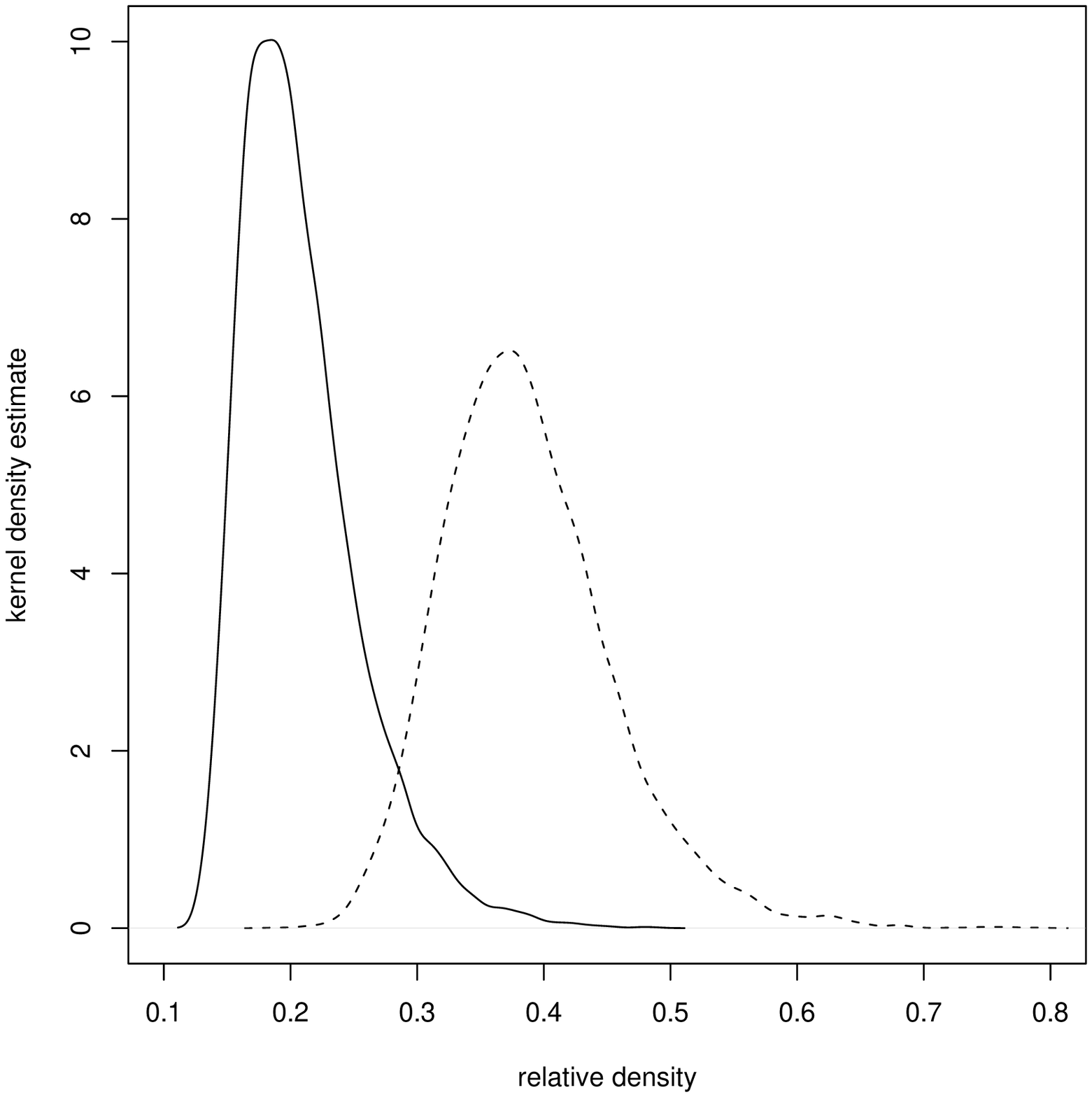}}}
\rotatebox{-90}{ \resizebox{1.7 in}{!}{ \includegraphics{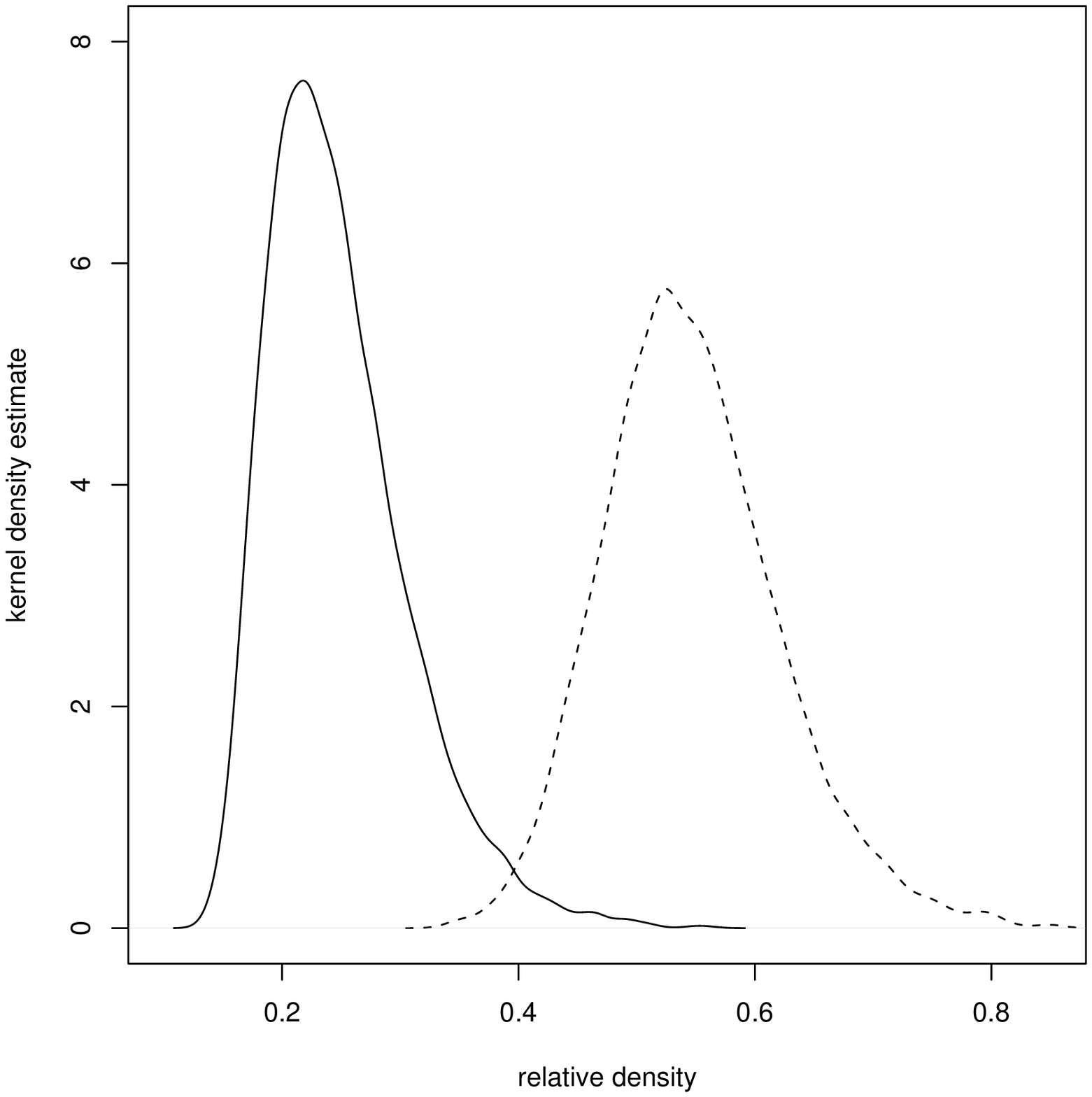}}}
\rotatebox{-90}{ \resizebox{1.7 in}{!}{ \includegraphics{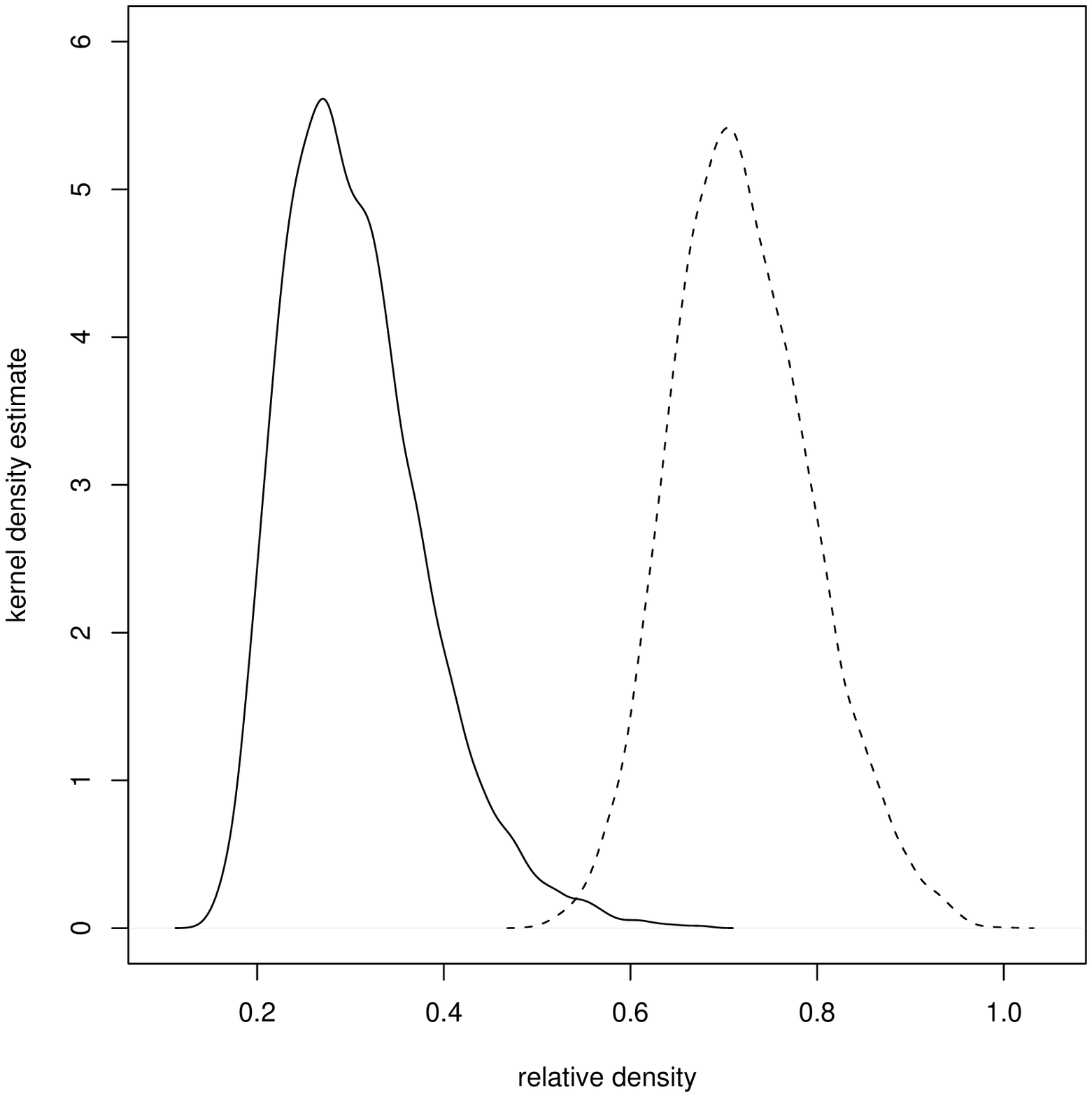}}}
\rotatebox{-90}{ \resizebox{1.7 in}{!}{ \includegraphics{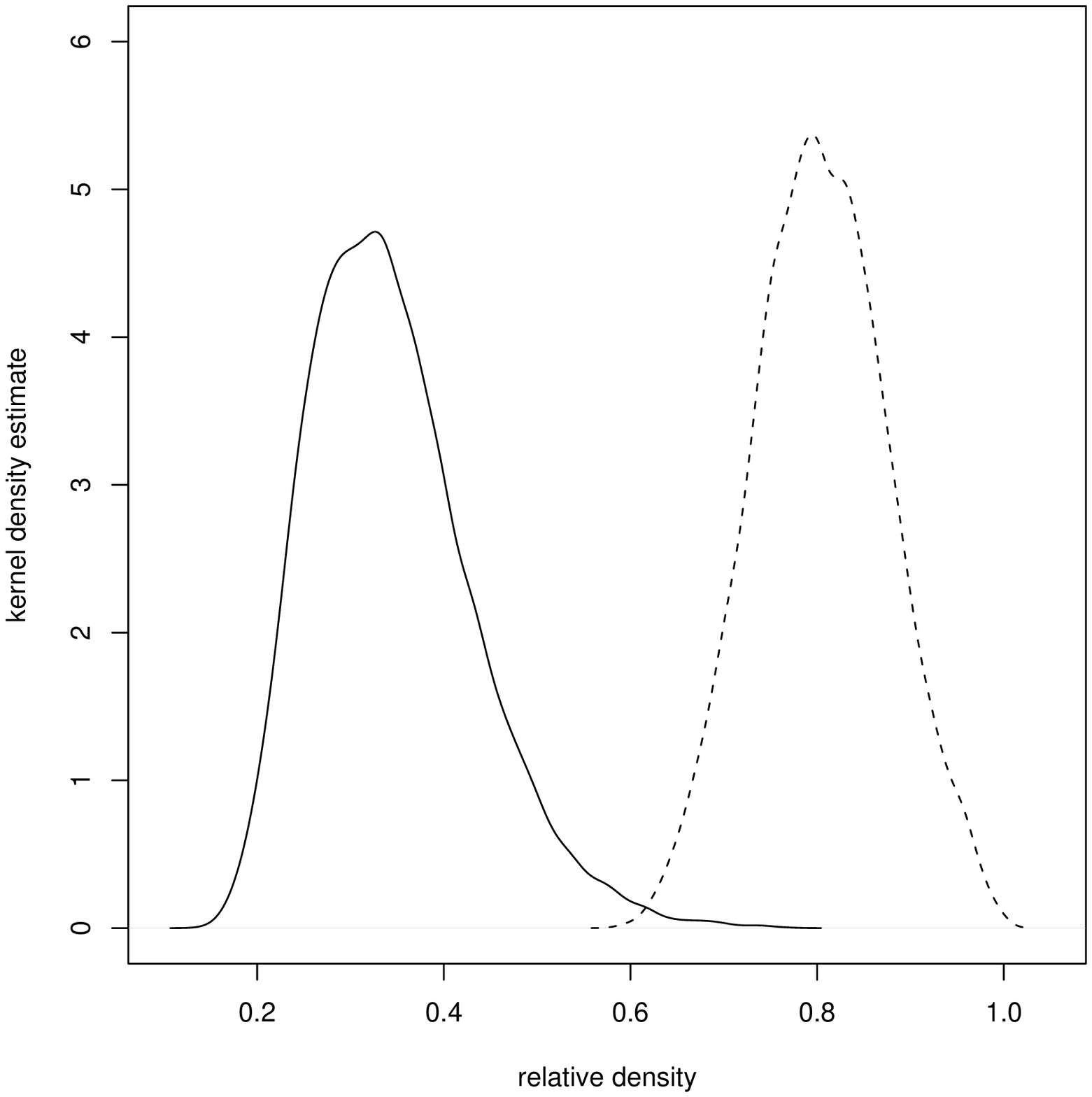}}}
\rotatebox{-90}{ \resizebox{1.7 in}{!}{ \includegraphics{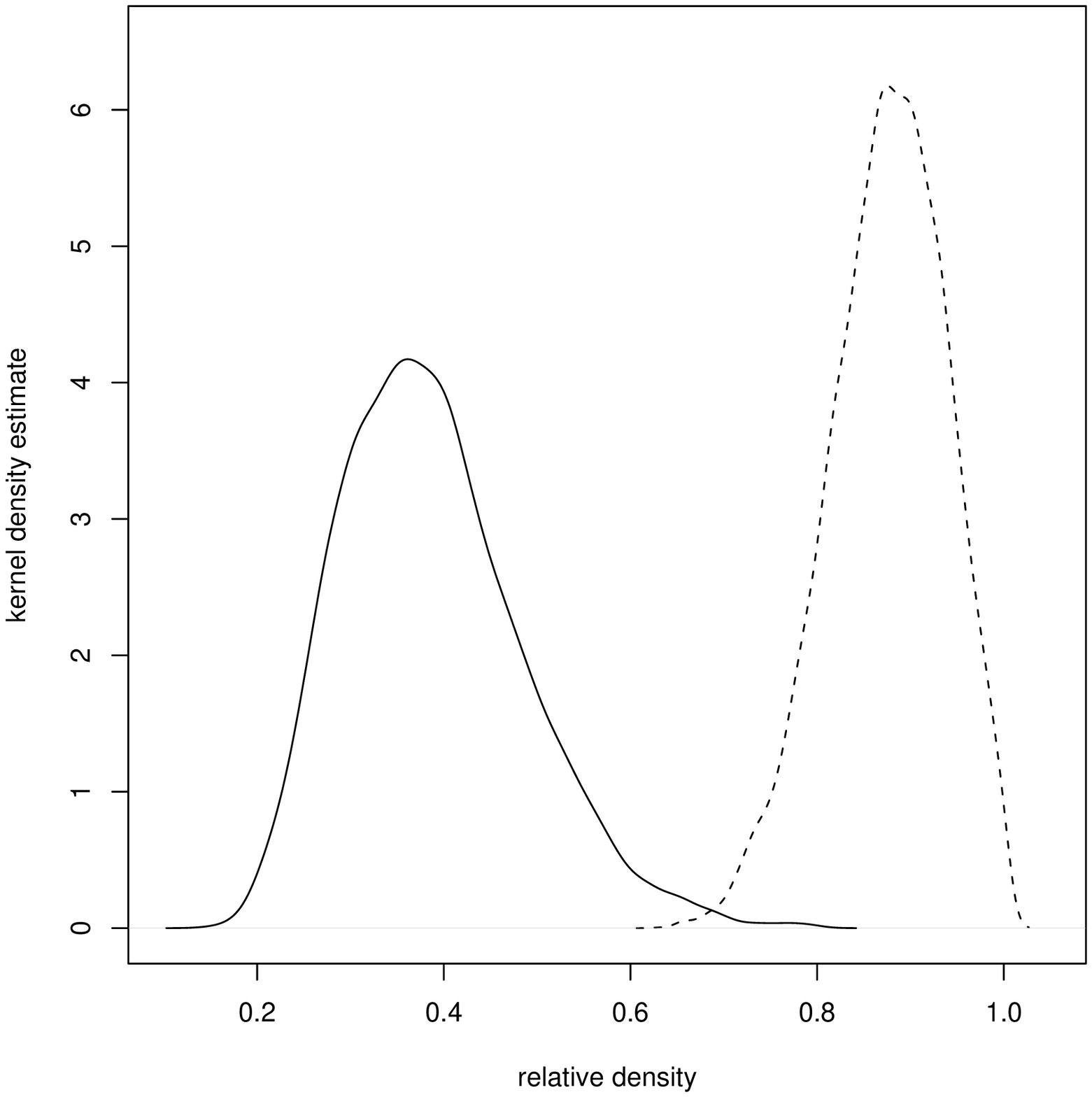}}}
\caption{\label{fig:seg2sim1-6}
Kernel density estimates for the null (solid) and the segregation alternative $H^S_{\sqrt{3}/4}$ (dashed) for $r=1,\,11/10,\,6/5,\,4/3,\,\sqrt{2},\,3/2$ (left-to-right).
}
\end{figure}

For segregation with $\epsilon=\sqrt{3}/4\,\approx .433$, we run the Monte Carlo experiments for six $r$ values, $1,\,11/10,\,6/5,\\
4/3,\,\sqrt{2},\,3/2$. In Figure \ref{fig:seg2sim1-6}, are the kernel density estimates for the null case and the segregation alternative with $\epsilon=\sqrt{3}/4$, for the six $r$ values with $n=10$ and $N=10,000$. Observe that under $H^S_{\sqrt{3}/4}$, kernel density estimate is skewed right for $r=1$ and kernel density estimates are almost symmetric for $r=11/10,\,6/5,\,4/3,\,\sqrt{2}$, with most symmetry occurring at $r=4/3$, kernel density estimate is skewed left for $r=3/2$.

The empirical critical values, empirical significance levels, and empirical power estimates under $H^S_{\sqrt{3}/4}$ are presented in Table \ref{tab:emp-val-S2}.
\begin{table}[ht]
\centering
\begin{tabular}{|c|c|c|c|c|c|c|}
\hline
$r$  & 1 & 11/10 & 6/5 & 4/3 & $\sqrt{2}$ & 3/2  \\
\hline
$\widehat{C}^S_n$         & $.2\bar{4}$ & .3 & $.3\bar{5}$ & $.\bar 4$ & .5 & $.\bar 5$ \\
\hline
$\widehat{\alpha}^S_{mc}(10)$      & .0318       & .0411 & .0479 & .0484 &  .0481 & .043 \\
\hline
$\widehat{\beta}^S_{mc}\left(10,\,\sqrt{3}/4 \right)$  & .1247        &  .9138 & .998 & 1.0  &  1.0 & 1.0 \\
\hline
\end{tabular}
\caption{
\label{tab:emp-val-S2}
The empirical critical values, empirical significance levels, and empirical power estimates under $H^S_{\sqrt{3}/4}$, $N=10,000$, and $n=10$ at $\alpha=.05$.}
\end{table}

\begin{figure}[]
\centering
\psfrag{kernel density estimate}{ \Huge{\bfseries{kernel density estimate}}}
\psfrag{relative density}{ \Huge{\bfseries{relative density}}}
\rotatebox{-90}{ \resizebox{1.7 in}{!}{ \includegraphics{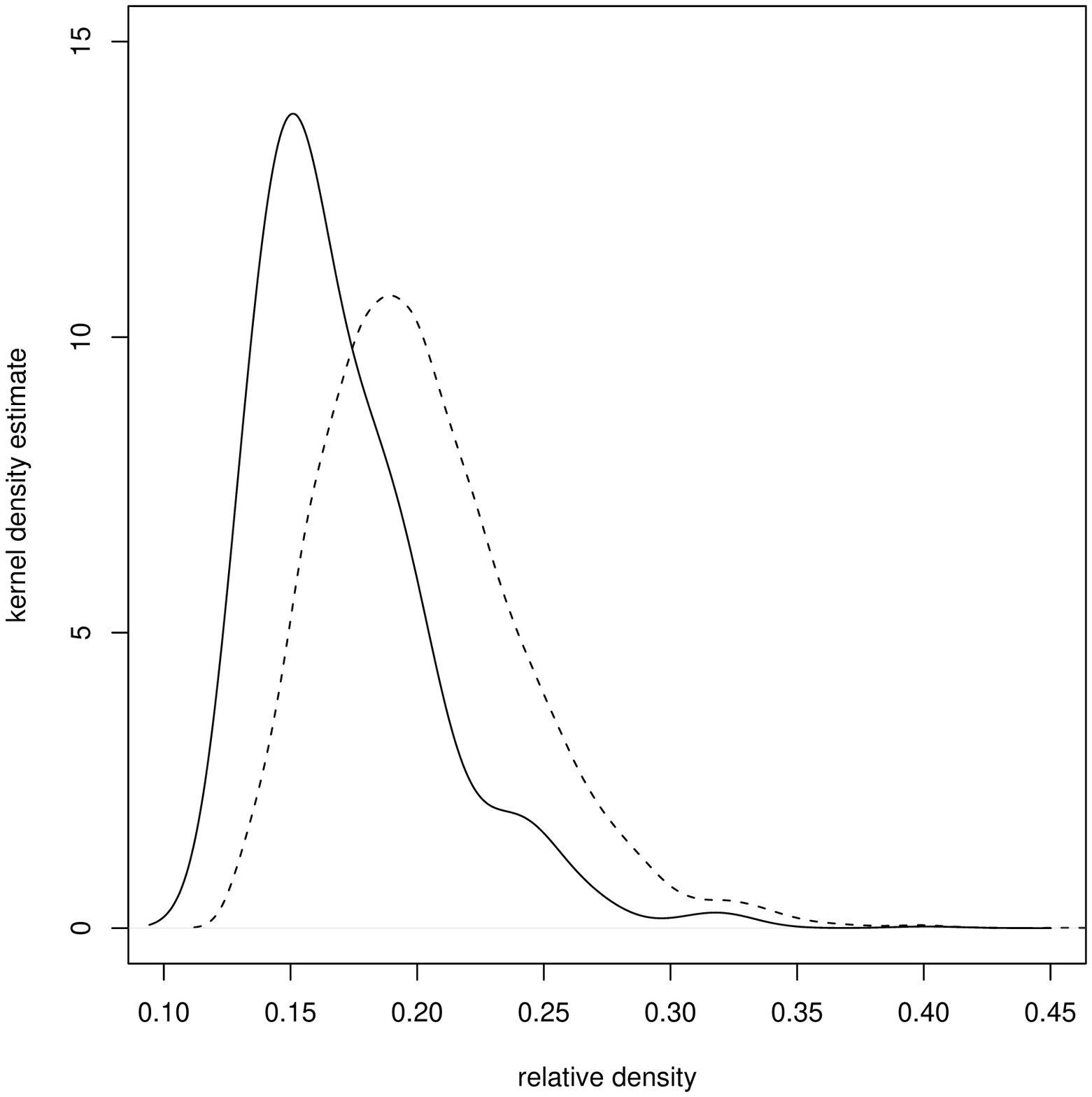}}}
\rotatebox{-90}{ \resizebox{1.7 in}{!}{ \includegraphics{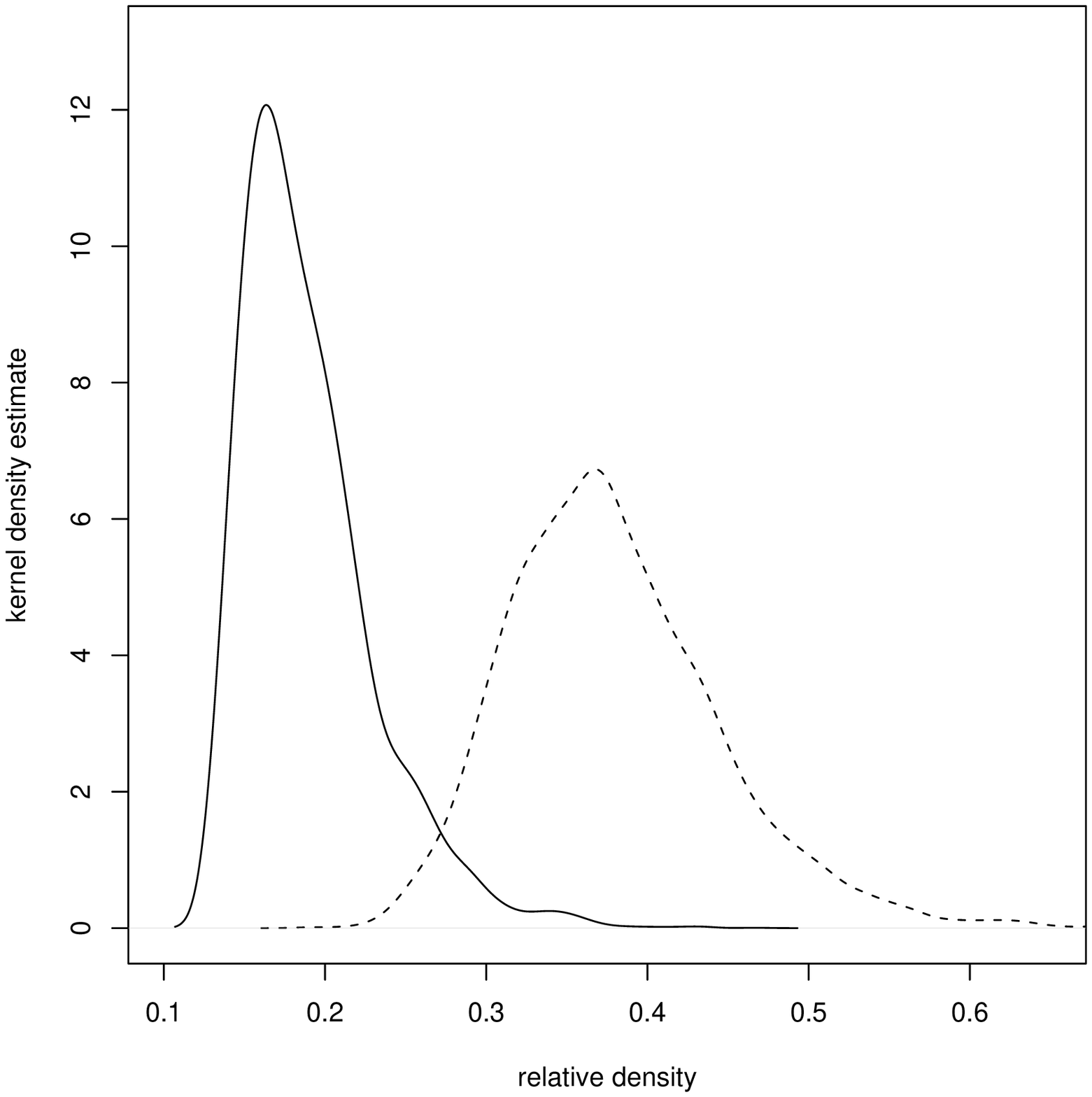}}}
\rotatebox{-90}{ \resizebox{1.7 in}{!}{ \includegraphics{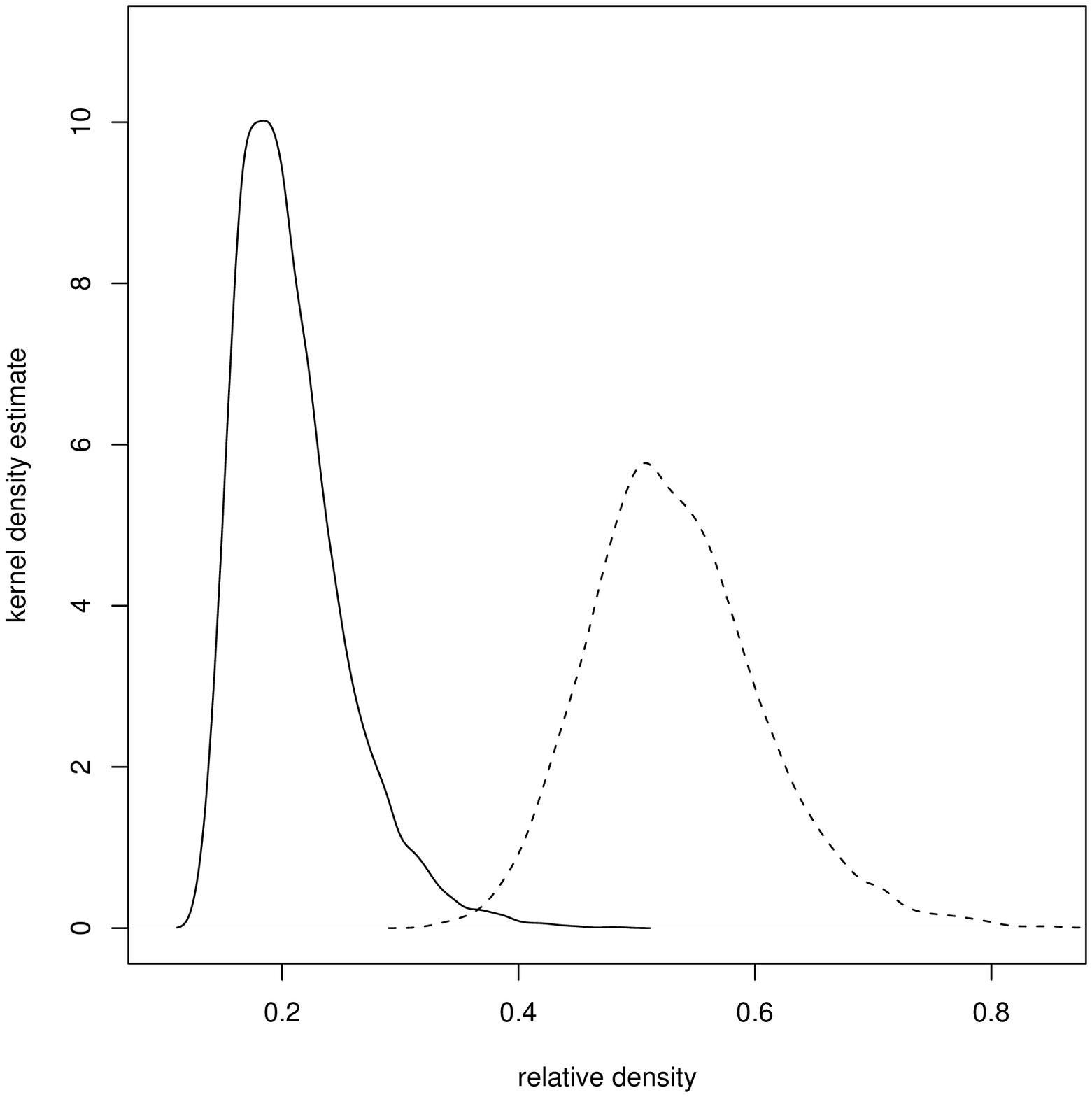}}}
\rotatebox{-90}{ \resizebox{1.7 in}{!}{ \includegraphics{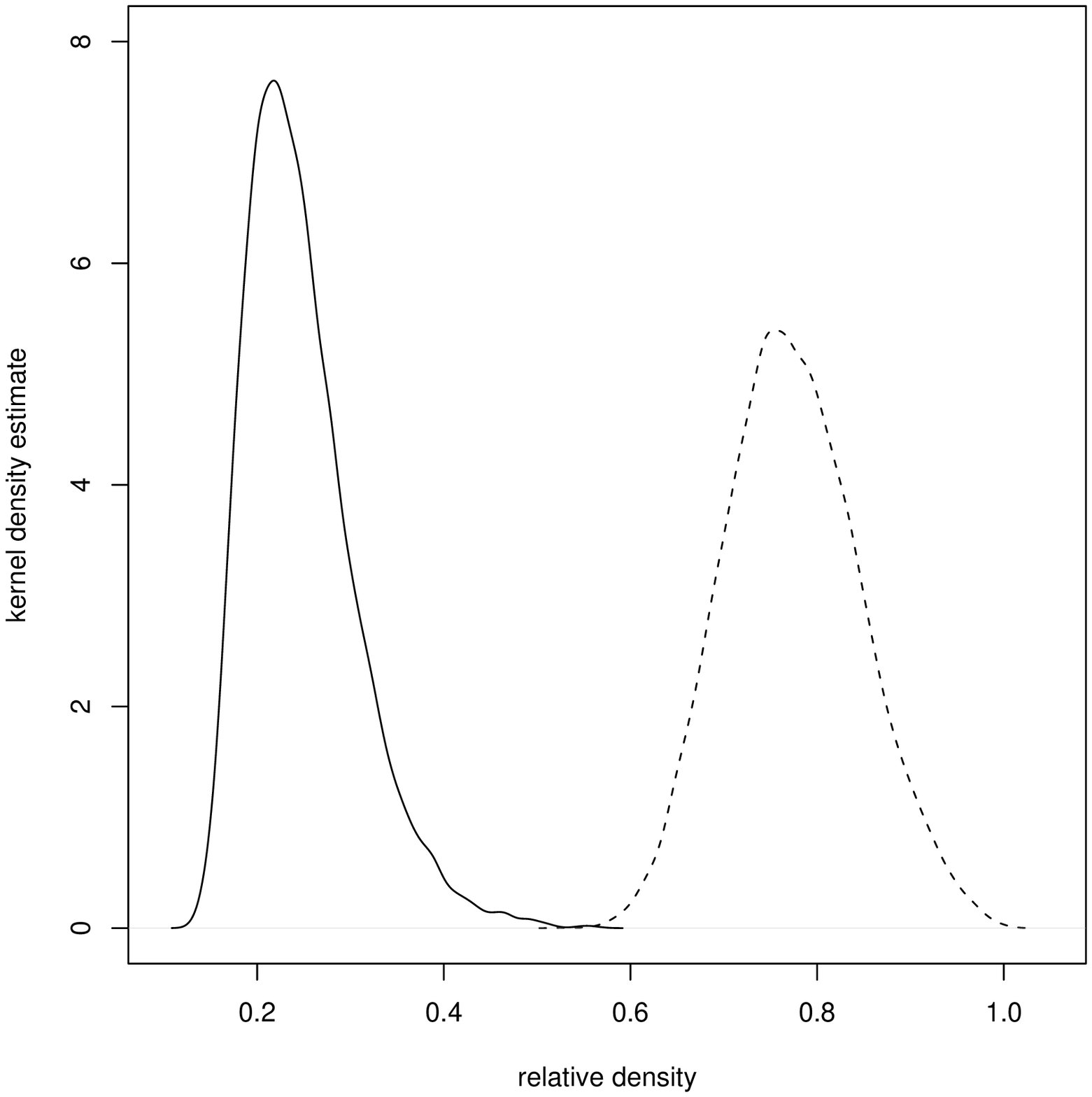}}}
\rotatebox{-90}{ \resizebox{1.7 in}{!}{ \includegraphics{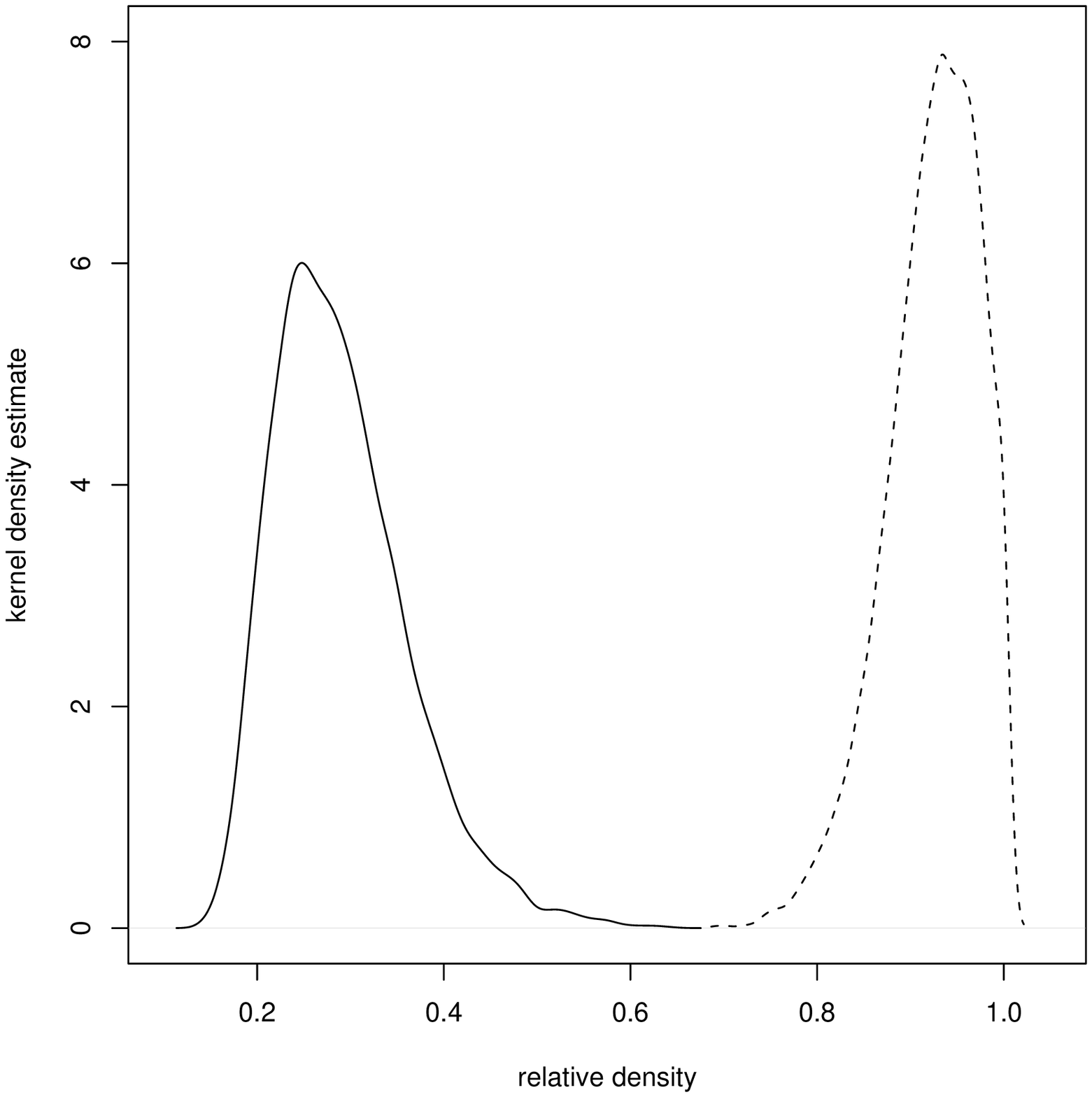}}}
\rotatebox{-90}{ \resizebox{1.7 in}{!}{ \includegraphics{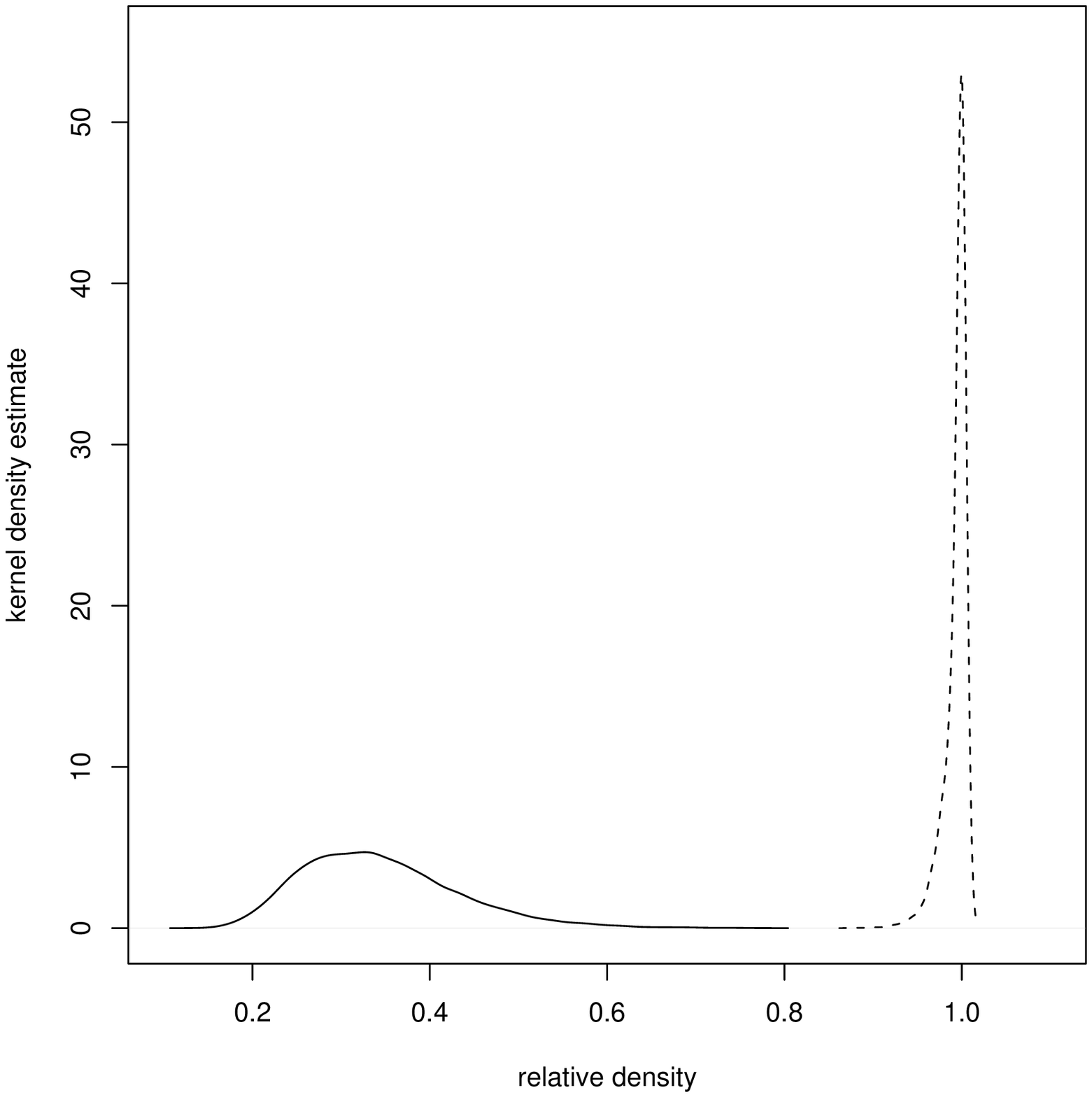}}}
\caption{\label{fig:seg3sim1-6}
Kernel density estimates for the null (solid) and the segregation alternative $H^S_{2\,\sqrt{3}/7}$ (dashed) for $r=1,\,21/20,\,11/10,\,6/5,\,4/3,\,\sqrt{2}$ (left-to-right).
}
\end{figure}

For segregation with $\epsilon=2\,\sqrt{3}/7 \approx .495$, we run the Monte Carlo experiments for six $r$ values, $1,\,21/20,\\
11/10,\,6/5,\,4/3,\,\sqrt{2}$. In Figure \ref{fig:seg3sim1-6}, are the kernel density estimates for the null case and the segregation alternative with $\epsilon=2\,\sqrt{3}/7$, for the six $r$ values with $n=10$ and $N=10,000$. Observe that under $H^S_{2\,\sqrt{3}/7}$, kernel density estimate is skewed right for $r=1$ and kernel density estimates are almost symmetric for $r=21/20,\,11/10,\,6/5$, with most symmetry occurring at $r=6/5$, kernel density estimate is skewed left for $r=\sqrt{2}$.

The empirical critical values, empirical significance levels, and empirical power estimates under $H^S_{2\,\sqrt{3}/7}$ are presented in Table \ref{tab:emp-val-S3}.

\begin{table}[]
\centering
\begin{tabular}{|c|c|c|c|c|c|c|}
\hline
$r$  & 1 & 21/20 & 11/10 &6/5 & 4/3 & $\sqrt{2}$ \\
\hline
$\widehat{C}^S_n$         & $.2\bar{4}$ & $.2\bar{8}$ & .3 & $.3\bar{5}$ & $.4\,\bar 2$ & .5 \\
\hline
$\widehat{\alpha}^S_{mc}(10)$      & .0318       & .0447   & .0411 & .0479 & .0477 &  .0481 \\
\hline
$\widehat{\beta}^S_{mc}\left( 10,\, 2\,\sqrt{3}/7 \right)$  & .1247        &  .9728 & 1.0 & 1.0  &  1.0 & 1.0 \\
\hline
\end{tabular}
\caption{
\label{tab:emp-val-S3}
The empirical critical values, empirical significance levels, and empirical power estimates under $H^S_{2\,\sqrt{3}/7}$, $N=10,000$, and $n=10$ at $\alpha=.05$.}
\end{table}

\begin{figure}[]
\centering
\psfrag{power}{ \Huge{\bfseries{power}}}
\psfrag{r}{\Huge{$r$}}
\rotatebox{-90}{ \resizebox{1.8 in}{!}{ \includegraphics{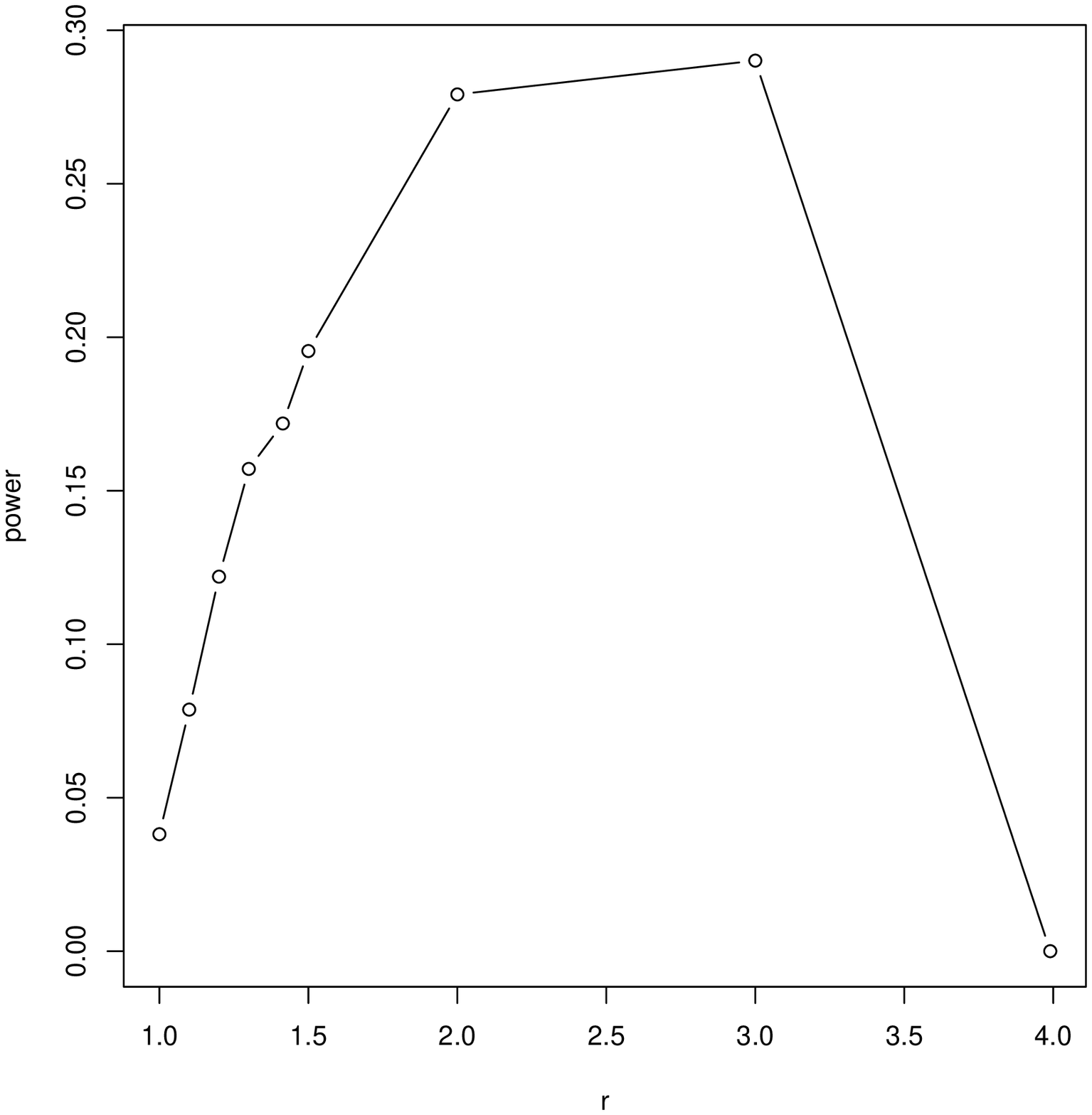}}}
\rotatebox{-90}{ \resizebox{1.8 in}{!}{ \includegraphics{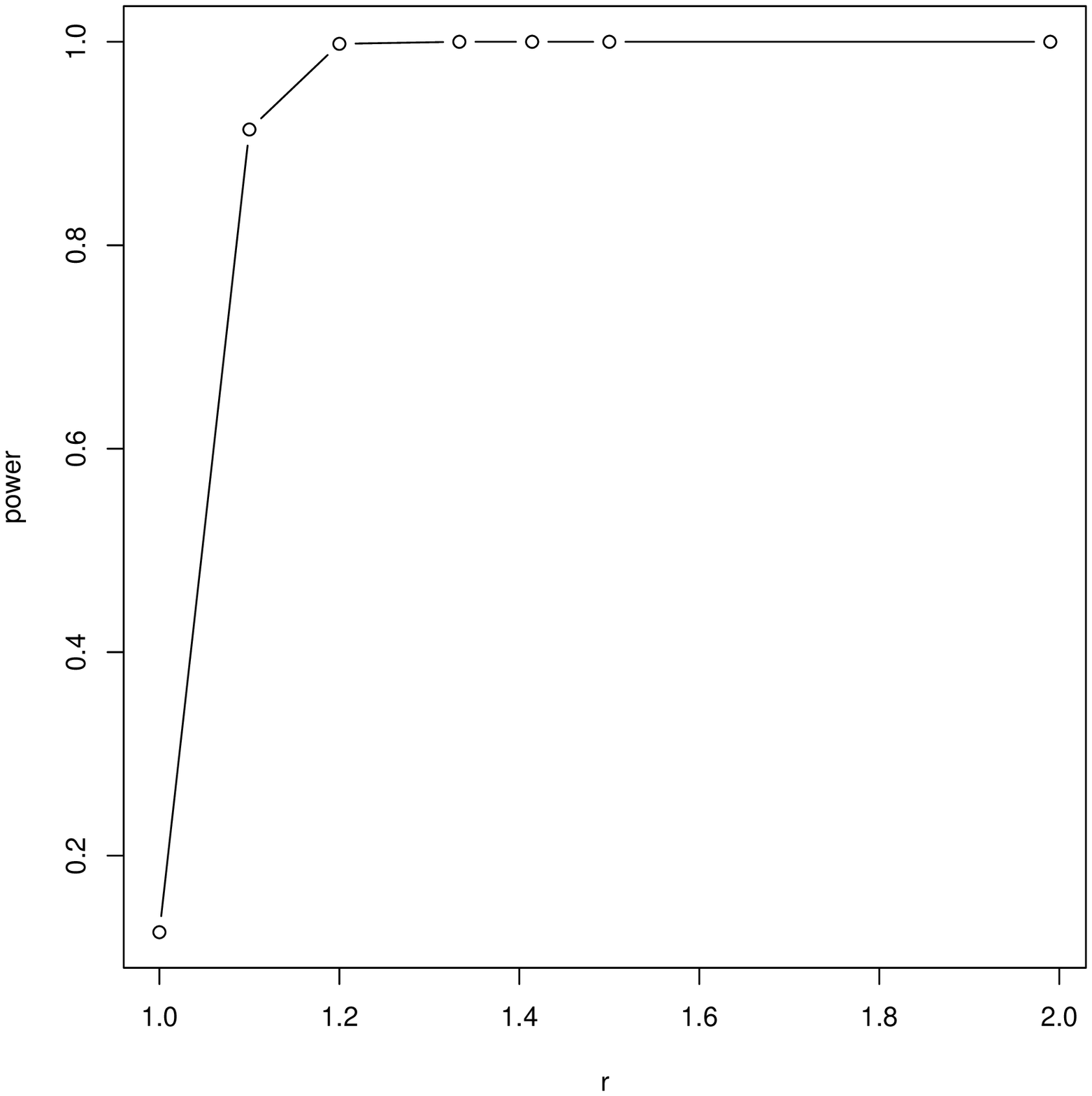}}}
\rotatebox{-90}{ \resizebox{1.8 in}{!}{ \includegraphics{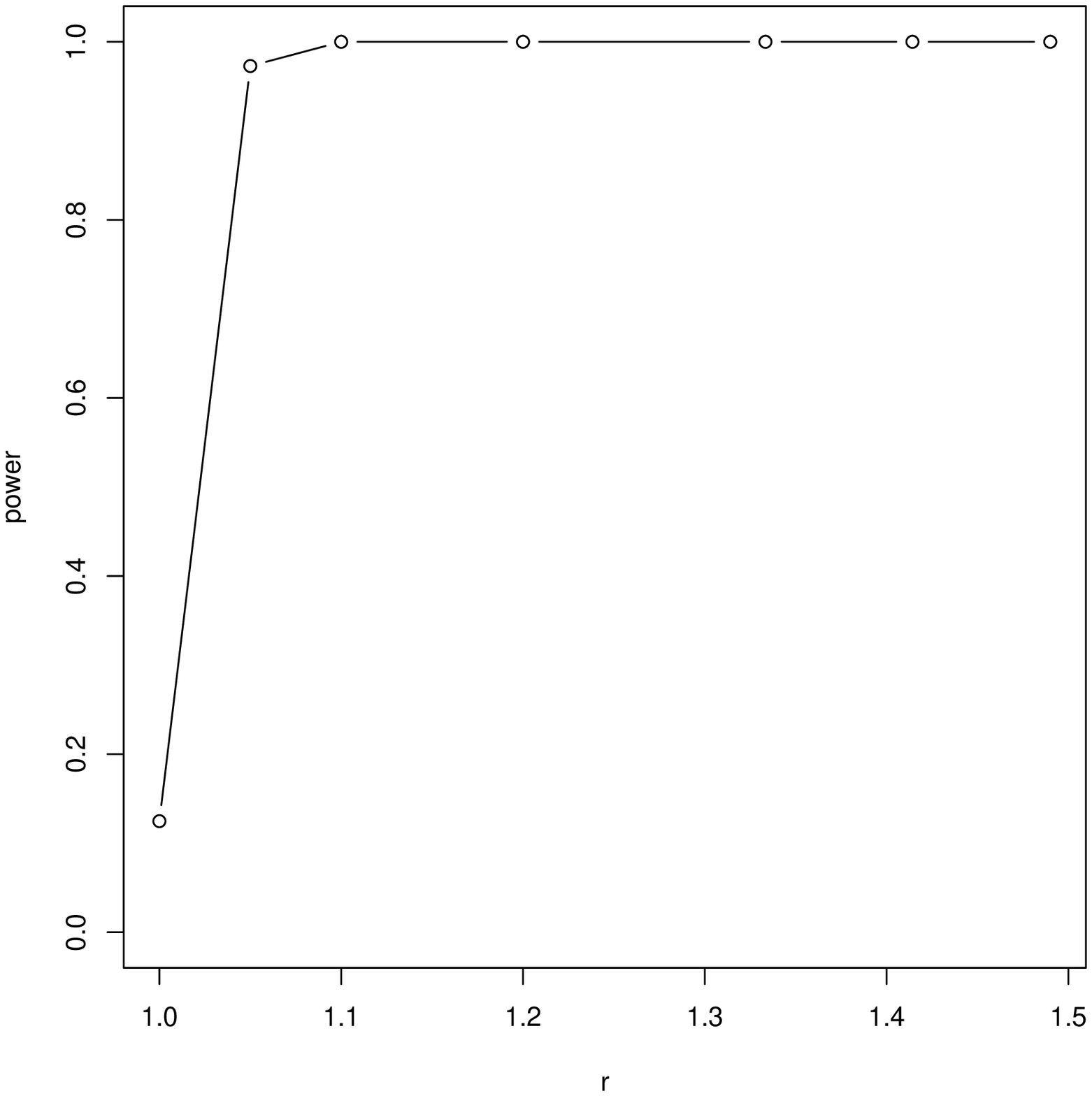}}}
\caption{ \label{fig:segpow1-3}
Monte Carlo power using the empirical critical value against segregation alternatives
$H^S_{\sqrt{3}/8}$ (left),
$H^S_{\sqrt{3}/4}$ (middle)
and
$H^S_{2\,\sqrt{3}/7}$ (right)
as a function of $r$, for $n=10$.}
\end{figure}

We also plot the empirical power as a function of $r$ in Figure \ref{fig:segpow1-3}. Let $r^*_S(\epsilon)$ be the value of $r$ at which maximum Monte Carlo power estimate occurs, then $r^*_S(\sqrt{3}/8)=3$.  Furthermore, Monte Carlo power estimate increases as $r$ gets larger and then decreases, due to the magnitude of $r$ and $n$. Because for small $n$ and large $r$, the critical value is approximately 1 under $H_0$, as we get a complete digraph with high probability.

Furthermore, $r^*_S\bigl( \sqrt{3}/4\bigr) \in \bigl\{ 4/3,\,\sqrt{2},\,3/2 \bigr\}$ and $r^*_S\bigl( 2\,\sqrt{3}/7 \bigr) \in  \bigl\{ 11/10,\,6/5,\,4/3,\,\sqrt{2} \bigr\}$.  Monte Carlo power estimates increase as $r$ gets larger. The phenomenon happened above for $\epsilon=\sqrt{3}/8$ does not occur, because $r$ values are not large enough to yield complete digraphs under $H_0$ with high probability.

\begin{figure}[]
\centering
\psfrag{power}{ \Huge{\bfseries{power}}}
\psfrag{r}{\Huge{$r$}}
\rotatebox{-90}{ \resizebox{1.8 in}{!}{ \includegraphics{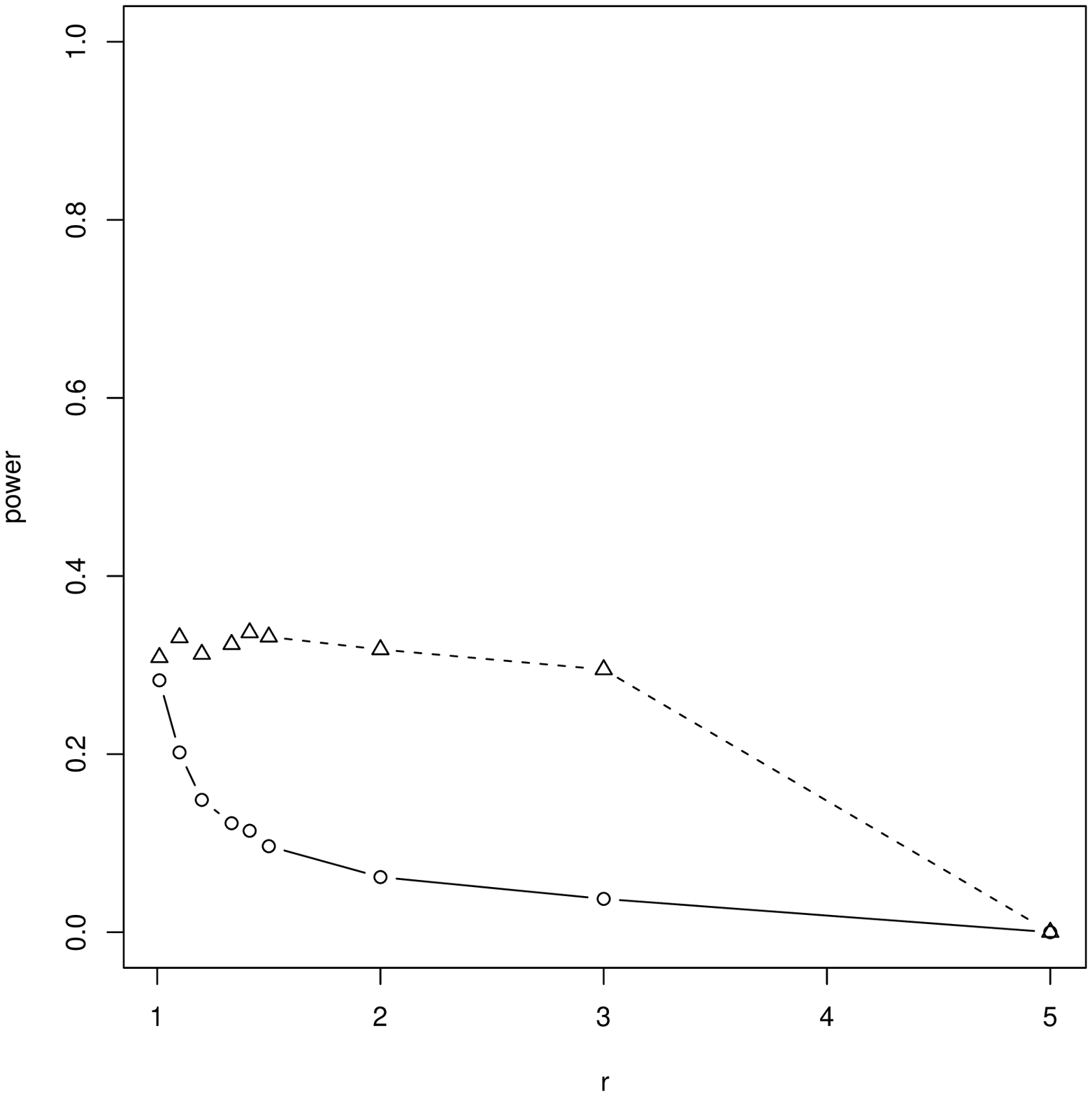}}}
\rotatebox{-90}{ \resizebox{1.8 in}{!}{ \includegraphics{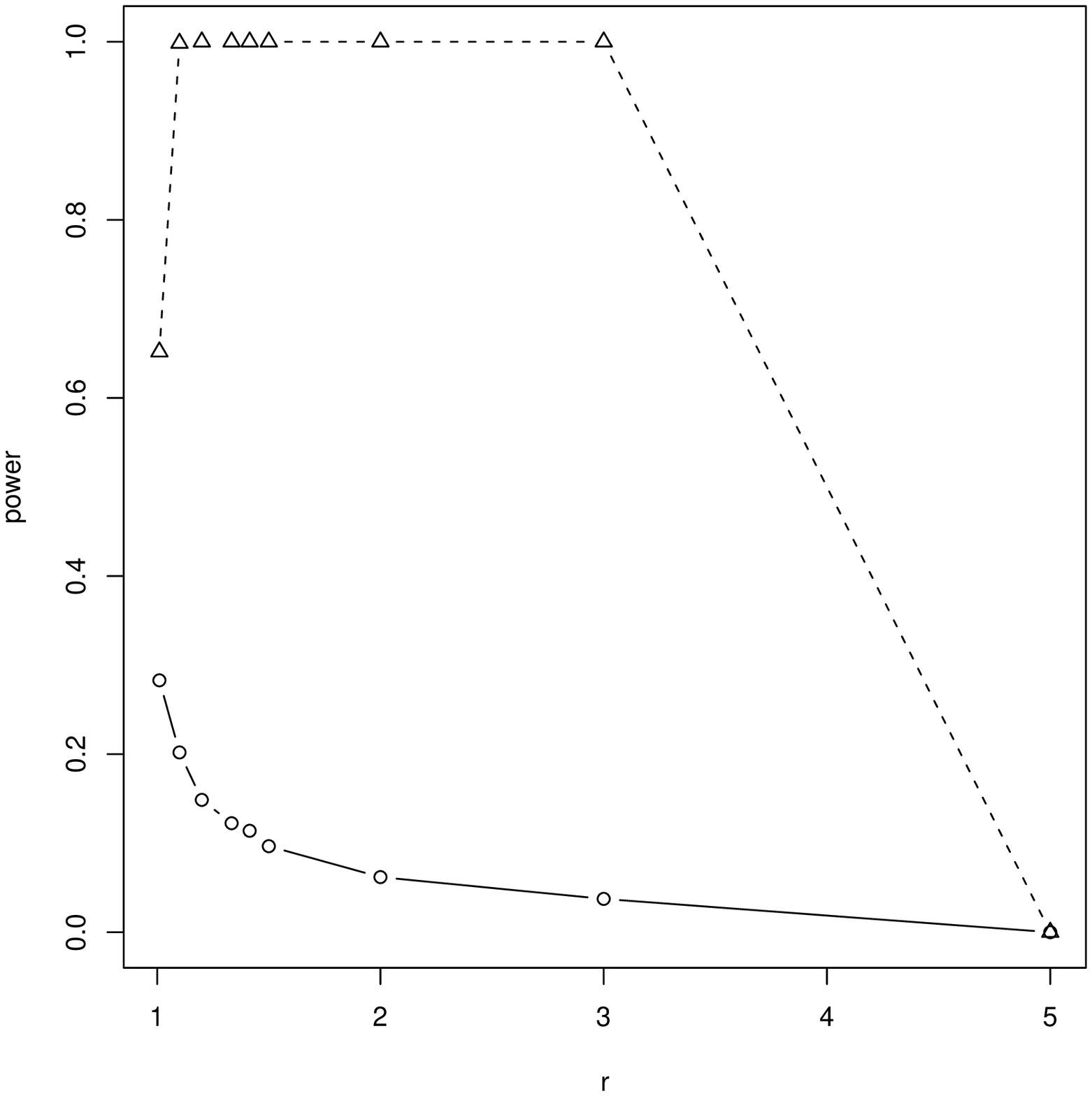}}}
\rotatebox{-90}{ \resizebox{1.8 in}{!}{ \includegraphics{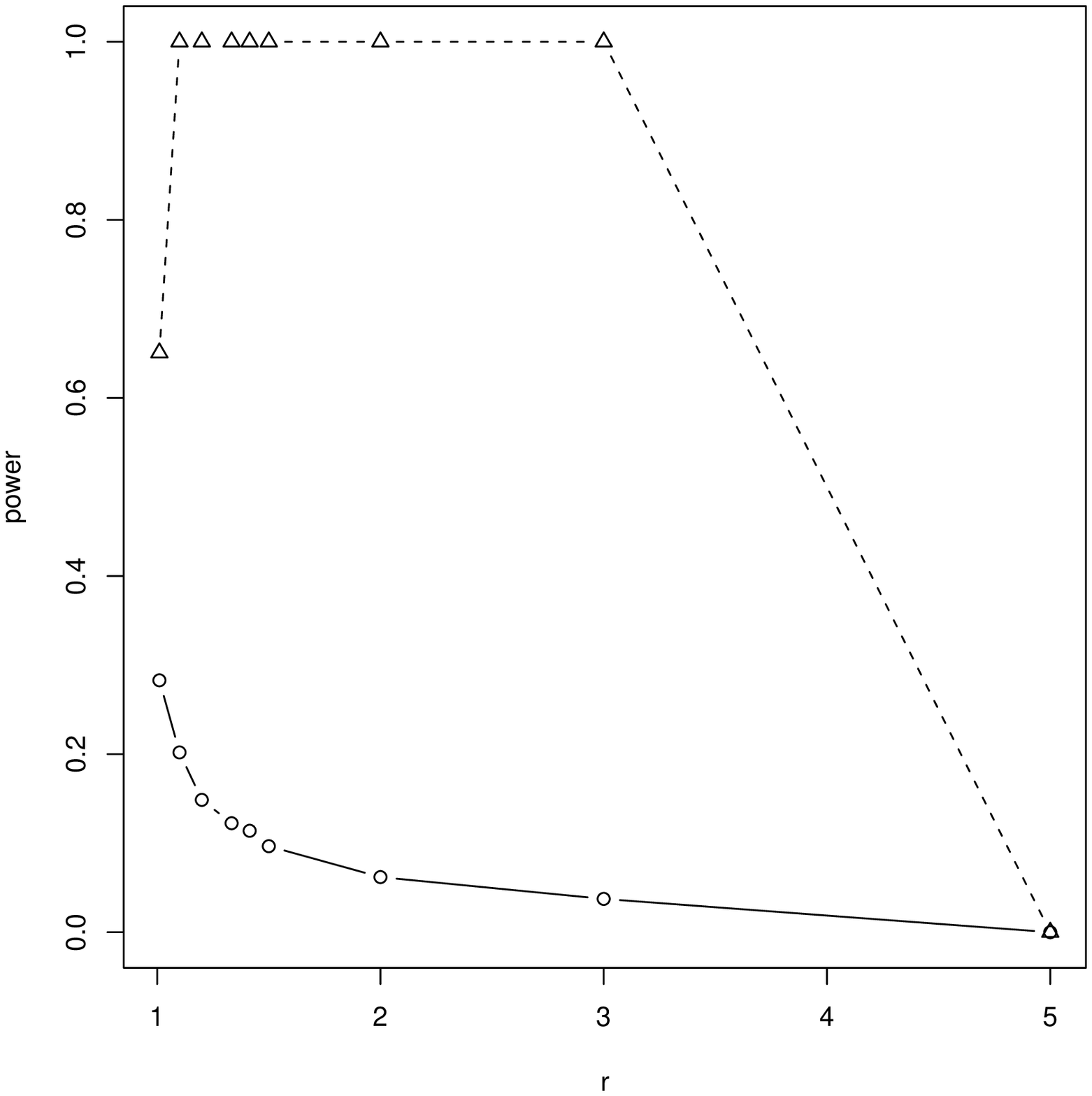}}}
\caption{ \label{fig:SegSimPowerCurve}
Monte Carlo power using the asymptotic critical value against segregation alternatives
$H^S_{\sqrt{3}/8}$ (left),
$H^S_{\sqrt{3}/4}$ (middle),
 and
$H^S_{2\,\sqrt{3}/7}$ (right)
as a function of $r$, for $n=10$.  The circles represent the empirical significance levels while triangles represent the empirical power values.
}
\end{figure}

For a given alternative and sample size,
we may consider analyzing the power of the test --- using the asymptotic critical value---
as a function of the proximity factor $r$.
Let $R_j:=\frac{\sqrt{n}\,\bigl( \rho_j(n)-\mu(r) \bigr)}{\sqrt{\nu(r)}}$ be the standardized relative density for experiment $j$ with sample size $n$ for $j=1,2,\ldots,N$. For each $r$ value, the level $\alpha$ asymptotic critical value is  $\mu(r)+z_{(1-\alpha)} \cdot \sqrt{\nu(r)/n}$.  We estimate the empirical power as $\widehat{\beta}^S_n(r,\epsilon):=\frac{1}{N}\sum_{j=1}^{N}\I \left( R_j > z_{1-\alpha} \right)$.  In Figure \ref{fig:SegSimPowerCurve}, we present
a Monte Carlo investigation of power $\widehat{\beta}^S_n(r,\epsilon)$ against
$H^S_{\sqrt{3}/8}$, $H^S_{\sqrt{3}/4}$, and $H^S_{2\,\sqrt{3}/7}$
as a function of $r$ for $n=10$.  The empirical significance level is $\widehat{\alpha}^S_{mc}(n):=\frac{1}{N}\sum_{j=1}^{N}\I\left( R_j > z_{1-\alpha}|H_0 \right)$. Then $\widehat{\alpha}_S(10)$, is about $.05$ for $r=2,\,3$ which have the empirical power $\widehat{\beta}^S_{10}\left( r,\sqrt{3}/8 \right) \approx .35$, and $\widehat{\beta}^S_{10}(r,\epsilon) =1$ for $\epsilon=\sqrt{3}/4,2\,\sqrt{3}/7$.  So, for small sample sizes, moderate values of $r$ are more appropriate for normal approximation, as they yield the desired significance level and the more severe the segregation, the higher the power estimate at each $r$.

The empirical significance level, and empirical power $\widehat{\beta}^S_n(r,\epsilon)$ values under $H^S_{\epsilon}$ for $\epsilon=\sqrt{3}/8,\,\sqrt{3}/4,\,2\,\sqrt{3}/7$ are presented in Table \ref{tab:asy-emp-val-S1}.
\begin{table}[]
\centering
\begin{tabular}{|c|c|c|c|c|c|c|c|c|c|}
\hline
$r$  & 1 & 11/10 &6/5 & 4/3 & $\sqrt{2}$ &3/2 & 2 & 3 & 5\\
\hline
$\widehat{\alpha}_S(n)$ & .2829 & .2019 & .1486 & .1224 & .1139 & .0966 & .0619 & .0374 & .000 \\
\hline
$\widehat{\beta}^S_{n}\left( r,\sqrt{3}/8 \right)$ & .3086 & .3309 & .3123 & .3233 & .3365 & .3317 & .3175 & .2950 & .0000 \\
\hline
$\widehat{\beta}^S_{n}\left( r,\sqrt{3}/4 \right)$ & .6519 & .9985 & 1.0000 & 1.0000 & 1.0000 & 1.0000 & 1.0000 & 1.0000 & .0000 \\
\hline
$\widehat{\beta}^S_{n}\left( r,2\,\sqrt{3}/7 \right)$ & .6508 & 1.0000 & 1.0000 & 1.0000 & 1.0000 & 1.0000 & 1.0000 & 1.0000 & .0000 \\
\hline
\end{tabular}
\caption{
\label{tab:asy-emp-val-S1}
The empirical significance level and empirical power values under $H^S_{\epsilon}$ for $\epsilon=\sqrt{3}/8,\,\sqrt{3}/4,\,2\,\sqrt{3}/7$, $N=10,000$, and $n=10$ at $\alpha=.05$.}
\end{table}
Note that even for $n=10$, the plots of the empirical power $\widehat{\beta}^S_n(r,\epsilon)$ resemble the curves of the asymptotic power function $\Pi_S(r)$ in Section \ref{sec:APF}.

\subsection{Monte Carlo Power Analysis Under Association}

In association alternatives with $\epsilon>0$, we implement the Monte Carlo experiment for $r \in \bigl\{ 1,\,11/10,\,6/5,\,4/3,\\
\sqrt{2},\,3/2,\,2,\,3,\,5,\,10 \bigr\}$.    Then for each $r$ value, we estimate the empirical critical value $\widehat{C}^A_n:=\rho_{(\lfloor \alpha\,N \rfloor)}$ and the empirical significance level $\widehat{\alpha}^A_{mc}(n):=\frac{1}{N}\sum_{j=1}^{N}\I\left( \rho_{j}(n) < \widehat{C}^A_n \right)$ under $H_0$ and the empirical power $\widehat{\beta}^A_{mc}(n,\epsilon):=\frac{1}{N}\sum_{j=1}^{N}\I\left( \rho_{j} < \widehat{C}^A_n \right)$.  We implement the Monte Carlo simulation for three $\epsilon$ values; $5\,\sqrt{3}/24,\\
\sqrt{3}/12,\,\sqrt{3}/21$.

\begin{figure}[]
\centering
\psfrag{kernel density estimate}{ \Huge{\bfseries{kernel density estimate}}}
\psfrag{relative density}{ \Huge{\bfseries{relative density}}}
\rotatebox{-90}{ \resizebox{1.7 in}{!}{ \includegraphics{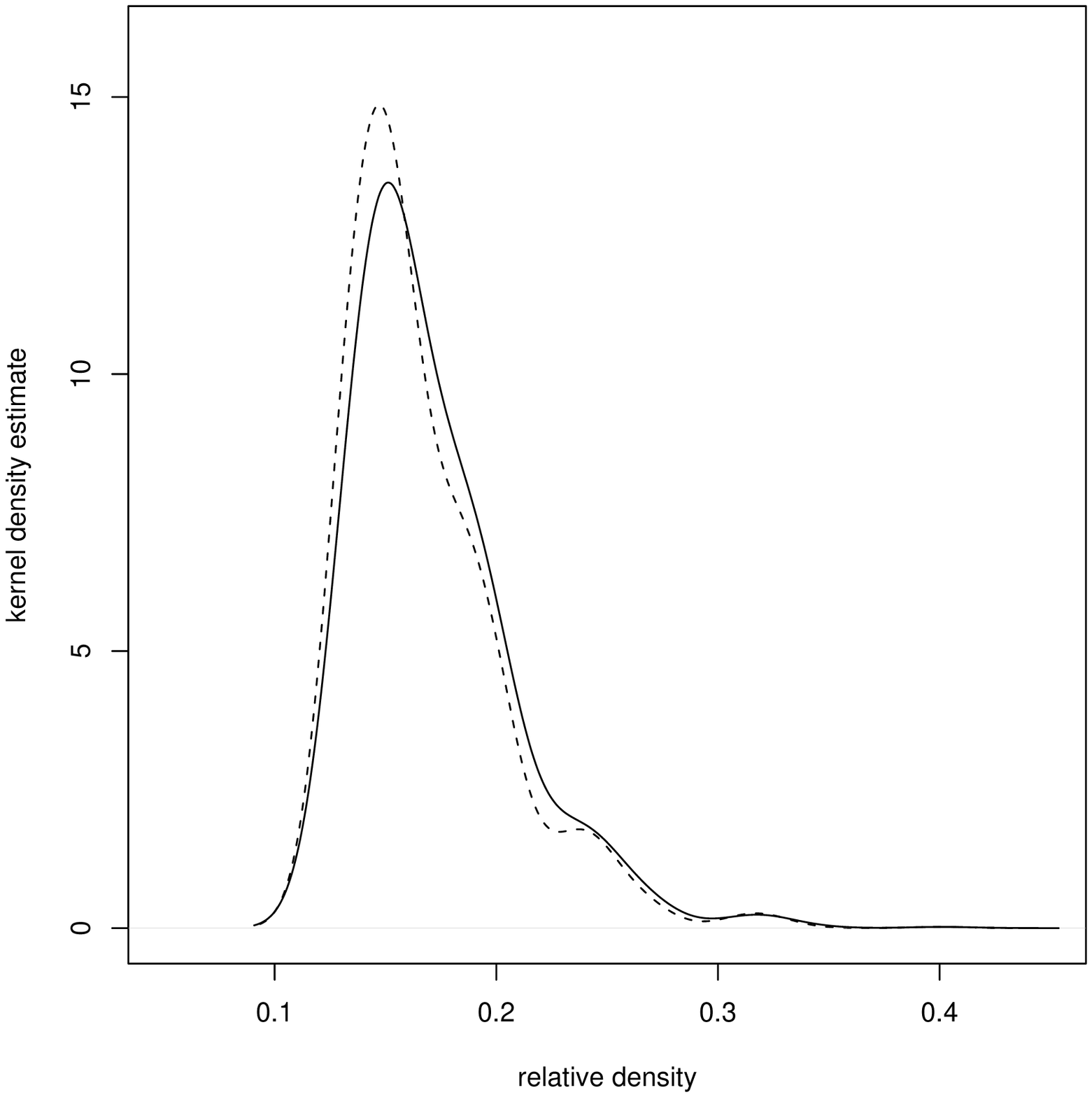}}}
\rotatebox{-90}{ \resizebox{1.7 in}{!}{ \includegraphics{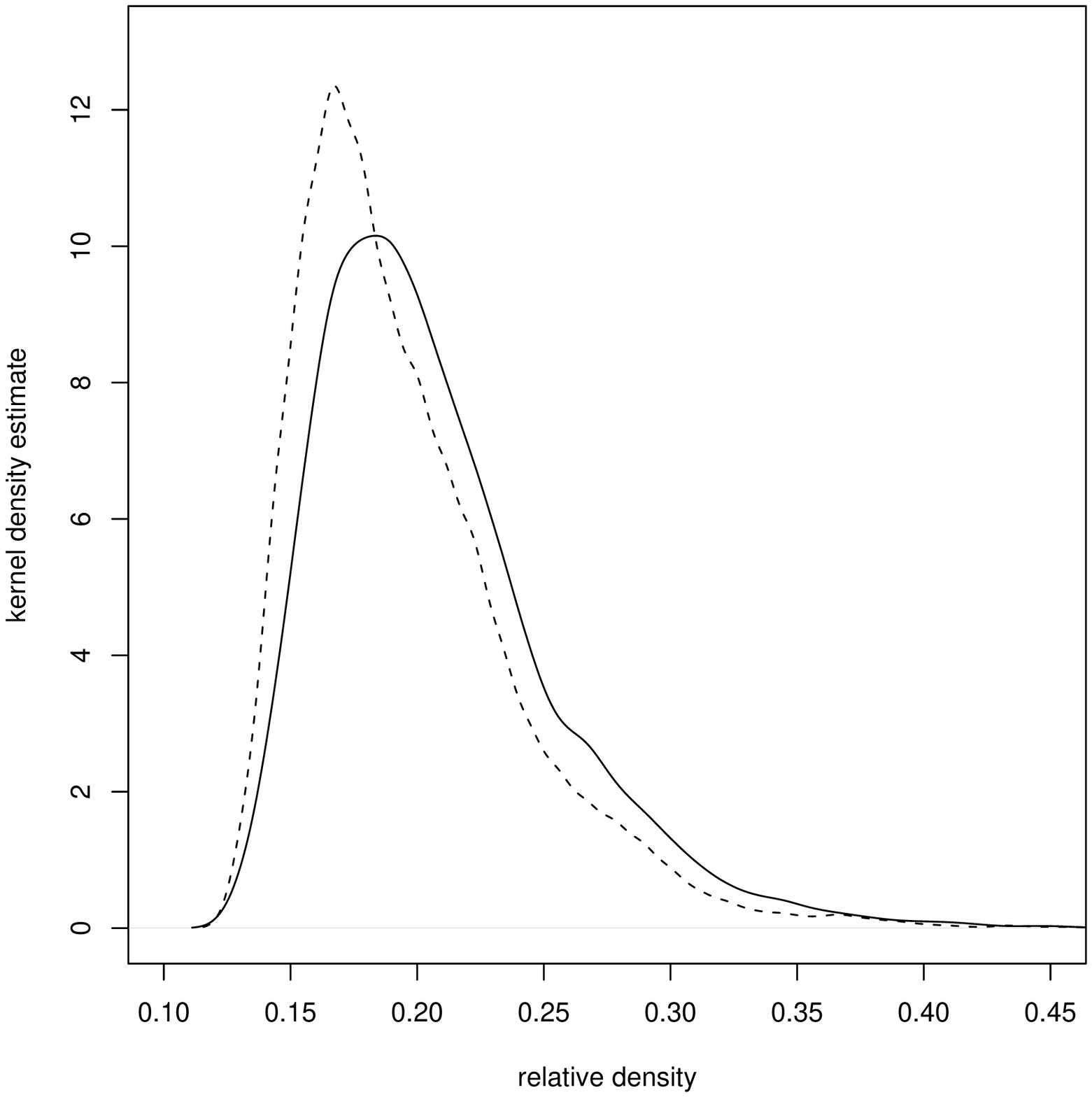}}}
\rotatebox{-90}{ \resizebox{1.7 in}{!}{ \includegraphics{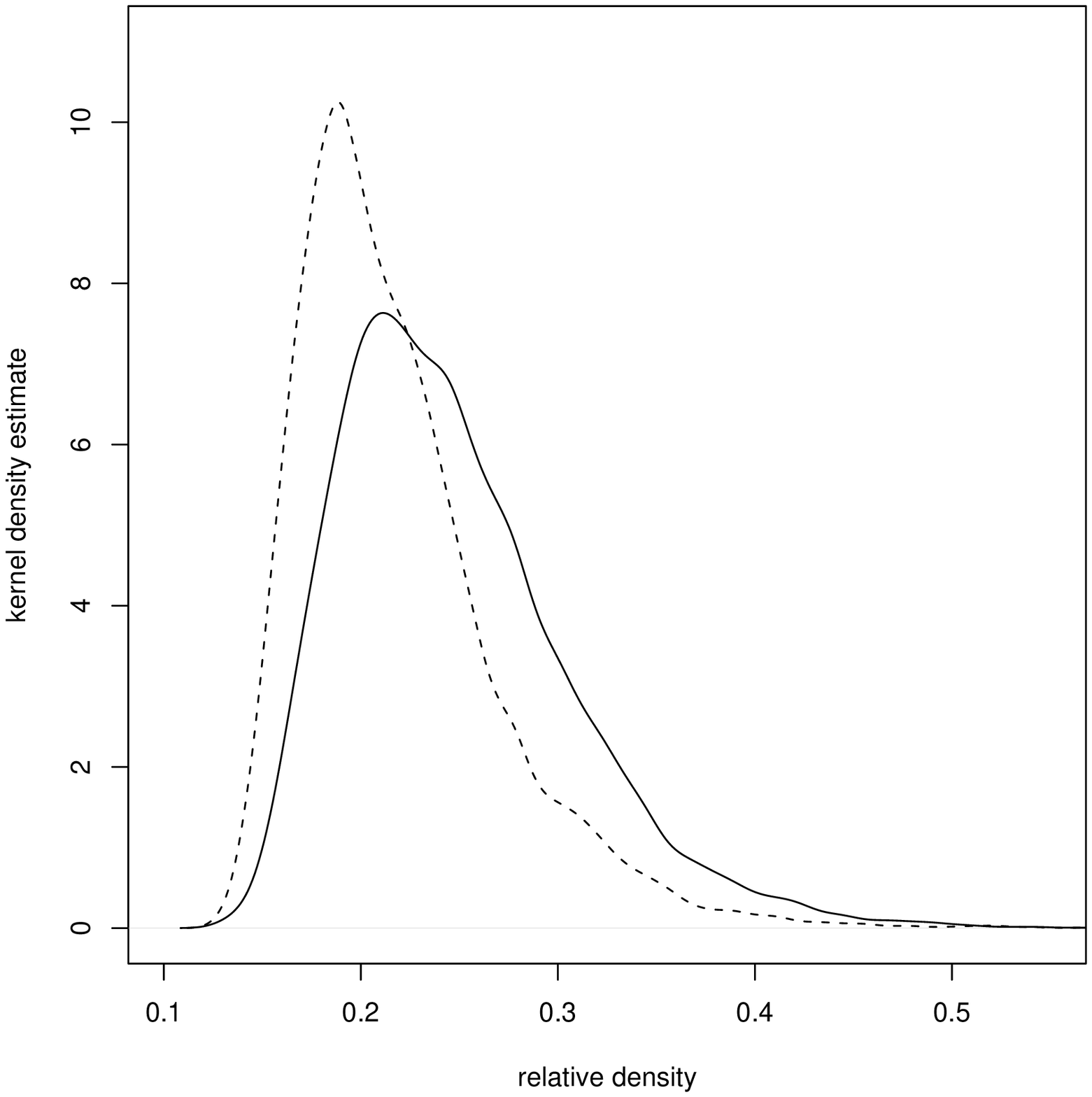}}}
\rotatebox{-90}{ \resizebox{1.7 in}{!}{ \includegraphics{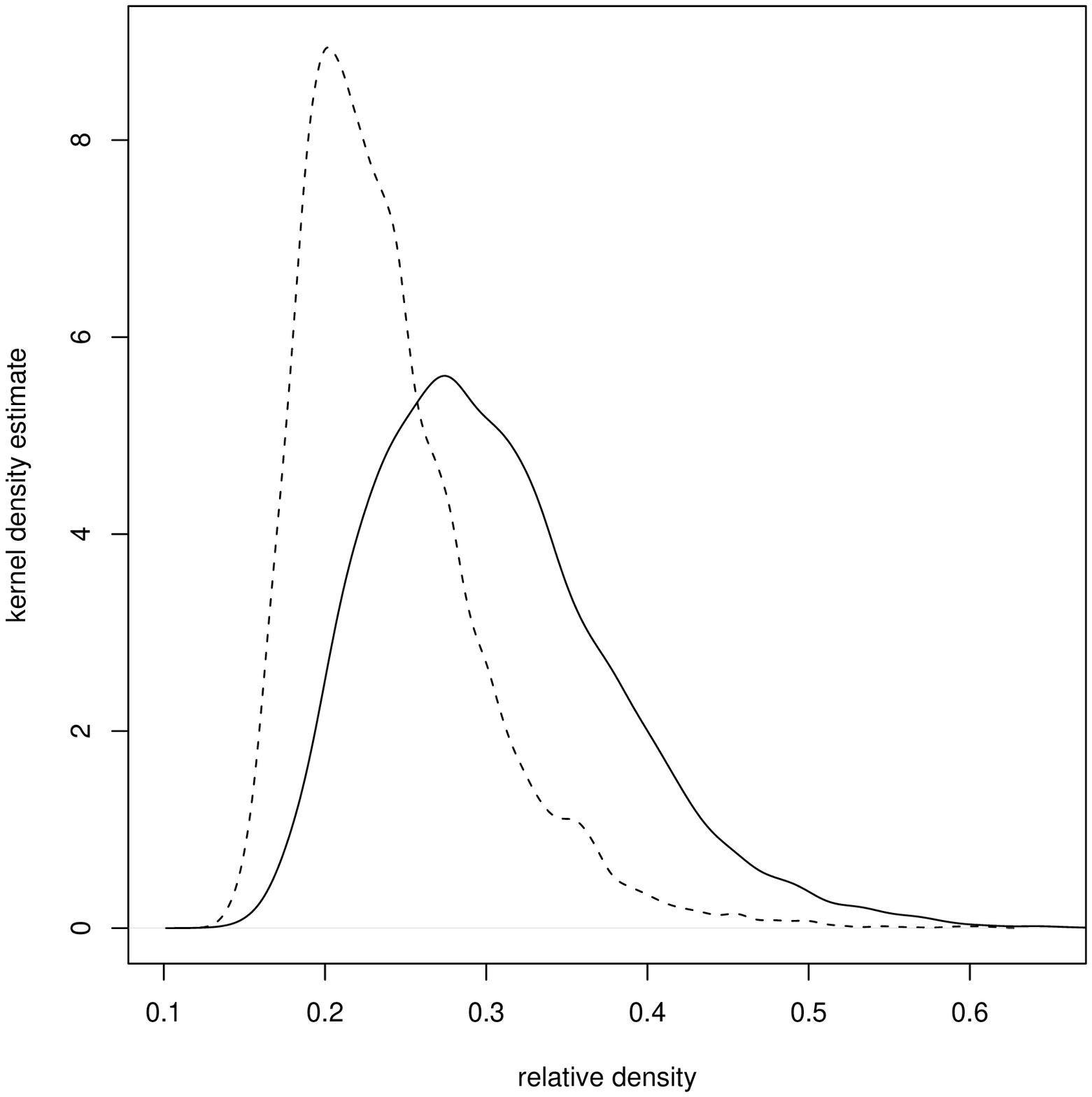}}}
\rotatebox{-90}{ \resizebox{1.7 in}{!}{ \includegraphics{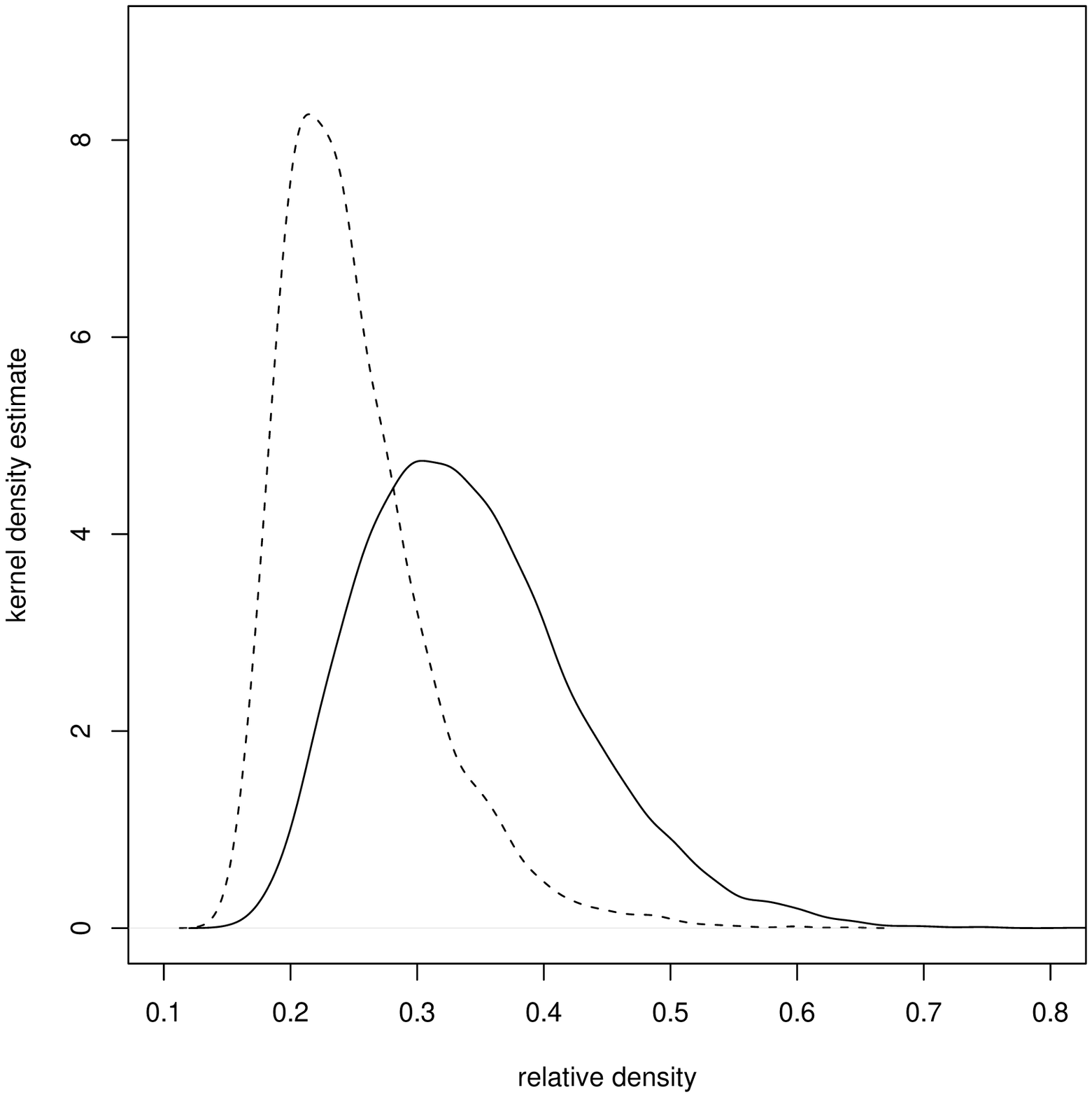}}}
\rotatebox{-90}{ \resizebox{1.7 in}{!}{ \includegraphics{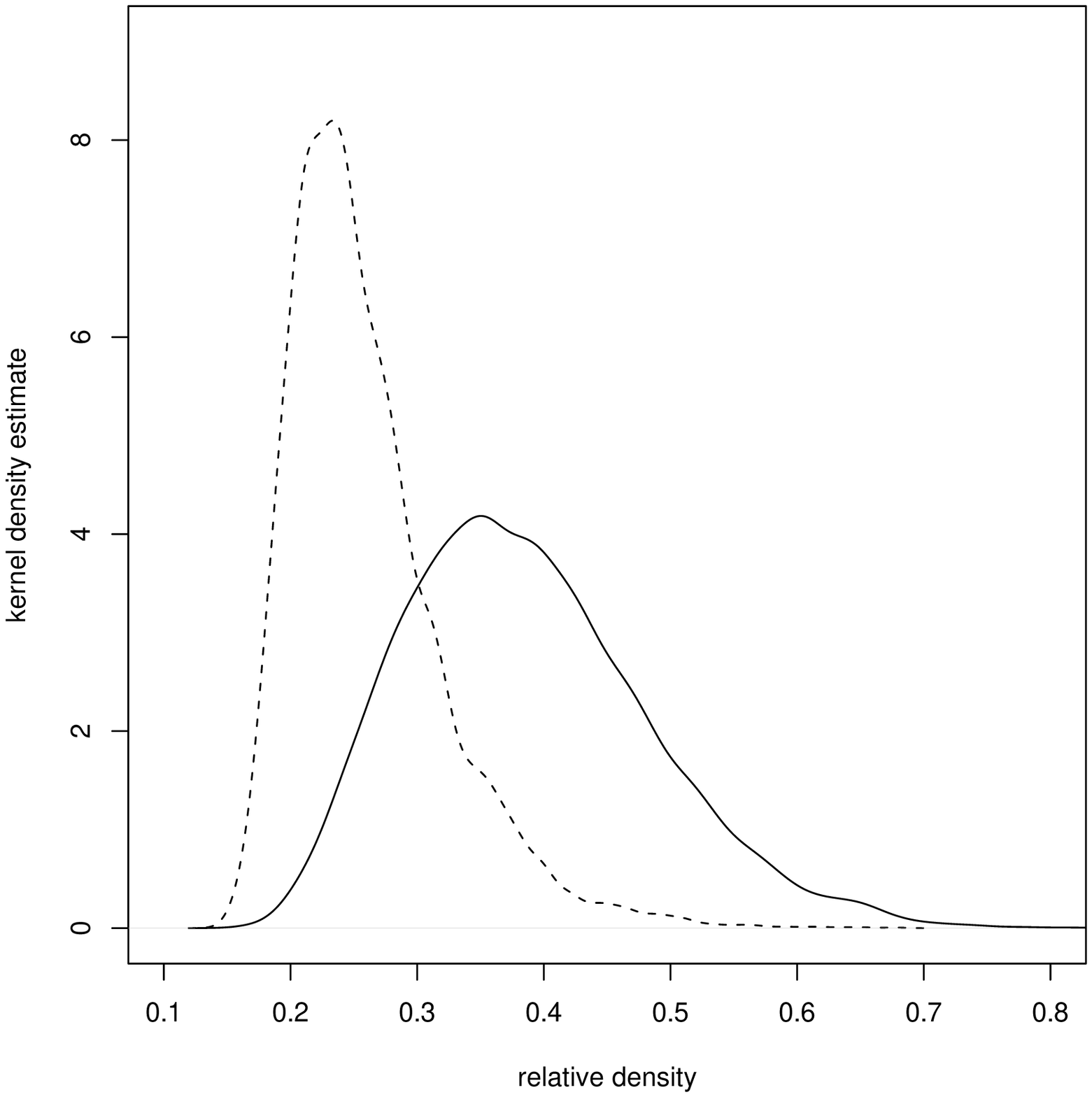}}}
\rotatebox{-90}{ \resizebox{1.7 in}{!}{ \includegraphics{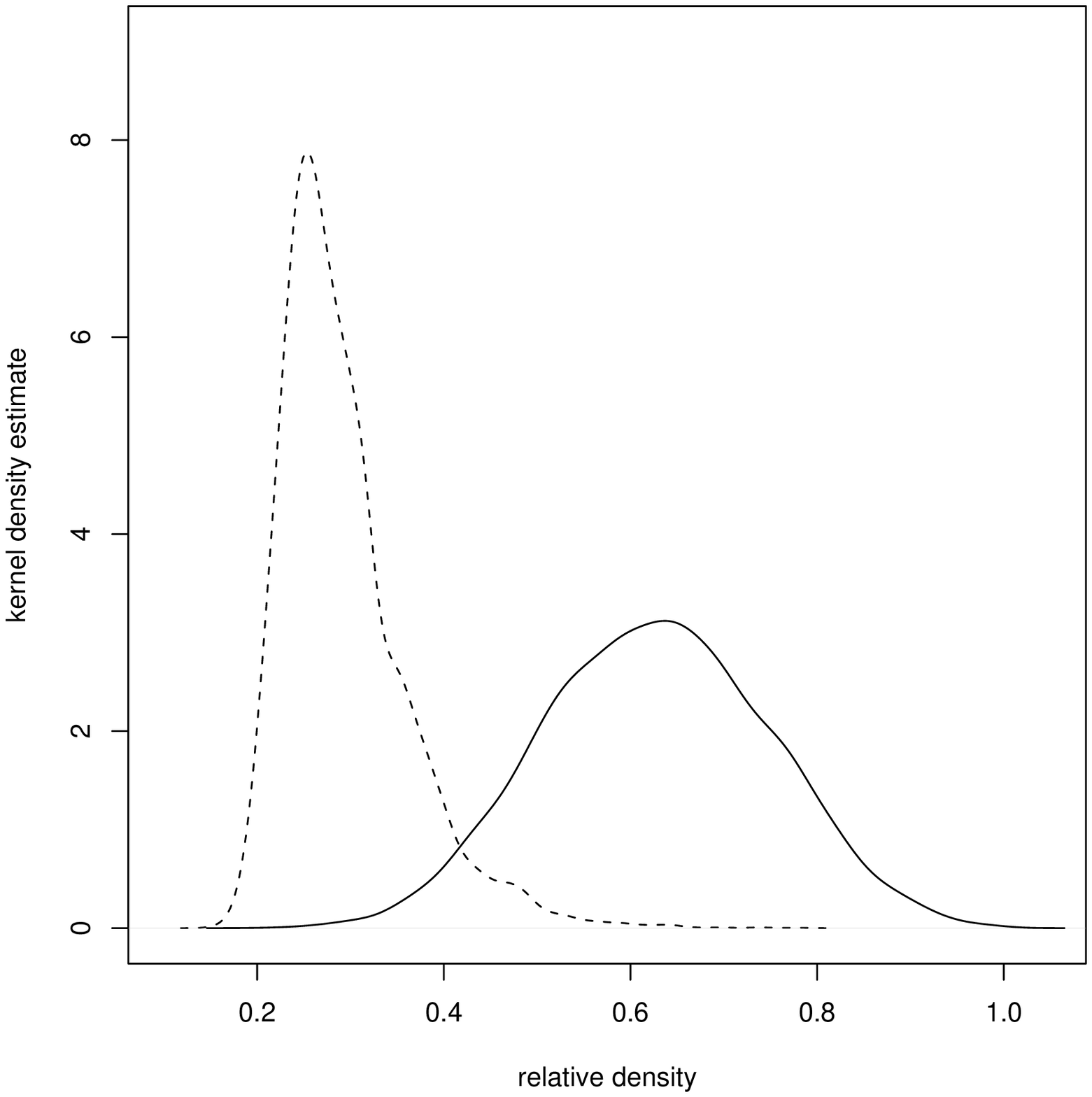}}}
\rotatebox{-90}{ \resizebox{1.7 in}{!}{ \includegraphics{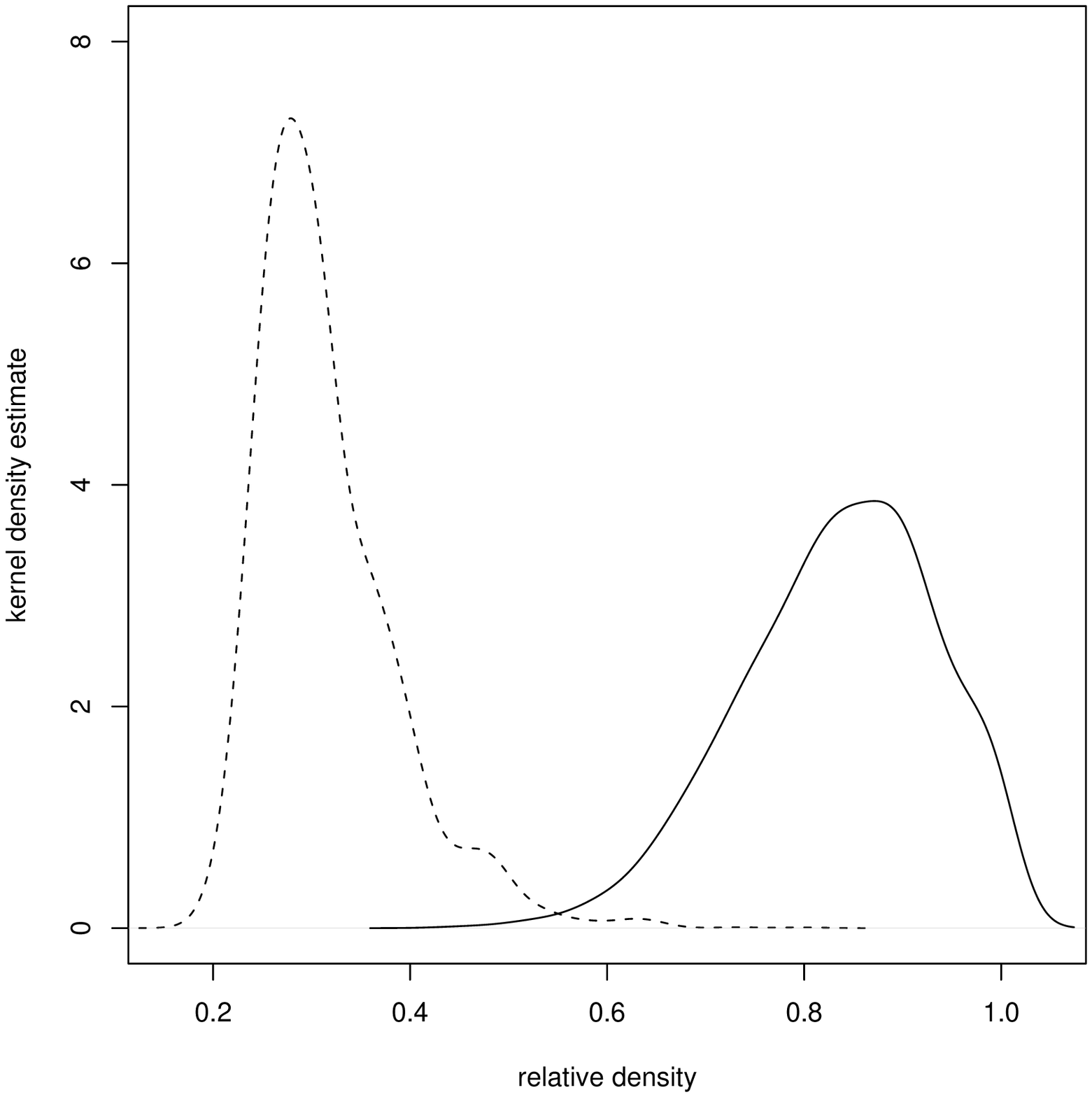}}}
\rotatebox{-90}{ \resizebox{1.7 in}{!}{ \includegraphics{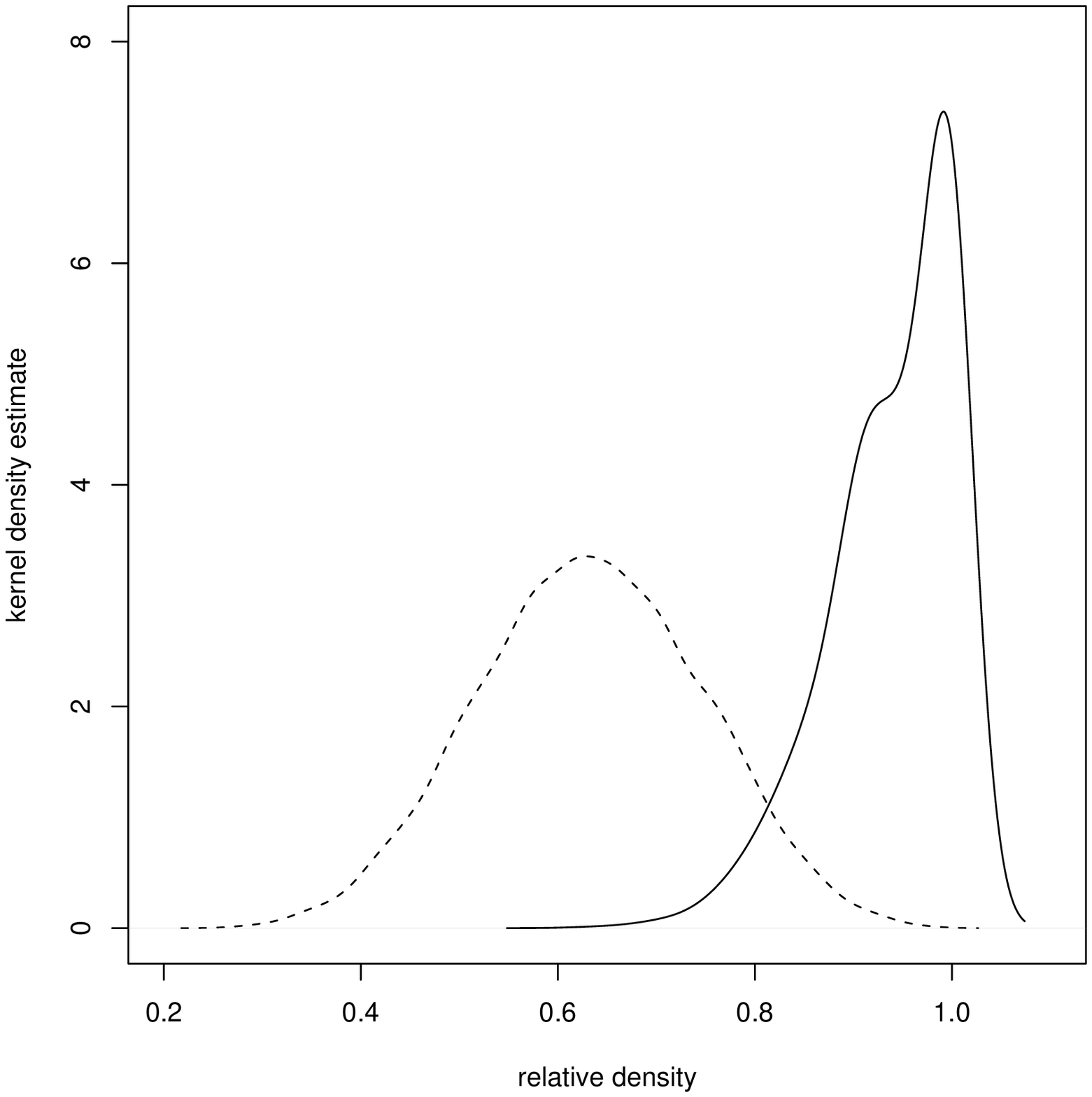}}}
\rotatebox{-90}{ \resizebox{1.7 in}{!}{ \includegraphics{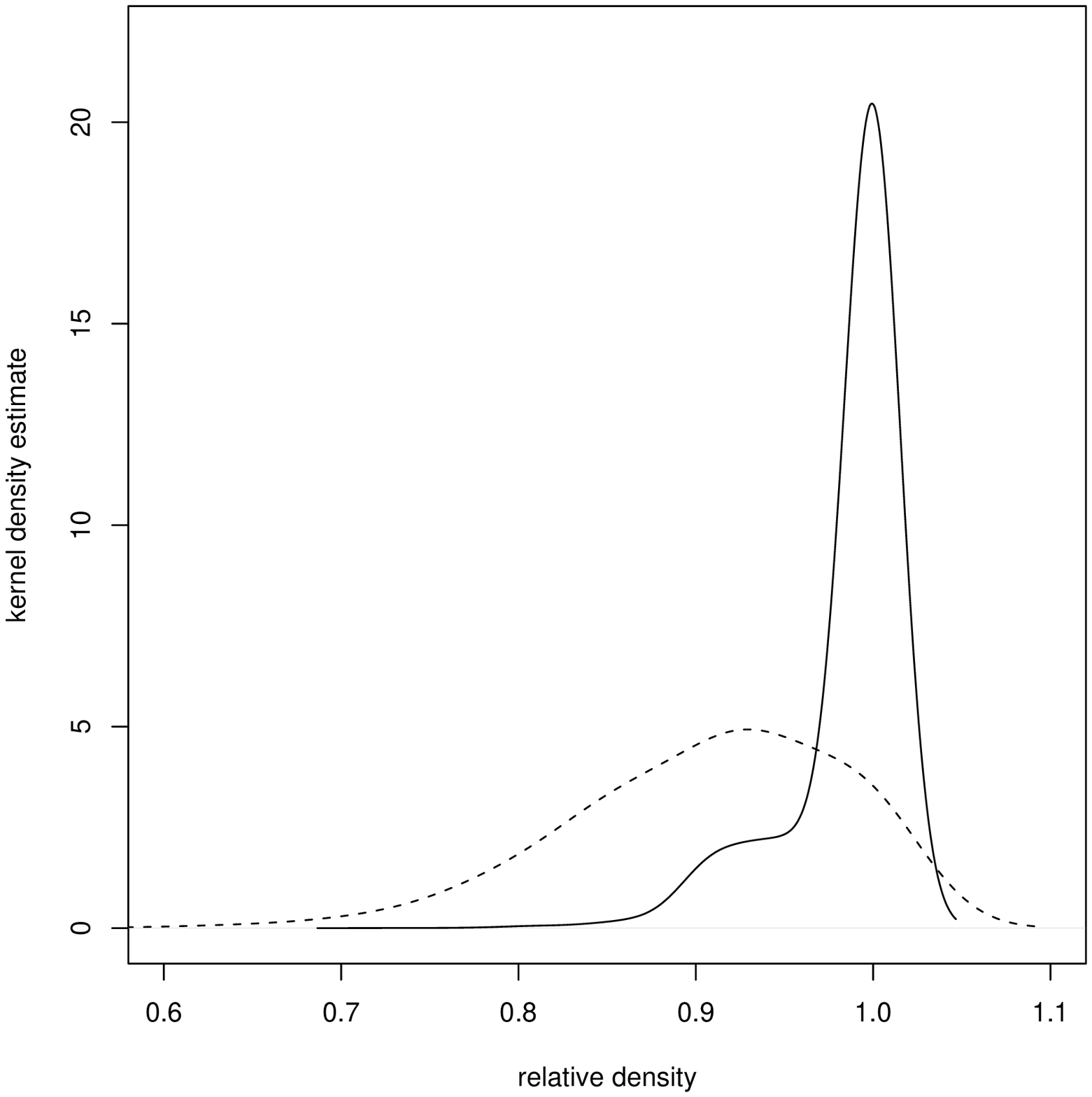}}}
\caption{\label{fig:aggsim1-10}
Kernel density estimates for the null (solid) and the association alternative $H^A_{5\,\sqrt{3}/24}$ (dashed) for $r= 1,\, 11/10,\, 6/5,\, 4/3,\, \sqrt{2},\, 3/2,\, 2,\, 3,\, 5,\, 10$ (left-to-right).
}
\end{figure}

The empirical critical values, empirical significance levels, and empirical power estimates under $H^A_{\epsilon}$ are presented in Table \ref{tab:emp-val-A}.

For association with $\epsilon=5\,\sqrt{3}/24\,\approx.36$, in Figure \ref{fig:aggsim1-10}, are the kernel density estimates for the null case and the segregation alternative for the ten $r$ values with $n=10,\,N=10,000$. Observe that, under $H_0$, kernel density estimates are skewed right for $r=1,\,11/10$, (with skewness increasing as $r$ gets smaller) and  kernel density estimates are almost symmetric for $r=6/5,\,4/3,\,\sqrt{2}, 3/2,\, 2 $, with most symmetry occurring at $r=3/2$,  kernel density estimates are skewed left for $r=3,\,5,\,10$, (with skewness increasing as $r$ gets larger). Under $H^A_{5\,\sqrt{3}/24}$, kernel density estimates are skewed right for $r=1,\,11/10,\,6/5,\,4/3,\,3/2,\,2,\,3$, (with skewness increasing as $r$ gets smaller) and kernel density estimate is almost symmetric for $r=5$, kernel density estimate is skewed left for $r=10$.

\begin{figure}[]
\centering
\psfrag{kernel density estimate}{ \Huge{\bfseries{kernel density estimate}}}
\psfrag{relative density}{ \Huge{\bfseries{relative density}}}
\rotatebox{-90}{ \resizebox{1.7 in}{!}{ \includegraphics{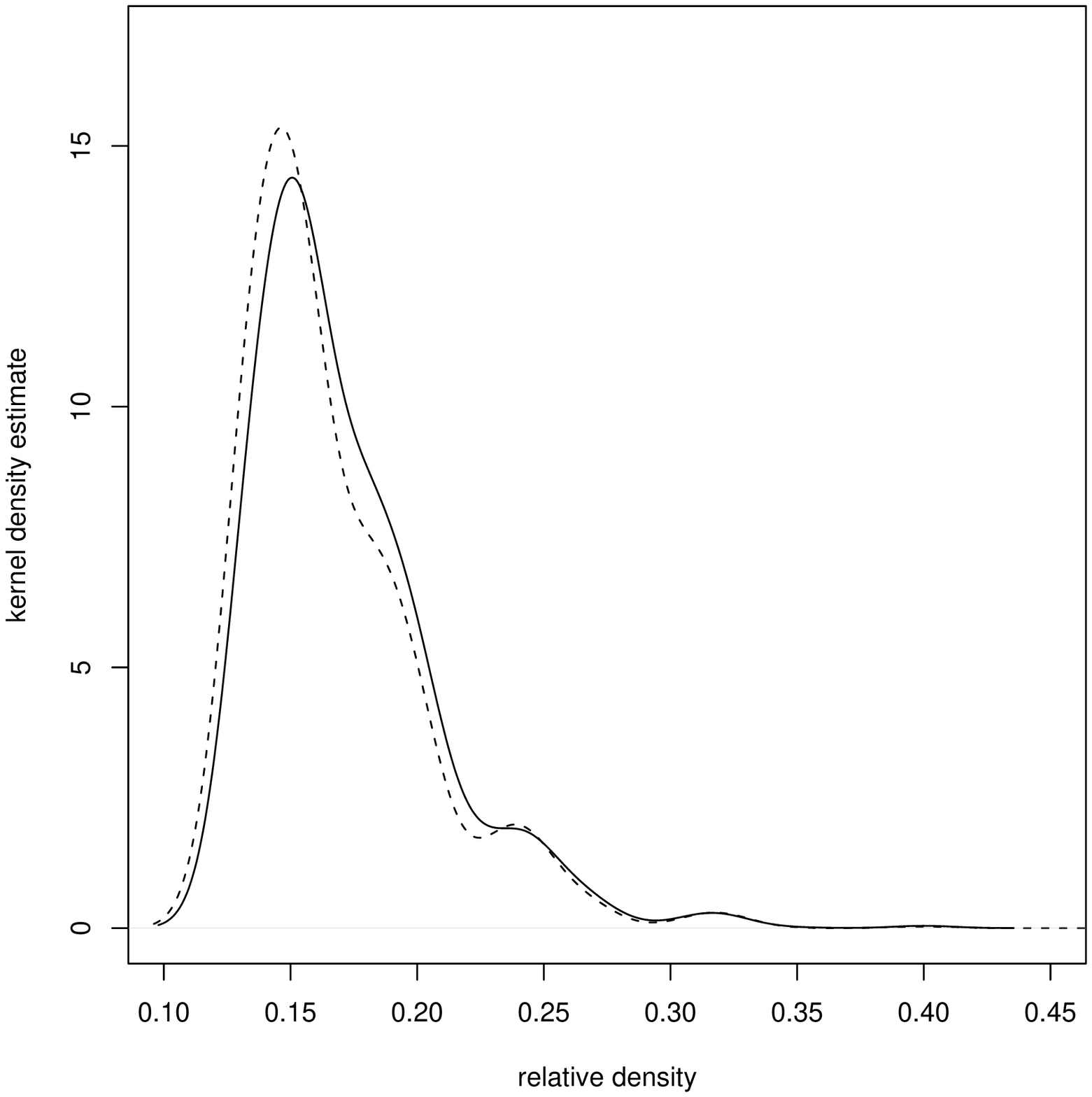}}}
\rotatebox{-90}{ \resizebox{1.7 in}{!}{ \includegraphics{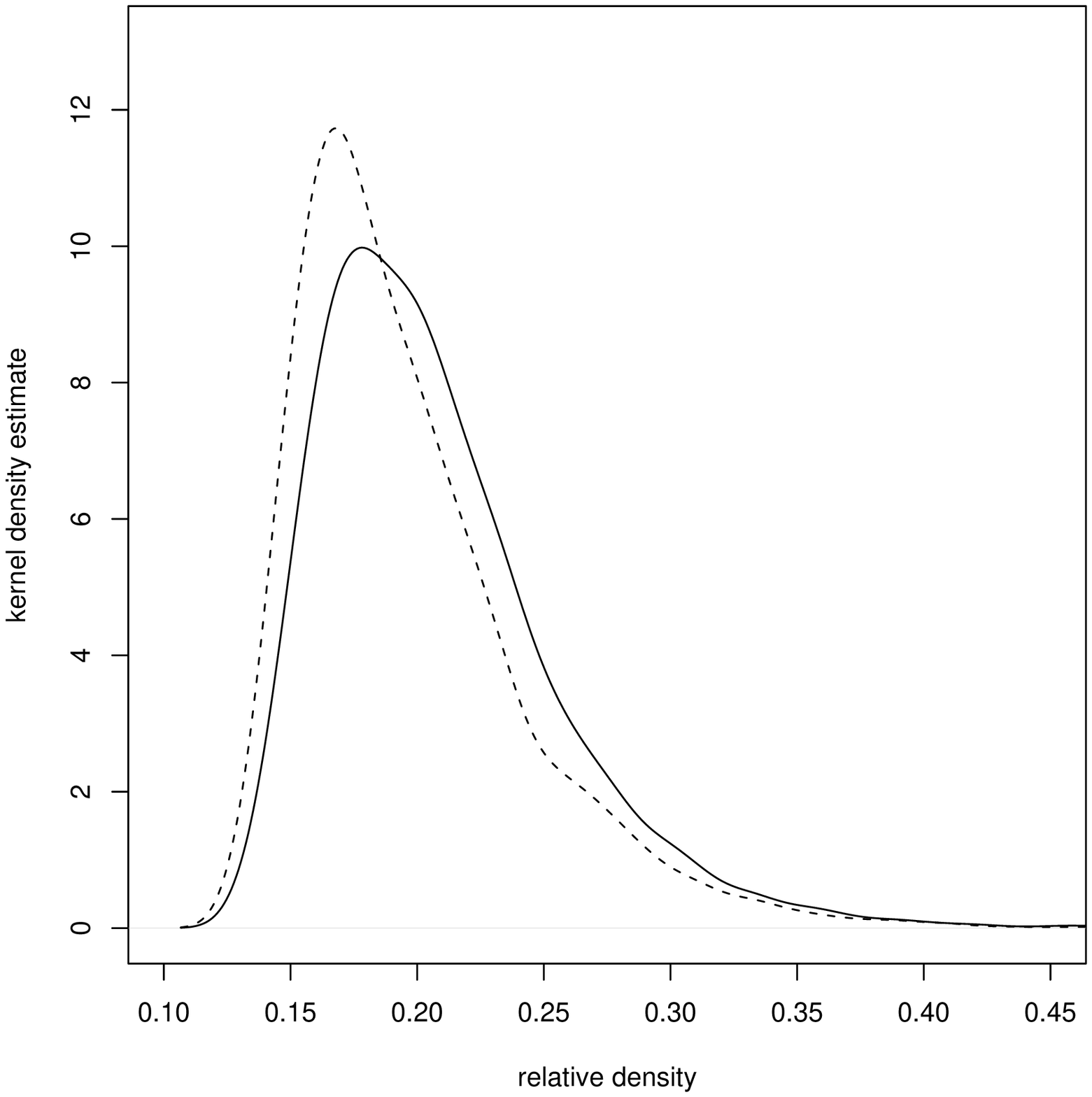}}}
\rotatebox{-90}{ \resizebox{1.7 in}{!}{ \includegraphics{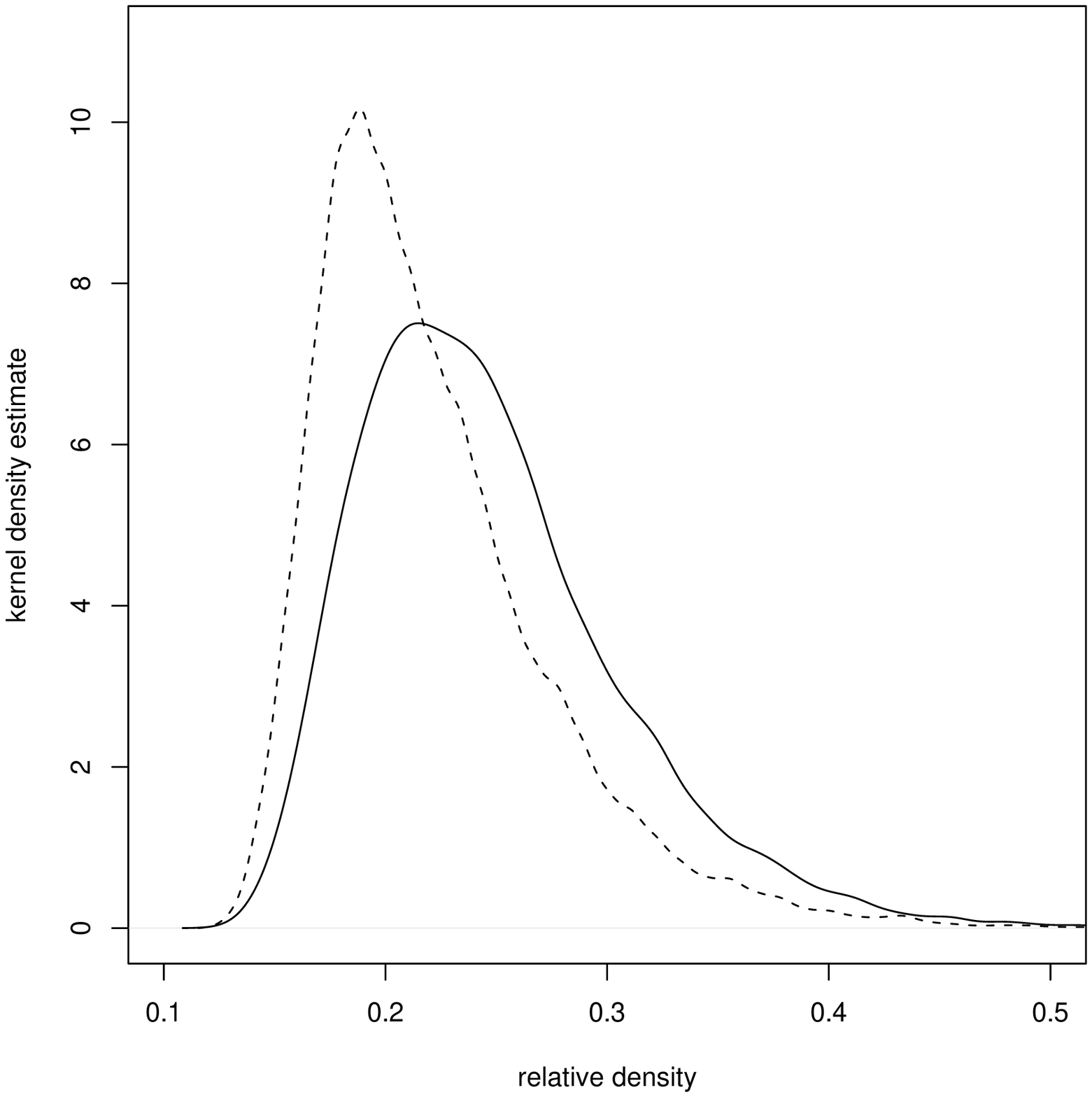}}}
\rotatebox{-90}{ \resizebox{1.7 in}{!}{ \includegraphics{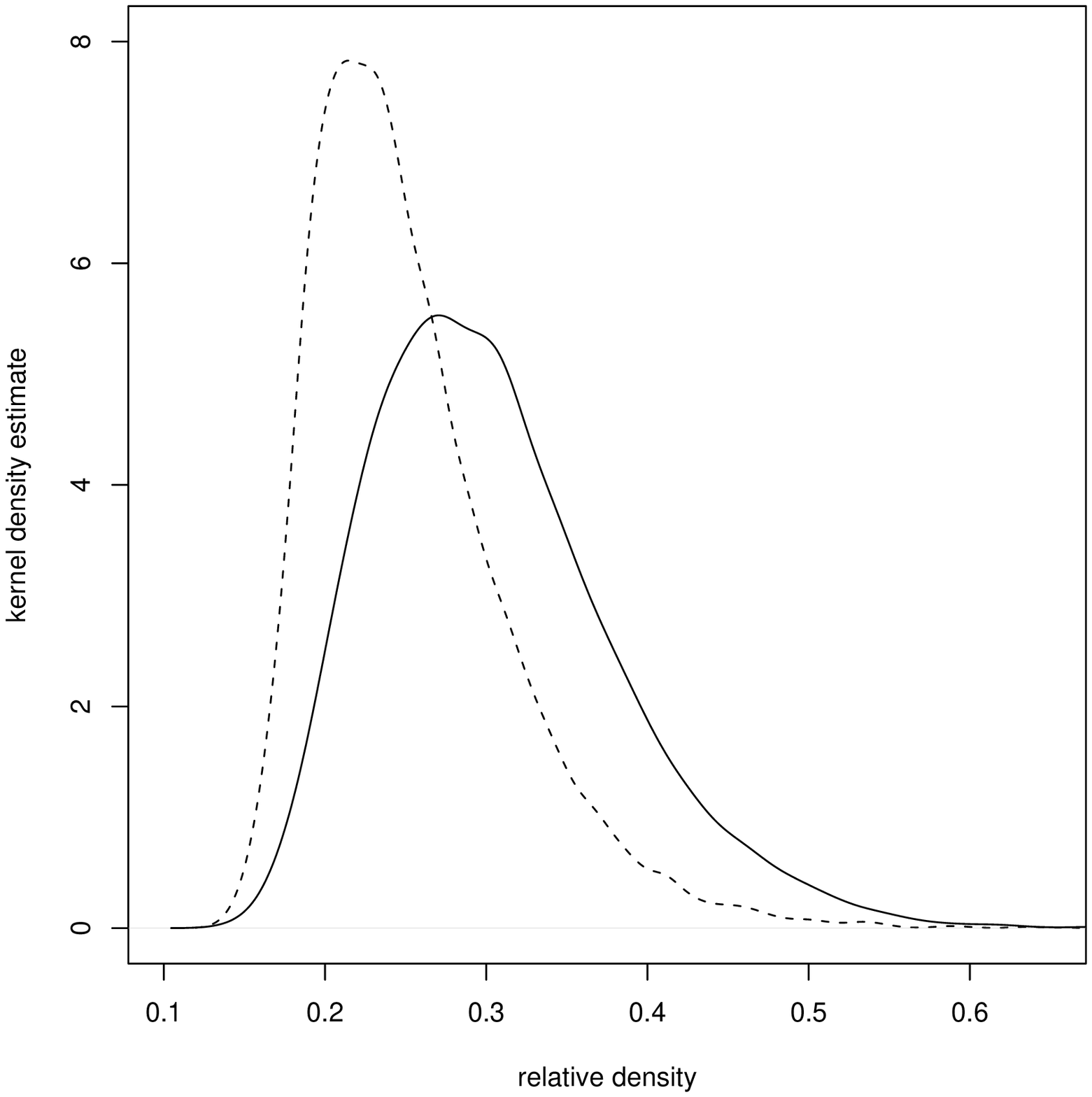}}}
\rotatebox{-90}{ \resizebox{1.7 in}{!}{ \includegraphics{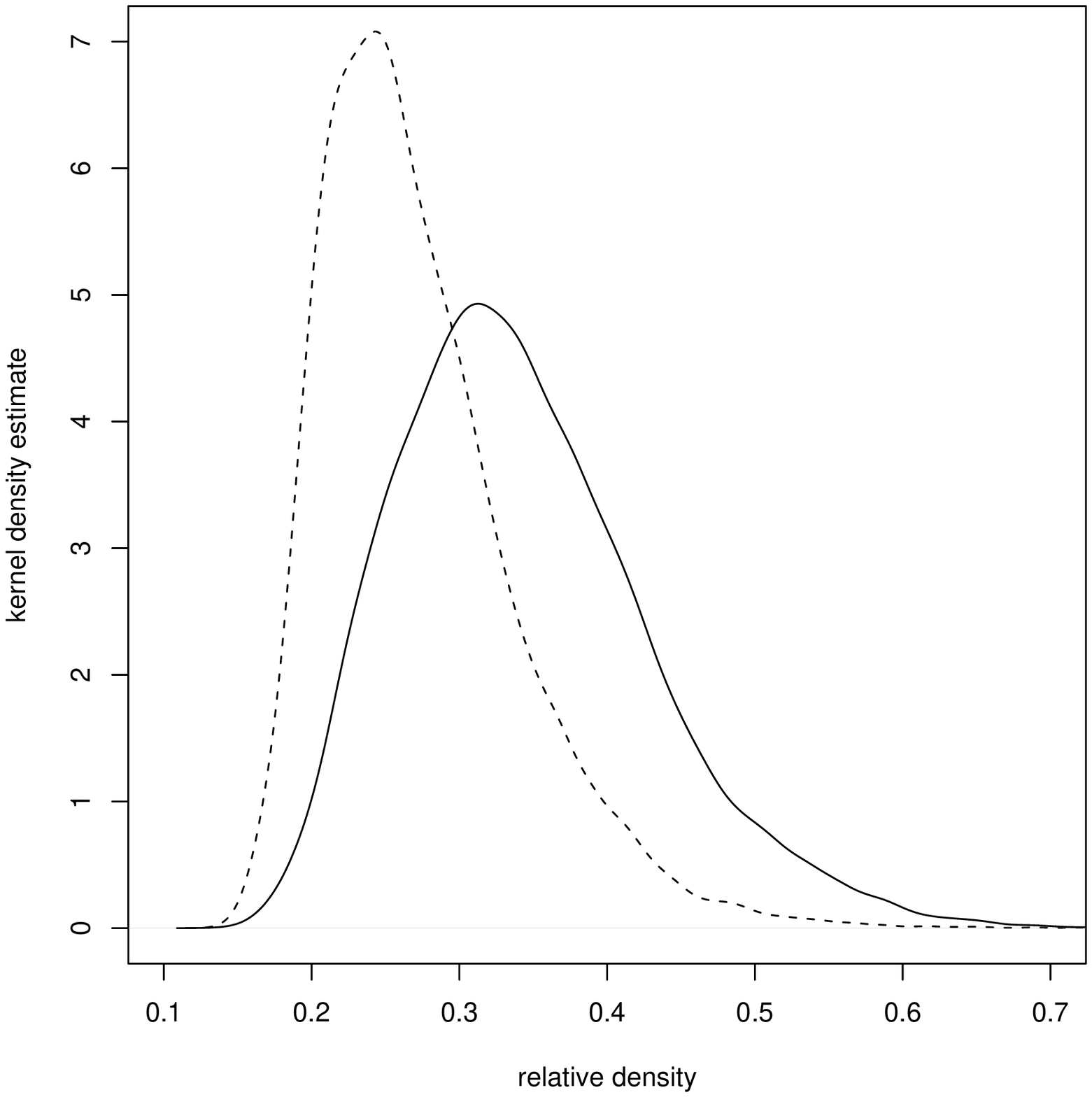}}}
\rotatebox{-90}{ \resizebox{1.7 in}{!}{ \includegraphics{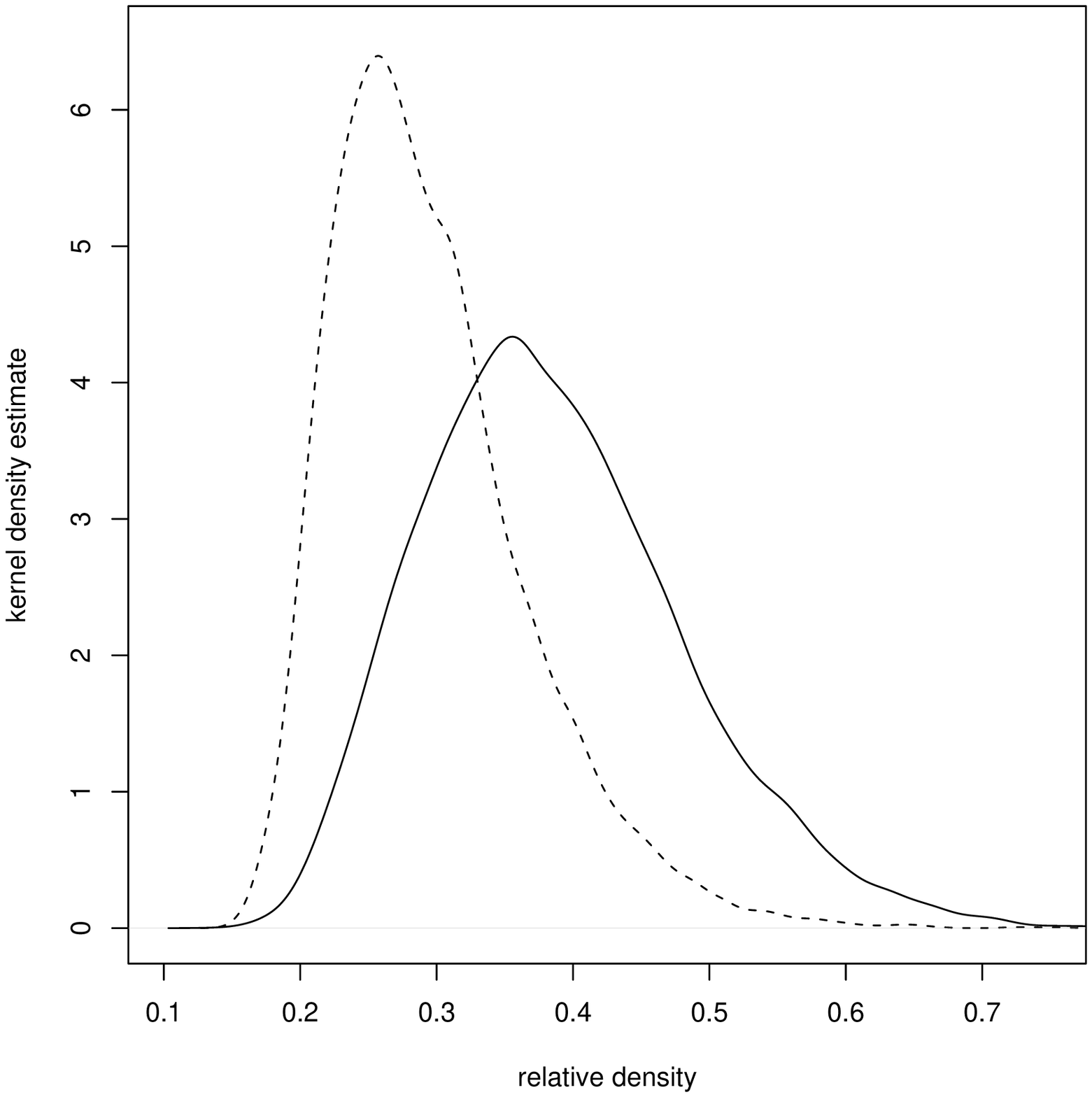}}}
\rotatebox{-90}{ \resizebox{1.7 in}{!}{ \includegraphics{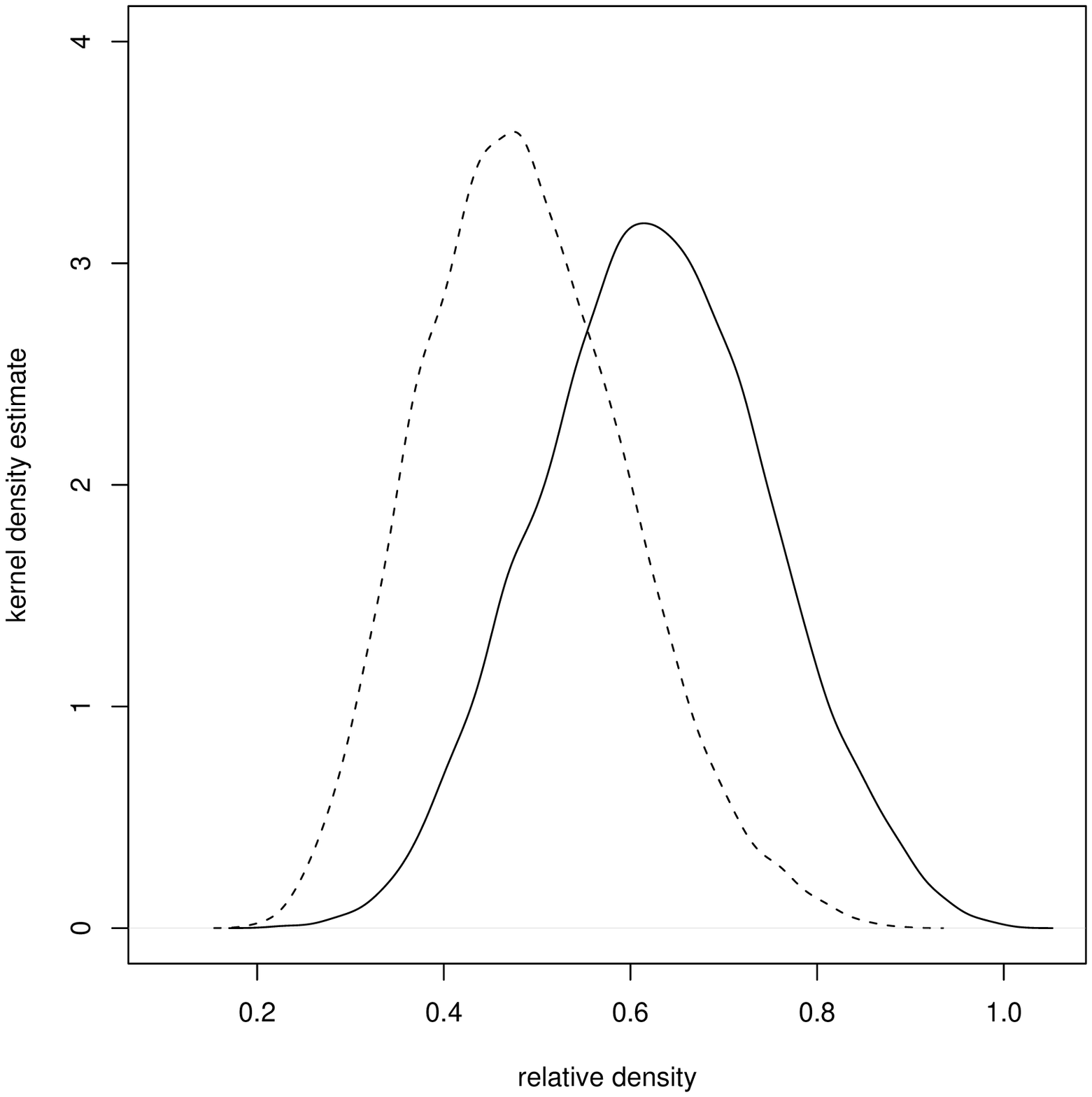}}}
\rotatebox{-90}{ \resizebox{1.7 in}{!}{ \includegraphics{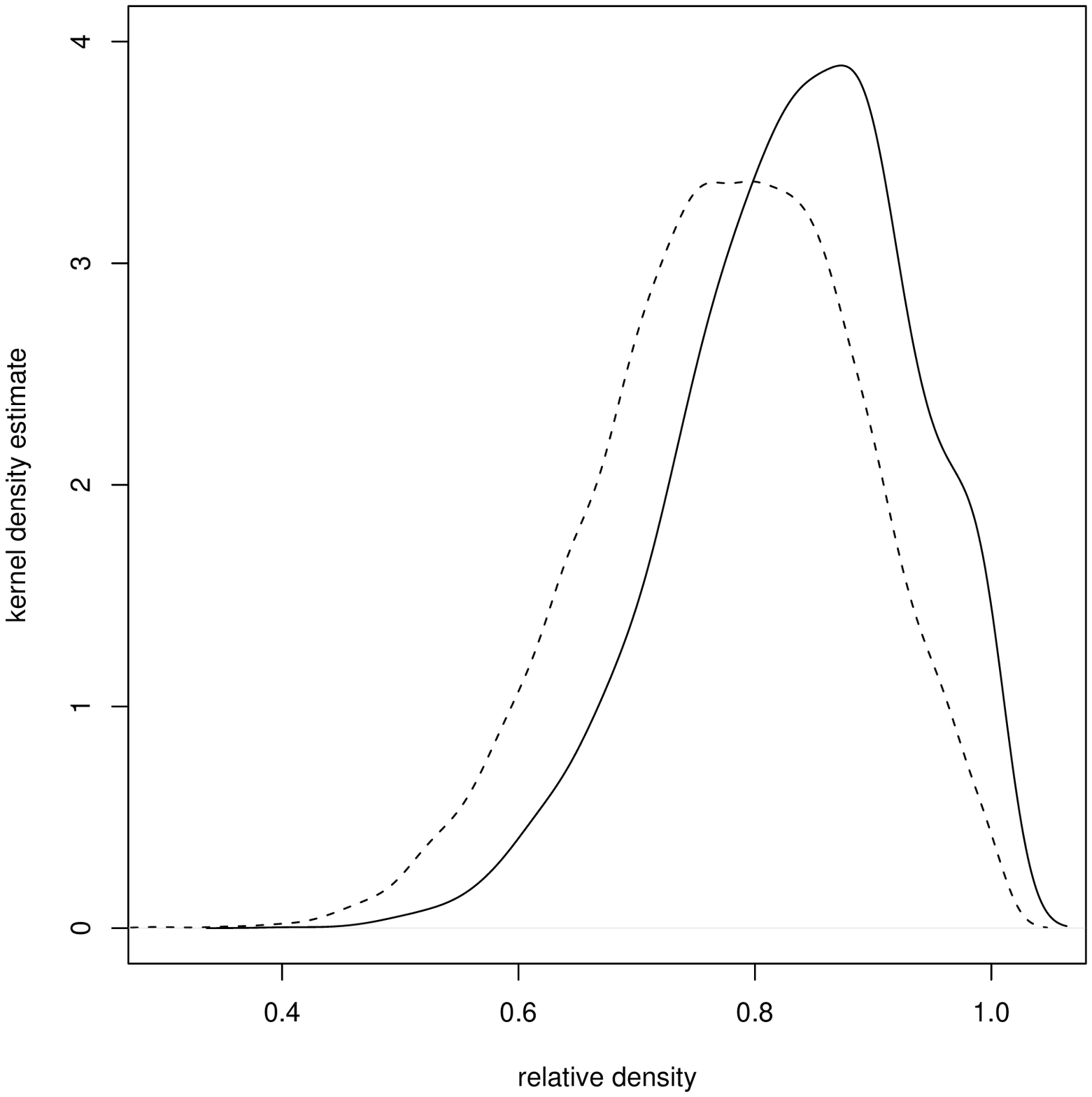}}}
\rotatebox{-90}{ \resizebox{1.7 in}{!}{ \includegraphics{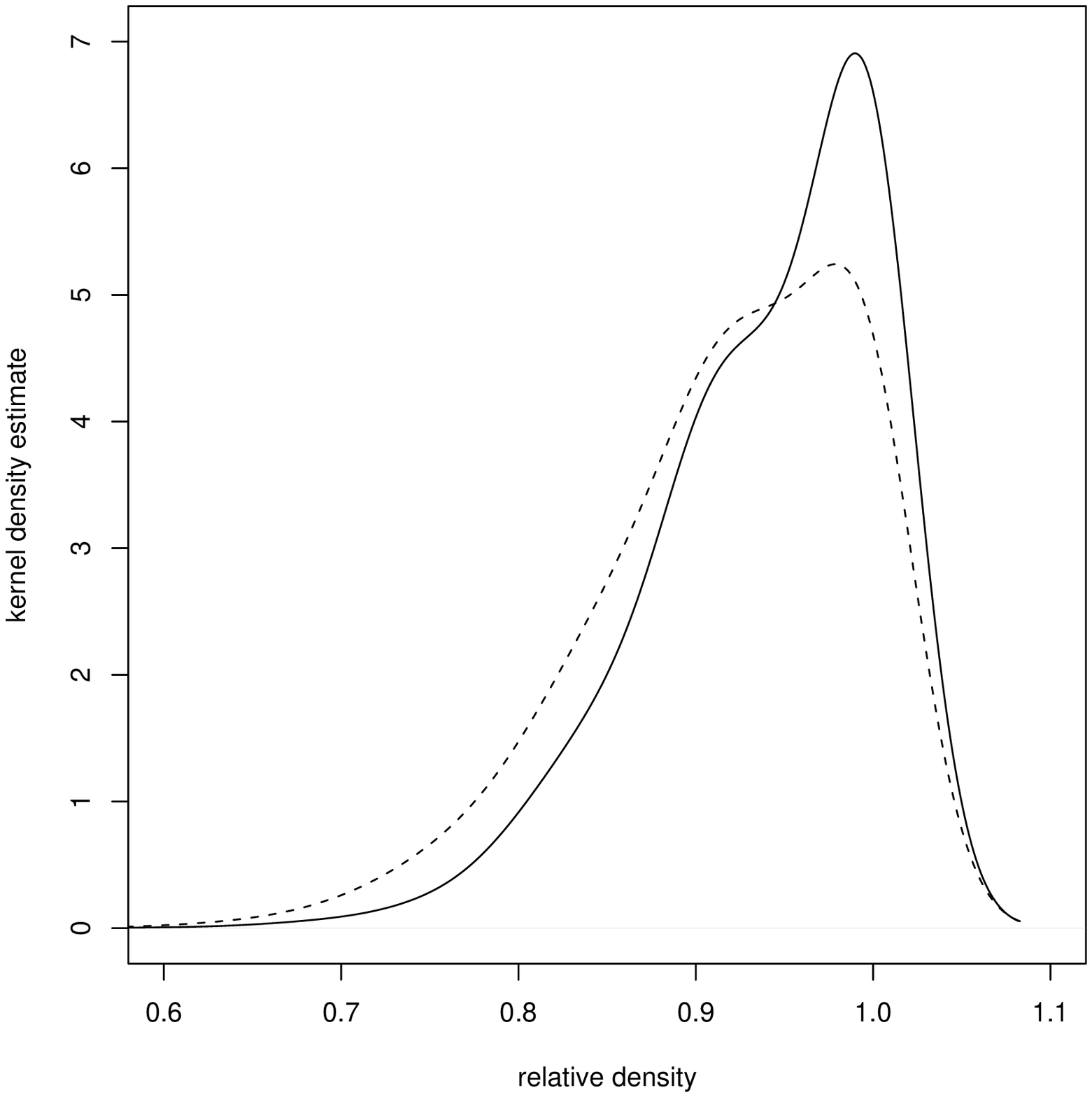}}}
\rotatebox{-90}{ \resizebox{1.7 in}{!}{ \includegraphics{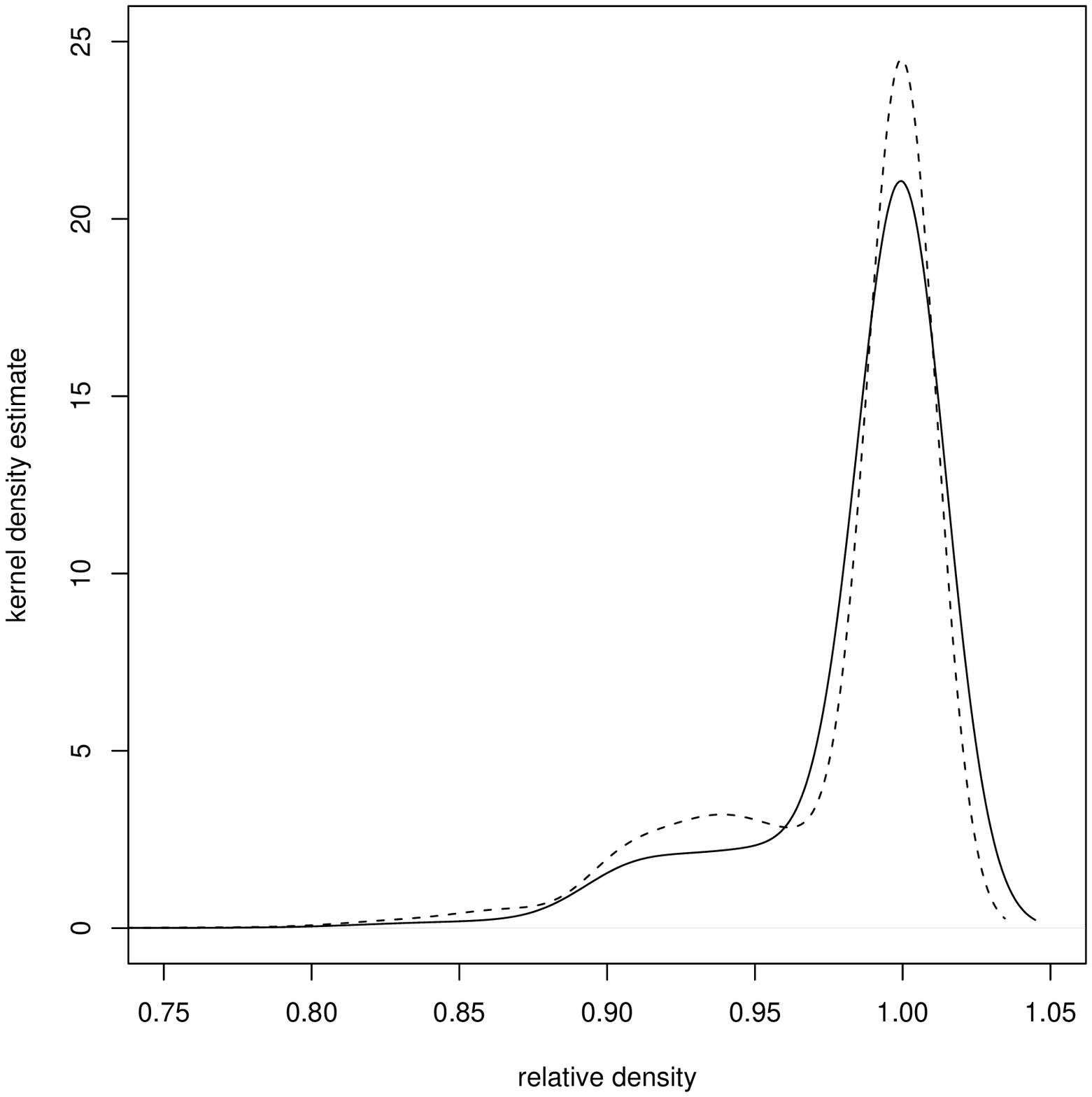}}}
\caption{\label{fig:agg2sim1-10}
Kernel density estimates for the null (solid) and the association alternative $H^A_{\sqrt{3}/12}$ (dashed) for $r= 1,\, 11/10,\, 6/5,\, 4/3,\, \sqrt{2},\, 3/2,\, 2,\, 3,\, 5,\, 10$ (left-to-right).
}
\end{figure}

For association with $\epsilon=\sqrt{3}/12\approx .144$, in Figure \ref{fig:agg2sim1-10}, are the kernel density estimates for the null case and the segregation alternative for the ten $r$ values with  $n=10,\,N=10,000$. Observe that under  $H^A_{\sqrt{3}/12}$, kernel density estimates are skewed right for $r=1,\,11/10,\,6/5,\,4/3$, (with skewness increasing as $r$ gets smaller) and kernel density estimates are almost symmetric for $r=\sqrt{2}, 3/2,\, 2 $, with most symmetry occurring at  $r=2 $, kernel density estimates are skewed left for $r=3,\,5,\,10$, (with skewness increasing as  $r$ gets larger).

\begin{figure}[]
\centering
\psfrag{kernel density estimate}{ \Huge{\bfseries{kernel density estimate}}}
\psfrag{relative density}{ \Huge{\bfseries{relative density}}}
\rotatebox{-90}{ \resizebox{2.2 in}{!}{ \includegraphics{agg2sim2.ps}}}
\rotatebox{-90}{ \resizebox{2.2 in}{!}{ \includegraphics{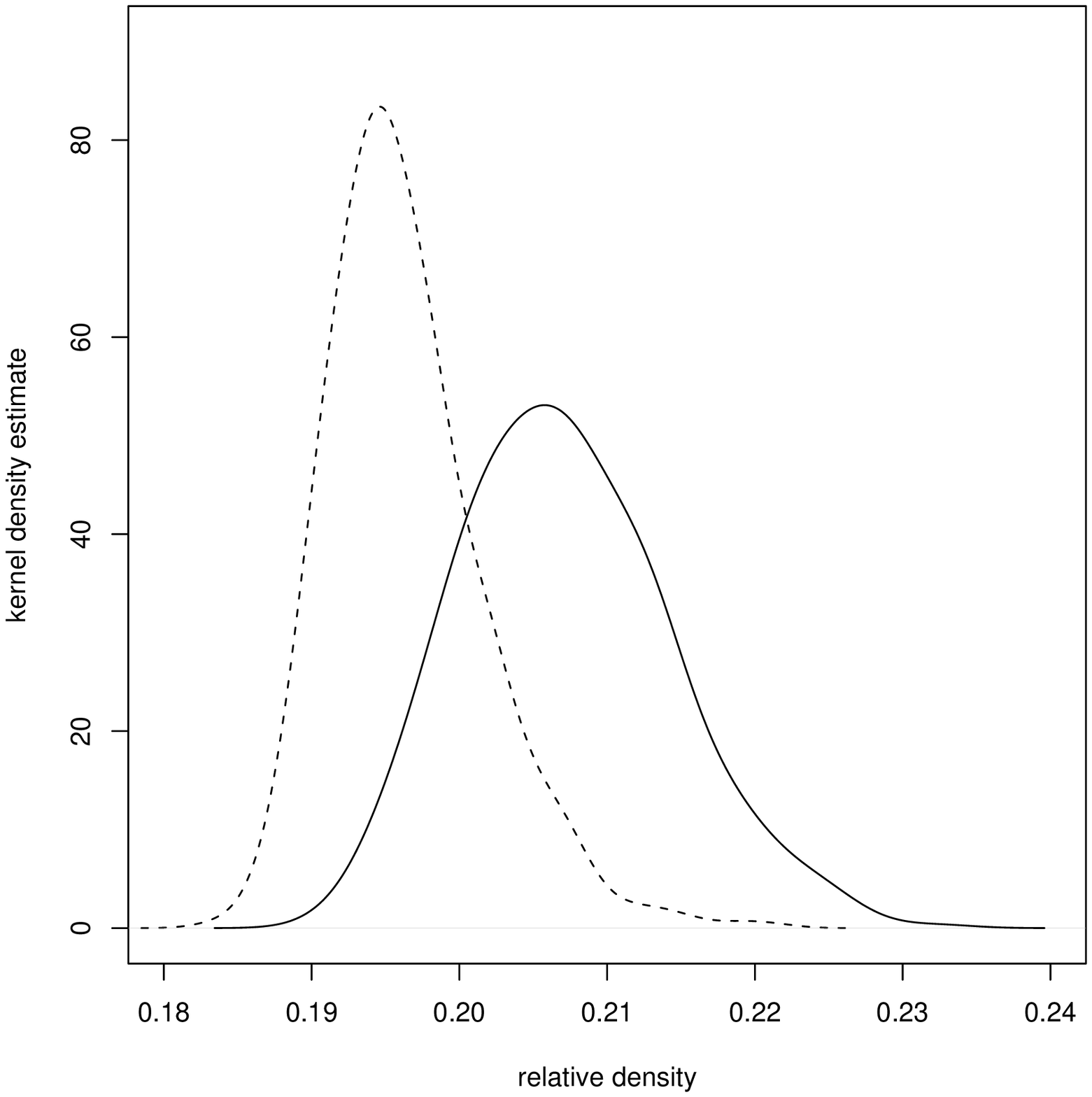}}}
\caption{ \label{fig:AggSimPowerPlots}
Two Monte Carlo experiments against the association alternative $H^A_{\sqrt{3}/12}$.
Depicted are kernel density estimates for $\rho_n(11/10)$ for
$n=10$ (left) and $n=100$ (right) under the null (solid) and alternative (dashed).
}
\end{figure}

Note also that for $r=11/10$ with $n=10,\,N=1000$, the kernel density estimates are very similar, implying small power. With $N=10,000$ $\widehat{\beta}^A_{mc}\left(10,\, \sqrt{3}/12 \right)= 0.0921$, $\widehat{C}^A_n = 0.1\bar{5}$, and $\widehat{\alpha}^A_{mc}(10)= 0.0484$.  See Figure \ref{fig:AggSimPowerPlots}. Note that for large $n$, there is more separation between null and alternative kernel densities, which implies higher power. With $n=100,\,N=1000$ and get $\widehat{C}^A_n = 0.1963$, $\widehat{\alpha}^A_{mc}(100)= 0.049$, and $\widehat{\beta}^A_{mc}\left(100,\sqrt{3}/12 \right) = 0.56$.

\begin{figure}[]
\centering
\psfrag{kernel density estimate}{ \Huge{\bfseries{kernel density estimate}}}
\psfrag{relative density}{ \Huge{\bfseries{relative density}}}
\rotatebox{-90}{ \resizebox{1.7 in}{!}{ \includegraphics{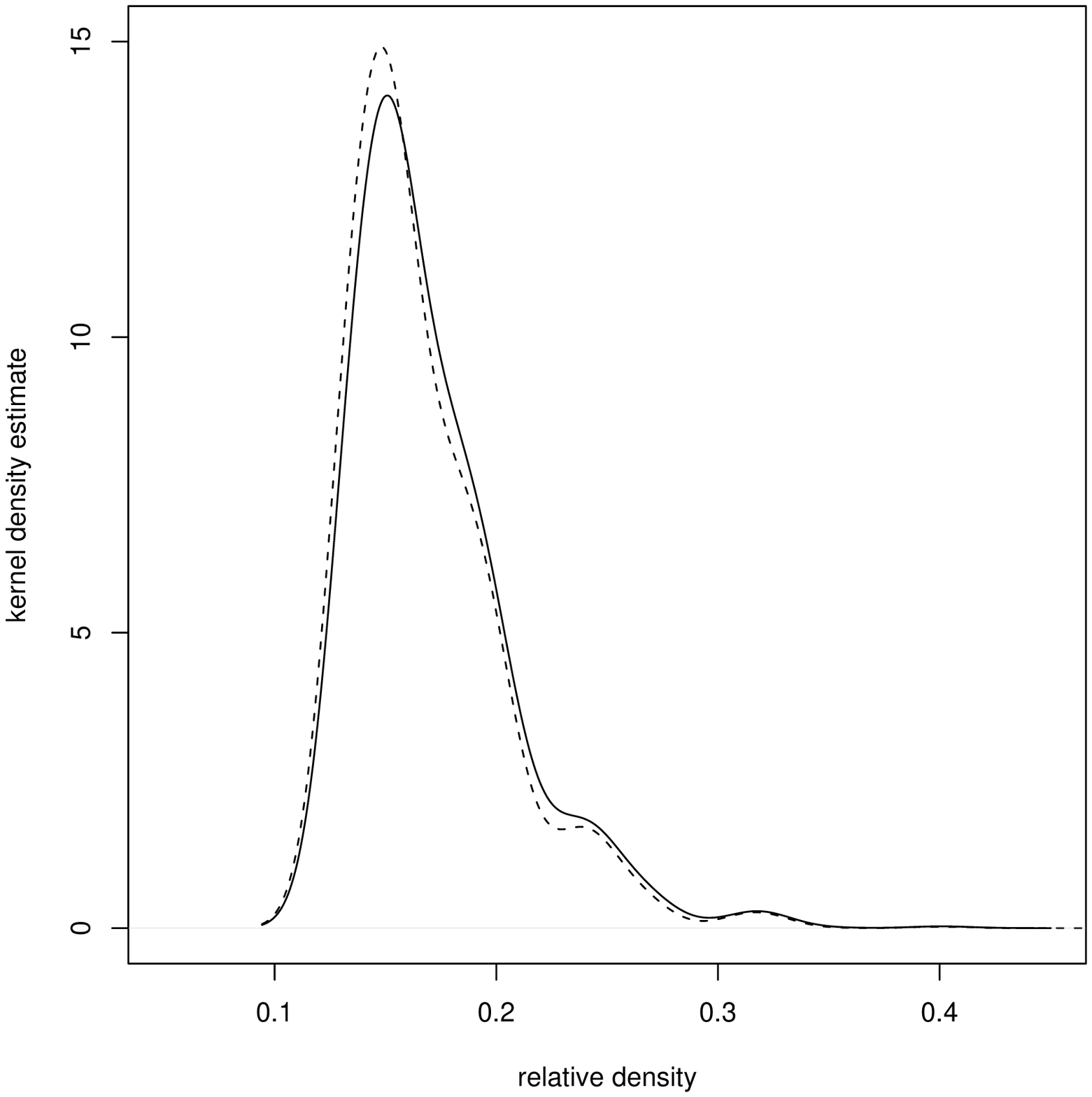}}}
\rotatebox{-90}{ \resizebox{1.7 in}{!}{ \includegraphics{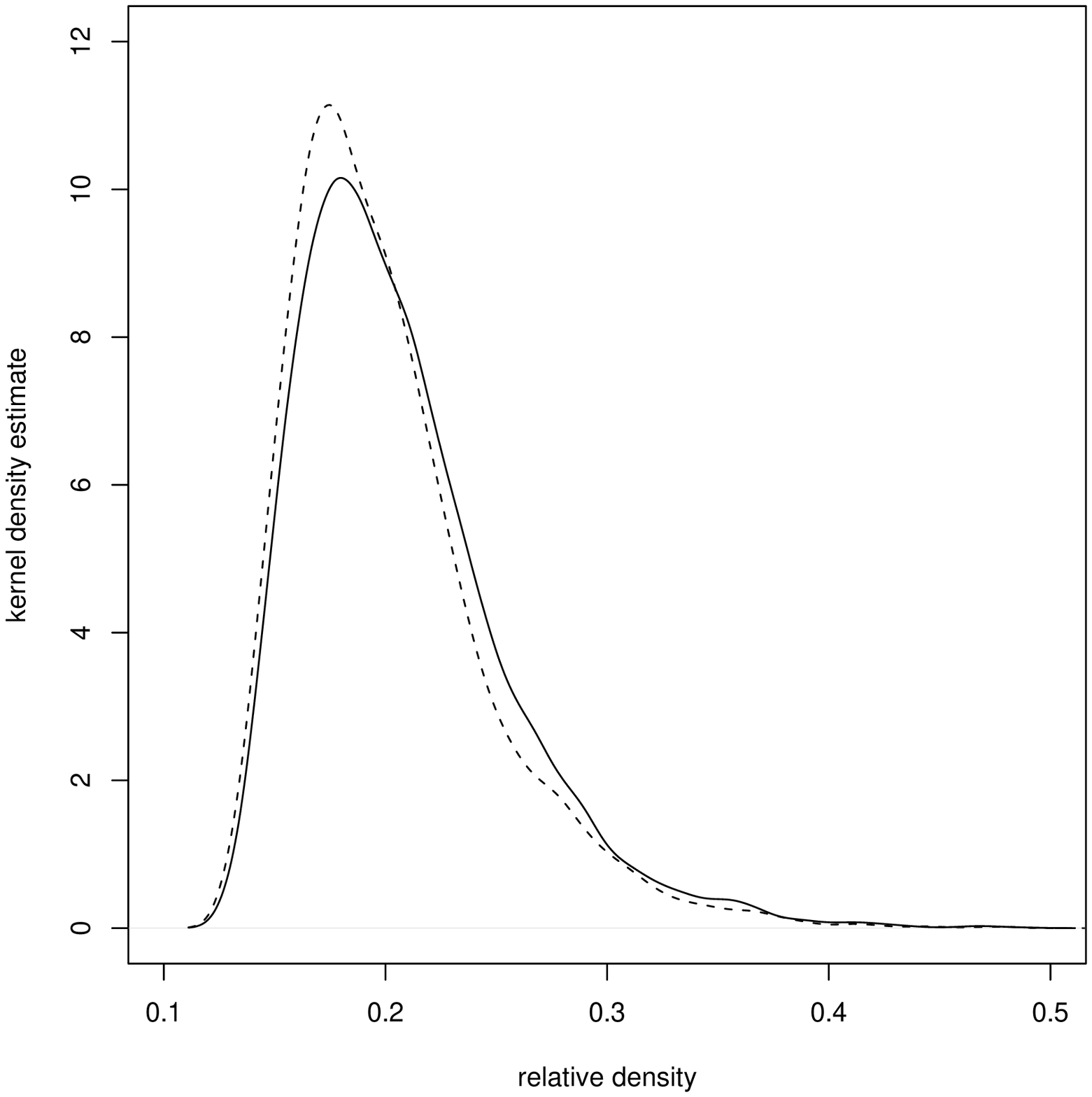}}}
\rotatebox{-90}{ \resizebox{1.7 in}{!}{ \includegraphics{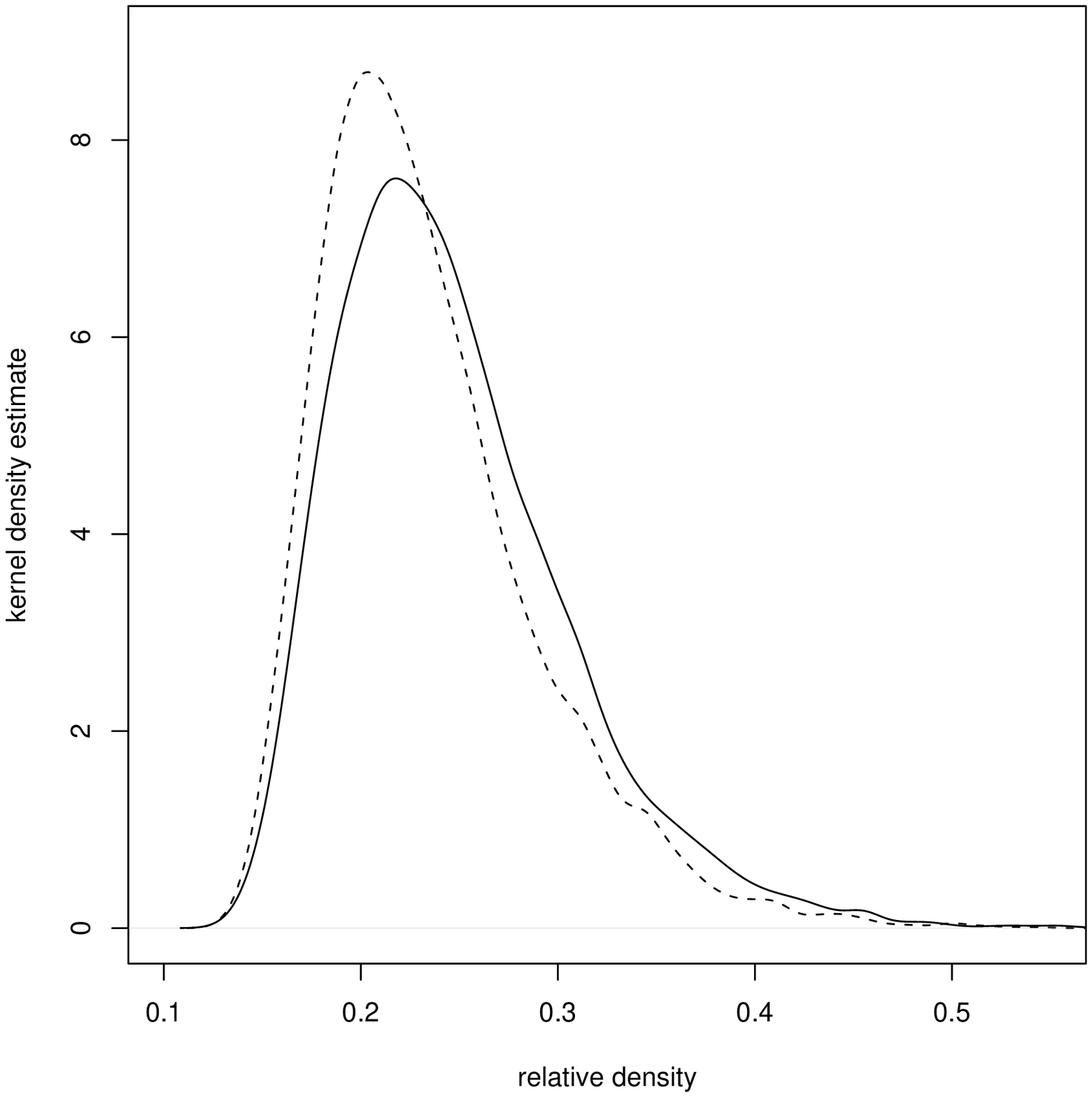}}}
\rotatebox{-90}{ \resizebox{1.7 in}{!}{ \includegraphics{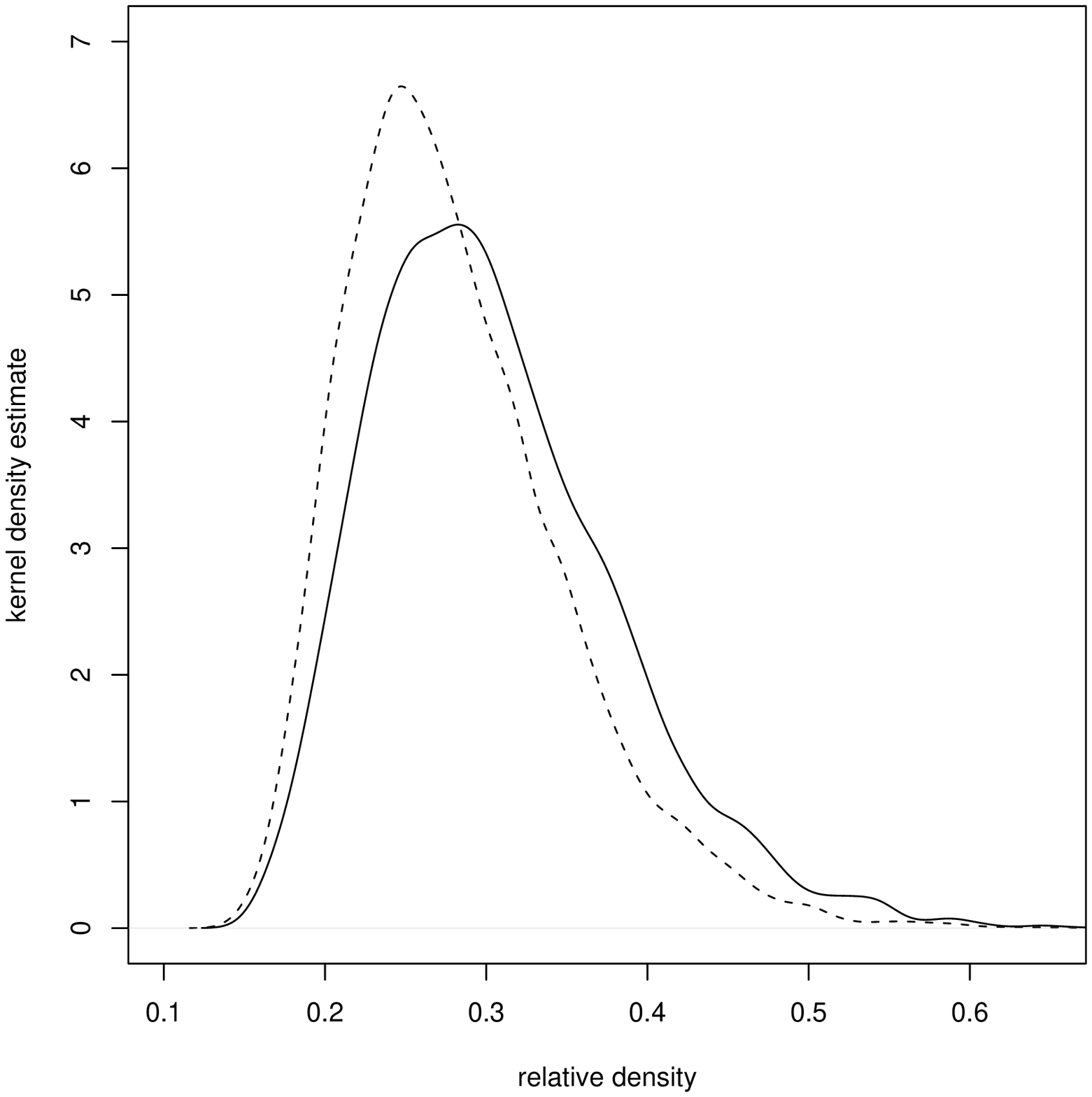}}}
\rotatebox{-90}{ \resizebox{1.7 in}{!}{ \includegraphics{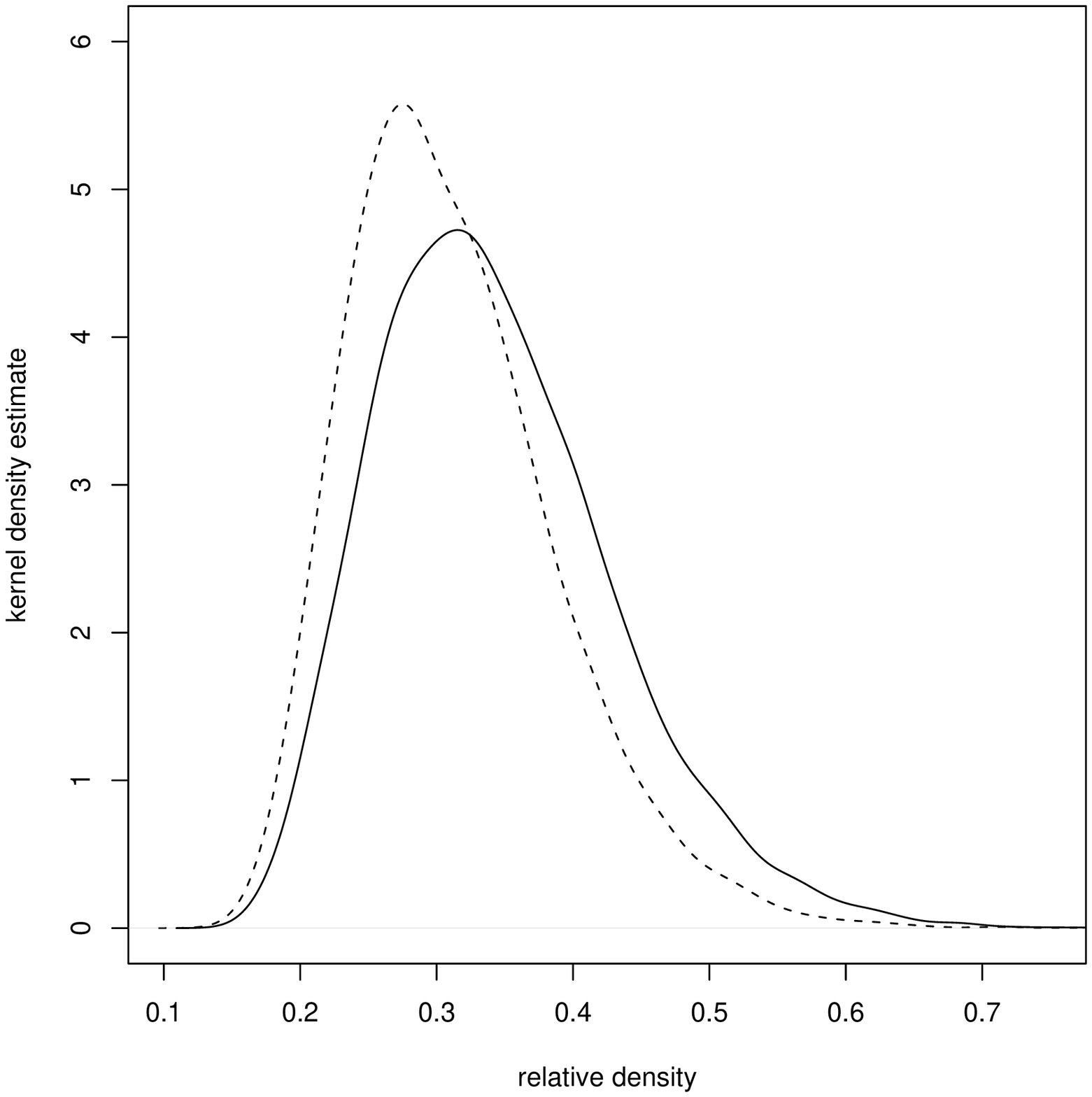}}}
\rotatebox{-90}{ \resizebox{1.7 in}{!}{ \includegraphics{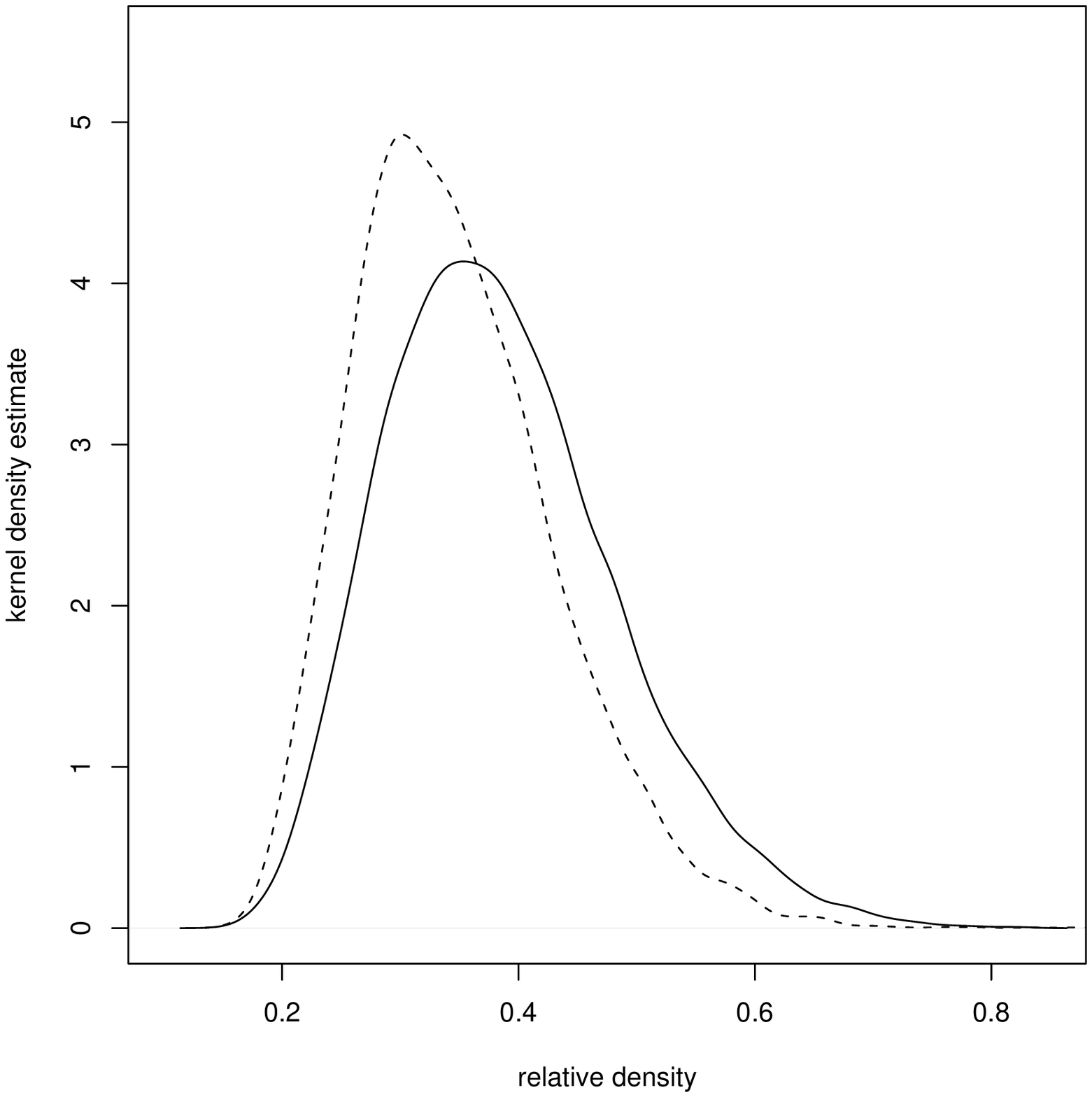}}}
\rotatebox{-90}{ \resizebox{1.7 in}{!}{ \includegraphics{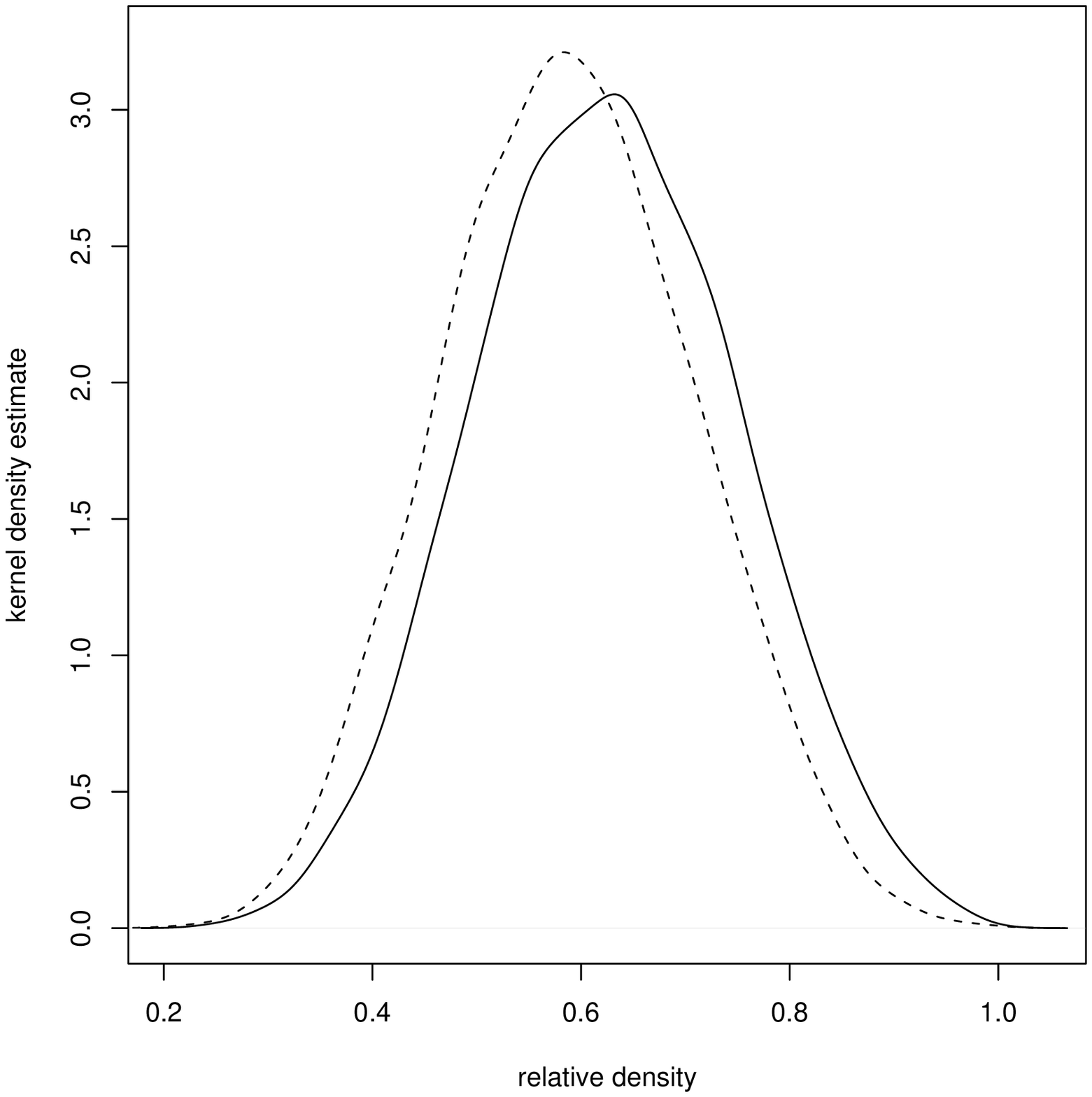}}}
\rotatebox{-90}{ \resizebox{1.7 in}{!}{ \includegraphics{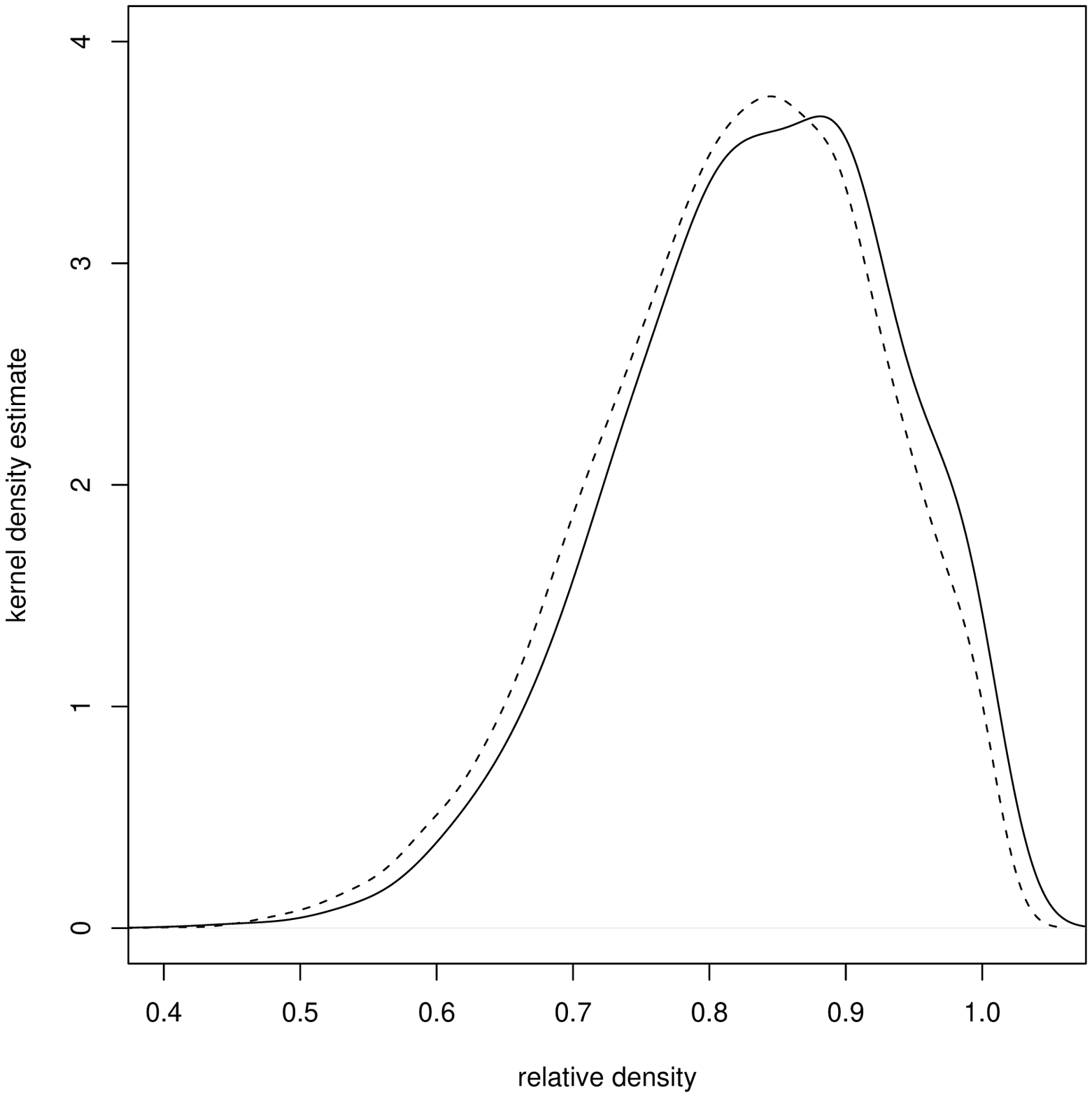}}}
\rotatebox{-90}{ \resizebox{1.7 in}{!}{ \includegraphics{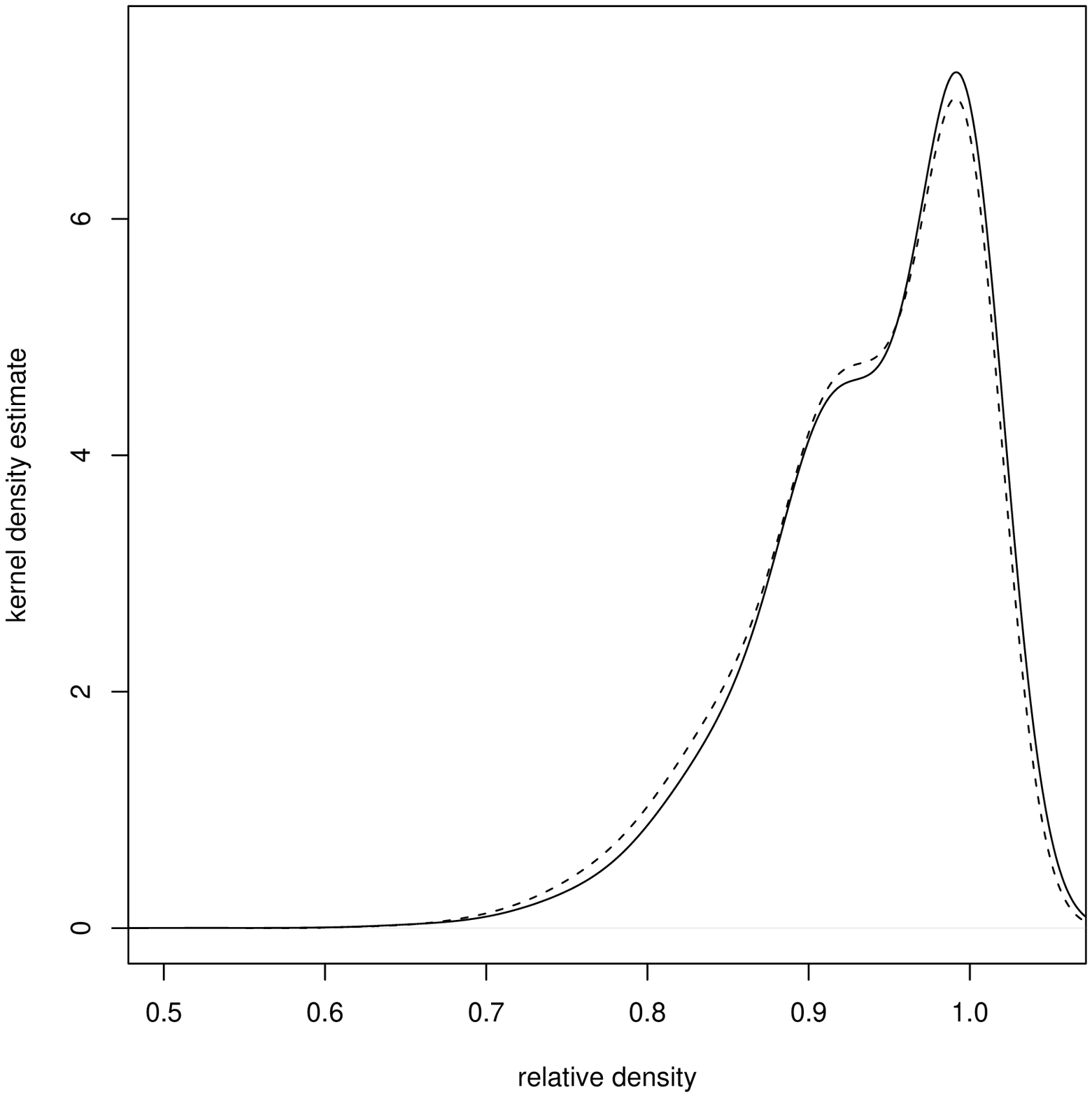}}}
\rotatebox{-90}{ \resizebox{1.7 in}{!}{ \includegraphics{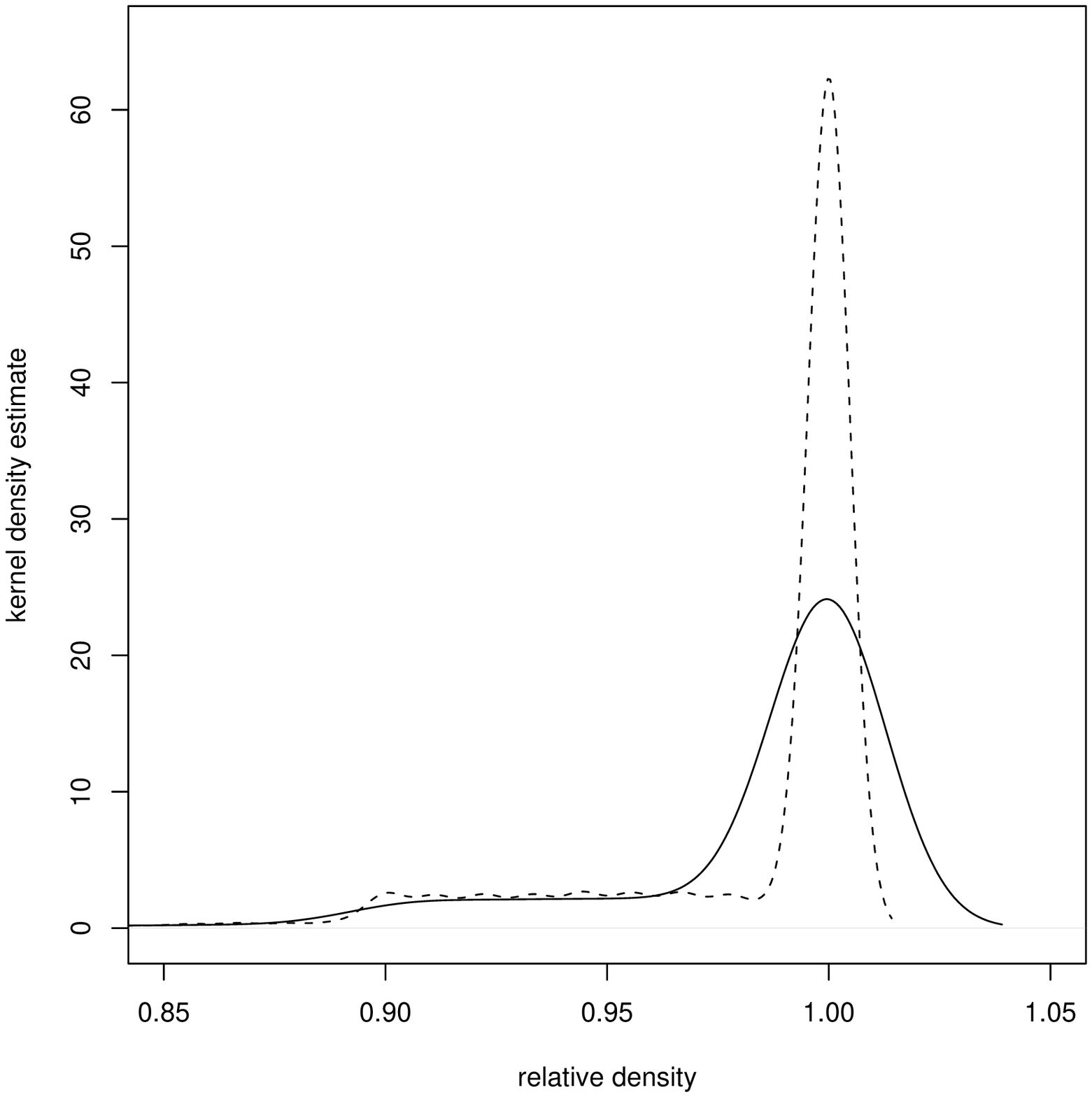}}}
\caption{\label{fig:agg3sim1-10}
Kernel density estimates for the null (solid) and the association alternative $H^A_{\sqrt{3}/21}$ (dashed) for $r= 1,\, 11/10,\, 6/5,\, 4/3,\, \sqrt{2},\, 3/2,\, 2,\, 3,\, 5,\, 10$ (left-to-right).
}
\end{figure}

For association with $\epsilon=\sqrt{3}/21\approx .0825$, in Figure \ref{fig:agg3sim1-10}, are the kernel density estimates for the null case and the segregation alternative for the ten $r$ values with  $n=10,\,N=10,000$. Observe that under  $H^A_{\sqrt{3}/21}$, kernel density estimates are skewed right for $r=1,\,11/10,\,6/5,\,4/3$, (with skewness increasing as $r$ gets smaller) and kernel density estimates are almost symmetric for $r=\sqrt{2}, 3/2,\, 2 $, with most symmetry occurring at $r=3/2$, kernel density estimates are skewed left for $r=3,\,5,\,10$, (with skewness increasing as  $r$ gets larger).

\begin{figure}[]
\centering
\psfrag{power}{ \Huge{\bfseries{power}}}
\psfrag{r}{\Huge{$r$}}
\rotatebox{-90}{ \resizebox{1.8 in}{!}{ \includegraphics{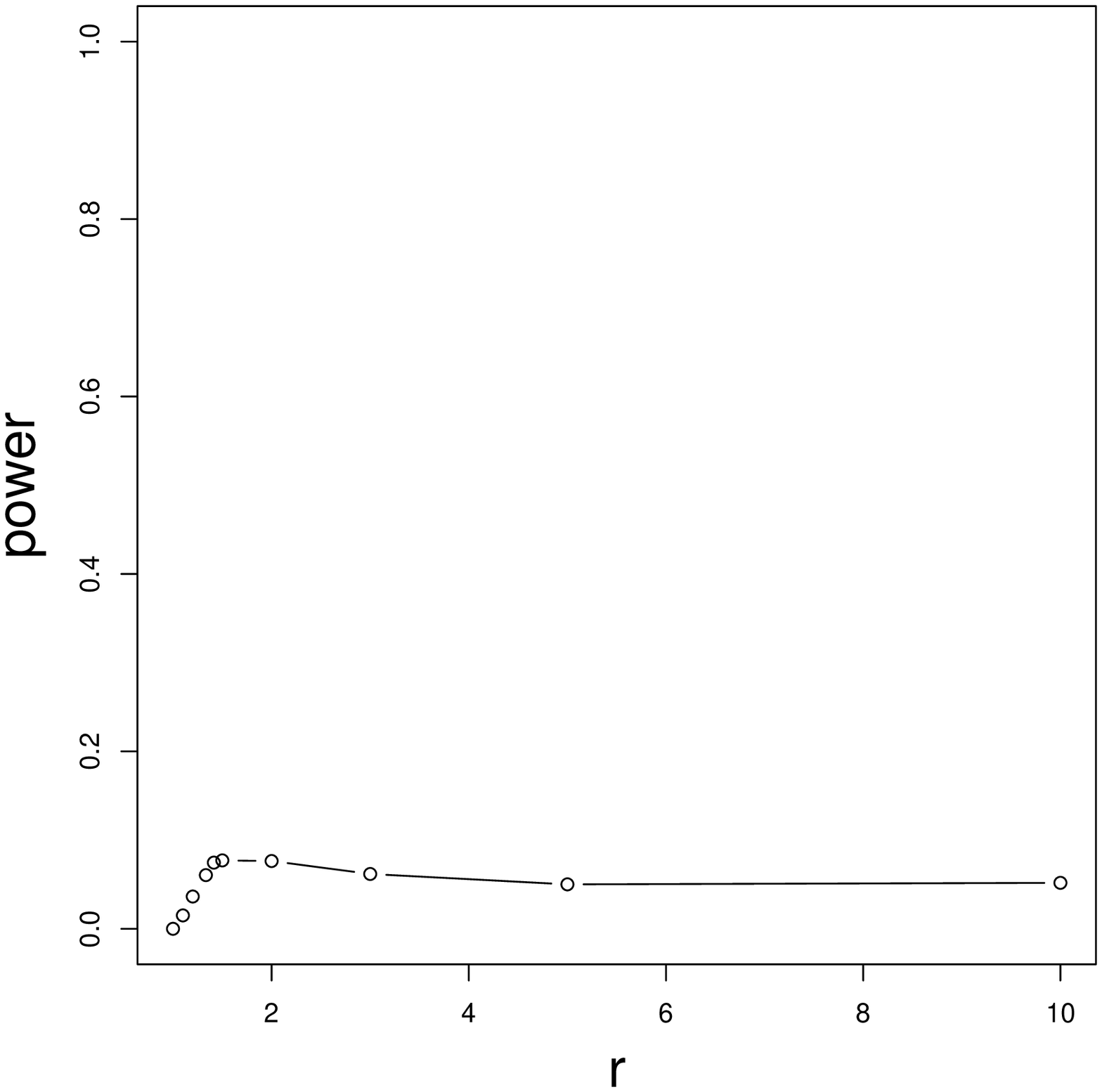}}}
\rotatebox{-90}{ \resizebox{1.8 in}{!}{ \includegraphics{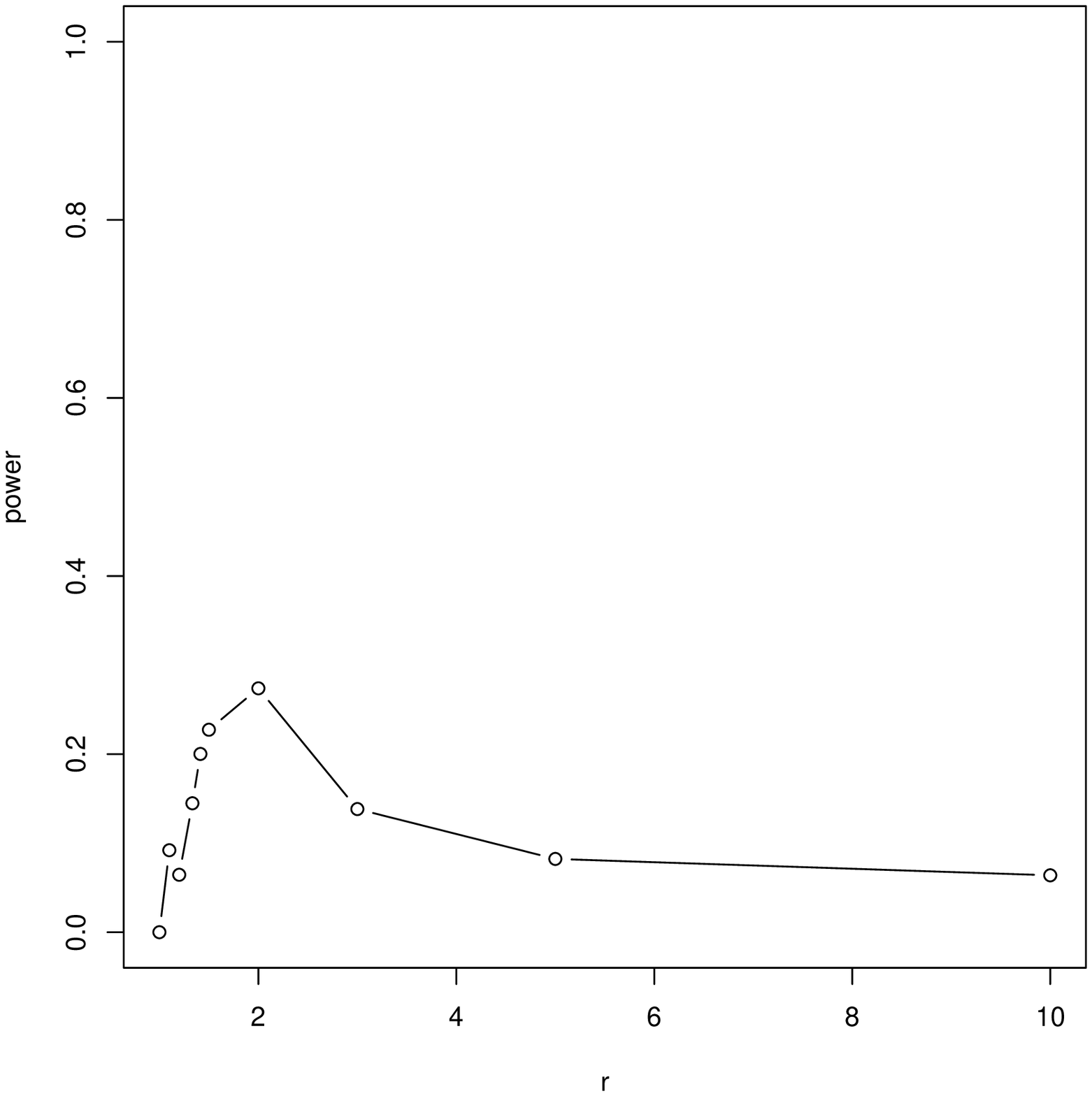}}}
\rotatebox{-90}{ \resizebox{1.8 in}{!}{ \includegraphics{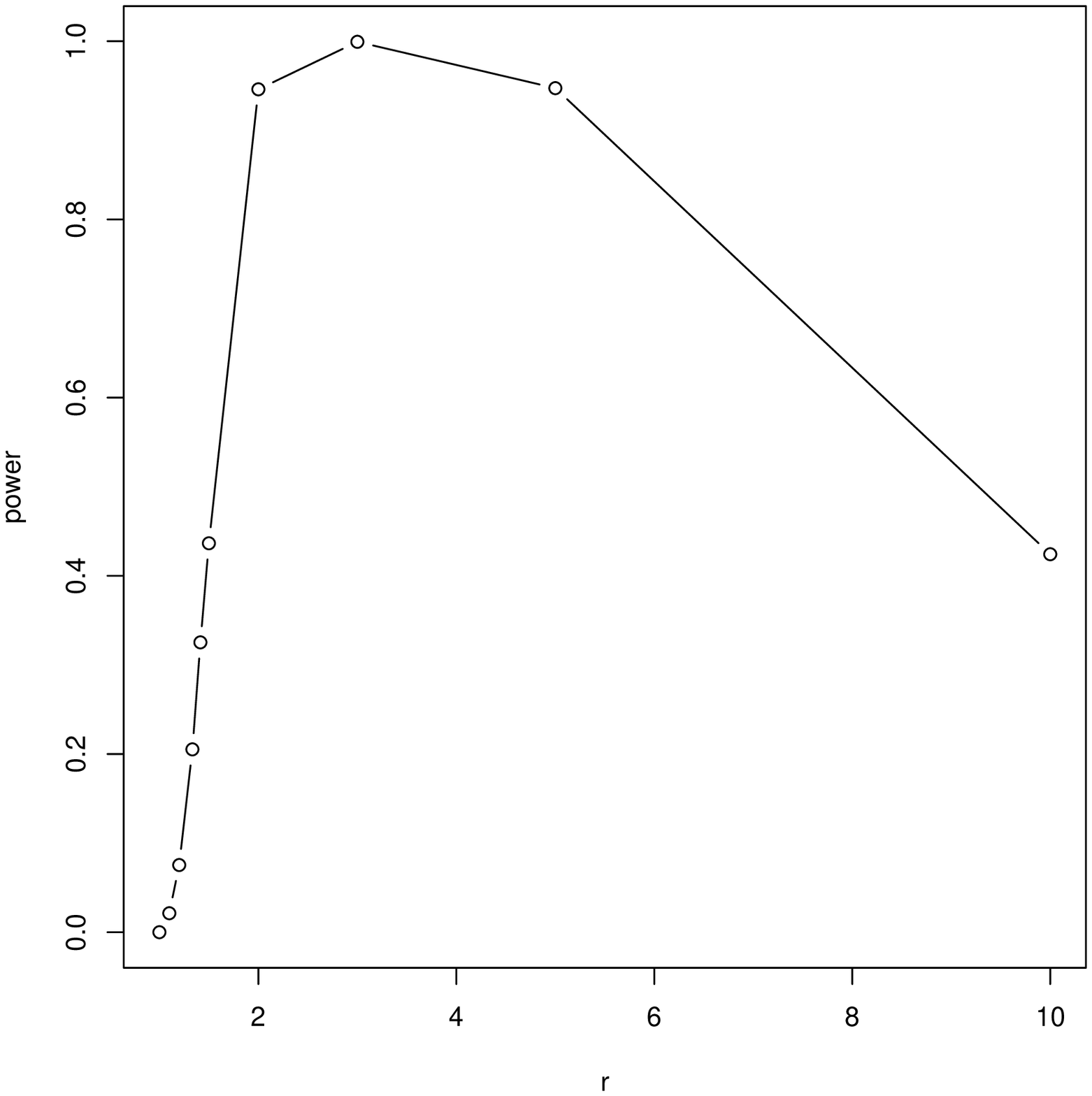}}}
\caption{ \label{fig:aggpow1-3}
Monte Carlo power using the empirical critical value against association alternatives
$ H^A_{\sqrt{3}/21}$ (left),
$H^A_{\sqrt{3}/12}$ (middle)
and
$H^A_{5\,\sqrt{3}/24}$ (right)
as a function of $r$, for $n=10$.}
\end{figure}

We also plot the empirical power as a function of $r$ in Figure \ref{fig:aggpow1-3}.  Let $r^*_A(\epsilon)$ be the value at which maximum Monte Carlo power estimate occurs.  Then  $r^*_A\left( 5\,\sqrt{3}/24 \right)=3$, $r^*_A\left( \sqrt{3}/12 \right)=2$, and for $r^*_A\left( \sqrt{3}/21 \right)=3/2$. Notice that the more severe the association the larger the value of $r^*_A$. Based on the analysis of the Monte Carlo power estimates, we suggest moderate $r$ values for moderate association.

\begin{table}[]
\centering
\begin{tabular}{|c|c|c|c|c|c|c|c|c|c|c|}
\hline
$r$  & 1 & 11/10 &6/5 & 4/3 & $\sqrt{2}$ & 3/2 & 2 & 3 & 5 & 10\\
\hline
$\widehat{C}^A_n$ & $.1\bar{3}$ & $.1\bar{4}$ & $.1\bar{6}$ & .2 & $.\bar 2$ & $.2\bar 4$ & $.4\,\bar 2$ & .65 & $.8\bar 2$ & $.9\bar 1$\\
\hline
$\widehat{\alpha}^A_{mc}(10)$& 0 & .0112  & .0208 & .0308 & .0363 &  .0359 & .0392 & .0413 & .0478 & .0398 \\
\hline
$\widehat{\beta}^A_{mc}(10,\,5\,\sqrt{3}/24)$ & 0 &  .0213 & .0754 & .2052  &  .3253 & .4365 & .946 & .9993 & .9473 & .4242 \\
\hline
$\widehat{\beta}^A_{mc}(10,\,\sqrt{3}/12)$ & 0 &  .0921 & .0645 & .1448  &  .2002 & .2274 & .2739 & .1383 & .0823 & .0639 \\
\hline
$\widehat{\beta}^A_{mc}(10,\,\sqrt{3}/21)$ & 0 &  .0151 & .0364 & .0605  &  .0746 & .0771 & .0764 & .0618 & .0501 & .0518\\
\hline
\end{tabular}
\caption{
\label{tab:emp-val-A}
The empirical critical values, empirical significance levels, and empirical power estimates under $H^A_{\epsilon}$ for $\epsilon=5\,\sqrt{3}/24,\,\sqrt{3}/12,\,\sqrt{3}/21$ and $n=10$ at $\alpha=.05$.}
\end{table}

We also estimate the power using the asymptotic critical value in association alternatives for various values of $r$. For each $r$ value, the level $\alpha$ asymptotic critical value is $\mu(r)+z_{\alpha} \cdot \sqrt{\nu(r)/n}$.  We estimate the empirical power as $\widehat{\beta}^A_n(r,\epsilon):=\frac{1}{N}\sum_{j=1}^{N}\I(R_j< z_{\alpha})$.

\begin{figure}[]
\centering
\psfrag{power}{ \Huge{\bfseries{power}}}
\psfrag{r}{\Huge{$r$}}
\rotatebox{-90}{ \resizebox{1.7 in}{!}{ \includegraphics{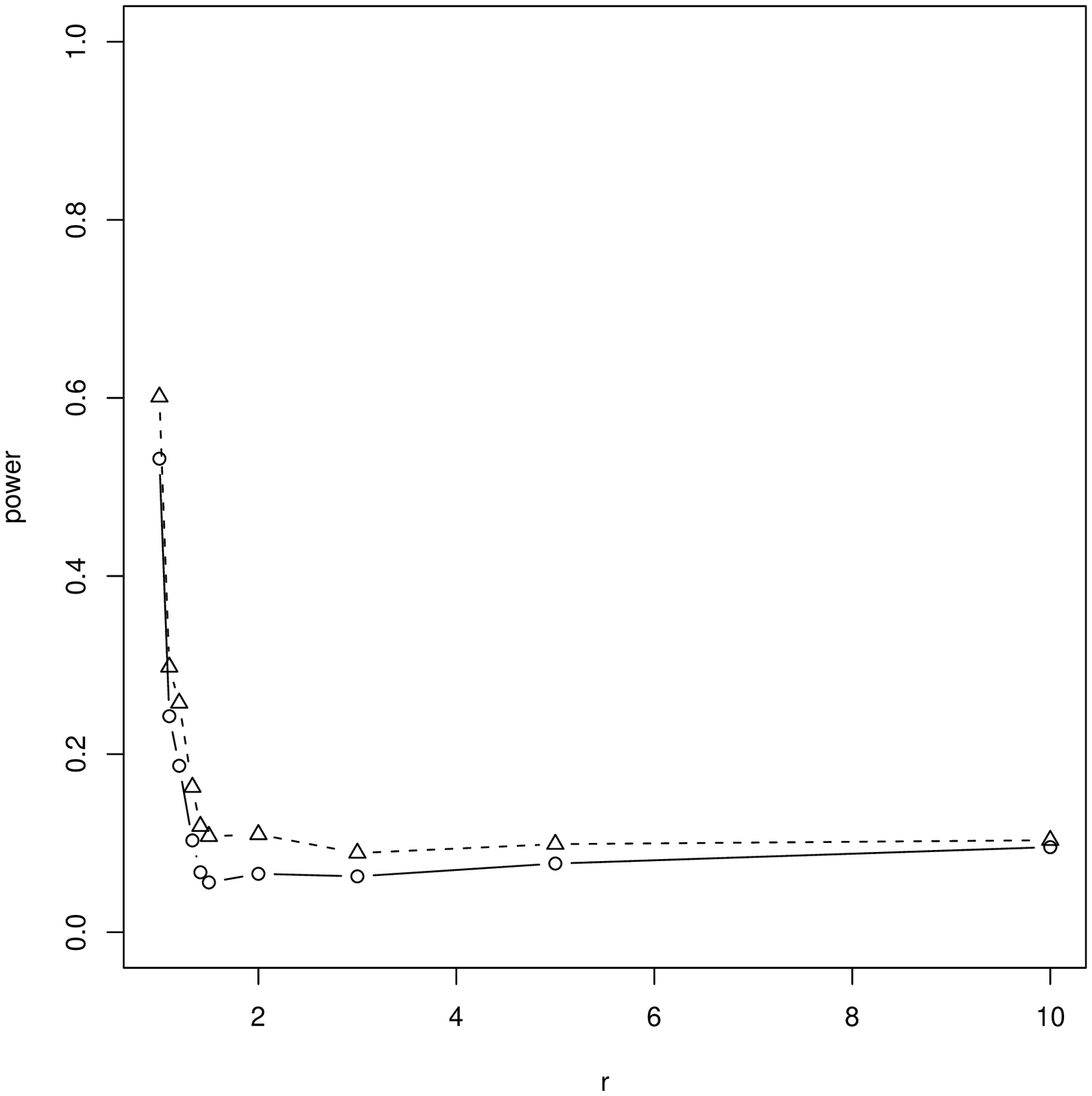}}}
\rotatebox{-90}{ \resizebox{1.7 in}{!}{ \includegraphics{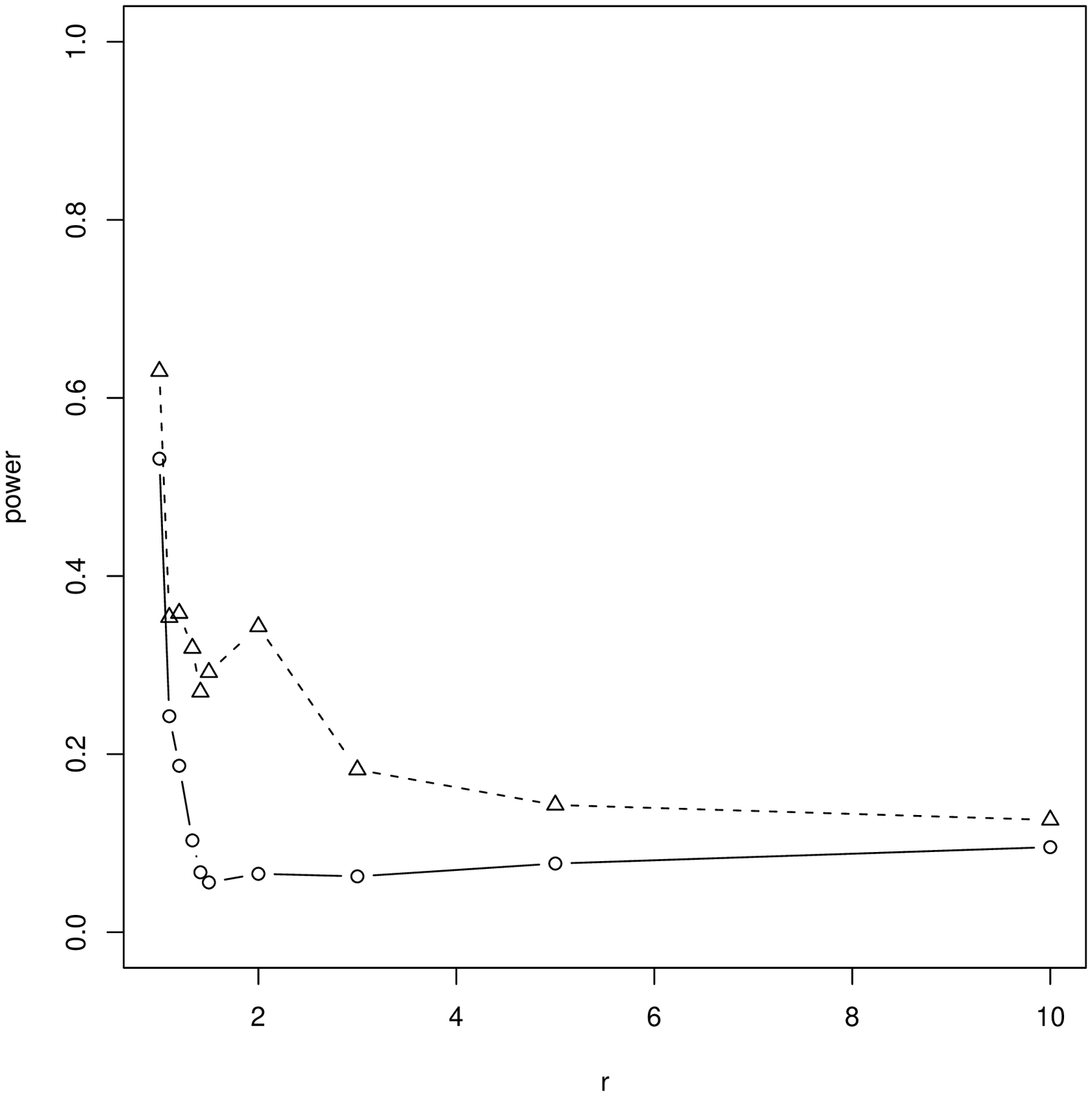}}}
\rotatebox{-90}{ \resizebox{1.7 in}{!}{ \includegraphics{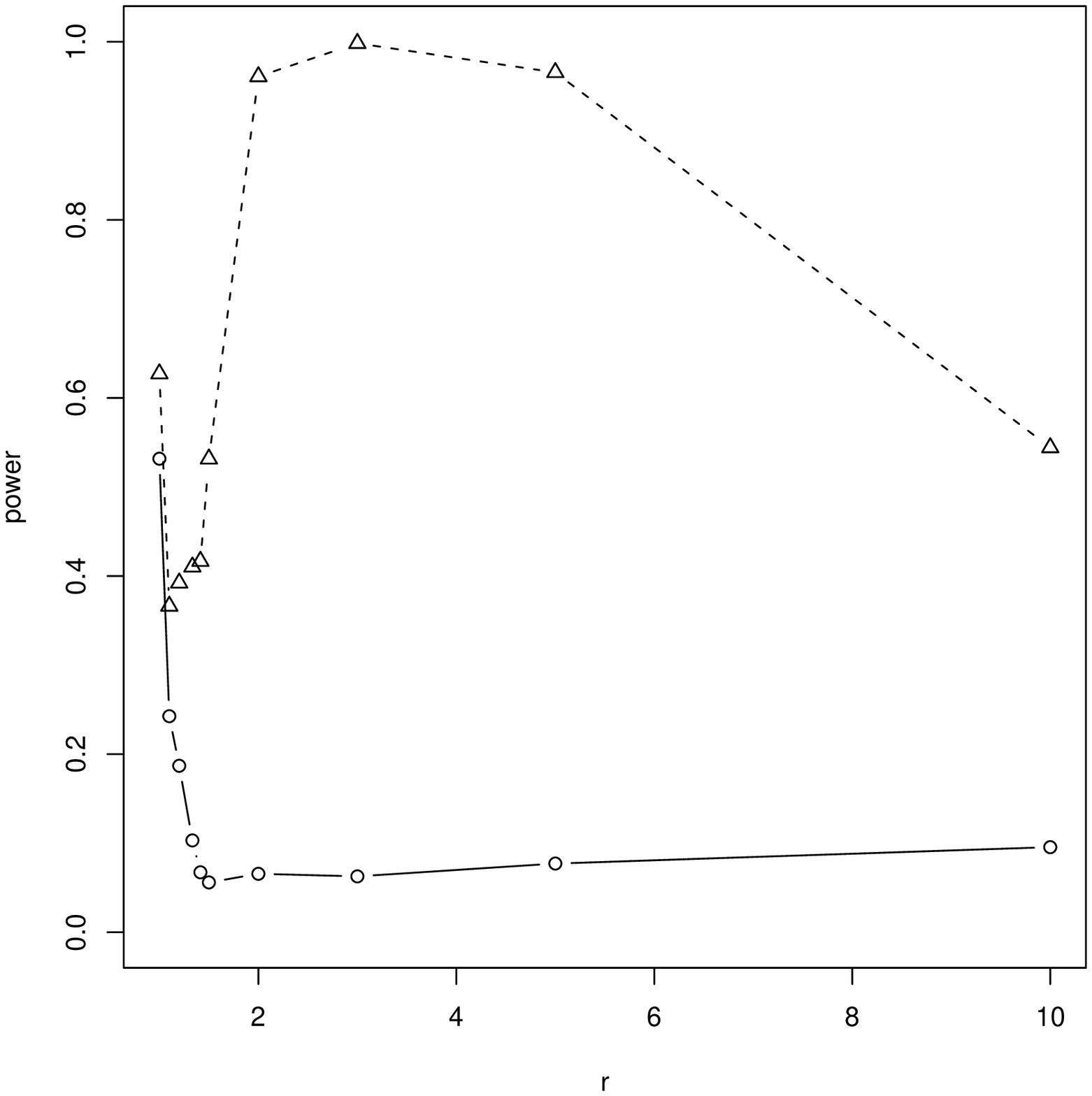}}}
\caption{ \label{fig:AggSimPowerCurve}
Monte Carlo power using the asymptotic critical value against association alternatives
$H^A_{\sqrt{3}/21}$ (left),
$H^A_{\sqrt{3}/12}$ (middle),
and
$H^A_{5\,\sqrt{3}/24}$ (right)
as a function of $r$, for $n=10$.  The circles represent the empirical significance levels while triangles represent the empirical power values.}
\end{figure}

In Figure \ref{fig:AggSimPowerCurve}, we present
a Monte Carlo investigation of power against
$H^A_{\sqrt{3}/21}$, $H^A_{\sqrt{3}/12}$,
and
$H^A_{5\,\sqrt{3}/24}$
as a function of $r$ for $n=10$.  The empirical significance level is $\widehat{\alpha}_A(n):=\frac{1}{N}\sum_{j=1}^{N}\I(R_j< z_{\alpha}|H_0)$. Then $\widehat{\alpha}_A(10)$, is about $.05$ for $r=\sqrt{2},3/2,\,2,\,3,\,5$ which have the empirical power $\widehat{\beta}^A_{10}\left( r,\sqrt{3}/12 \right) \le .35$ with maximum power at $r=2$, and $\widehat{\beta}^A_{10}\left( r=3,5\,\sqrt{3}/24 \right)=1$.  So, for small sample sizes, moderate values of $r$ are more appropriate for normal approximation, as they yield the desired significance level, and the more severe the association, the higher the power estimate.

The empirical significance levels and empirical power $\widehat{\beta}^S_n(r,\epsilon)$ values under $H^A_{\epsilon}$ for $\epsilon=5\,\sqrt{3}/24,\,\sqrt{3}/12,\sqrt{3}/21$ are presented in Table \ref{tab:asy-emp-val-A}.
\begin{table}[]
\centering
\begin{tabular}{|c|c|c|c|c|c|c|c|c|c|c|}
\hline
$r$  & 1 & 11/10 &6/5 & 4/3 & $\sqrt{2}$ &3/2 & 2 & 3 & 5 &10\\
\hline
$\widehat{\alpha}_A(10)$ & .5318 & .2426 & .1869 & .1031 & .0673 & .0559 & .0656 & .0627 & .0771 & .0955\\
\hline
$\widehat{\beta}^A_{10}\left( r,5\,\sqrt{3}/24 \right)$ & .6273 & .3663 & .3923 & .4103 & .4167 & .5316 & .9610 & .9983 & .9656 & .5443 \\
\hline
$\widehat{\beta}^A_{10}\left( r,\sqrt{3}/12 \right)$ & .6300 & .3537 & .3583 & .3190 & .2698 & .2919 & .3433 & .1825 & .1429 & .1261 \\
\hline
$\widehat{\beta}^A_{10}\left( r,\sqrt{3}/21 \right)$ & .6012 & .2979 & .2574 & .1629 & .1190 & .1077 & .1098 & .0889 & .0989 & .1033 \\
\hline
\end{tabular}
\caption{
\label{tab:asy-emp-val-A}
The empirical significance level and empirical power values under $H^A_{\epsilon}$ for $\epsilon=5\,\sqrt{3}/24,\,\sqrt{3}/12,\sqrt{3}/21$ with $N=10,000$, and $n=10$ at $\alpha=.05$.}
\end{table}

Note that even for $n=10$, the plots of the empirical power $\widehat{\beta}^A_{10}(r)$ resembles the curves of the asymptotic power function $\Pi_A(r)$ in Section \ref{sec:APF}.

\subsection{Pitman Asymptotic Efficacy}
\label{sec:Pitman}
Suppose that the distribution $F $ under consideration may be indexed by a set $\Theta \subset \R$ and consider $H_0:\theta=\theta_0 $ versus $H_a: \theta>\theta_0 $.

Pitman asymptotic efficacy (PAE)
provides for an investigation of ``local asymptotic power''
--- local around $H_0$.
This involves the limit as $n \rightarrow \infty$ as well as
the limit as $\epsilon \rightarrow 0$.

Consider the comparison of test sequences $S=\bigl\{ S_n \bigr\}$ satisfying the following conditions in a neighborhood $\theta \in [\theta_0,\theta_0+\delta] $ of the null parameter for some $\delta>0$.

\noindent
{\bfseries Pitman's Conditions:}
\begin{itemize}
\item[(PC1)] For some functions $\mu_n(\theta)$ and $\sigma_n(\theta)$, the distribution $F_{\theta}$ of $\bigl[S_n-\mu_n(\theta)\bigr]/\sigma_n(\theta)$ converges to $Z \sim \N(0,1)$ uniformly on $\bigl[ \theta_0,\theta_0+\delta \bigr] $, i.e.,
 $$\sup_{~~~~~~~~~\theta_0\le \theta \le \theta_0+\delta}\;\;\sup_{t \in \R} \left| P \left( \frac{S_n-\mu_n(\theta)}{\sigma_n(\theta)}\le t \right)-\Phi(t) \right|\rightarrow 0 \text{ as } n \rightarrow \infty.$$
\item[(PC2)] For $\theta \in [\theta_0,\theta_0+\delta] $, $\mu_n(\theta)$ is differentiable with $\mu_n'(\theta_0)>0$,
\item[(PC3)] For $\theta_n=\theta_0 + O\left( n^{-1/2} \right)$,   $\lim_{n\rightarrow \infty}\frac{\mu_n'(\theta_n)}{\mu_n'(\theta_0)}=1$,
\item[(PC4)] For $\theta_n=\theta_0 + O\left( n^{-1/2} \right)$,   $\lim_{n\rightarrow \infty}\frac{\sigma_n(\theta_n)}{\sigma_n(\theta_0)}=1$.
\item[(PC5)] For some constant $c>0$,
 $$\lim_{n\rightarrow \infty}\frac{\mu_n'(\theta_0)}{\sqrt{n}\,\sigma_n(\theta_0)}=c, $$
\end{itemize}
Condition (PC1) is equivalent to
\begin{itemize}
\item[(PC1)$^{\prime}$] For some functions $\mu_n(\theta)$ and $\sigma_n(\theta)$, the distribution $F_{\theta}$ of $\bigl[S_n-\mu_n(\theta_n)\bigr]/\sigma_n(\theta_n)$
converges to a standard normal distribution (see \cite{eeden:1963}).
\end{itemize}
Note that if $\mu_n^{(k)}(\theta_0)>0$ and $\mu_n^{(l)}(\theta_0)=0$, for all $l=1,\,2,\ldots,k-1$, then $\mu_n^{\prime}(\theta_0)$ in (PC2), (PC3), and (PC5) can be replaced by $\mu_n^{(k)}(\theta_0)>0$ and
$\mu_n^{\prime}(\theta_n)$ in (PC3) can be replaced by $\mu_n^{(k)}(\theta_n)$ (see \cite{kendall:1979}).

{\bf Lemma 1:}
(Pitman-Noether)
\begin{itemize}
\item[(i)] Let $S=\bigl\{ S_n \bigr\}$ satisfy (PC1)-(PC5). Consider testing $H_0 $ by the critical regions $S_n > u_{\alpha_n}$ with $\alpha_n=P_{\theta_0}\bigl(S_n >u_{\alpha_n}\bigr)\rightarrow \alpha $ as $n \rightarrow \infty $ where $\alpha \in (0,1)$.  For $\beta \in (0,1-\alpha)$ and $\theta_n=\theta_0 + O\left( n^{-1/2} \right)$, we have
 $$\beta_n(\theta_n)=P_{\theta_n}\bigl( T_n > u_{\alpha_n} \bigr)\rightarrow \beta \text{ iff } c\,\sqrt{n}\bigl( \theta_n-\theta \bigr)\rightarrow \Phi^{-1}(1-\alpha)-\Phi^{-1}(\beta).$$
\item[(ii)] Let $S=\bigl\{ S_n \bigr\}$ and $Q=\bigl\{ Q_n \bigr\}$ each satisfy satisfy {\em (PC1)-(PC5)}. Then the asymptotic relative efficiency of $S $ relative to $Q $ is given by $ARE(S,Q)=\left(c_S/c_Q \right)^2$.
\end{itemize}

Thus, to evaluate $ARE(S,Q)$ under the conditions (PC1)-(PC5), we need only calculate the quantities $c_S $ and $c_Q $, where
 $$c_S=\lim_{n \rightarrow \infty}\frac{\mu'_{S_n}(\theta_0)}{\sqrt{n}\cdot \sigma_{S_n}(\theta_0)} \text{ and } c_Q=\lim_{n \rightarrow \infty}\frac{\mu'_{Q_n}(\theta_0)}{\sqrt{n}\cdot \sigma_{Q_n}(\theta_0)}$$
 $PAE(S)=c_S^2$ is called the {\em Pitman Asymptotic Efficacy} (PAE) of the test based on $S_n $. Using similar notation and terminology for $Q_n $,
 $$ARE(S,Q)=\frac{\PAE(S)}{\PAE(Q)}.$$

For segregation or association alternatives
the PAE of $\rho_n(r)$ is given by $\PAE(r) = \frac{\left( \mu^{(k)}(r,\epsilon=0) \right)^2}{\nu(r)}$ where $k$ is the minimum order of the derivative with respect to $\epsilon$ for which $\mu^{(k)}(r,\epsilon=0) \not= 0$.  That is, $\mu^{(k)}(r,\epsilon=0) \not=0$ but $\mu^{(l)}(r,\epsilon=0)=0$ for $l=1,2,\ldots,k-1$.

\subsubsection{Pitman Asymptotic Efficacy Under Segregation Alternatives}

Consider the test sequences $\rho(r)=\bigl\{ \rho_n(r) \bigr\}$ for sufficiently small $\epsilon>0$ and $r \in \bigl[ 1,\sqrt{3}/(2\,\epsilon) \bigr)$.

In the PAE framework above, $\theta=\epsilon$ and $\theta_0=0$.  Suppose, $\mu_n(\epsilon)=E^S_{\epsilon}[\rho_n(r)]=\mu_S(r,\epsilon)$. For $\epsilon \in \bigl[ 0,\sqrt{3}/8 \bigr)$,
$$\mu_S(r,\epsilon)=\sum_{j=1}^5 \varpi_{1,j}(r,\epsilon)\,\I(r \in \mI_j)$$
with the corresponding intervals $\mI_1=\Bigl[ 1,3/2-\sqrt{3}\,\epsilon \Bigr)$, $\mI_2=\Bigl[3/2-\sqrt{3}\,\epsilon,3/2 \Bigr)$, $\mI_3=\Bigl[ 3/2,2-4\,\,\epsilon/\sqrt{3} \Bigr)$, $\mI_4=\Bigl[ 2-4\,\,\epsilon/\sqrt{3},2 \Bigr)$, $\mI_5=\Bigl[ 2,\sqrt{3}/(2\,\epsilon) \Bigr)$.  See Appendix 2 for the explicit form of $\mu(r,\epsilon)$ and Appendix 3 for derivation. Notice that as $\epsilon \rightarrow 0$, only $\mI_1=\Bigl[ 1,3/2-\sqrt{3}\,\epsilon \Bigr)$, $\mI_3=\Bigl[ 3/2,2-4\,\,\epsilon/\sqrt{3} \Bigr)$, $\mI_5=\Bigl[ 2,\sqrt{3}/(2\,\epsilon)\Bigr)$ do not vanish, so we only keep the components of $\mu_S(r,\epsilon)$ on these intervals.

Furthermore, $\sigma_S^2(n,\epsilon)=\Var^S_{\epsilon}(\rho_n(r))=\frac{1}{2\,n\,(n-1)}\Var^S_{\epsilon}[h_{12}]+\frac{(n-2)}{n\,(n-1)}$, $\nu_S(r,\epsilon)=\Cov^S_{\epsilon}[h_{12},h_{13}].$  The explicit forms of $\Var^S_{\epsilon}[h_{12}]$ and $\Cov^S_{\epsilon}[h_{12},h_{13}]$  are not calculated, since we only need $\lim_{n\rightarrow \infty}\sigma_n^2(\epsilon=0)=\nu(r)$ which is given in Equation \ref{eq:Asyvar}.

Notice that $\E^S_{\epsilon}|h_{12}|^3 \le 8 < \infty $ and
$\E^S_{\epsilon}[h_{12}\,h_{13}]-\E^S_{\epsilon}[h_{12}]^2=\Cov^S_{\epsilon}[h_{12},h_{13}]>0$
then by \cite{callaert:1978}
 $$\sup\text{}_{t\in \R} \left| P_{\epsilon} \left( \sqrt{n}\frac{\bigl( \rho_n(r)-\mu_S(r,{\epsilon}) \bigr)}{\sqrt{\nu_S(r,\epsilon)}}\le t \right)-\Phi(t)\right| \le C\,\E^S_{\epsilon}\left| h_{12} \right|^3\, \left[ \nu_S(r,\epsilon) \right]^{-\frac{3}{2}}\,n^{-\frac{1}{2}}$$
where $C $ is an absolute constant and $\Phi(\cdot)$ is the standard normal distribution function.  Then (PC1) follows for each $r \in \Bigl[ 1,\sqrt{3}/(2\,\epsilon) \Bigr)$ and $\epsilon \in \Bigl[0,\sqrt{3}/4 \Bigr)$.

Differentiating $\mu_S(r,\epsilon)$ with respect to $\epsilon$ yields
\begin{multline*}
 \mu_S^{\prime}(r,\epsilon)= \varpi_{1,1}^{\prime}(r,\epsilon)\,\I\left(r \in \bigl[ 1,3/2-\sqrt{3}\,\epsilon \bigr)\right)+\varpi_{1,3}^{\prime}(r,\epsilon)\,\I\left( r \in [3/2,2-4\,\,\epsilon/\sqrt{3}) \right)\\
+\varpi_{1,5}^{\prime}(r,\epsilon)\,\I\left( r \in \bigl[ 2,\sqrt{3}/(2\,\epsilon)\bigr) \right)
\end{multline*}
where
\begin{align*}
\varpi_{1,1}^{\prime}(r,\epsilon)&=\frac{2\,\epsilon\,(144\,{\epsilon}^2\,(r^2-1)+36-37\,r^2)}{27\,(2\,\epsilon-1)^3(2\,\epsilon+1)^3}\\
\varpi_{1,3}^{\prime}(r,\epsilon)&= \Bigl[2\,\sqrt{3}\Bigl( (2\,r-3)\,64\,\,{\epsilon}^3+(7\,r^2+r^4-24\,\,r+20)\,16\,\sqrt{3}{\epsilon}^2+(r-3)\,48\,\epsilon+3\,\sqrt{3}\,r^4+96\,\sqrt{3}\,r\\
&-36\,\sqrt{3}-60\,\sqrt{3}\,r^2\Bigr)\epsilon\Bigr]/\Bigl[9\,(2\,\epsilon+1)^3(2\,\epsilon-1)^3r^2 \Bigr]\\
\varpi_{1,5}^{\prime}(r,\epsilon)&=\frac{8\,\sqrt{3}\,\epsilon\,\left( 48\,{\epsilon}^3+(3\,r^4+3\,r^2-20)\,4\,\,\sqrt{3}\,{\epsilon}^2+36\,\epsilon+9\,\sqrt{3}-9\,\sqrt{3}\,r^2 \right)}{27\,r^2(2\,\epsilon+1)^3(2\,\epsilon-1)^3}.
\end{align*}

Hence, $\mu_S^{\prime}(r,\epsilon=0)=0$, so we need higher order derivatives for (PC2).
A detailed discussion is available in \cite{kendall:1979}.

Differentiating  $\mu_S^{\prime}(r,\epsilon)$ with respect to $\epsilon$  yields
\begin{multline*}
 \mu_S^{\prime\prime}(r,\epsilon)= \varpi_{1,1}^{\prime\prime}(r,\epsilon)\,\I\left(r \in \bigl[ 1,3/2-\sqrt{3}\,\epsilon \bigr)\right)+\varpi_{1,3}^{\prime\prime}(r,\epsilon)\,\I\left( r \in [3/2,2-4\,\,\epsilon/\sqrt{3}) \right)\\
+\varpi_{1,5}^{\prime\prime}(r,\epsilon)\,\I\left( r \in \bigl[ 2,\sqrt{3}/(2\,\epsilon)\bigr) \right)
\end{multline*}
where
\begin{align*}
\varpi_{1,1}^{\prime\prime}(r,\epsilon)&=-\frac{2\,(r^2-1)\,1728\,{\epsilon}^4+(72-77\,r^2)\,4\,\,{\epsilon}^2+36-37\,r^2}{27\,(4\,\,{\epsilon}^2-1)^4}\\
\varpi_{1,3}^{\prime\prime}(r,\epsilon)&= -2\,\Bigl[(2\,r-3)\,512\,\sqrt{3}\,{\epsilon}^5+(20+r^4+7\,r^2-24\,\,r)\,576\,{\epsilon}^4+(2\,r-3)\,1024\,\,\sqrt{3}\,{\epsilon}^3+(20-108\,r^2\\
&+96\,r+9\,r^4)\,36\,{\epsilon}^2+(-3+2\,r)\,96\,\sqrt{3}\,\epsilon-108+9\,r^4-180\,r^2+288\,r\Bigr]/\Bigl[9\,r^2(2\,\epsilon+1)^4(2\,\epsilon-1)^4\,\Bigr]\\
\varpi_{1,5}^{\prime\prime}(r,\epsilon)&= -8\,\Bigl[128\,\sqrt{3}\,{\epsilon}^5+(-20+3\,r^4+3\,r^2)\,48\,{\epsilon}^4+256\,\sqrt{3}\,{\epsilon}^3+(-5-12\,r^2+3\,r^4)\,12\,{\epsilon}^2\\
&+24\,\,\epsilon\,\sqrt{3}+9-9\,r^2\Bigr]/\Bigl[9\,r^2(2\,\epsilon+1)^4(2\,\epsilon-1)^4\,\Bigr].
\end{align*}
Thus,
\begin{equation}
\label{eq:PAE-S-''}
 \mu_S^{\prime\prime}(r,\epsilon=0)=
\begin{cases}
          -\frac{8}{3}+{\frac{74}{27}}\,r^2 &\text{for} \quad r \in [1,3/2)\\
          -2\,{\frac{(r^2-4\,\,r+2)(r^2+4\,\,r-6)}{r^2}} &\text{for} \quad r \in [3/2,2)\\
          -\frac{8\,(1-r^2)}{r^2} &\text{for} \quad r \in [2,\sqrt{3}/(2\,\epsilon)).
\end{cases}
\end{equation}
Observe that $\mu_S^{\prime\prime}(r,\epsilon=0)>0$ for all $r \in \Bigl[ 1,\sqrt{3}/(2\,\epsilon)\Bigr)$, so (PC2) holds with the second derivative. (PC3) in the second derivative form follows from continuity of $\mu_S^{\prime\prime}(r,\epsilon)$ in $\epsilon$ and (PC4) follows from continuity of $\sigma_n^2(r,\epsilon)$ in $\epsilon$.

Next, we find $c_S(\rho(r))=\lim_{n\rightarrow \infty}\frac{\mu_S^{\prime\prime}(r,\epsilon=0)}{\sqrt{n}\,\sigma_n(r,\epsilon=0)}=\frac{\mu_S^{\prime\prime}(r,\epsilon=0)}{\sqrt{\nu(r)}}$, where numerator is given in Equation \ref{eq:PAE-S-''} and denominator is given in Equation \ref{eq:Asyvar}.  We can easily see that $c_S(\rho(r))>0$, since $c_S(\rho(r))$ is increasing in $r$ and $c_S(\rho(r=1))>0$. Then (PC5) follows.
So under segregation alternatives $H^S_{\epsilon}$, the PAE of $\rho_n(r)$ is given by
$$
\PAE^S(r) =c^2_S(\rho(r))=
   \frac{\left( \mu_S^{\prime\prime}(r,\epsilon=0) \right)^2}{\nu(r)}.
$$

\begin{figure}[ht]
\centering
\psfrag{r}{\scriptsize{$r$}}
\psfrag{pS}{\scriptsize{$\PAE^S(r)$}}
\epsfig{figure=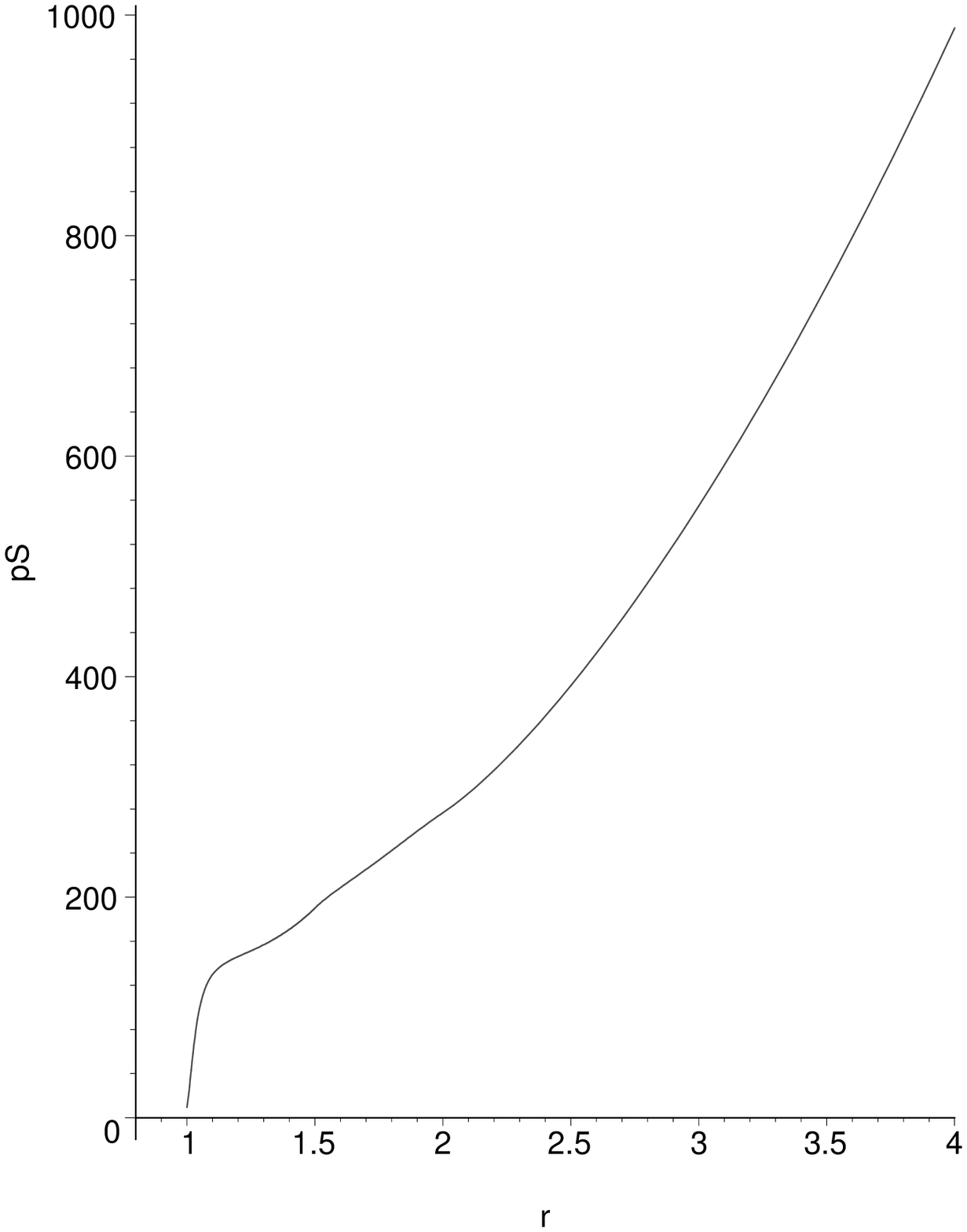, height=180pt, width=180pt}
\psfrag{pA}{\scriptsize{$\PAE^A(r)$}}
\epsfig{figure=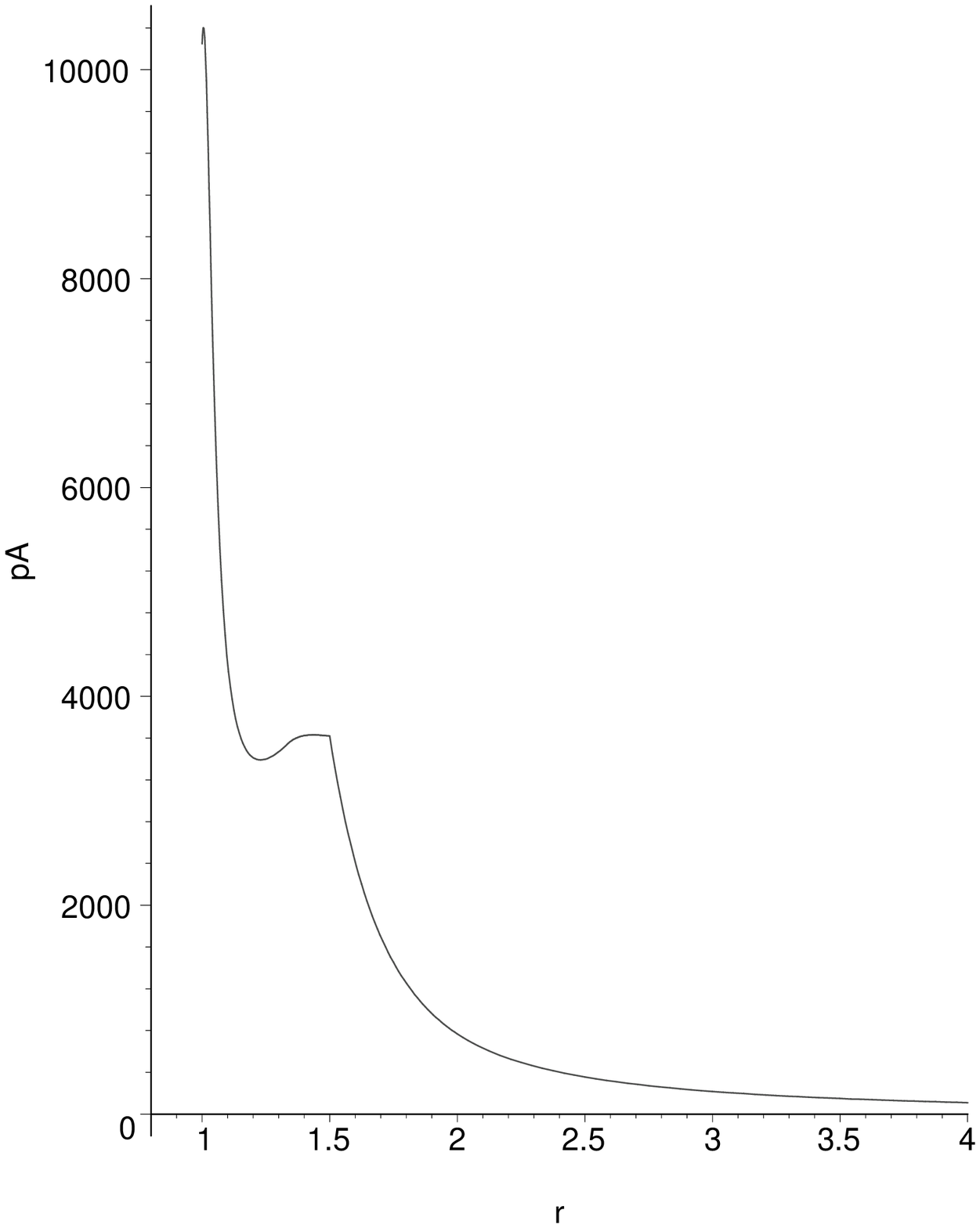, height=180pt, width=180pt}
\caption{
\label{fig:PAE-seg-agg}
Pitman asymptotic efficacy against segregation (left) and against association (right) as a function of $r$.
}
\end{figure}

In Figure \ref{fig:PAE-seg-agg} (left), we present the PAE as a function of $r$ for segregation. Notice that
$\PAE^S(r=1) = 160/7 \approx 22.8571$, $\lim_{r \rightarrow \infty} \PAE^S(r) = \infty$. Based on the PAE analysis, we suggest,
for large $n$ and small $\epsilon$,
choosing $r$ large for testing against segregation.  However, for small and moderate values of $n$, normal approximation is not appropriate due to the skewness in the density of $\rho_n(r)$, Therefore, for small $n$, we suggest moderate $r$ values.

PAE analysis is local (around $\epsilon=0$) and for arbitrarily large $n$. The comparison would hold in general provided that $\mu(r,\epsilon)$ is convex in $\epsilon$ for all $\epsilon \in \bigl[ 0,\sqrt{3}/3 \bigr)$.  As an alternative, we fix an $\epsilon$ and then compare the asymptotic behaviour of $\rho_n(r)$ with Hodges-Lehmann asymptotic efficacy in Section \ref{sec:Hodges-Lehmann-Seg}.

\subsubsection{Pitman Asymptotic Efficacy Under Association Alternatives}
Consider the test sequences $\rho(r)=\bigl\{ \rho_n(r) \bigr\}$ for sufficiently small $\epsilon >0$ and $r \in [1,\infty)$.

In the PAE framework above, $\theta=\epsilon$ and $\theta_0=0$.  Suppose, $\mu_n(\epsilon)=\E_{\epsilon}[\rho_n(r)]=\mu_A(r,\epsilon)$. For $\epsilon \in \Bigl[0,\left(7\,\sqrt{3}-3\,\sqrt{15}\right)/12 \approx .042 \Bigr)$,
$$\mu_A(r,\epsilon)=\sum_{j=1}^6 \varpi_{1,j}(r,\epsilon)\,\I(r \in \mI_j)$$
with the corresponding intervals $\mI_1=\Bigl[1,\left( 1+2\,\sqrt{3}\,\epsilon \right)/\left( 1-\sqrt{3}\,\epsilon \right)\Bigr)$, $\mI_2=\Bigl[\left( 1+2\,\sqrt{3}\,\epsilon \right)/\left( 1-\sqrt{3}\,\epsilon \right),\\
4\,\,\left( 1-\sqrt{3}\,\epsilon \right)/3 \Bigr)$, $\mI_3=\Bigl[ 4\,\,\left( 1-\sqrt{3}\,\epsilon \right)/3,4\,\,\left( 1+2\,\sqrt{3}\,\epsilon \right)/3 \Bigr)$, $\mI_4=\Bigl[ 4\,\,\left( 1+2\,\sqrt{3}\,\epsilon \right)/3,3/(2\,\left( 1-\sqrt{3}\,\epsilon \right))\Bigr)$, $\mI_5=\Bigl[3/(2\,\left( 1-\sqrt{3}\,\epsilon \right)),2 \Bigr)$ and $\mI_6=[2,\infty)$.  Notice that as $\epsilon \rightarrow 0$, only $\mI_j$ for $j=2,4,5,6$ do not vanish, so we only keep the components of $\mu_A(r,\epsilon)$ on these intervals. See Section \ref{sec:arc-prob-agg} for the explicit form of $\mu_A(r,\epsilon)$ and Section \ref{sec:derivation-par-eps} for derivation.

Furthermore, $\sigma_n^2(\epsilon)=\Var^A_{\epsilon}(\rho_n(r))=\frac{1}{2\,n\,(n-1)}\Var^A_{\epsilon}[h_{12}]+\frac{(n-2)}{n\,(n-1)} \, \Cov^A_{\epsilon}[h_{12},h_{13}] $ whose explicit form is not calculated, since we only need $\lim _{n \rightarrow \infty}\sqrt{n}\, \sigma_n(\epsilon=0)=\nu(r)$ which is given Equation \ref{eq:Asyvar}.

(PC1) follows for each $r \in [1,\infty)$ and $\epsilon \in \bigl[ 0,\sqrt{3}/3 \bigr)$ as in the segregation case.

Differentiating $\mu_A(r,\epsilon)$ with respect to $\epsilon$, then we get
\begin{multline*}
\mu_A^{\prime}(r,\epsilon)=\varpi_{1,2}^{\prime}(r,\epsilon)\,\I(r \in [1,4/3))+\varpi_{1,4}^{\prime}(r,\epsilon)\,\I(r \in [4/3,3/2))\\
+\varpi_{1,5}^{\prime}(r,\epsilon)\,\I(r \in [3/2,2))+\varpi_{1,6}^{\prime}(r,\epsilon)\,\I(r \in [2,\infty))
\end{multline*}
where
\begin{align*}
\varpi_{1,2}^{\prime}(r,\epsilon)&=-2\,\Bigl[\sqrt{3}\Bigl(-1152\,r^4{\epsilon}^3+720\,\sqrt{3}\,r^4{\epsilon}^2-288\,r^4\,\epsilon+11\,\sqrt{3}\,r^4+2592\,\sqrt{3}\,r^2{\epsilon}^2-10368\,\sqrt{3}\,r{\epsilon}^2\\
&+432\,\sqrt{3}\,r^2+6480\,\sqrt{3}{\epsilon}^2-864\,\,\sqrt{3}\,r+432\,\sqrt{3}\Bigr)\epsilon\Bigr]/\Bigl[\left( -6\,\epsilon+\sqrt{3} \right)^3\left( 6\,\epsilon+\sqrt{3} \right)^3r^2\Bigr],\\
\varpi_{1,4}^{\prime}(r,\epsilon)&= -2\,\Bigl[\sqrt{3}\Bigl(-1152\,r^4{\epsilon}^3+720\,\sqrt{3}\,r^4{\epsilon}^2-288\,r^4\,\epsilon+11\,\sqrt{3}\,r^4-1296\,\sqrt{3}\,r^2{\epsilon}^2+108\,\sqrt{3}\,r^2\\
&-2160\,\sqrt{3}{\epsilon}^2-144\,\sqrt{3}\Bigr)\epsilon\Bigr]/\Bigl[\left( -6\,\epsilon+\sqrt{3} \right)^3\left( 6\,\epsilon+\sqrt{3} \right)^3r^2\Bigr],\\
\varpi_{1,5}^{\prime}(r,\epsilon)&=\frac{2\,\epsilon\,(3\,r^4-72\,r^2-240\,{\epsilon}^2+192\,r-124)}{r^2(12\,{\epsilon}^2-1)^3},\\
\varpi_{1,6}^{\prime}(r,\epsilon)&=-\frac{40\,\epsilon}{r^2(12\,{\epsilon}^2-1)^2}.
\end{align*}
Hence $\mu_A^{\prime}(r,\epsilon=0)=0$, so we differentiate $\mu_A^{\prime}(r,\epsilon)$ with respect to $\epsilon$ and get
\begin{multline*}
\mu_A^{\prime\prime}(r,\epsilon)=\varpi_{1,2}^{\prime\prime}(r,\epsilon)\,\I(r \in [1,4/3))+\varpi_{1,4}^{\prime\prime}(r,\epsilon)\,\I(r \in [4/3,3/2))\\
+\varpi_{1,5}^{\prime\prime}(r,\epsilon)\,\I(r \in [3/2,2))+\varpi_{1,6}^{\prime\prime}(r,\epsilon)\,\I(r \in [2,\infty))
\end{multline*}
where
\begin{align*}
\varpi_{1,2}^{\prime\prime}(r,\epsilon)&=-6\,\Bigl[\sqrt{3}\Bigl(-27648\,r^4{\epsilon}^5+25920\,\sqrt{3}\,r^4{\epsilon}^4-18432\,r^4{\epsilon}^3+2820\,\sqrt{3}\,r^4{\epsilon}^2+93312\,\sqrt{3}\,r^2{\epsilon}^4\\
&-576\,r^4\,\epsilon-373248\,\sqrt{3}\,r{\epsilon}^4+11\,\sqrt{3}\,r^4+33696\,\sqrt{3}\,r^2{\epsilon}^2+233280\,\sqrt{3}{\epsilon}^4-82944\,\sqrt{3}\,r{\epsilon}^2\\
&+432\,\sqrt{3}\,r^2+45360\,\sqrt{3}{\epsilon}^2-864\,\,\sqrt{3}\,r+432\,\sqrt{3}\Bigr)\Bigr]/\Bigl[\left( 6\,\epsilon+\sqrt{3} \right)^4\,\left( -6\,\epsilon+\sqrt{3} \right)^4r^2\Bigr],\\
\varpi_{1,4}^{\prime\prime}(r,\epsilon)&= -6\,\Bigl[\sqrt{3}\Bigl(-27648\,r^4{\epsilon}^5+25920\,\sqrt{3}\,r^4{\epsilon}^4-18432\,r^4{\epsilon}^3+2820\,\sqrt{3}\,r^4{\epsilon}^2-46656\,\sqrt{3}\,r^2{\epsilon}^4\\
&-576\,r^4\,\epsilon+11\,\sqrt{3}\,r^4+2592\,\sqrt{3}\,r^2{\epsilon}^2-77760\,\sqrt{3}{\epsilon}^4+108\,\sqrt{3}\,r^2-15120\,\sqrt{3}{\epsilon}^2-144\,\sqrt{3}\Bigr)\Bigr]\\
&/\Bigl[\left( 6\,\epsilon+\sqrt{3} \right)^4\,\left( -6\,\epsilon+\sqrt{3} \right)^4r^2\Bigr],\\
\varpi_{1,5}^{\prime\prime}(r,\epsilon)&=-\frac{2\,(180\,r^4{\epsilon}^2+3\,r^4-4320\,r^2{\epsilon}^2-8640\,{\epsilon}^4+11520\,r{\epsilon}^2-72\,r^2-8160\,{\epsilon}^2+192\,r-124)}{r^2(12\,{\epsilon}^2-1)^4},\\
\varpi_{1,6}^{\prime\prime}(r,\epsilon)&=\frac{40\,(36\,{\epsilon}^2+1)}{r^2(12\,{\epsilon}^2-1)^3}.
\end{align*}
Thus,
\begin{equation}
\label{eq:PAE-A-''}
\mu_A^{\prime\prime}(r,\epsilon=0)=
\begin{cases}
           -\frac{22}{9}\,r^2+192\,r^{-1}-96\,r^{-2}-96 &\text{for} \quad r \in [1,4/3)\\
           -\frac{22}{9}\,r^2+32\,r^{-2}-24 &\text{for} \quad r \in [4/3,3/2)\\
           -6\,r^2-384\,\,r^{-1}+248\,r^{-2}+144 &\text{for} \quad r \in [3/2,2)\\
            -40\,r^{-2} &\text{for} \quad r \in [2,\infty).
\end{cases}
\end{equation}
Note that $\mu_A^{\prime\prime}(r,\epsilon=0)>0$ for all $r \in [1,\infty)$, so (PC2) follows with the second derivative. (PC3) and (PC4) follow from continuity of $\mu_A^{\prime\prime}(r,\epsilon)$ and $\sigma_n^2(r,\epsilon)$ in $\epsilon$.

Next, we find $c_A(\rho(r))=\lim_{n\rightarrow \infty}\frac{\mu_A^{\prime\prime}(r,\epsilon=0)}{\sqrt{n}\,\sigma_n(r,\epsilon=0)}=\frac{\mu_A^{\prime\prime}(r,0)}{\sqrt{\nu(r)}}$, by substituting the numerator from Equation \ref{eq:PAE-A-''} and denominator from Equation \ref{eq:Asyvar}.  We can easily see that $c_A(\rho(r))<0$, for all $r \ge 1 $. Then (PC5) holds, so under association alternatives $H^A_{\epsilon}$, the PAE of $\rho_n(r)$ is
$$
\PAE^A(r) = c^2_A(\rho(r))=
   \frac{\bigl( \mu_A^{\prime\prime}(r,\epsilon=0) \bigr)^2}{\nu(r)}.
$$

In Figure \ref{fig:PAE-seg-agg} (right), we present the PAE as a function of $r$ for association. Notice that
$\PAE^A(r=1) = 174240/17 \approx 10249.4118$, $\lim_{r \rightarrow \infty} \PAE^A(r) = 0$,
$\argsup_{r \in [1,\infty)} \PAE^A(r) \approx 1.006$ with supremum $\approx 10399.7726$. $\PAE^A(r)$ has also a local supremum at $r_l\approx 1.4356$ with local supremum $\approx 3630.8932$.
Based on the Pitman asymptotic efficacy analysis, we suggest,
for large $n$ and small $\epsilon$,
choosing $r$ small for testing against association.  However, for small and moderate values of $n$ normal approximation is not appropriate due to the skewness in the density of $\rho_n(r)$. Therefore, for small $n$, we suggest moderate $r$ values.

We also calculate  Hodges-Lehmann asymptotic efficacy for fixed alternatives in Section \ref{sec:Hodges-Lehmann-Agg}.

\subsection{Hodges-Lehmann Asymptotic Efficacy}
\label{sec:Hodges-Lehmann}

Unlike PAE, HLAE does not involve the limit as $\epsilon \rightarrow 0$.
Since this requires the mean and, especially,
the asymptotic variance of $\rho_n(r)$ \emph{under an alternative},
we investigate HLAE for specific values of $\epsilon$. See Appendix 4 for a sample derivation of $\mu_S(r,\varepsilon)$ and $\nu_S(r,\varepsilon)$.

\subsubsection{Hodges-Lehmann Asymptotic Efficacy Under Segregation Alternatives}
\label{sec:Hodges-Lehmann-Seg}
In the HLAE framework, $\theta=\E^S_{\epsilon}[\rho_n(r)]=\mu_S(r,\epsilon)$ and $\theta_0=\mu(r)$.  Then testing $H_0:\epsilon=0$ versus $H^S_{\epsilon}:\epsilon>0$ is equivalent to $H_0:\E[\rho_n]=\mu(r)$ versus $H^S_{\epsilon}:\E^S_{\epsilon}[\rho_n]=\mu_S(r,\epsilon)>\mu(r)$.  Let $\delta=\frac{\mu_S(r,\epsilon)-\mu(r)}{\sqrt{\nu_S(r,\epsilon)}}$ and $\widetilde{R}=\frac{\rho_n(r)-\mu(r)}{\sqrt{\nu_S(r,\epsilon)}}$, then $\widetilde{R}\stackrel{\mathcal L}{\rightarrow}\N(\delta,1)$.

Then HLAE of $\rho_n(r)$ is given by
$$
\HLAE^S(r,\epsilon):=\frac{\bigl( \mu_S(r,\epsilon)-\mu(r))^2}{\nu_S(r,\epsilon\bigr)}.
$$
We calculate HLAE of $\rho_n(r)$ under $H^S_{\epsilon}$ for $\epsilon=\sqrt{3}/8$, $\epsilon=\sqrt{3}/4$, and $\epsilon=2\,\sqrt{3}/7$.

\begin{figure}[]
\centering
\scalebox{.3}{\input{ArcProbEps.pstex_t}}
\scalebox{.3}{\input{AsyVarEps.pstex_t}}
\caption{
\label{fig:mean-var-seg}
The mean $\mu_S(r,\epsilon)$ (left)
and asymptotic variance $\nu_S(r,\epsilon)$ (right) as a function of $r$ under segregation with $\epsilon=0,\sqrt{3}/8, \sqrt{3}/4,2\,\sqrt{3}/7$.}
\end{figure}
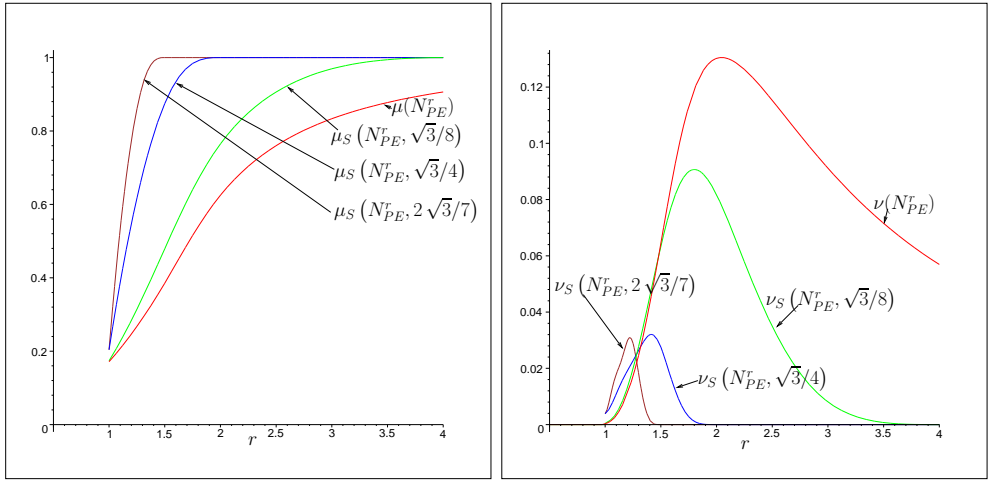

With $\epsilon=\sqrt{3}/8$, $\rho_n(r)$ is non-degenerate for $r \in [1,4)$
 $$ \mu_S\left( r,\sqrt{3}/8 \right)=
\begin{cases}
          {\frac{2287}{9126}}\,r^2-\frac{1}{13} &\text{for} \quad r \in [1,9/8)\\
          -\frac{5905\,r^4-36864\,\,r^3+62910\,r^2-46656\,r+13122}{9126\,r^2} &\text{for} \quad r \in  [9/8,3/2)\\
          \frac{61\,r^4-768\,r^3+3494\,\,r^2-5120\,r+2466}{338\,r^2} &\text{for} \quad r \in [3/2,2)\\
          -\frac{3\,r^4-422\,r^2+606}{338\,r^2} &\text{for} \quad r \in [2,3)\\
          \frac{3\,r^4-48\,r^3+530\,r^2-768}{338\,r^2} &\text{for} \quad r \in [2,4]
\end{cases}$$
and
$$\nu_S\left( r,\sqrt{3}/8 \right)=\sum_{j=1}^{12}\nu_j\left( r,\sqrt{3}/8 \right)\,\I(\mI_j)$$
where
{\small
\begin{align*}
\nu_1\left( r,\sqrt{3}/8 \right)&=\Bigl[9959911\,r^{10}-46006272\,r^9-430526\,r^8+258785280\,r^7-385799609\,r^6+162699264\,\,r^5\\
&-83976048\,r^4+201277440\,r^3-129392640\,r^2+12939264\Bigr]/\Bigl[104104845\,r^4\Bigr],\\
\nu_2\left( r,\sqrt{3}/8 \right)&=\Bigl[9959911\,r^{10}-46006272\,r^9-430526\,r^8+258785280\,r^7-415110891\,r^6 +272331072\,r^5\\
&-158725008\,r^4-16174080\,r^3+315394560\,r^2-310542336\,r+90574848\Bigr]/\Bigl[104104845\,r^4\Bigr],\\
\nu_3\left( r,\sqrt{3}/8 \right)&=\Bigl[3144167\,r^{12}+15335424\,\,r^{11}-378655166\,r^{10}+2750459904\,\,r^9-11800111467\,r^8\\
&+31878202752\,r^7-54792387144\,r^6+60339341664\,\,r^5-42745183272\,r^4+19903426272\,r^3\\
&-6790168926\,r^2+1989715104\,\,r-373071582\Bigr]/\Bigl[104104845\,r^6\Bigr],
\end{align*}
\begin{align*}
\nu_4\,\left( r,\sqrt{3}/8 \right)&=-\Bigl[8177689\,r^{12}-54153216\,r^{11}+320428478\,r^{10}-2459326464\,\,r^9+11854698987\,r^8\\
&-32751603072\,r^7+55010737224\,\,r^6-59029241184\,\,r^5+42131073672\,r^4-20886001632\,r^3\\
&+7379714142\,r^2-1694942496\,r+170415414\Bigr]/\Bigl[104104845\,r^6\Bigr],\\
\nu_5\left( r,\sqrt{3}/8 \right)&=-\Bigl[8177689\,r^{12}-54153216\,r^{11}+320428478\,r^{10}-2459326464\,\,r^9+12509010411\,r^8\\
&-37904305536\,r^7+71918042184\,\,r^6-88617024864\,\,r^5+71256548232\,r^4-36176875776\,r^3\\
&+10724592861\,r^2-1694942496\,r+170415414\Bigr]/\Bigl[104104845\,r^6\Bigr].\\
\nu_6\left( r,\sqrt{3}/8 \right)&=-\Bigl[2718937\,r^{12}-39596544\,r^{11}+434455742\,r^{10}-3154811904\,\,r^9+14086429683\,r^8\\
&-39680803584\,\,r^7+72881433288\,r^6-88893062496\,r^5+71547681672\,r^4\,\\
&-36487418112\,r^3+10828106973\,r^2-1694942496\,r+170415414\Bigr]/\Bigl[104104845\,r^6\Bigr],\\
\nu_7\left( r,\sqrt{3}/8 \right)&=-\Bigl[1027\,r^{12}-19968\,r^{11}+295626\,r^{10}-3265792\,r^9+23210081\,r^8\\
&-103077696\,r^7+289042360\,r^6-511170304\,\,r^5+553668600\,r^4-343186304\,\,r^3\\
&+109133095\,r^2-20431008\,r+5845554\Bigr]/\Bigl[428415\,r^6\Bigr],\\
\nu_8\left( r,\sqrt{3}/8 \right)&=-\Bigl[637\,r^{12}-19968\,r^{11}+299370\,r^{10}-3265792\,r^9+23199551\,r^8\\
&-103077696\,r^7+289042360\,r^6-511170304\,\,r^5+553700190\,r^4-343186304\,\,r^3\\
&+109133095\,r^2-20431008\,r+5788692\Bigr]/\Bigl[428415\,r^6\Bigr],\\
\nu_9\left( r,\sqrt{3}/8 \right)&=-\Bigl[637\,r^{12}-19968\,r^{11}+299370\,r^{10}-3265792\,r^9+24051519\,r^8\\
&-112023360\,r^7+328179640\,r^6-602490624\,\,r^5+673558110\,r^4-427086848\,r^3\\
&+133604087\,r^2-20431008\,r+5788692\Bigr]/\Bigl[428415\,r^6\Bigr],\\
\nu_{10}\left( r,\sqrt{3}/8 \right)&=\Bigl[130\,r^{12}-2496\,r^{11}+22134\,\,r^{10}-122720\,r^9+452225\,r^8-1010880\,r^7+1075400\,r^6\\
&+26624\,\,r^5-1993566\,r^4+5324800\,r^3-5083895\,r^2+303264\,\,r-37908\Bigr]/\Bigl[428415\,r^6\Bigr],\\
\nu_{11}\left( r,\sqrt{3}/8 \right)&=-\Bigl[330\,r^8-8896\,r^7+85445\,r^6-342624\,\,r^5+332000\,r^4+1148560\,r^3\\
&-1180986\,r^2-5324800\,r+6678947\Bigr]/\Bigl[428415\,r^4\Bigr],\\
\nu_{12}\left( r,\sqrt{3}/8 \right)&=-\frac{(330\,r^5-4936\,r^4+12453\,r^3+47388\,r^2-12992\,r-128256)(r-4)^3}{428415\,r^4},
\end{align*}
}
and the corresponding intervals are $\mI_1=[1,12/11),\; \mI_2=[12/11,9/8),\; \mI_3=\Bigl[9/8,\sqrt{6}/2\Bigr),\; \mI_4=\Bigl[ \sqrt{6}/2,21/16 \Bigr),\; \mI_5=[21/16,4/3),\; \mI_6=[4/3,3/2),\; \mI_7=\Bigl[3/2,\sqrt{3}\Bigr),\; \mI_8=\Bigl[\sqrt{3},7/4\Bigr),\; \mI_9=[7/4,2),\; \mI_{10}=[2,3),\; \mI_{11}=[3,7/2),\; \mI_{12}=[7/2,4)$. See Section \ref{sec:derivation-par-cov-1} for derivation and Figure \ref{fig:mean-var-seg} for the graph of $\mu\left( r,\sqrt{3}/8 \right)$ and $\nu\left( r,\sqrt{3}/8 \right)$.

\begin{figure}[]
\centering
\psfrag{r}{\scriptsize{$r$}}
\epsfig{figure=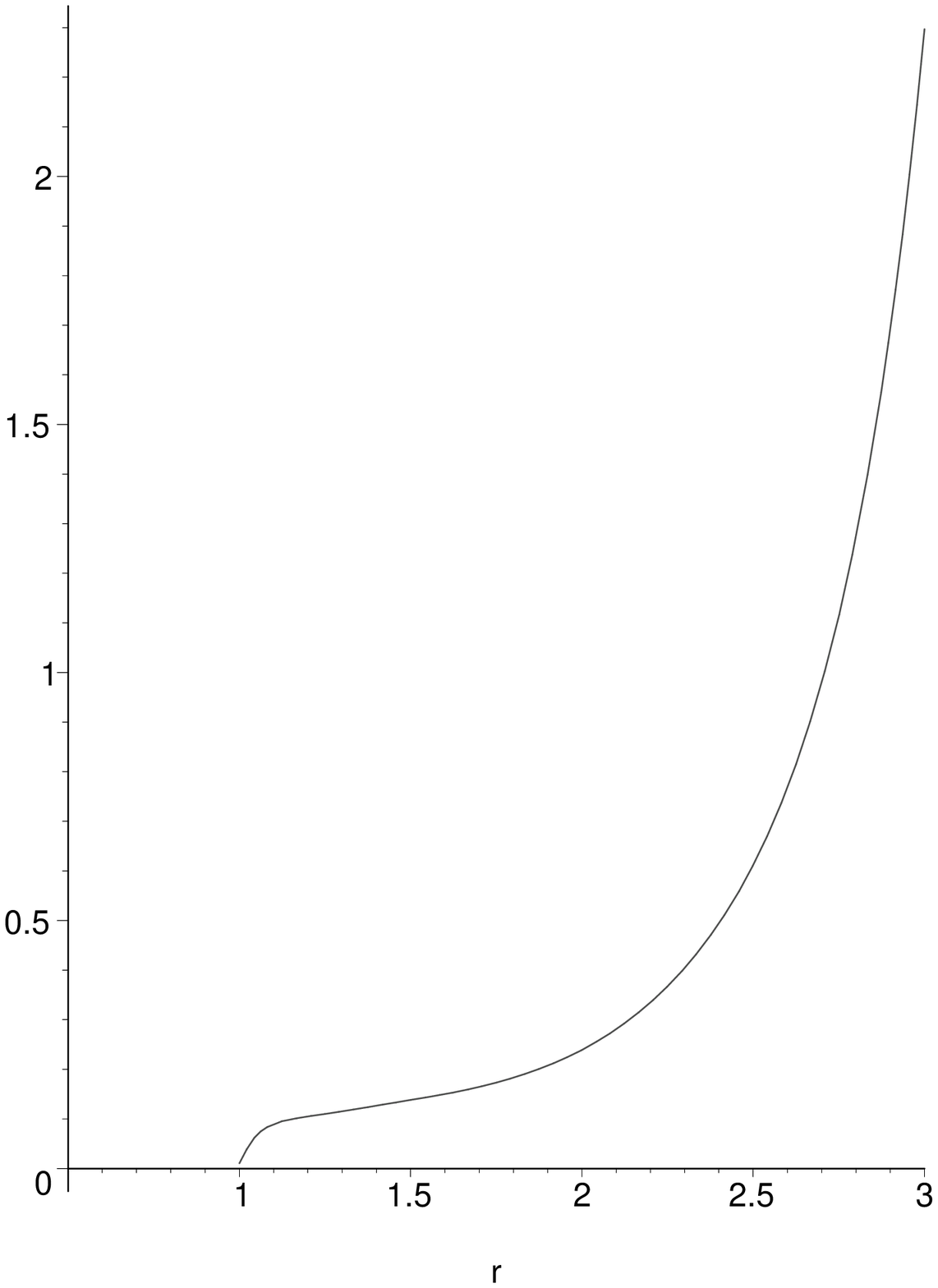, height=125pt, width=125pt}
\epsfig{figure=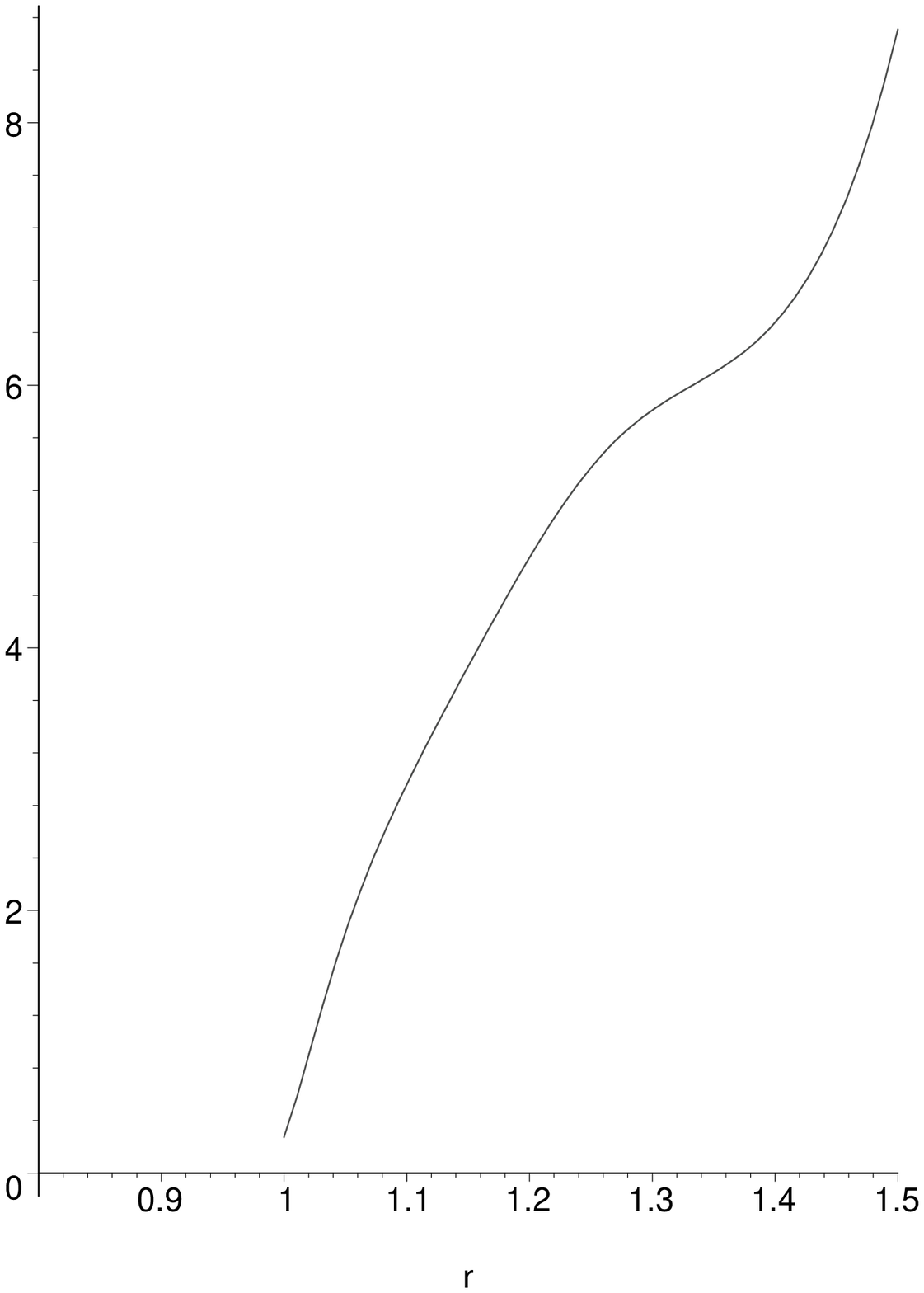, height=125pt, width=125pt}
\epsfig{figure=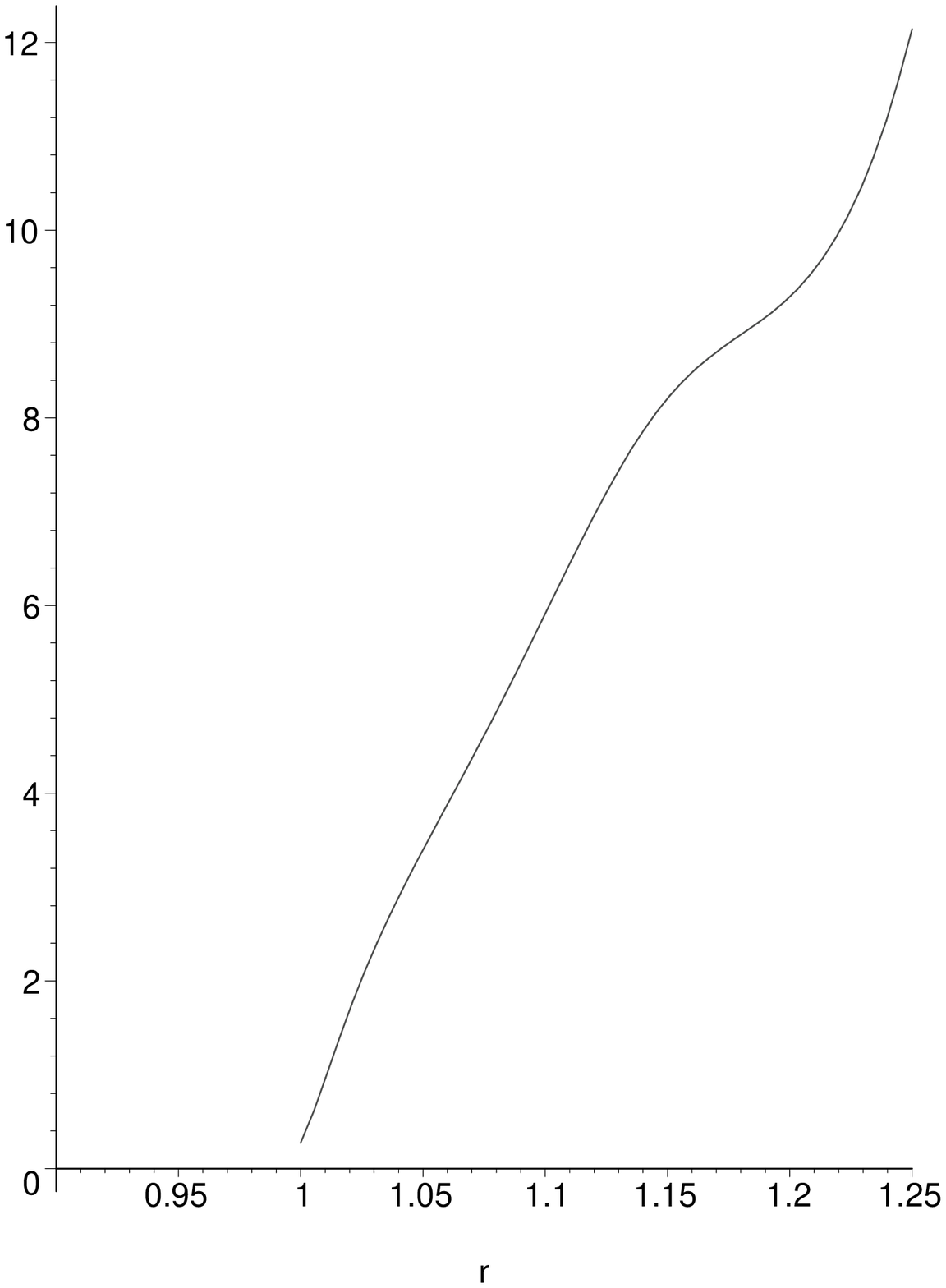, height=125pt, width=125pt}
\caption{ \label{fig:HLAEPlots for Seg}
Hodges-Lehmann asymptotic efficacy
against segregation alternative $H^S_{\epsilon}$
as a function of $r$
for $\epsilon= \sqrt{3}/8, \sqrt{3}/4, 2\,\sqrt{3}/7$ (left to right).}
\end{figure}

Then we get $\HLAE^S\left( r,\sqrt{3}/8 \right)=\frac{\left( \mu_S\left( r,\sqrt{3}/8 \right)-\mu(r)\right)^2}{\nu_S\left( r,\sqrt{3}/8 \right)}$ by substituting the relevant terms. See Figure \ref{fig:HLAEPlots for Seg}.

With $\epsilon=\sqrt{3}/4$, $r\in [1,2)$
 $$ \mu_S\left( r,\sqrt{3}/4 \right)=
\begin{cases}
          -{\frac{67}{54}}\,r^2+{\frac{40}{9}}\,r-3 &\text{for} \quad r \in [1,3/2)\\
           {\frac{7\,r^4-48\,r^3+122\,r^2-128\,r+48}{2\,r^2}} &\text{for} \quad r \in [3/2,2)
\end{cases}$$
and
$$\nu_S\left( r,\sqrt{3}/4 \right)=\sum_{j=1}^5\nu_j\left( r,\sqrt{3}/4 \right)\,\I(\mI_j)$$
where
{\small
\begin{align*}
\nu_1\left( r,\sqrt{3}/4 \right)&=-\Bigl[14285\,r^7-28224\,\,r^6-233266\,r^5+1106688\,r^4-2021199\,r^3+1876608\,r^2\\
&-880794\,\,r+165888\Bigr]/\Bigl[3645\,r\Bigr],\\
\nu_2\left( r,\sqrt{3}/4 \right)&=-\Bigl[14285\,r^{10}-28224\,\,r^9-233266\,r^8+1106688\,r^7-1234767\,r^6-3431808\,r^5\\
&+14049126\,r^4-22228992\,r^3+18895680\,r^2-8503056\,r+1594323\Bigr]/\Bigl[3645\,r^4\Bigr],\\
\nu_3\left( r,\sqrt{3}/4 \right)&=-\Bigl[14285\,r^{10}-28224\,\,r^9-233266\,r^8+1106688\,r^7-2545713\,r^6+5903280\,r^5\\
&-13456044\,r^4+20636208\,r^3-18305190\,r^2+8503056\,r-1594323\Bigr]/\Bigl[3645\,r^4\Bigr],\\
\nu_4\,\left( r,\sqrt{3}/4 \right)&=\Bigl[104920\,r^8-111072\,r^7+1992132\,r^6-15844032\,r^5+50174640\,r^4+6377292\\
&-34012224\,\,r+73220760\,r^2-81881280\,r^3+1909\,r^{10}-27072\,r^9\Bigr]/\Bigl[14580\,r^4\Bigr],\\
\nu_5\left( r,\sqrt{3}/4 \right)&=-\Bigl[-1187904\,\,r^5+1331492\,r^6+433304\,\,r^2+611163\,r^{10}-850240\,r^9-198144\,r\\
&+955392\,r^4-705536\,r^3-387680\,r^{11}+1118472\,r^8-1308960\,r^7+175984\,\,r^{12}\\
&-46176\,r^{13}+5120\,r^{14}+56016\Bigr]/\Bigl[20\,r^4\Bigr],
\end{align*}
}
and the corresponding intervals are $\mI_1=[1,9/8),\; \mI_2=[9/8,9/7),\; \mI_3=[9/7,4/3),\; \mI_4=[4/3,3/2),\; \mI_5=[3/2,2)$. See Figure \ref{fig:mean-var-seg} for the graph of $\mu_S\left( r,\sqrt{3}/4 \right)$ and $\nu_S\left( r,\sqrt{3}/4 \right)$.

Then we get $\HLAE\left( r,\sqrt{3}/4 \right)=\frac{\left( \mu_S\left( r,\sqrt{3}/4 \right)-\mu(r)\right)^2}{\nu_S\left( r,\sqrt{3}/4 \right)}$ by substituting the relevant terms. See Figure \ref{fig:HLAEPlots for Seg}.

With $\epsilon=2\,\sqrt{3}/7$,   $r \in [1,3/2)$
 $$ \mu_S\left( r,2\,\sqrt{3}/7 \right)=
\begin{cases}
          -{\frac{241}{54}}\,r^2+{\frac{38}{3}}\,r-8 &\text{for} \quad r \in [1,9/7)\\
           \frac{80\,r^4-432\,r^3+866\,r^2-756\,r+243}{2\,r^2} &\text{for} \quad r \in [9/7,3/2)
\end{cases}$$
and
$$\nu_S\left( r,2\,\sqrt{3}/7 \right)=\sum_{j=1}^6\nu_j\left( r,2\,\sqrt{3}/7 \right)\,\I(\mI_j)$$
where
\begin{align*}
\nu_1\left( r,2\,\sqrt{3}/7 \right)&=-\Bigl[2495087\,r^7-5067342\,r^6-29145379\,r^5+134149248\,r^4-230713503\,r^3\\
&+202262778\,r^2-90317349\,r+16336404\Bigr]/\Bigl[14580\,r\Bigr],\\
\nu_2\left( r,2\,\sqrt{3}/7 \right)&=-\Bigl[2495087\,r^{10}-5067342\,r^9-29145379\,r^8+134149248\,r^7-140359071\,r^6\\
&-378587142\,r^5+1465530651\,r^4-2206303596\,r^3+1786050000\,r^2-765450000\,r\\
&+136687500\Bigr]/\Bigl[14580\,r^4\Bigr],\\
\nu_3\left( r,2\,\sqrt{3}/7 \right)&=-\Bigl[2495087\,r^{10}-5067342\,r^9-29145379\,r^8+134149248\,r^7-309668679\,r^6\\
&+731864538\,r^5-1559738349\,r^4+2174176404\,\,r^3-1767825000\,r^2+765450000\,r\\
&-136687500\Bigr]/\Bigl[14580\,r^4\Bigr],\\
\nu_4\,\left( r,2\,\sqrt{3}/7 \right)&=\Bigl[1000147\,r^8-654768\,r^7+77561559\,r^6-527363136\,r^5+1468526760\,r^4\\
&+1767825000\,r^2-765450000\,r-2157840000\,r^3+136687500+24337\,r^{10}\\
&-321426\,r^9\Bigr]/\Bigl[14580\,r^4\Bigr],\\
\nu_5\left( r,2\,\sqrt{3}/7 \right)&={\frac{24337}{14580}}\,r^6-{\frac{17857}{810}}\,r^5+{\frac{1000147}{14580}}\,r^4-{\frac{18188}{405}}\,r^3-{\frac{174113}{1620}}\,r^2+{\frac{8176}{45}}\,r-78, \\
\nu_6\left( r,2\,\sqrt{3}/7 \right)&=-\frac{ \left( 8\,r^6-106\,r^5+8709\,r^4-39684\,\,r^3+68000\,r^2-51192\,r+14256 \right)  \left( 2\,r-3 \right) ^4}{20\,r^4}.
\end{align*}
The corresponding intervals are $\mI_1=[1,15/14),\; \mI_2=[15/14,15/13),\; \mI_3=[15/13,7/6),\; \mI_4=[7/6,5/4),\; \mI_5=[5/4,9/7),\; \mI_6=[9/7,3/2)$.  See Figure \ref{fig:mean-var-seg} for the graph of $\mu_S\left( r,2\,\sqrt{3}/7 \right)$ and $\nu_S\left( r,2\,\sqrt{3}/7 \right)$.

Then we get $\HLAE^S\left( r,2\,\sqrt{3}/7 \right)=\frac{\left( \mu_S\left( r,2\,\sqrt{3}/7 \right)-\mu_S\left( r,2\,\sqrt{3}/7 \right)\right)^2}{\nu_S\left( r,2\,\sqrt{3}/7 \right)}$ by substituting the relevant terms. In Figure \ref{fig:HLAEPlots for Seg} are the graphs of $\HLAE^S(r,\epsilon)$ for $\epsilon=\sqrt{3}/8,\,\sqrt{3}/4,\,2\,\sqrt{3}/7$.

From Figure \ref{fig:HLAEPlots for Seg}, we see that,
under $H^S_{\epsilon}$,
$\HLAE^S(r,\epsilon)$
appears to be an increasing function, dependent on $\epsilon$, of $r$.
Let $r_{\delta}(\epsilon)$ be the minimum $r$ such that
$\rho_n(r)$ becomes degenerate under the alternative $H^S_{\epsilon}$.
Then
$r_{\delta}\left( \sqrt{3}/8 \right)=4$,
$r_{\delta}\left( \sqrt{3}/4 \right)=2$, and
$r_{\delta}\left( 2\,\sqrt{3}/7 \right)=3/2$.  In fact, for $\epsilon \in \bigl( 0,\sqrt{3}/4 \bigr]$, $r_{\delta}(\epsilon)=\sqrt{3}/(2\,\epsilon)$ and for $\epsilon \in \left( \sqrt{3}/4,\sqrt{3}/3 \right)$, $r_{\delta}(\epsilon)=\sqrt{3}/\epsilon -2$.
Notice that $\lim_{r \rightarrow r_{\delta}(\epsilon)}\HLAE^S(\rho_n(r),\epsilon)=\infty$,
which is in agreement with PAE analysis because as $\epsilon \rightarrow 0$ HLAE becomes PAE, and as $\epsilon \rightarrow 0$, $r_{\delta}(\epsilon) \rightarrow \infty$
and under $H_0$, $\rho_n(r)$ is degenerate for $r=\infty$. The above result for HLAE can also be generalized for arbitrary $\epsilon$ as follows.

{\bf Proposition 1}
Let $\widetilde r:=\argsup_{r \in [1,r_{\delta}(\epsilon)]}\HLAE^S(r,\epsilon)$ where $r_{\delta}(\epsilon)$ is the value of $r$ at which $\rho_n(r)$ becomes degenerate under $H^S_{\epsilon}$. Then $\widetilde r=r_{\delta}(\epsilon)$. In particular, for $\epsilon \in \bigl[0,\sqrt{3}/4 \bigr]$, $r_{\delta}=\sqrt{3}/(2\,\epsilon)$ and  for $\epsilon \in \bigl(\sqrt{3}/4,\sqrt{3}/3\bigr]$,  $r_{\delta}=\sqrt{3}/\epsilon -2$.

{\bfseries Proof:}
Recall that   $\HLAE^S(r,\epsilon)=\frac{\bigl( \mu_S(r,\epsilon)-\mu(r) \bigr)^2}{\nu_S(r,\epsilon)}$.  For $\epsilon \in \left[ 0,\sqrt{3}/4 \right] $, $\mu_S(r,\epsilon)\rightarrow 1 $ and $\nu(r,\epsilon)\rightarrow 0$ as $r \rightarrow r_{\delta}(\epsilon)=\sqrt{3}/(2\,\epsilon)$.  Hence $\HLAE^S(r,\epsilon)\rightarrow \infty $ as $r \rightarrow r_{\delta}(\epsilon)=\sqrt{3}/(2\,\epsilon)$.  So for $\epsilon \in \left[ 0,\sqrt{3}/4 \right] $, the $\widetilde r =\sqrt{3}/(2\,\epsilon)$.  For $\epsilon \in \Bigl(\sqrt{3}/4,\sqrt{3}/3\Bigr] $, the result follows similarly. $\blacksquare$

So HLAE suggests choosing $r$ larger as the segregation gets more severe,
but choosing $r$ too large will reduce power since
$r \geq r_{\delta}(\epsilon)$ guarantees the complete digraph
under the alternative and, as $r$ increases therefrom,
provides an ever greater probability of seeing the complete digraph
under the null.

In Figure \ref{fig:mean-var-seg}, we plot the graphs of mean and asymptotic variance for $r \in [1,4]$ under segregation with $\epsilon=0,\sqrt{3}/8, \sqrt{3}/4,2\,\sqrt{3}/7$. Notice that $\mu_S(r,\epsilon)$ gets larger as $\epsilon$ gets larger at each $r$ which is in agreement with the $\mu_S(r,\epsilon)$ expressions in Section \ref{sec:arc-prob-seg}.  However, the same ordering holds for $\nu_S(r,\epsilon)$ at each $r$ only for large $r$, but for small $r$ the ordering is reversed. Furthermore, both the $\sup_{r \in [1,\infty]} \nu_S(r,\epsilon)$ and $\argsup_{r \in [1,\infty]} \nu_S(r,\epsilon)$ seem to decrease as $\epsilon$ increases.

\subsubsection{Hodges-Lehmann Asymptotic Efficacy Under Association Alternatives}
\label{sec:Hodges-Lehmann-Agg}
In the HLAE framework, $\theta=\E^A_{\epsilon}[\rho_n(r)]=\mu_A(r,\epsilon)$ and $\theta_0=\mu_A(r,\epsilon=0)$.  Then testing $H_0:\epsilon=0$ versus $H_a:\epsilon>0$ is equivalent to $H_0:\E[\rho_n]=\mu(r)$ versus $H_a:\E^A_{\epsilon}[\rho_n]=\mu_A(r,\epsilon)<\mu(r)$.  Let $\delta=\frac{\mu(r,\epsilon)-\mu(r)}{\sigma_n(r,\epsilon)}$ and $\widetilde{R}=\frac{\rho_n(r)-\mu(r)}{\sigma_n(r,\epsilon)}$, then $\widetilde{R}\stackrel{\mathcal L}{\rightarrow}\N(\delta,1)$.

Hodges-Lehmann asymptotic efficacy (HLAE)
(\cite{hodges:1956})
is given by
$$
\HLAE^A(r,\epsilon):=\frac{\bigl( \mu_A(r,\epsilon)-\mu(r) \bigr)^2}{\nu_A(r,\epsilon)}.
$$

Rather than an arbitrary $\epsilon$ we pick specific values: $5\,\sqrt{3}/24,\,\sqrt{3}/12$, and $\sqrt{3}/21$.  Recall that $\rho_n(r)$ is degenerate as $r\rightarrow \infty $. Furthermore $\rho_n(r)$ is degenerate when $r=1$.

\begin{figure}[]
\centering
\scalebox{.3}{\input{ArcProbEpsAgg.pstex_t}}
\scalebox{.3}{\input{AsyVarEpsAgg.pstex_t}}
\caption{
\label{fig:mean-var-agg}
The mean $\mu_A(r,\epsilon)$ (left)
and asymptotic variance $\nu_A(r,\epsilon)$ (right) as a function of $r$ under association with $\epsilon=0,\sqrt{3}/21, \sqrt{3}/12,5\,\sqrt{3}/24$.}
\end{figure}
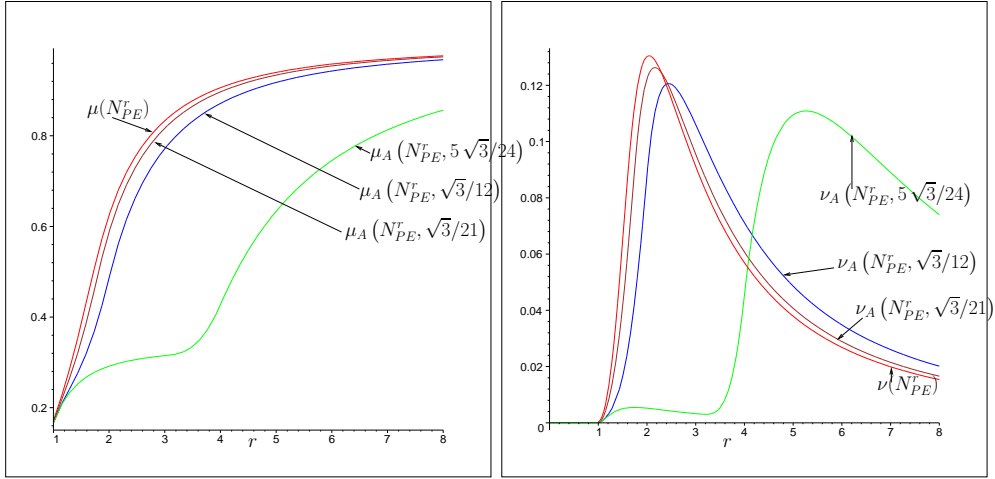

With $\epsilon=5\,\sqrt{3}/24$,
 $$ \mu_A\left( r,5\,\sqrt{3}/24 \right)=
\begin{cases}
          -\frac{1}{6}\,r^{-2}+\frac{1}{3} &\text{for} \quad r \in [1,3)\\
          \frac{1}{3}\,r^2-\frac{8}{3}\,r-{\frac{55}{6}}\,r^{-2}+{\frac{19}{3}} &\text{for} \quad r \in  [3,4)\\
          -{\frac{55}{6}}\,r^{-2}+1 &\text{for} \quad r \in [4,\infty)
\end{cases}$$
and
$$\nu_A\left( r,5\,\sqrt{3}/24 \right)=\sum_{j=1}^5\nu_j\left( r,5\,\sqrt{3}/24 \right)\,\I(\mI_j)$$
where
\begin{align*}
\nu_1\left( r,5\,\sqrt{3}/24 \right)&=\frac{r^4-2\,r^2+1}{27\,r^6},\\
\nu_2\left( r,5\,\sqrt{3}/24 \right)&=-\Bigl[120\,r^{10}-2176\,r^9+15340\,r^8-50304\,\,r^7+58754\,\,r^6+74880\,r^5-248577\,r^4\\
&+138240\,r^3+47172\,r^2+23328\,r-7305\Bigr]/\Bigl[405\,r^6\Bigr],\\
\nu_3\left( r,5\,\sqrt{3}/24 \right)&=-\Bigl[120\,r^{10}-2176\,r^9+15180\,r^8-48960\,r^7+58754\,\,r^6+47440\,r^5-176547\,r^4\\
&+138240\,r^3-70477\,r^2+23328\,r-7305\Bigr]/\Bigl[405\,r^6\Bigr],\\
\nu_4\,\left( r,5\,\sqrt{3}/24 \right)&=\Bigl[10\,r^{12}-192\,r^{11}+1320\,r^{10}-2944\,r^9-7590\,r^8+49920\,r^7-69986\,r^6-46480\,r^5\\
&+184137\,r^4-143360\,r^3+71917\,r^2-23520\,r+7315\Bigr]/\Bigl[405\,r^6\Bigr],\\
\nu_5\left( r,5\,\sqrt{3}/24 \right)&=\frac{787\,r^4-7601\,r^2-16032\,r+9265}{135\,r^6}.
\end{align*}
The corresponding intervals are $\mI_1=(1,3),\; \mI_2=[3,7/2),\; \mI_3=\Bigl[7/2,2+\sqrt{3}\Bigr),\; \mI_4=\Bigl[ 2+\sqrt{3},4 \Bigr),\; \mI_5=[4,\infty)$. See Figure \ref{fig:mean-var-agg} for the graph of $\mu_A\left( r,5\,\sqrt{3}/24 \right)$ and $\nu_A\left( r,5\,\sqrt{3}/24 \right)$.

\begin{figure}[]
\centering
\scalebox{.22}{\input{hlae_agg3.pstex_t}}
\scalebox{.22}{\input{hlae_agg2.pstex_t}}
\scalebox{.22}{\input{hlae_agg1.pstex_t}}
\caption{ \label{fig:HLAEPlots for Agg}
Hodges-Lehmann asymptotic efficacy
against association alternative $H^A_{\epsilon}$
as a function of $r$
for $\epsilon= \sqrt{3}/21,\, \sqrt{3}/12,\,  5\,\sqrt{3}/24$ (left to right).}
\end{figure}
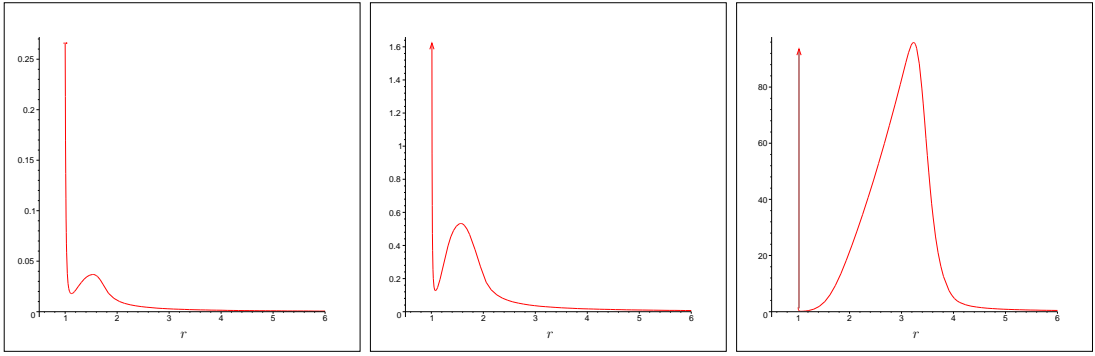

Then we get $\HLAE^A\left( r,5\,\sqrt{3}/24 \right)=\frac{\left( \mu_A\left( r,5\,\sqrt{3}/24 \right)-\mu(r) \right)^2}{\nu_A\left( r,5\,\sqrt{3}/24 \right)}$ by substituting the relevant terms. See Figure \ref{fig:HLAEPlots for Agg}.

With $\epsilon=\sqrt{3}/12$,
 $$ \mu_A\left( r,\sqrt{3}/12 \right)=
\begin{cases}
          \frac{6\,r^4-16\,r^3+18\,r^2-5}{18\,r^2} &\text{for} \quad r \in [1,2)\\
           -{\frac{37}{18}}\,r^{-2}+1 &\text{for} \quad r \in [2,\infty)
\end{cases}$$
and
$$\nu_A\left( r,\sqrt{3}/12 \right)=\sum_{j=1}^3\nu_j\left( r,\sqrt{3}/12 \right)\,\I(\mI_j)$$
where
\begin{align*}
\nu_1\left( r,\sqrt{3}/12 \right)&=\Bigl[10\,r^{12}-96\,r^{11}+240\,r^{10}+192\,r^9-1830\,r^8+3360\,r^7-2650\,r^6+240\,r^5+1383\,r^4\\
&-1280\,r^3+540\,r^2-144\,r+35\Bigr]/\Bigl[405\,r^6\Bigr],\\
\nu_2\left( r,\sqrt{3}/12 \right)&=\Bigl[10\,r^{12}-96\,r^{11}+240\,r^{10}+192\,r^9-1670\,r^8+2784\,\,r^7-2650\,r^6+2400\,r^5-1047\,r^4\\
&-1280\,r^3+1269\,r^2-144\,r+35\Bigr]/\Bigl[405\,r^6\Bigr],\\
\nu_3\left( r,\sqrt{3}/12 \right)&=\frac{537\,r^4-683\,r^2-2448\,r+1315}{405\,r^6}.
\end{align*}
The corresponding intervals are $\mI_1=[1,3/2),\; \mI_2=[3/2,2),\; \mI_3=[2,\infty)$. See Figure \ref{fig:mean-var-agg} for the graph of $\mu\left( r,\sqrt{3}/12 \right)$ and $\nu\left( r,\sqrt{3}/12 \right)$.

Then we get $\HLAE^A\left( r,\sqrt{3}/12 \right)=\frac{\left( \mu_A\left( r,\sqrt{3}/12 \right)-\mu(r)\right)^2}{\nu_A\left( r,\sqrt{3}/12 \right)}$ by substituting the relevant terms. See Figure \ref{fig:HLAEPlots for Agg}.

With $\epsilon=\sqrt{3}/21$,
 $$ \mu_A\left( r,\sqrt{3}/21 \right)=
\begin{cases}
          \frac{7839\,r^4-27648\,r^3+49152\,r^2-35840\,r+9216}{16200\,r^2} &\text{for} \quad r \in [1,8/7)\\
          \frac{2719\,r^4-5592\,r^3+5760\,r^2-1536}{8100\,r^2} &\text{for} \quad r \in (8/7,3/2)\\
          \frac{53\,r^4+2744\,r^3-7296\,r^2+8064\,\,r-3104}{2700\,r^2} &\text{for} \quad r \in (3/2,12/7)\\
          \frac{2719\,r^4-1440\,r^2+2112}{16200\,r^2} &\text{for} \quad r \in (12/7,7/4)\\
          -\frac{2401\,r^4-73824\,\,r^2+153664\,\,r-88548}{16200\,r^2} &\text{for} \quad r \in (7/4,2)\\
          1-{\frac{89}{54}}\,r^{-2} &\text{for} \quad r \in [2,\infty)
\end{cases}$$
and
$$\nu_A\left( r,\sqrt{3}/21 \right)=\sum_{j=1}^{10}\nu_j\left( r,\sqrt{3}/21 \right)\,\I(\mI_j)$$
where
{\small
\begin{align*}
\nu_1\left( r,\sqrt{3}/21 \right)&=\Bigl[4124031\,r^{12}-22708224\,\,r^{11}-389826\,r^{10}+369129408\,r^9-1592672721\,r^8\\
&+3532359672\,r^7-4721848374\,\,r^6+4050858048\,r^5-2387433568\,r^4+995033088\,r^3\\
&-209048784\,\,r^2-43352064\,\,r+25952256\Bigr]/\Bigl[65610000\,r^6\Bigr],\\
\nu_2\left( r,\sqrt{3}/21 \right)&=\Bigl[6594660\,r^{12}-31178952\,r^{11}-14911074\,\,r^{10}+441735648\,r^9-1578842961\,r^8\\
&+3311083512\,r^7-4669163574\,\,r^6+4366966848\,r^5-2522908768\,r^4+778272768\,r^3\\
&-93443280\,r^2+14450688\,r-8650752\Bigr]/\Bigl[65610000\,r^6\Bigr],\\
\nu_3\left( r,\sqrt{3}/21 \right)&=\Bigl[826701\,r^{12}-7118748\,r^{11}+14155864\,\,r^{10}+18467640\,r^9-104968680\,r^8\\
&+165877272\,r^7-128355690\,r^6+27338184\,\,r^5+47304144\,r^4-52684800\,r^3\\
&+24413592\,r^2-7225344\,r+1966080\Bigr]/\Bigl[32805000\,r^6\Bigr],\\
\nu_4\,\left( r,\sqrt{3}/21 \right)&=\Bigl[826701\,r^{12}-7118748\,r^{11}+14155864\,\,r^{10}+18467640\,r^9+20074008\,r^8\\
&-671672808\,r^7+2194076310\,r^6-3382581816\,r^5+2840904144\,r^4-1262284800\,r^3\\
&+240413592\,r^2-7225344\,r+1966080\Bigr]/\Bigl[32805000\,r^6\Bigr],
\end{align*}
\begin{align*}
\nu_5\left( r,\sqrt{3}/21 \right)&=\Bigl[826701\,r^{12}-7118748\,r^{11}+14155864\,\,r^{10}+18467640\,r^9-137116617\,r^8\\
&+512952192\,r^7-1511673690\,r^6+2773418184\,\,r^5-2883095856\,r^4+1560115200\,r^3\\
&-335586408\,r^2-7225344\,r+1966080\Bigr]/\Bigl[32805000\,r^6\Bigr],\\
\nu_6\left( r,\sqrt{3}/21 \right)&=\Bigl[826701\,r^{12}-7118748\,r^{11}+14155864\,\,r^{10}+18467640\,r^9-91939401\,r^8\\
&+125718912\,r^7-128697690\,r^6+139178184\,\,r^5-60695856\,r^4-52684800\,r^3\\
&+48413592\,r^2-7225344\,r+1966080\Bigr]/\Bigl[32805000\,r^6\Bigr],\\
\nu_7\left( r,\sqrt{3}/21 \right)&=\Bigl[226415\,r^{12}-1426740\,r^{11}+334536\,r^{10}+17196648\,r^9-87678147\,r^8\\
&+311364480\,r^7-711864862\,r^6+944809880\,r^5-684036240\,r^4+238099456\,r^3\\
&-24048504\,\,r^2-7633920\,r+4761088\Bigr]/\Bigl[10935000\,r^6\Bigr],\\
\nu_8\left( r,\sqrt{3}/21 \right)&=\Bigl[5786907\,r^{12}-42712488\,r^{11}+76274888\,r^{10}+51865788\,r^8-300043296\,r^7\\
&+132202536\,r^6+171413760\,r^5-93614976\,r^4+147517440\,r^3-194460480\,r^2\\
&+67608576\,r-29061120\Bigr]/\Bigl[262440000\,r^6\Bigr],\\
\nu_9\left( r,\sqrt{3}/21 \right)&=-\Bigl[2470629\,r^{12}-25412184\,\,r^{11}+112001848\,r^{10}-1958438076\,r^8\\
&+5449924256\,r^7+6150612888\,r^6-55820599296\,r^5+109663683136\,r^4\,\\
&-97335694848\,r^3+40552466112\,r^2-9825887232\,r+3078523200\Bigr]/\Bigl[262440000\,r^6\Bigr],\\
\nu_{10}\left( r,\sqrt{3}/21 \right)&=\frac{493829\,r^4-433645\,r^2-1765008\,r+929955}{455625\,r^6}.
\end{align*}
}

The corresponding intervals are $\mI_1=\Bigl( 1,2\,\sqrt{14}/7 \Bigr),\; \mI_2=\Bigl[2\,\sqrt{14}/7,8/7\Bigr),\; \mI_3=[8/7,5/4),\; \mI_4=[5/4,4/3),\; \mI_5=[4/3,10/7),\; \mI_6=[10/7,3/2),\; \mI_7=[3/2,12/7),\; \mI_8=[12/7,7/4),\; \mI_9=[7/4,2),\; \mI_{10}=[2,\infty)$.  See Figure \ref{fig:mean-var-agg} for the graph of $\mu_A\left( r,\sqrt{3}/21 \right)$ and $\nu_A\left( r,\sqrt{3}/21 \right)$.

Then we get $\HLAE^A\left( r,\sqrt{3}/21 \right)=\frac{\left( \mu_A\left( r,\sqrt{3}/21 \right)-\mu(r) \right)^2}{\nu_A\left( r,\sqrt{3}/21 \right))}$ by substituting the relevant terms. See Figure \ref{fig:HLAEPlots for Agg}.

Notice that for $\epsilon=5\,\sqrt{3}/24$ and $\epsilon=\sqrt{3}/12$, $\argsup _{r \ge 1}\HLAE^A(r,\epsilon)=1$.  This result for HLAE can be generalized for arbitrary $\epsilon$ as follows.

{\bf Proposition 2}
Let $r^*:=\argsup_{r \ge 1}\HLAE^A(r,\epsilon)$ and $\epsilon \ge \sqrt{3}/12$. Then $r^*=1$.

{\bfseries Proof:}
Recall that $\HLAE^A(r,\epsilon)=\frac{\left( \mu_A(r,\epsilon)-\mu(r)\right)^2}{\nu_A(r,\epsilon)}$.  For $\epsilon \in \Bigl[\sqrt{3}/12,\sqrt{3}/3 \Bigr)$, $\mu_A(r=1,\epsilon)\rightarrow 1 $ and $\nu_A(r=1,\epsilon)\rightarrow 0$.  Hence $\HLAE^A(r,\epsilon)\rightarrow \infty $ as $r \rightarrow 1$.  So the desired result follows. $\blacksquare$

For $\epsilon \in \Bigl[0,\sqrt{3}/12 \Bigr] $, it seems that for a while $r^*=1$ with respect to HLAE, e.g. for $\epsilon=\sqrt{3}/21$.  But for sufficiently small $\epsilon$, $r^*>1 $ holds.  This can also be seen as $\epsilon \rightarrow 0$ in which case HLAE becomes PAE and the optimal value is about $1.006 $ with respect to PAE.  Furthermore, observe that the argsup for HLAE gets closer to 1 as $\epsilon \rightarrow 0$ and $\nu_A(r,\epsilon)>0$ for $\epsilon \in \Bigl(0,\sqrt{3}/12\Bigr)$ and $\nu(r,\epsilon)$ gets larger as $\epsilon \rightarrow 0$.

Figure \ref{fig:HLAEPlots for Agg}
contains a graph of HLAE
against association as a function of $r$
for $\epsilon= 5\,\sqrt{3}/24,\, \sqrt{3}/12, \, \sqrt{3}/21$.  Notice that since $\nu(r=1,\epsilon)=0$ for $\epsilon \ge \sqrt{3}/12$, $\HLAE^A(r=1,\epsilon)=\infty$ for $\epsilon \ge \sqrt{3}/12$ and $\lim_{r\rightarrow \infty}\HLAE^A(r,\epsilon)=0$.
In Figure \ref{fig:HLAEPlots for Agg} we see that,
against $H^A_{\epsilon}$,
$\HLAE^A(r,\epsilon)$
has a local supremum for some $r>1$.
Let $\widetilde r_l$ be the value at which this local supremum is attained.  Then
$\widetilde r_l\left( 5\,\sqrt{3}/24 \right)\approx 3.2323$,
$\widetilde r_l\left( \sqrt{3}/12 \right) \approx 1.5676$, and
$\widetilde r_l\left( \sqrt{3}/21 \right) \approx 1.533$.
Note that, as $\epsilon$ gets smaller, $\widetilde r$ gets smaller.  Furthermore, $\HLAE^A\left(r=1,\sqrt{3}/21\right) < \infty$ and as $\epsilon \rightarrow 0$, so $\widetilde r$ becomes the global supremum, and $\PAE^A(r=1)=0$ and $\argsup_{r \ge 1}\PAE^A(r=1)\approx 1.006$. So HLAE suggests choosing moderate $r$ when testing against association, whereas PAE suggests choosing small $r$.

Derivation of $\mu_A(r,\epsilon)$ and $\nu_A(r,\epsilon)$ for association with $\epsilon=5\,\sqrt{3}/24,\,\sqrt{3}/12$, and $\sqrt{3}/21$ are similar to --- with the supports being the complements of--- the corresponding segregation cases.

In Figure \ref{fig:mean-var-agg}, we plot the graphs of mean and asymptotic variance for $r \in [1,8]$ under association with $\epsilon=0,\sqrt{3}/21, \sqrt{3}/12,5\,\sqrt{3}/24$. Notice that $\mu_A(r,\epsilon)$ gets smaller as $\epsilon$ gets larger at each $r$ which is in agreement with the $\mu_A(r,\epsilon)$ expressions in Appendix 2.  However, the same ordering does not hold for $\nu_A(r,\epsilon)$ at each $r$. For small $r$ the ordering is same as in $\mu_A(r,\epsilon)$, but for large $r$ the ordering is reversed. Furthermore, $\sup_{r \in [1,\infty]} \nu_A(r,\epsilon)$ seems to decrease as $\epsilon$ increases while $\argsup_{r \in [1,\infty]} \nu_A(r,\epsilon)$ seems to increase as $\epsilon$ increases.

\subsection{Asymptotic Power Function Analysis}
\label{sec:APF}

The asymptotic power function (see, e.g., \cite{kendall:1979}) can also be investigated as a function of $r$, $n$, and $\epsilon$
using the asymptotic critical value and an appeal to normality.

\subsubsection{Asymptotic Power Function Analysis Under Segregation}
\begin{figure}[]
\centering
\psfrag{r}{\scriptsize{$r$}}
\epsfig{figure=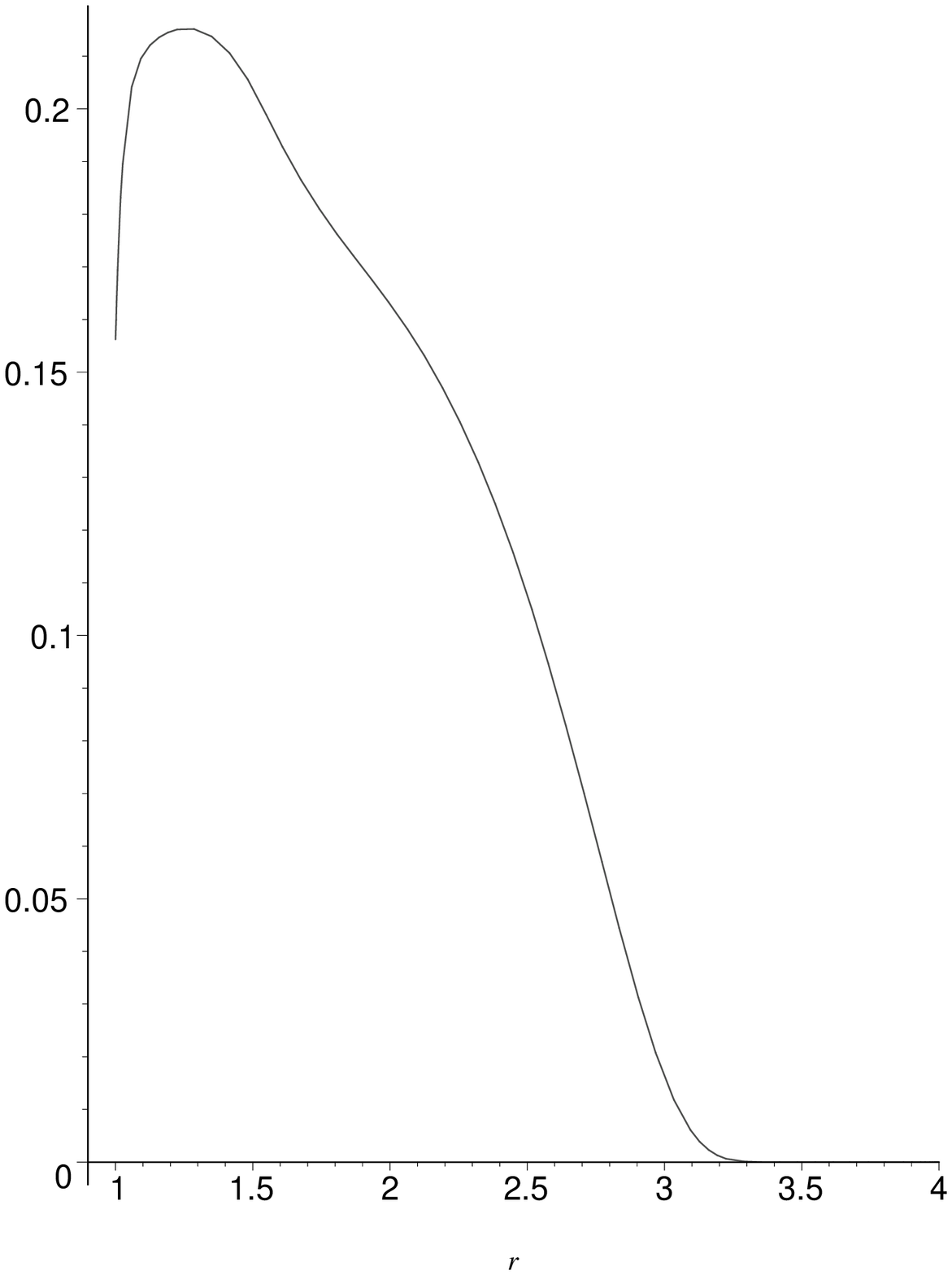, height=125pt, width=125pt}
\epsfig{figure=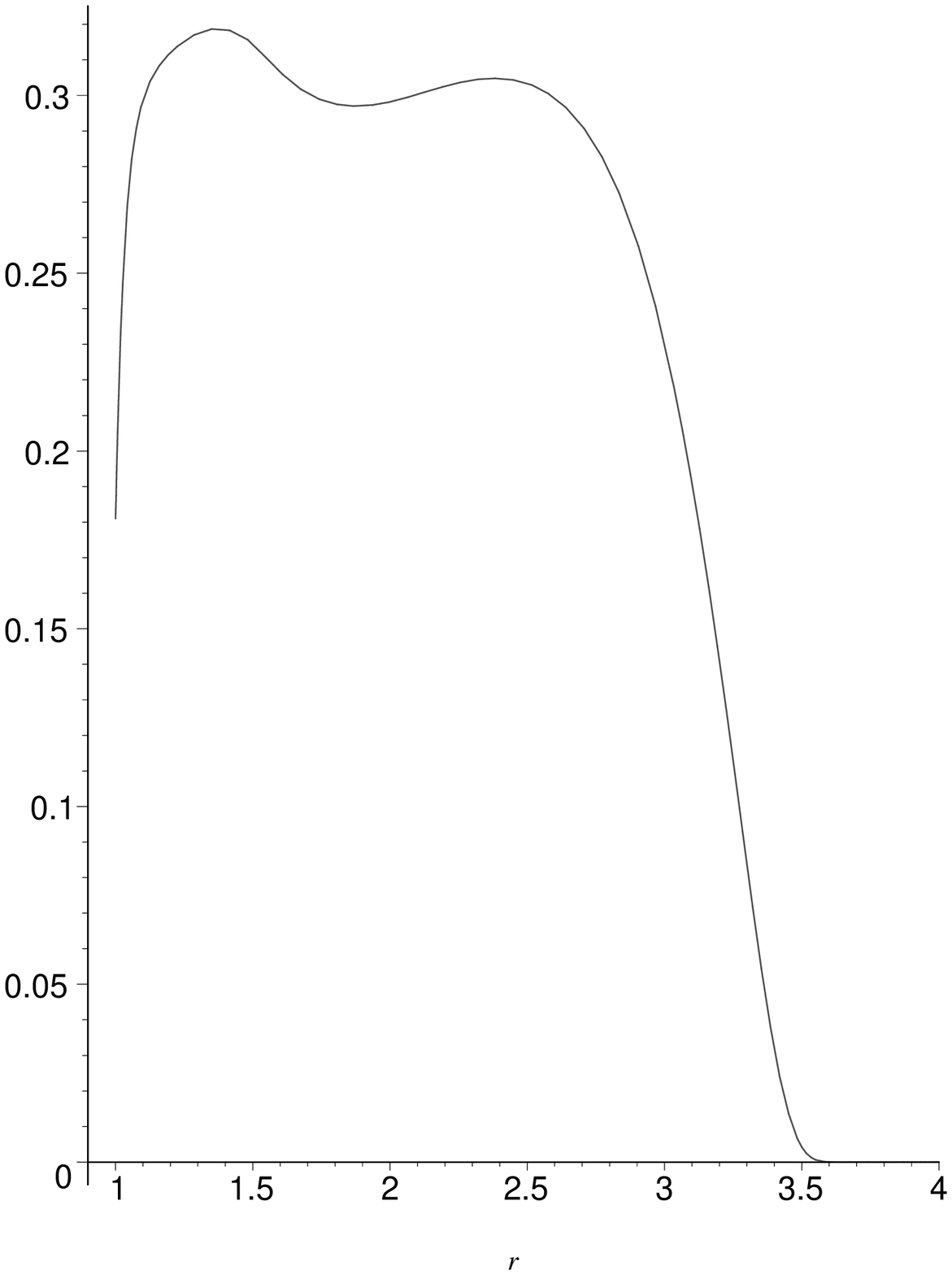, height=125pt, width=125pt}
\epsfig{figure=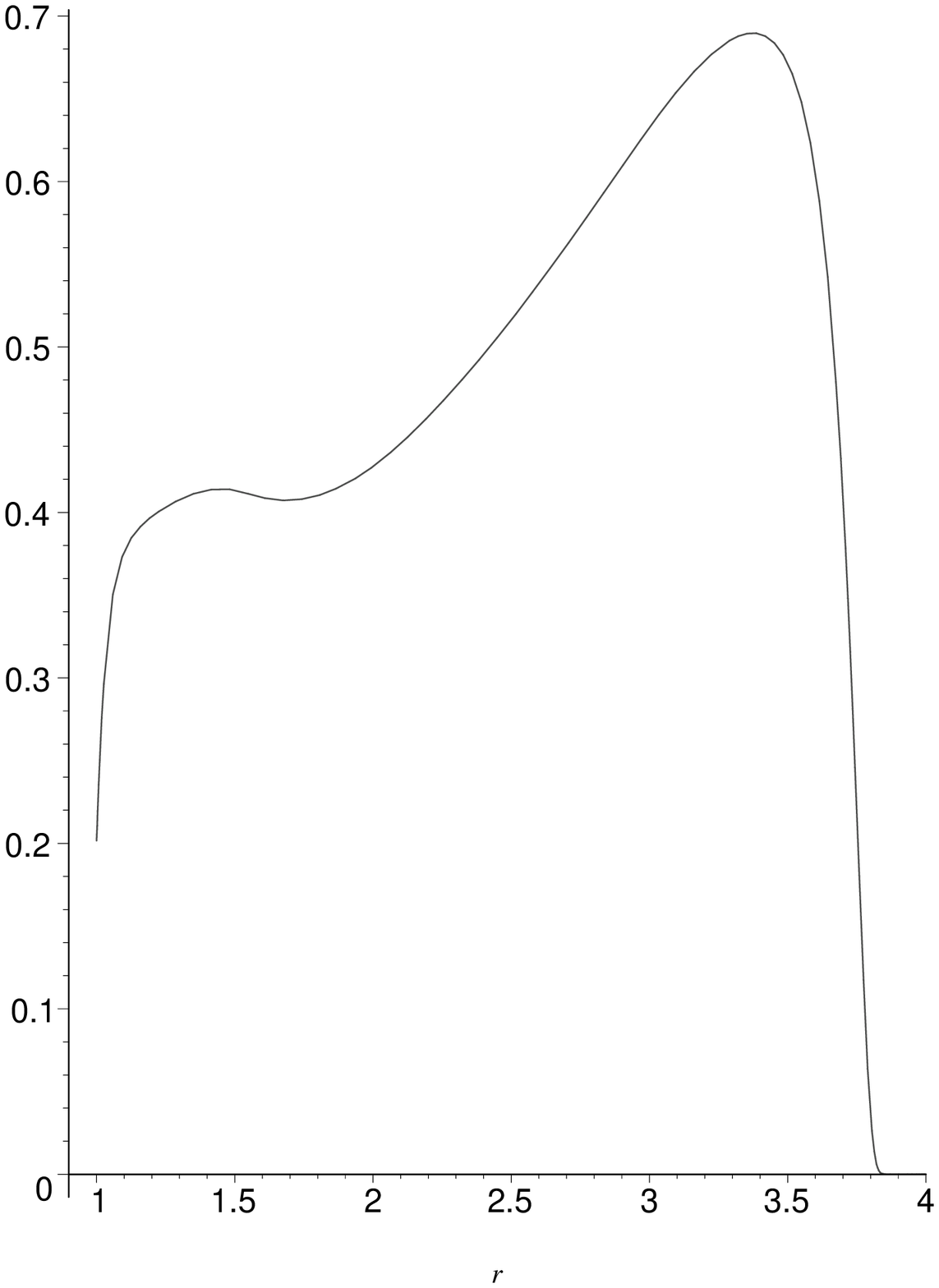, height=125pt, width=125pt}
\epsfig{figure=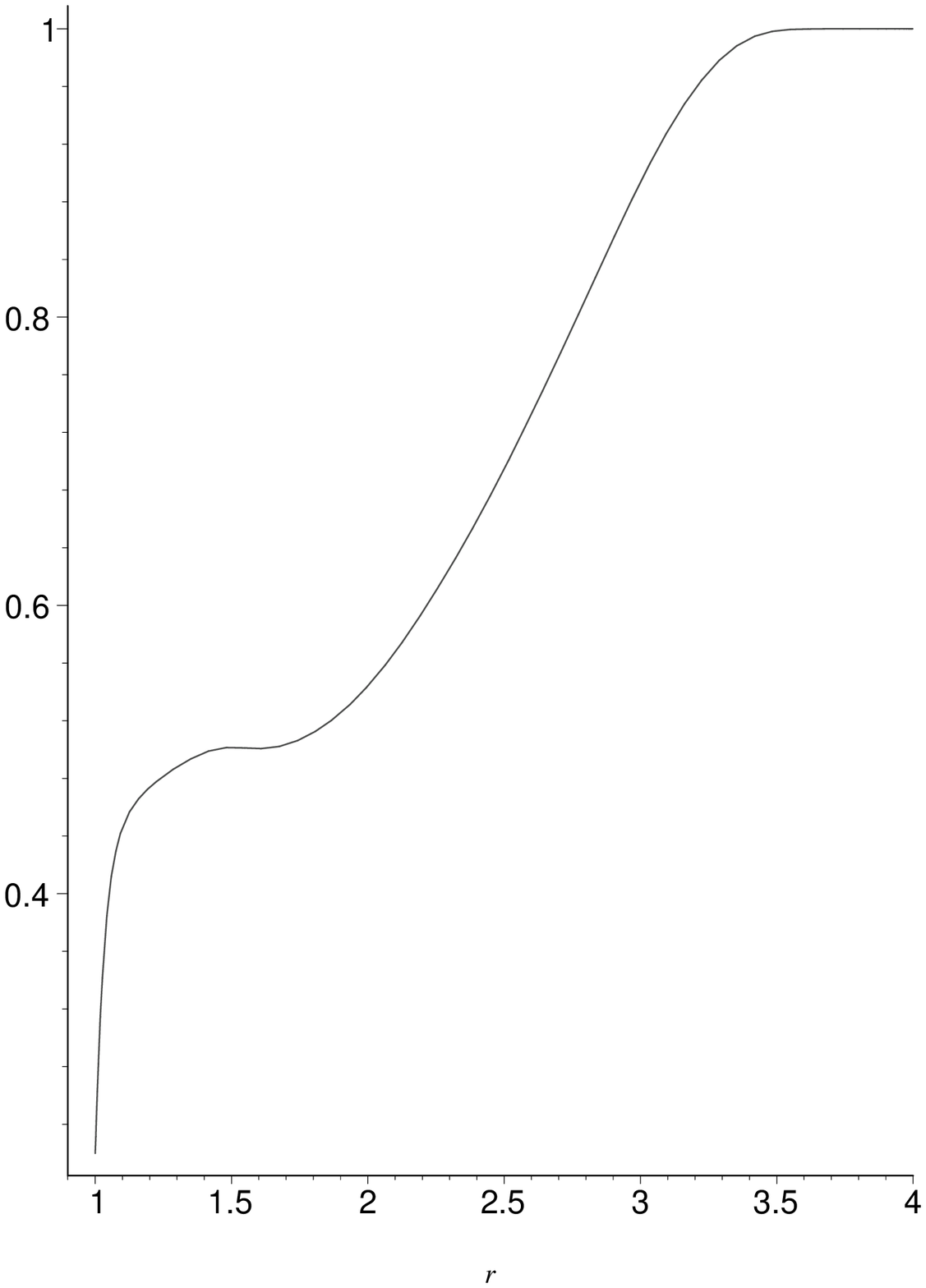, height=125pt, width=125pt}
\epsfig{figure=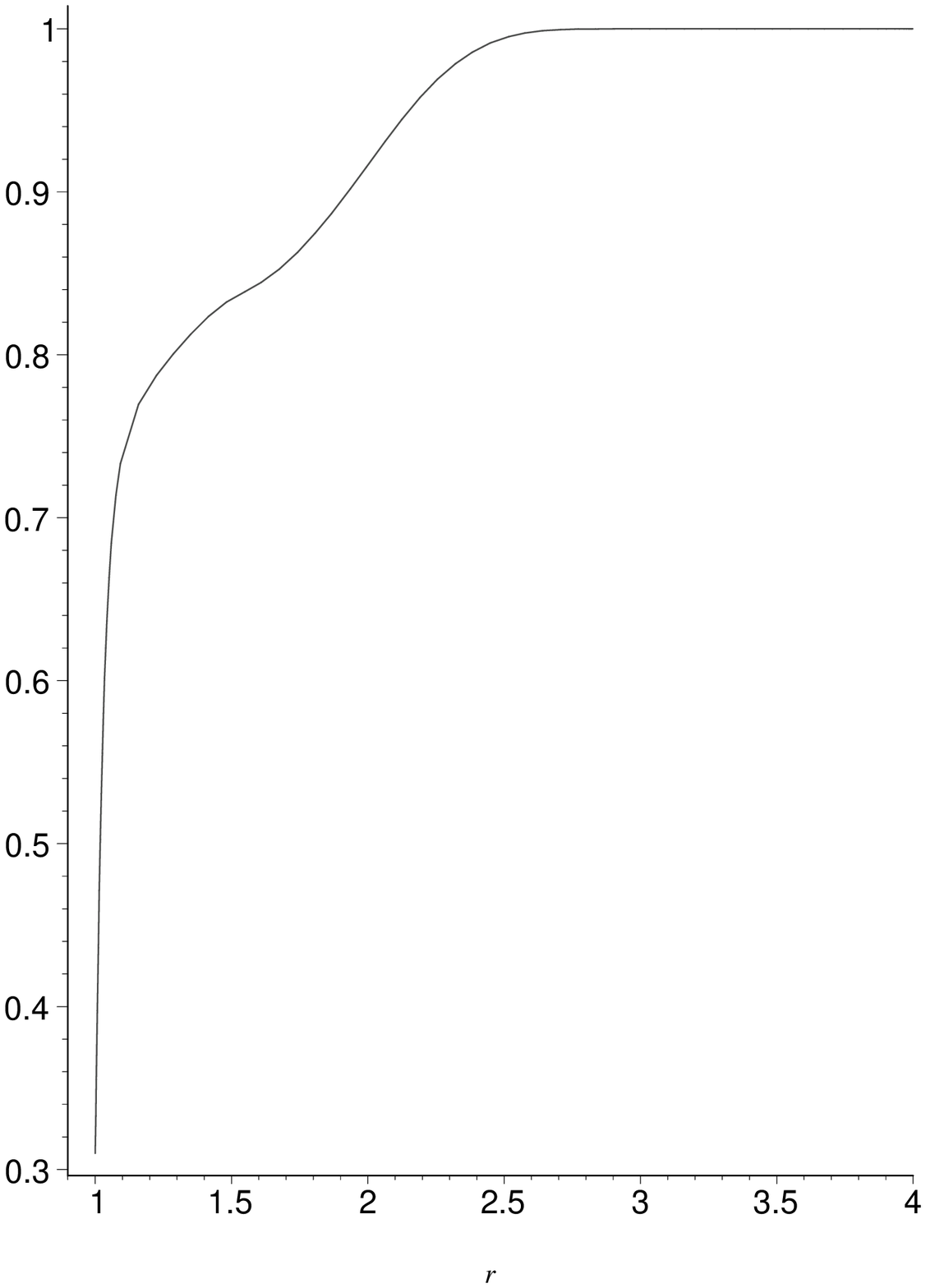, height=125pt, width=125pt}
\epsfig{figure=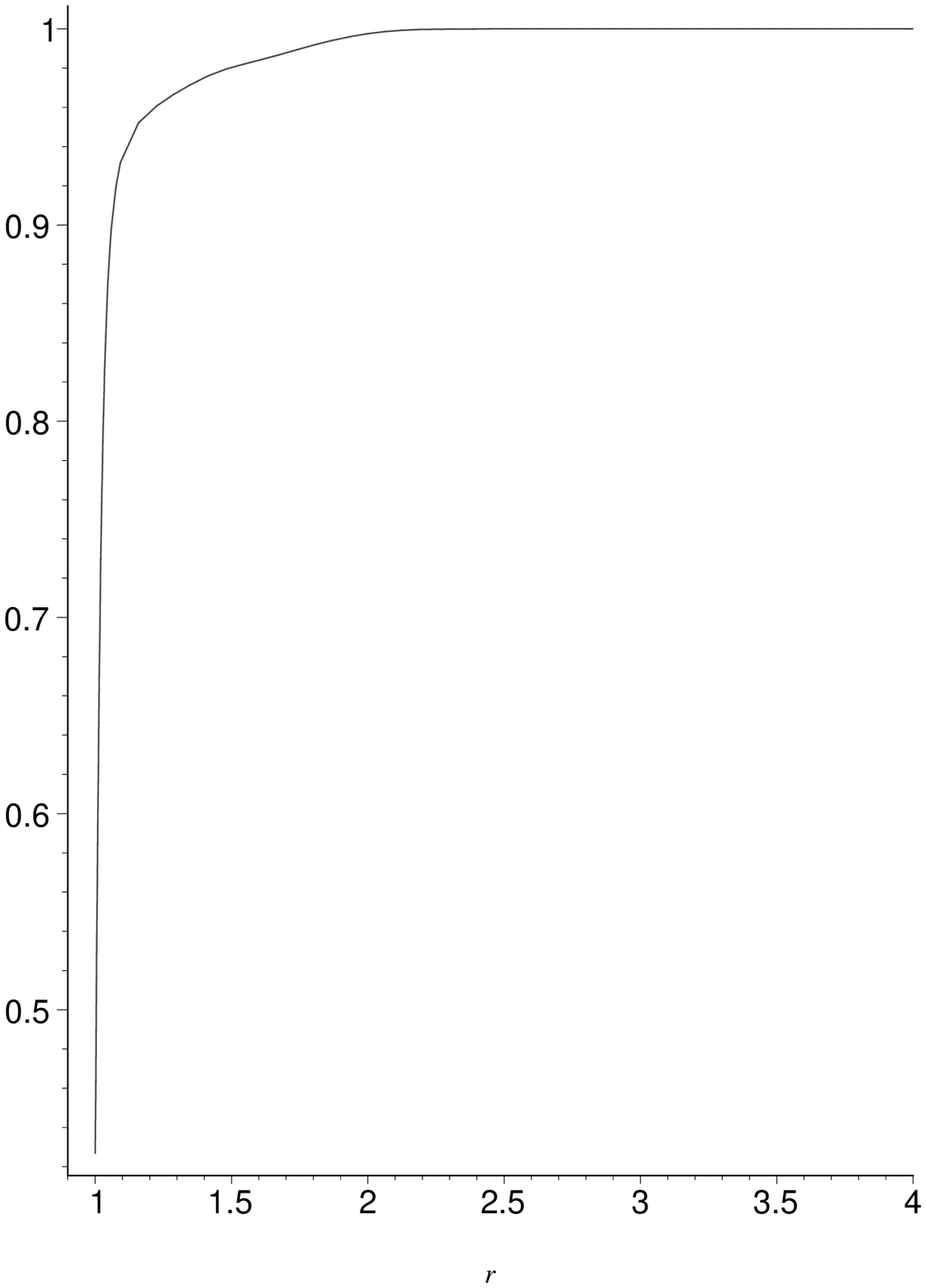, height=125pt, width=125pt}
\caption{
\label{fig:APF-seg-plots}
Asymptotic power function
against segregation alternative $H^S_{\sqrt{3}/8}$
as a function of $r$
for $n=5,10,15,20,50,100$ .
}
\end{figure}

Under segregation, for sufficiently large $n$, we reject $H_0$ when $\sqrt{n}\,\Bigl( \frac{\rho_n(r)-\mu_S(r,\epsilon)}{\sqrt{\nu_S(r,\epsilon)}} \Bigr) > z_{(1-\alpha)}$  where $z_{(1-\alpha)}$ is the $(1-\alpha)\times 100 $ percentile of the standard normal distribution, e.g. with $\alpha=.05 $, $z_{.95} \approx 1.645 $.  Then size $\alpha $ critical region for large samples is
 $$\rho_n(r)>\mu(r) + z_{(1-\alpha)} \cdot \sqrt{\nu(r)/n}.$$
Under a specific segregation alternative $H^S_{\epsilon}$,
the asymptotic power function is given by

\begin{eqnarray*}
\Pi_S(r,n,\epsilon)&:=&P\Bigl(\rho_n(r)>\mu(r) + z_{(1-\alpha)} \cdot \sqrt{\nu(r)/n}\Bigr)\\
&=&1-\Phi\Biggl(\frac{z_{(1-\alpha)}\,\sqrt{\nu(r)}}{\sqrt{\nu_S(r,\epsilon)}}+\frac{\sqrt{n}\,(\mu(r)-\mu_S(r,\epsilon))}{\sqrt{\nu_S(r,\epsilon)}}\Biggr).
\end{eqnarray*}
With $\epsilon=\sqrt{3}/8$, $\Pi_S(r,n,\epsilon)$ at level $\alpha=.05 $ is plotted in Figure \ref{fig:APF-seg-plots}.  Observe that $\Pi_S(r,n,\sqrt{3}/8)\rightarrow 0$ as $r \rightarrow 4$ for $n=5,10,15$. Let $r^*_g(n,\epsilon)$ be the the value at which $\Pi_S(r,n,\epsilon)$ attains its global supremum and $r^*_l(n,\epsilon)$ be the the value at which $\Pi_S(r,n,\epsilon)$ attains its local supremum. Then $r^*_g(5,\sqrt{3}/8)\approx 1.260$, $r^*_g(10,\sqrt{3}/8) \approx 1.3741$ and $r^*_l(10,\sqrt{3}/8)\approx 2.3818$, $r^*_g(15,\sqrt{3}/8) \approx 3.3724$ and $r^*_l(15,\sqrt{3}/8)\approx 1.45$, $r^*_g(20,\sqrt{3}/8)=4$ and $r^*_l(20,\sqrt{3}/8)\approx 1.5$. Finally, $r^*_g(n,\sqrt{3}/8)=4$ for $n=20,50,100$ and $\Pi_S(r,n,\sqrt{3}/8)$ has a hump for $n=10$ and $n=15$.

\begin{figure}[]
\centering
\psfrag{r}{\scriptsize{$r$}}
\epsfig{figure=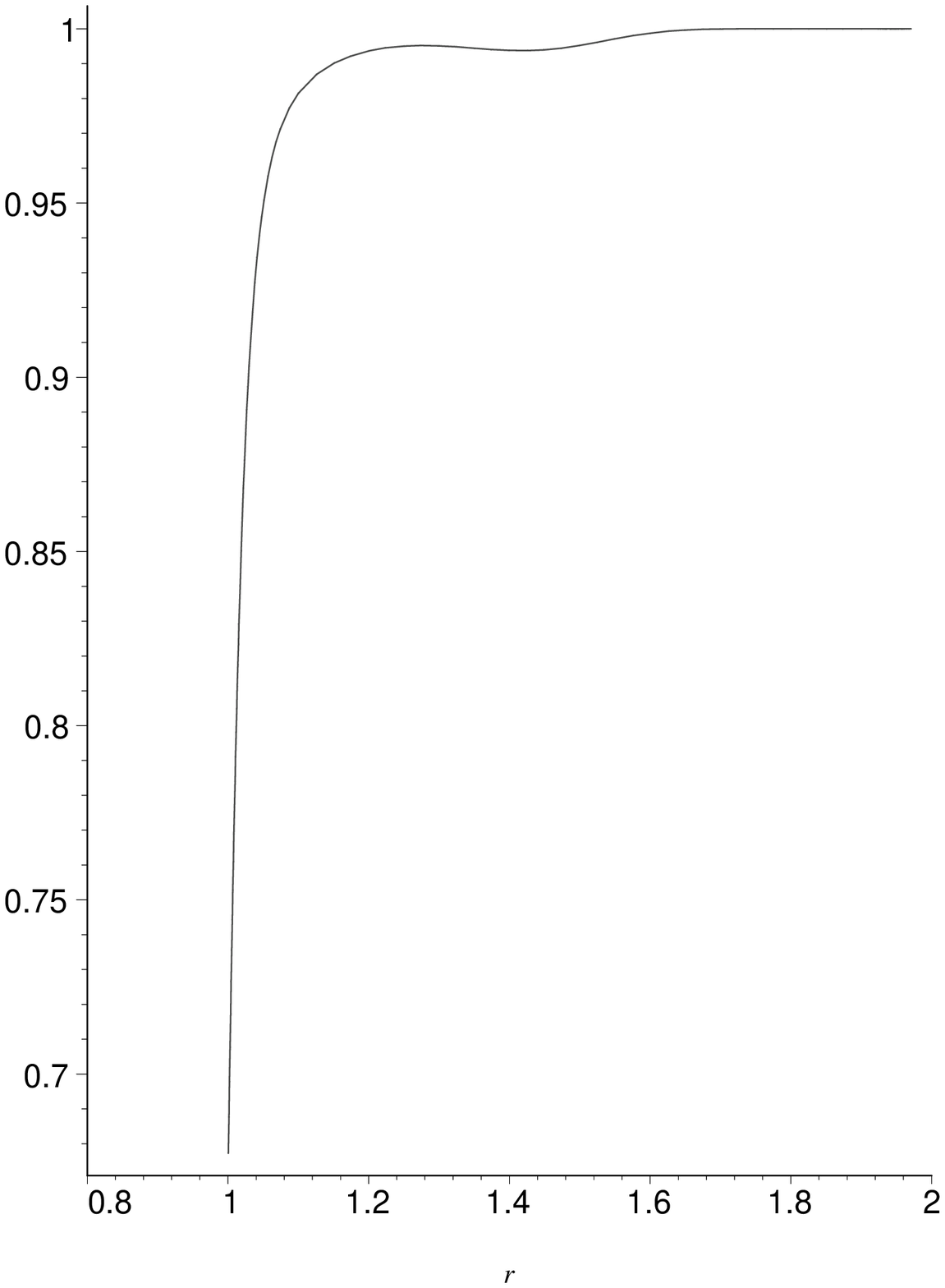, height=125pt, width=125pt}
\epsfig{figure=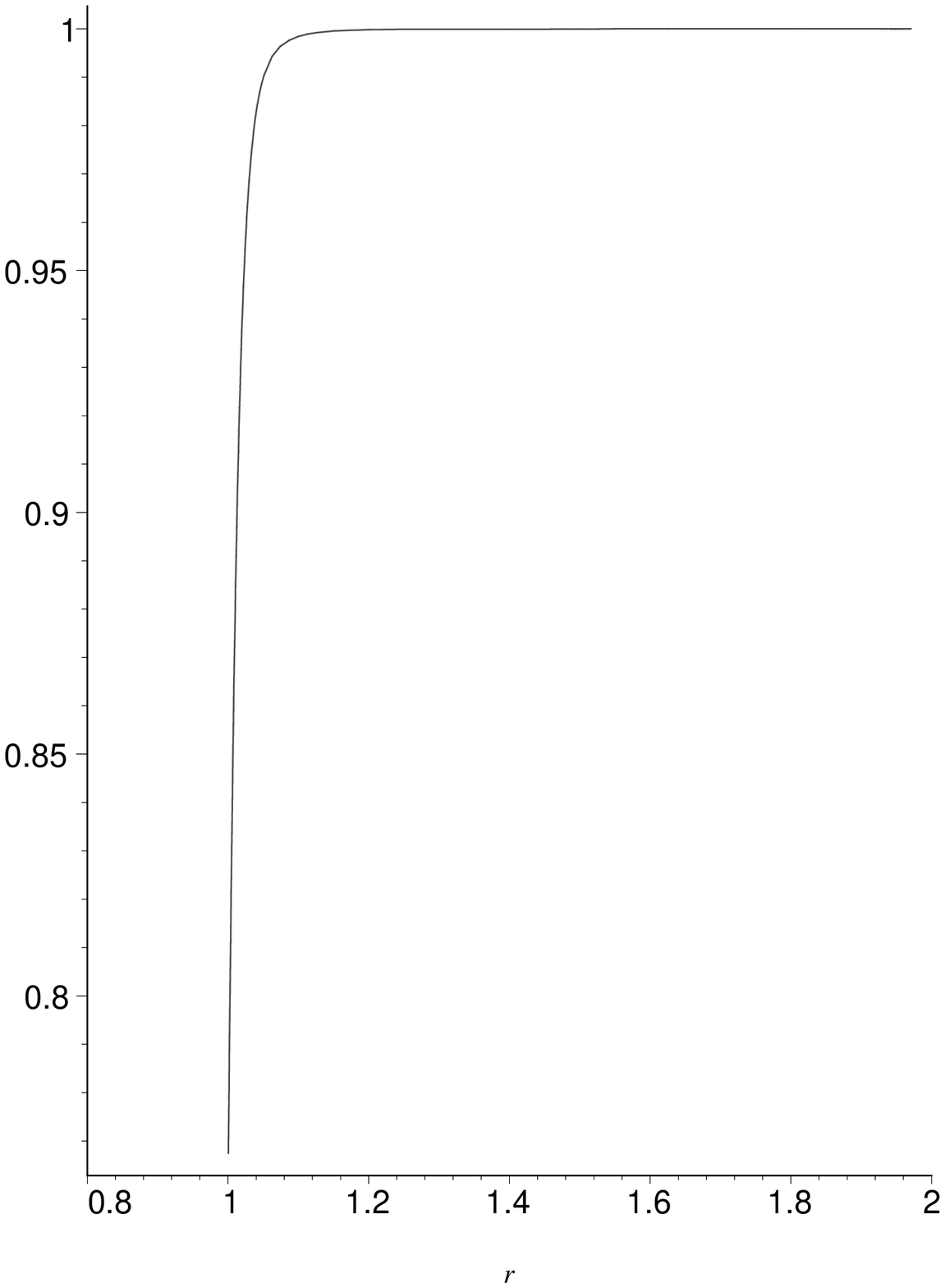, height=125pt, width=125pt}
\caption{
\label{fig:APF-seg2-plots}
Asymptotic power function
against segregation alternative $H^S_{\sqrt{3}/4}$
as a function of $r$
for $n=3$ (left) and $n=5$ (right).
}
\end{figure}

With $\epsilon=\sqrt{3}/4$, $\Pi_S(r,n,\epsilon)$ at level $\alpha=.05 $ is plotted in Figure \ref{fig:APF-seg2-plots}.  Observe that $\Pi_S(r,n,\sqrt{3}/4)\rightarrow 1$ as $r \rightarrow 2$ for $n=3,5$.  Moreover, $r^*_g(n,\sqrt{3}/4)=2$, for $n=3,\,5$.

\begin{figure}[]
\centering
\psfrag{r}{\scriptsize{$r$}}
\epsfig{figure=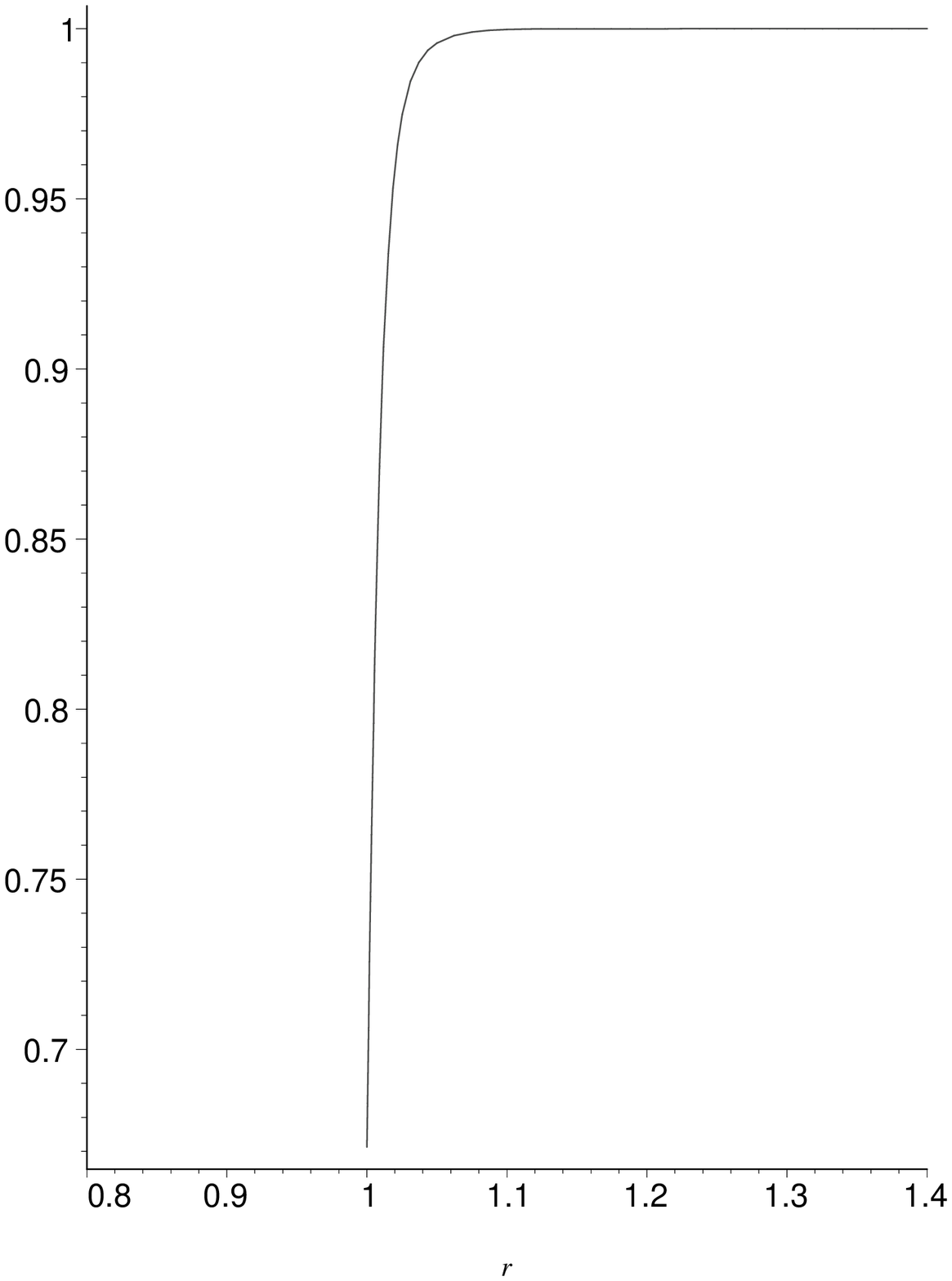, height=125pt, width=125pt}
\epsfig{figure=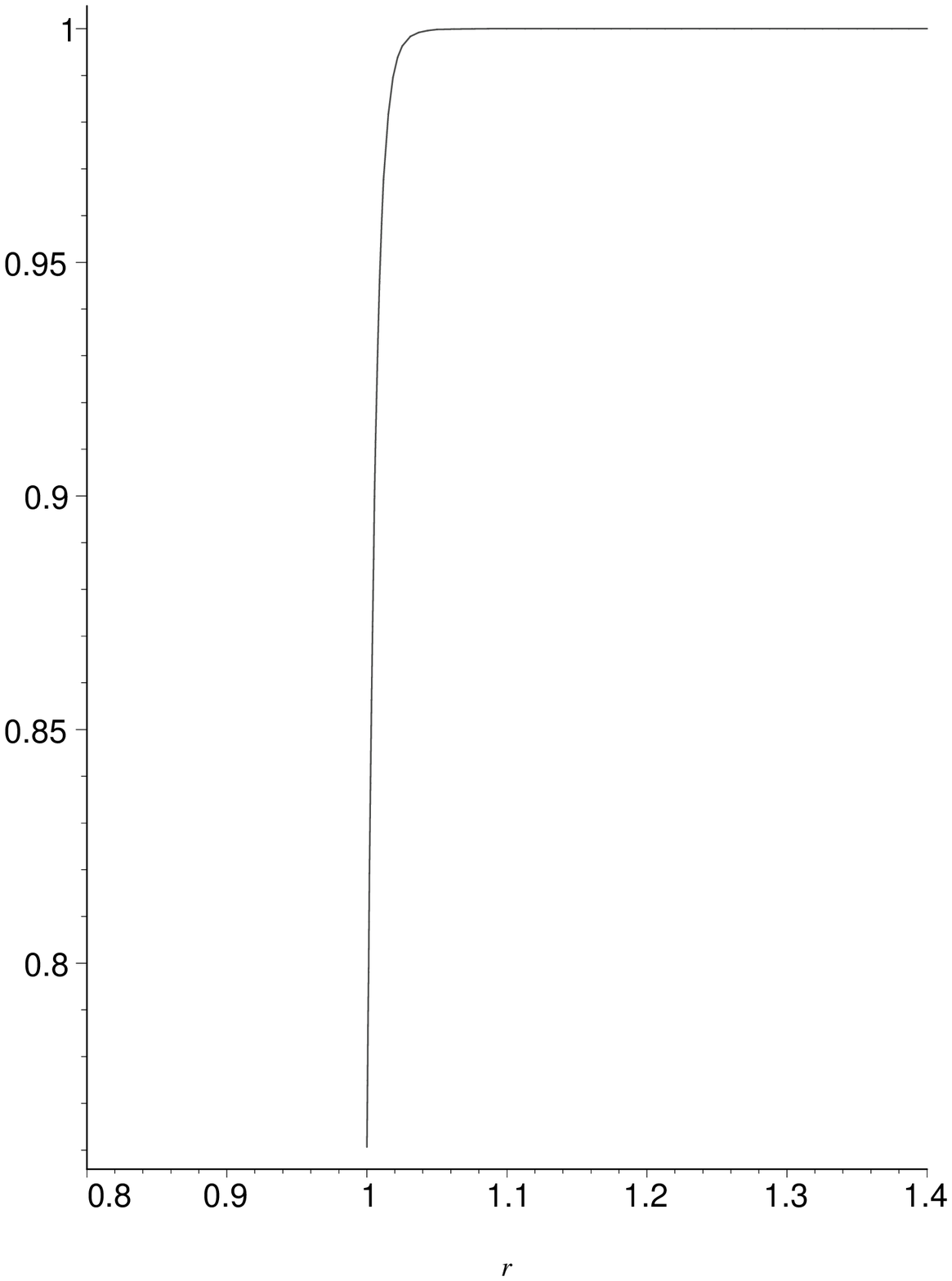, height=125pt, width=125pt}
\caption{
\label{fig:APF-seg3-plots}
Asymptotic power function
against segregation alternative $H^S_{\sqrt{3}/7}$
as a function of $r$
for $n=3$ (left) and $n=5$ (right).
}
\end{figure}

With $\epsilon=2\,\sqrt{3}/7$, $\Pi_S(r,n,\epsilon)$ at level $\alpha=.05 $ is plotted in Figure \ref{fig:APF-seg3-plots}.  Observe that $\Pi_S(r,n,2\,\sqrt{3}/7)\rightarrow 1$ as $r \rightarrow 3/2$ for $n=3,\,5$ and $r^*_g(n,2\,\sqrt{3}/7)=2$, for $n=3,\,5$.

\subsubsection{Asymptotic Power Function Analysis Under Association}
Under association, for sufficiently large $n$, we reject $H_0$ when $\sqrt{n}\,\Bigl( \frac{\rho_n(r)-\mu_A(r,\epsilon)}{\sqrt{\nu_A(r,\epsilon)}} \Bigr) < z_{\alpha}$  where $z_{\alpha}$ is the $\alpha \times 100 $ percentile of the standard normal distribution, e.g. with $\alpha=.05 $, $z_{.05} \approx -1.645$.  Then size $\alpha $ critical region for large samples is
 $$\rho_n(r)<\mu(r) + z_{\alpha} \cdot \sqrt{\nu(r)/n}.$$
 Under $H^A_{\epsilon}$,
we have
\begin{eqnarray*}
\Pi_A(r,n,\epsilon)&:=&P\Bigl(\rho_n(r)<\mu(r) + z_{\alpha} \cdot \sqrt{\nu(r)/n}\Bigr)\\
&=&\Phi\Biggl(\frac{z_{\alpha}\,\sqrt{\nu(r)}}{\sqrt{\nu_A(r,\epsilon)}}+\frac{\sqrt{n}\,\bigl( \mu(r)-\mu_A(r,\epsilon) \bigr)}{\sqrt{\nu_A(r,\epsilon)}}\Biggr).
\end{eqnarray*}

\begin{figure}[ht]
\centering
\psfrag{r}{\scriptsize{$r$}}
\epsfig{figure=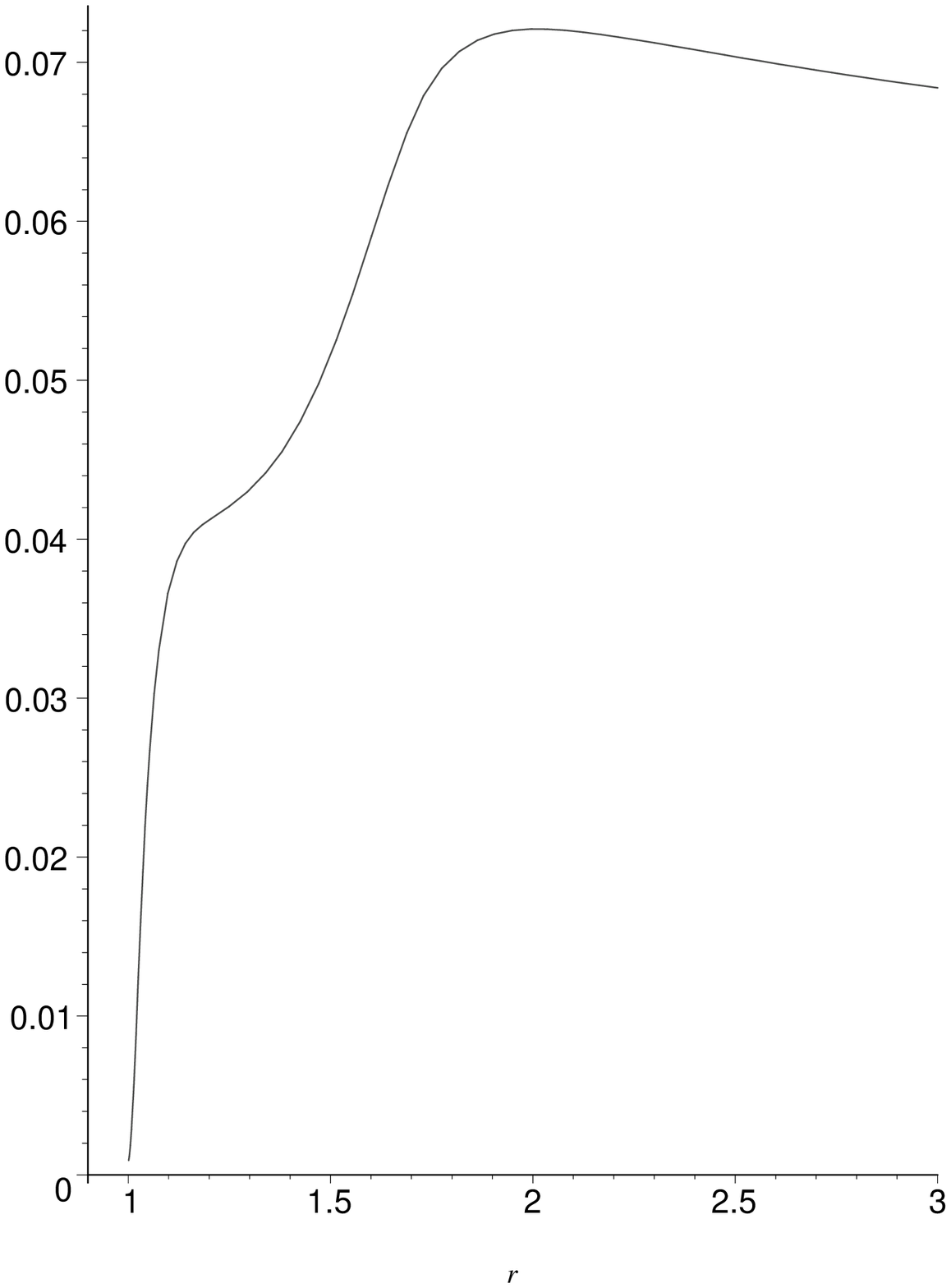, height=125pt, width=125pt}
\epsfig{figure=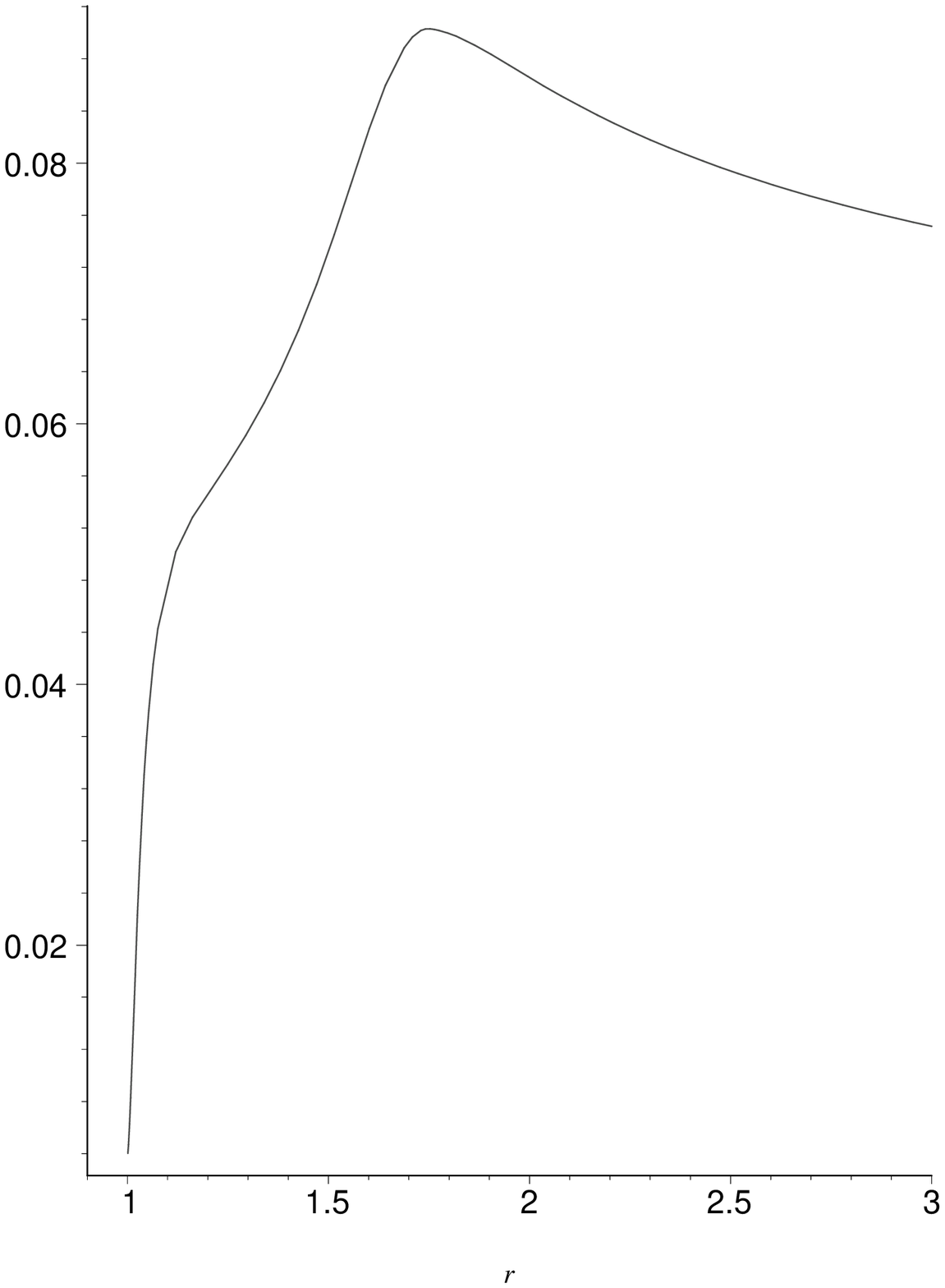, height=125pt, width=125pt}
\epsfig{figure=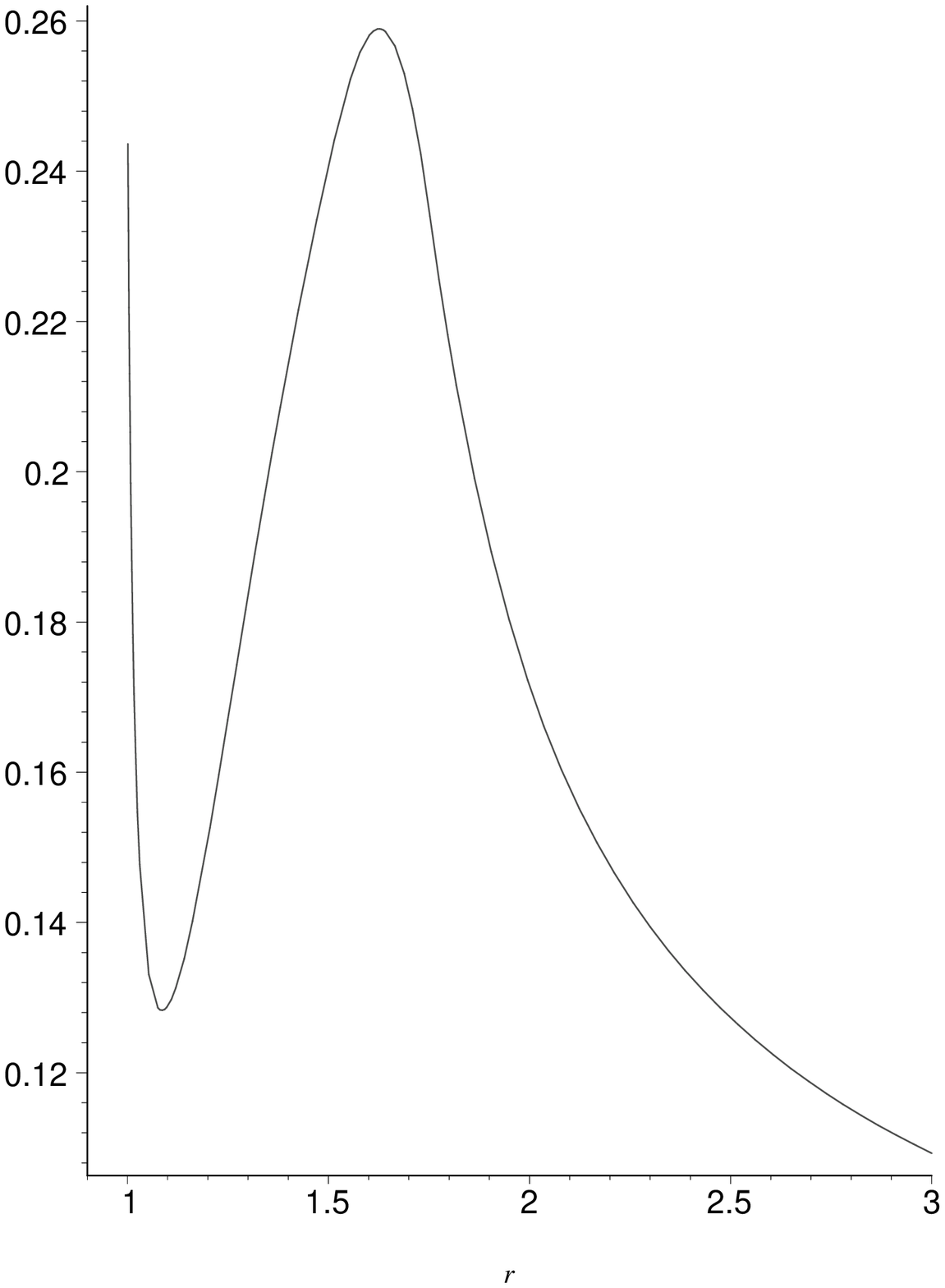, height=125pt, width=125pt}
\caption{
\label{fig:APF-agg3-plots}
Asymptotic power function
against association alternative $H^A_{\sqrt{3}/21}$
as a function of $r$
for $n=5,10,100$}
\end{figure}

With $\epsilon=\sqrt{3}/21$, $\Pi_A(r,n,\epsilon)$ at level $\alpha=.05 $ is plotted in Figure \ref{fig:APF-agg3-plots}.  Observe that $\Pi_A\left( r,n,\sqrt{3}/21 \right)\rightarrow .057$ as $r \rightarrow \infty$ for $n=5,\,10,\,100$.  Let $\widehat{r}\left( n,\epsilon \right)$ be the value at which $\Pi_A(r,n,\epsilon)$ attains its supremum. Then, $\widehat{r}\left( 5,\sqrt{3}/21 \right)\approx 2.01$, and $\widehat{r}\left( 10,\sqrt{3}/21 \right)\approx 1.875$, and $\widehat{r}\left( 100,\sqrt{3}/21 \right)\approx 1.645$. Moreover, $\Pi_A\left( r,100,\sqrt{3}/21 \right)$ attains a local infimum at $\approx 1.065$.

\begin{figure}[]
\centering
\psfrag{r}{\scriptsize{$r$}}
\epsfig{figure=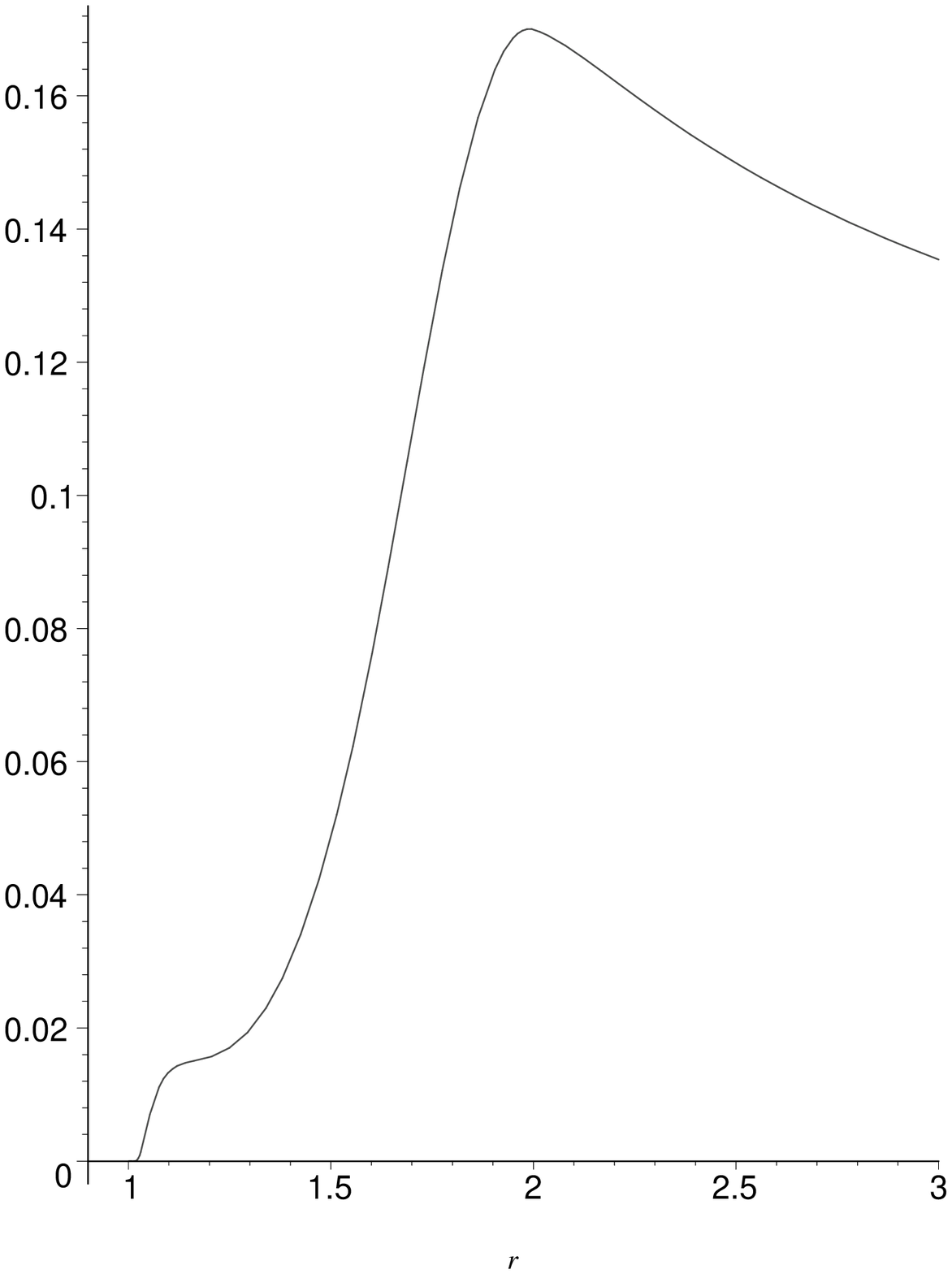, height=125pt, width=125pt}
\epsfig{figure=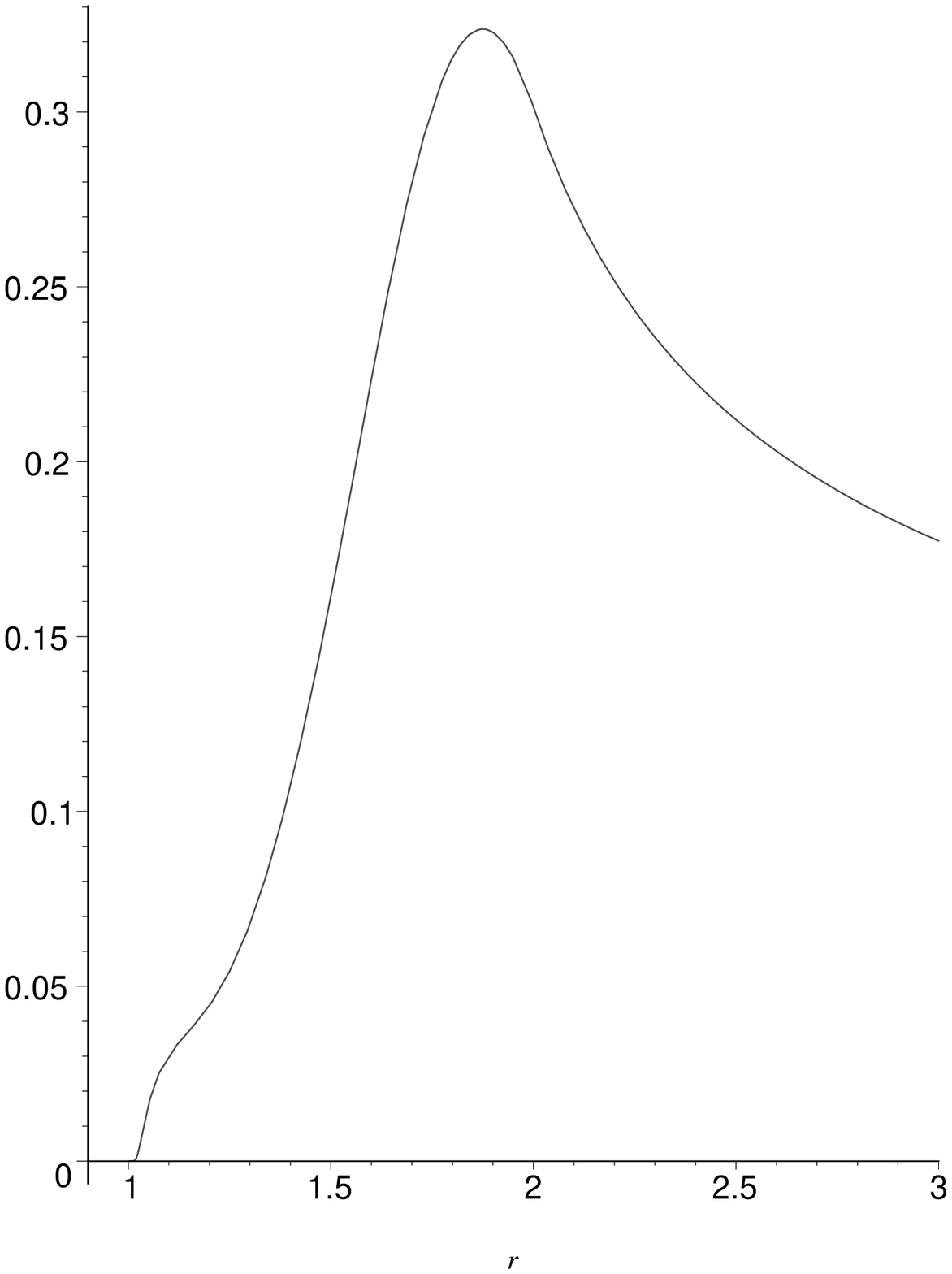, height=125pt, width=125pt}
\epsfig{figure=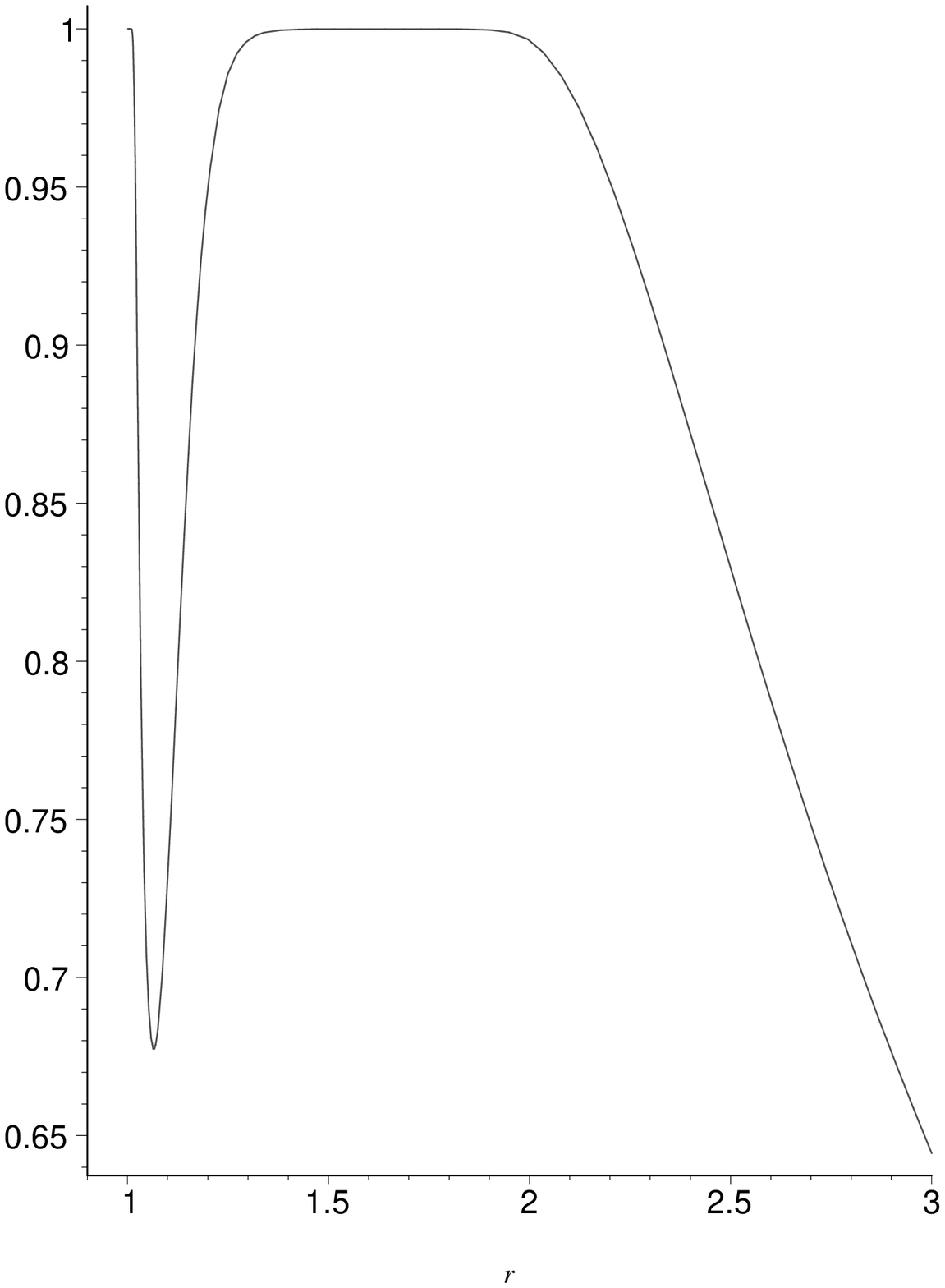, height=125pt, width=125pt}
\caption{
\label{fig:APF-agg2-plots}
Asymptotic power function
against association alternative $H^A_{\sqrt{3}/12}$
as a function of $r$
for $n=5,10,100$}
\end{figure}

With $\epsilon=\sqrt{3}/12$, $\Pi_A(r,n,\epsilon)$ at level $\alpha=.05 $ is plotted in Figure \ref{fig:APF-agg2-plots}.  Observe that $\Pi_A\left( r,n,\sqrt{3}/12 \right)\rightarrow .0766$ as $r \rightarrow \infty$ for $n=5,\,10,\,100$. Moreover, $\widehat{r}\left( 5,\sqrt{3}/12 \right)\approx 1.99$, $\widehat{r}\left( 10,\sqrt{3}/12 \right)\approx 1.75$, and $\widehat{r}\left( 100,\sqrt{3}/12 \right)\approx 1.60$. Moreover, $\Pi_A\left( r,100,\sqrt{3}/21 \right)$ attains a local infimum at $\approx 1.105$.

\begin{figure}[]
\centering
\psfrag{r}{\scriptsize{$r$}}
\epsfig{figure=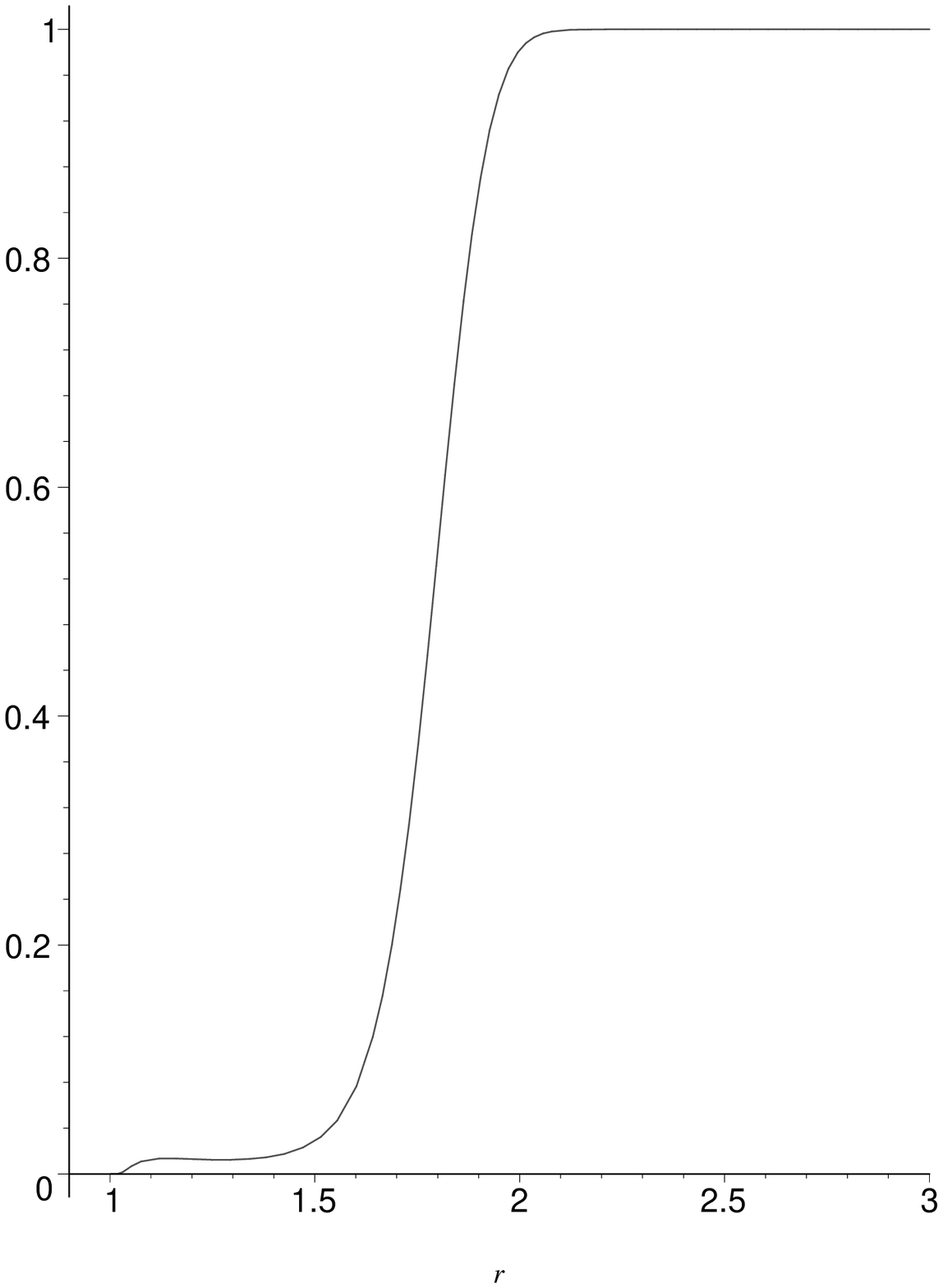, height=125pt, width=125pt}
\epsfig{figure=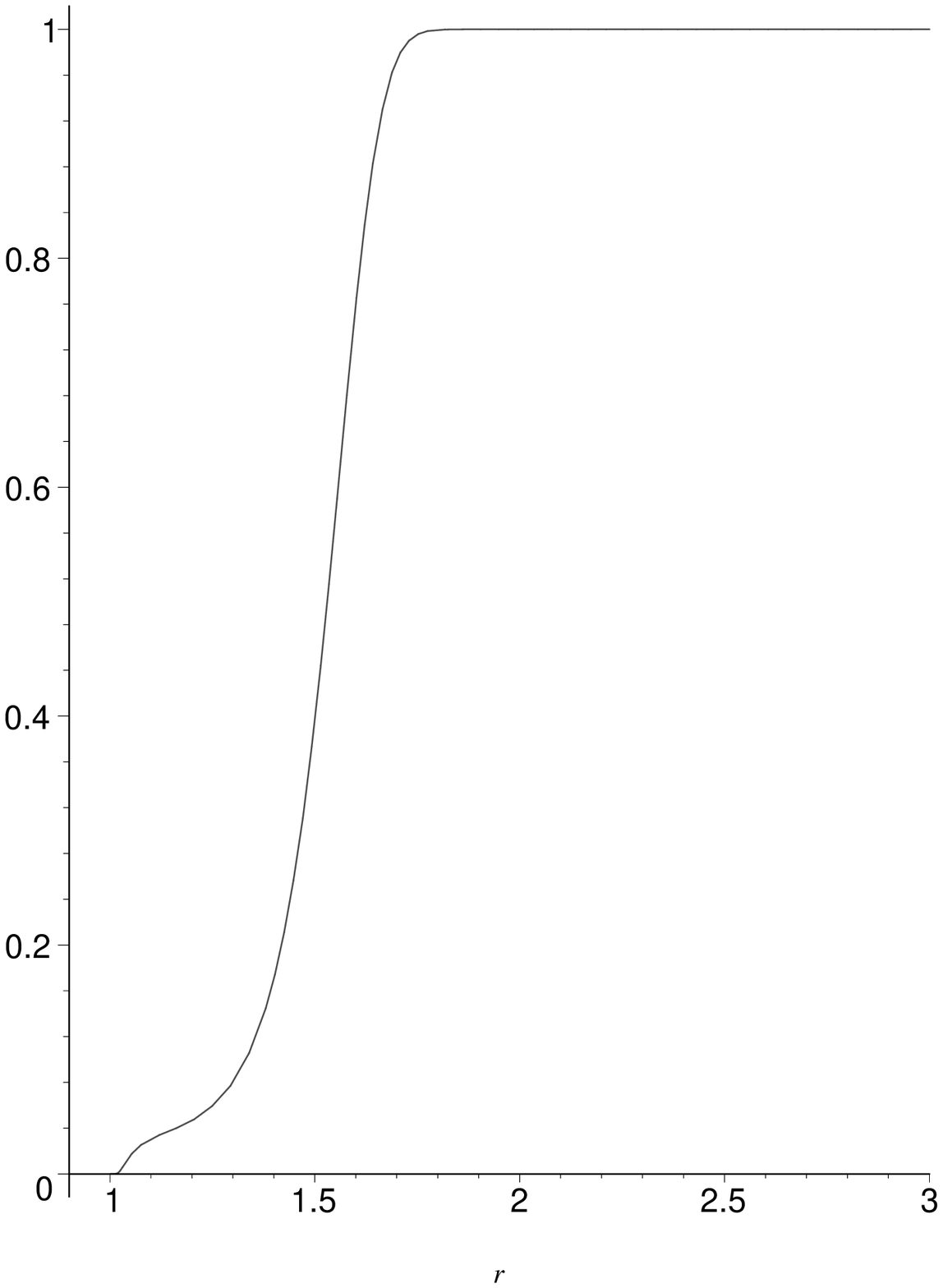, height=125pt, width=125pt}
\caption{
\label{fig:APF-agg-plots}
Asymptotic power function
against association alternative $H^A_{5\,\sqrt{3}/24}$
as a function of $r$ for $n=5,10$.
}
\end{figure}

With $\epsilon=5\,\sqrt{3}/24$, $\Pi_A(r,n,\epsilon)$ at level $\alpha=.05 $ is plotted in Figure \ref{fig:APF-agg-plots}.  Observe that $\Pi\left( r,n,5\,\sqrt{3}/24 \right)\rightarrow 1$ as $r \rightarrow 2$ for $n=5,\,10$.

\subsection{Multiple Triangle Case}

Suppose $\Y$ is a finite collection of points in $\R^2$ with $|\Y| \ge 3$.
Consider the Delaunay triangulation (assumed to exist) of $\Y$,
where $T_j$ denotes the $j^{th}$ Delaunay triangle,
$J$ denotes the number of triangles, and
$C_H(\Y)$ denotes the convex hull of $\Y$.
We wish to test
$$H_0: X_i \stackrel{iid}{\sim} \U(C_H(\Y))$$
against segregation and association alternatives.

\begin{figure}[]
\centering
\epsfig{figure=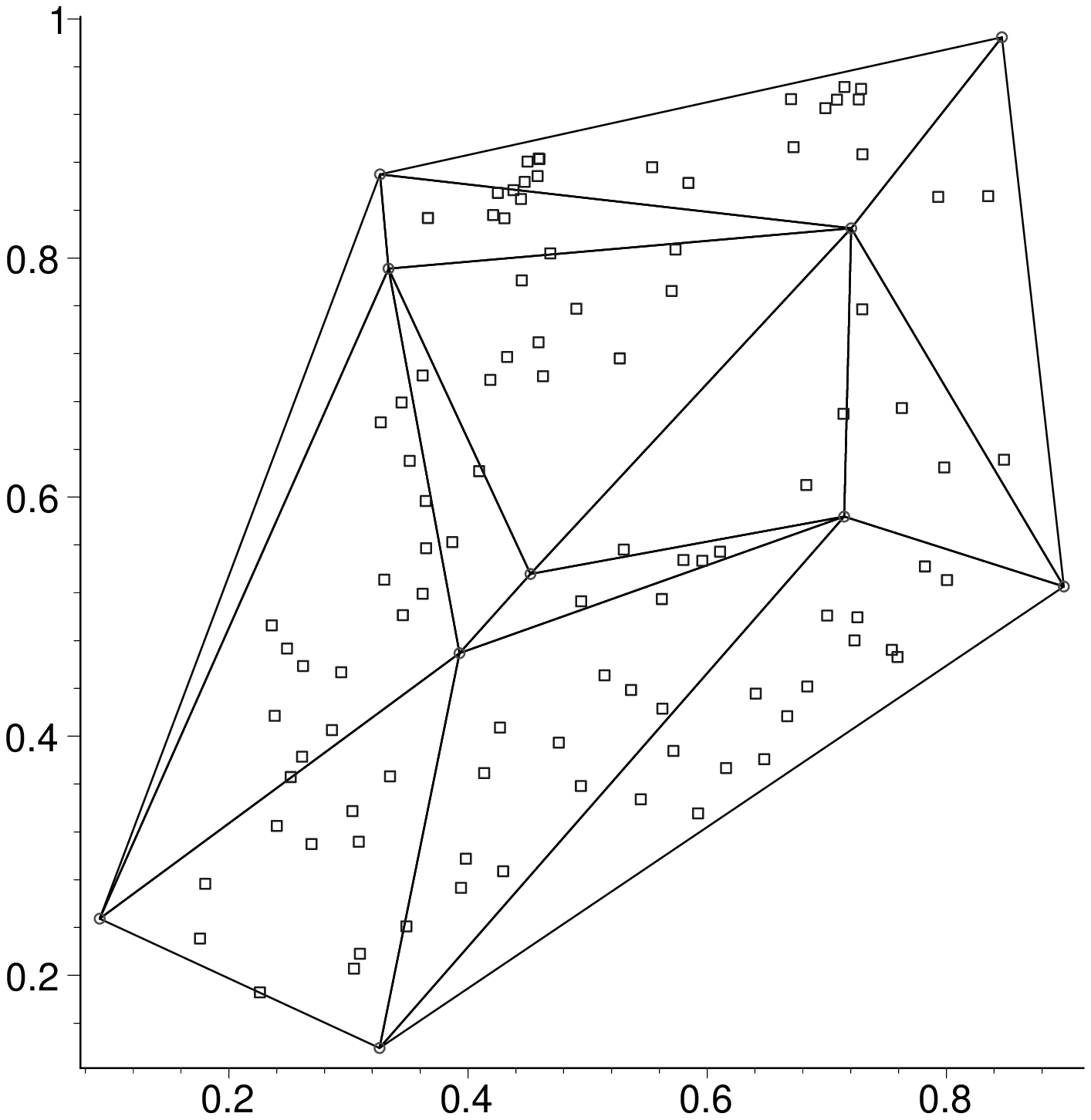, height=125pt, width=125pt}
\epsfig{figure=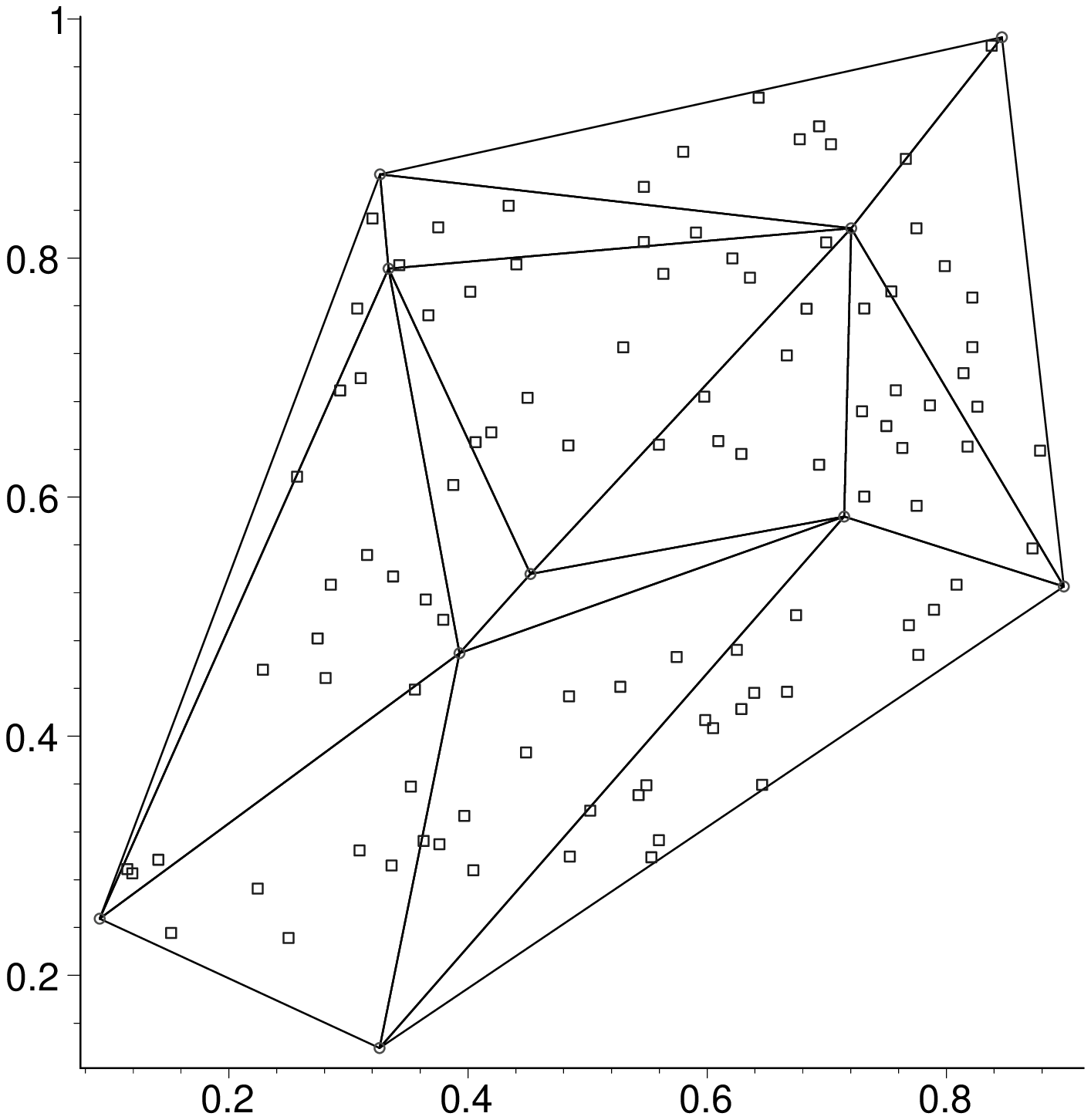, height=125pt, width=125pt}
\epsfig{figure=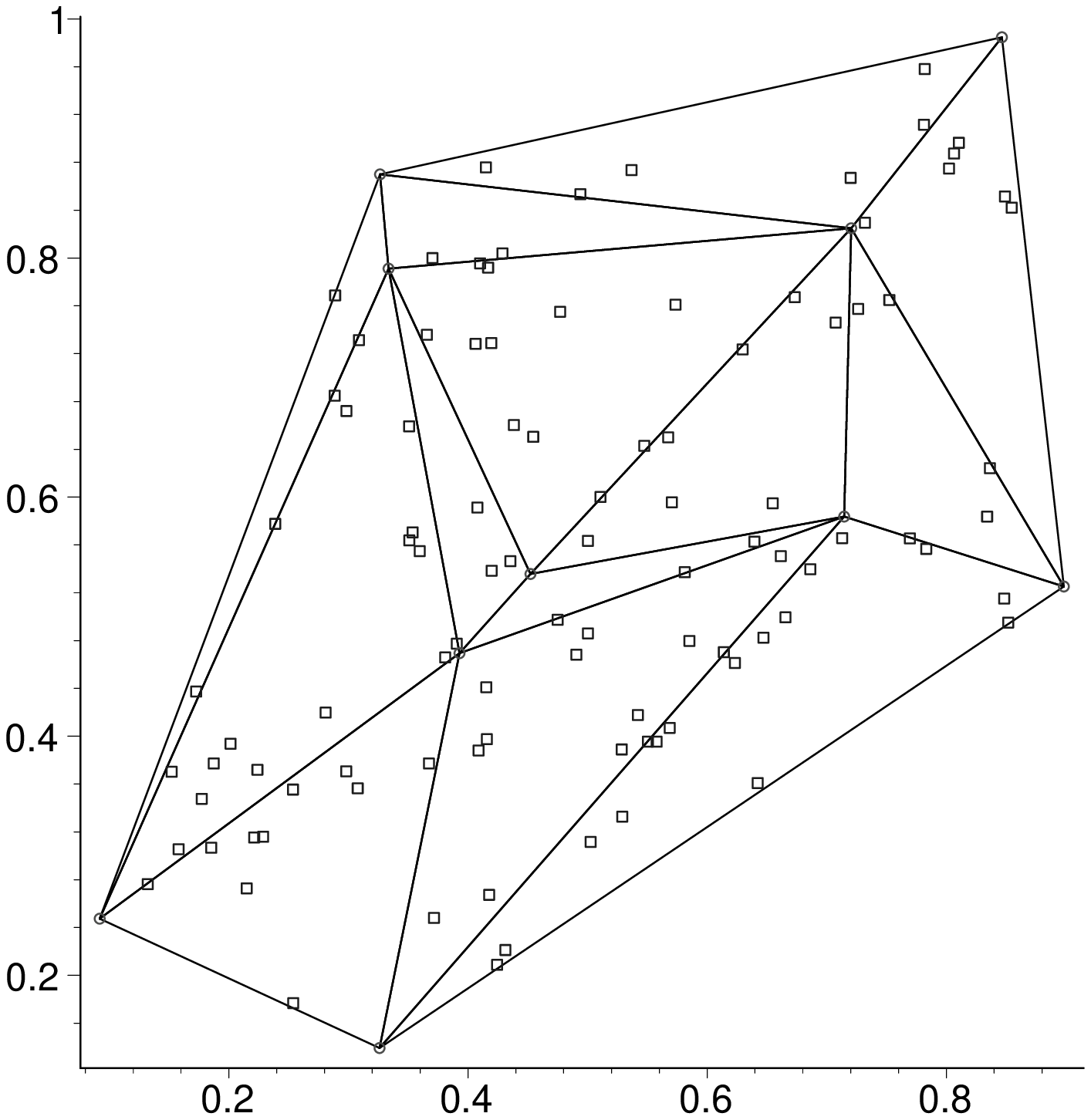, height=125pt, width=125pt}
\caption{\label{fig:deldata1}
Realization of segregation (left), $H_0$ (middle), and association (right) for $|\Y|=10$, $J=13$, and $n=100$.
}
\end{figure}

\begin{figure}[ht]
\centering
\epsfig{figure=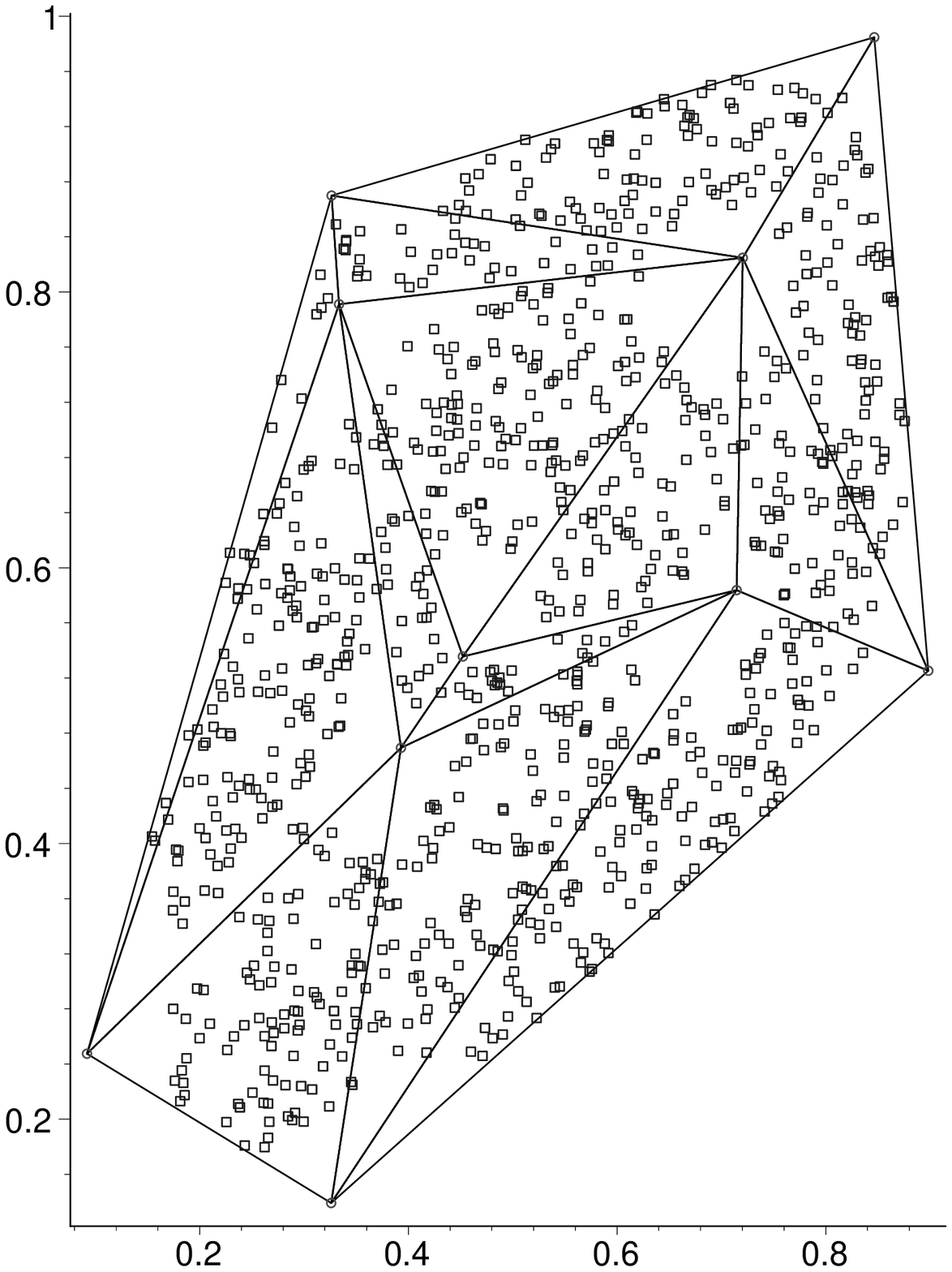, height=125pt, width=125pt}
\epsfig{figure=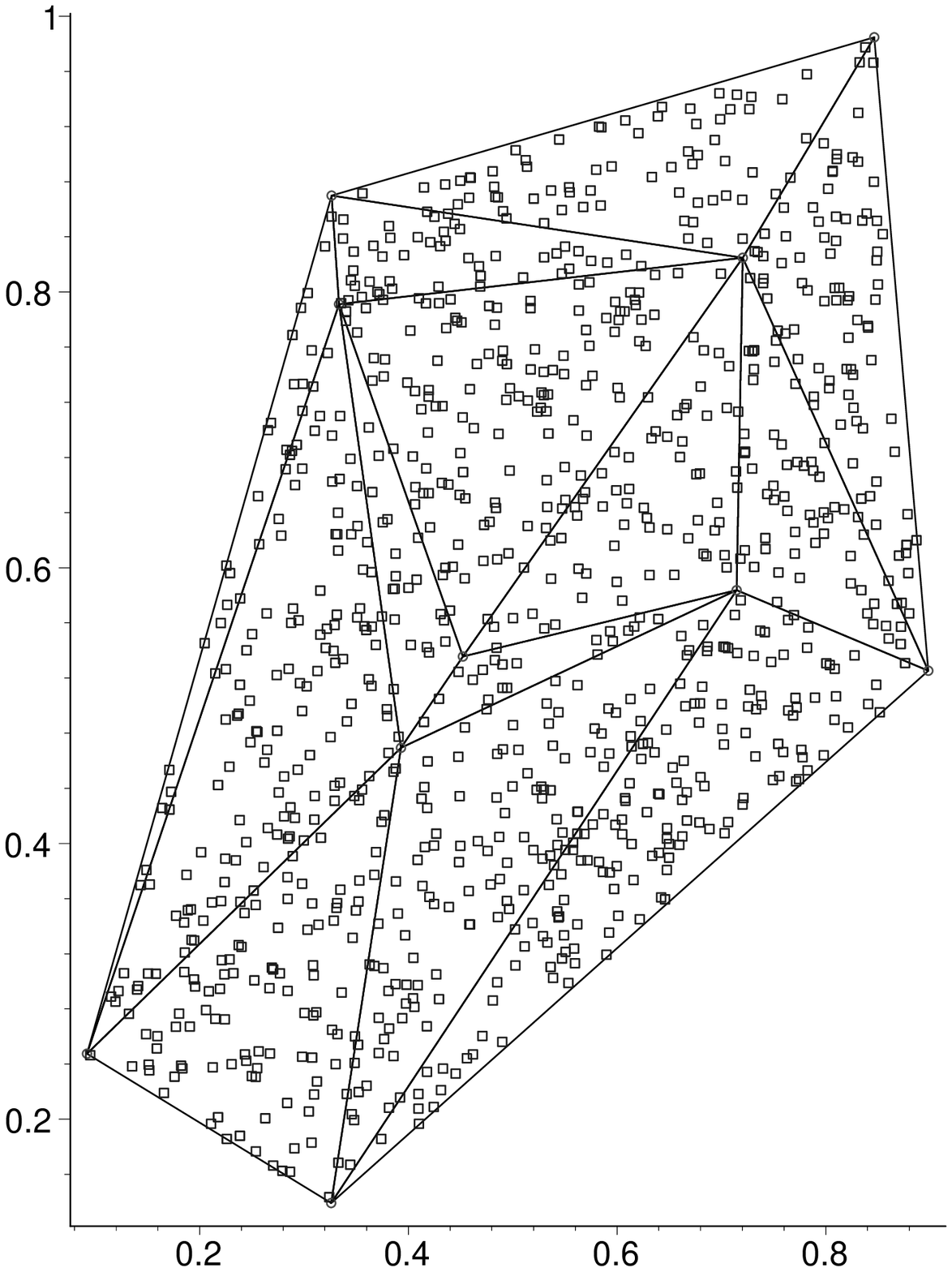,      height=125pt, width=125pt}
\epsfig{figure=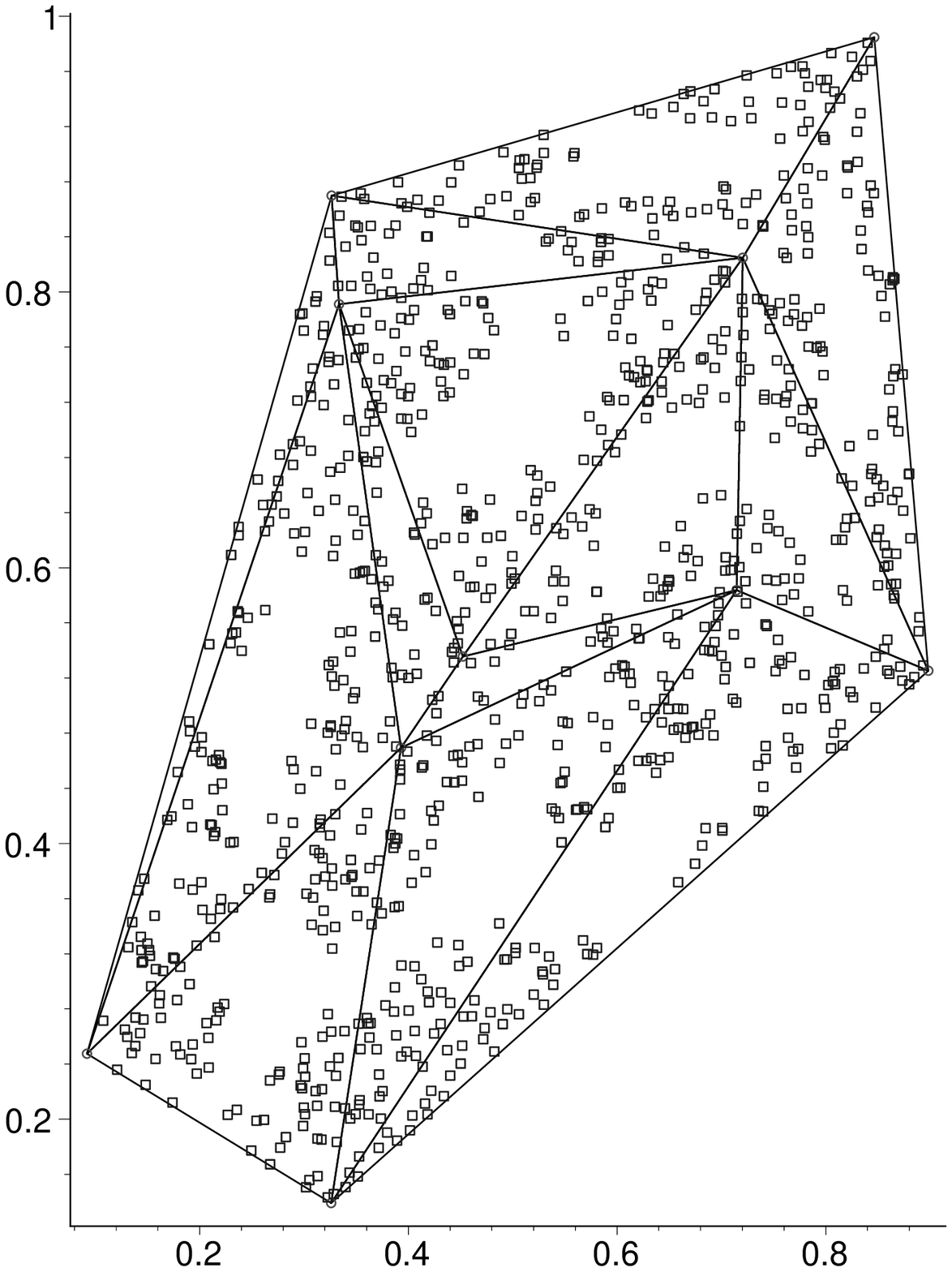, height=125pt, width=125pt}
\caption{\label{fig:deldata}
Realization of segregation (left), $H_0$ (middle), and association (right) for $|\Y|=10$, $J=13$, and $n=1000$.
}
\end{figure}

Figure \ref{fig:deldata1} and Figure \ref{fig:deldata} are graphs of
realizations of $n=100$ and $n=1000$ observations which are
independent and identically distributed
according to $\U(C_H(\Y))$
for $|\Y|=10$ and $J=13$ and realizations of $n=100$ and $n=1000$ observations
under segregation and association
for the same $\Y$, respectively.

The digraph $D$ is constructed using $N_{\Y_j}^r(\cdot)$ as described in Section \ref{sec:data-random-PCD}, here for $X_i \in T_j$ the three points in $\Y$ defining the
Delaunay triangle $T_j$ are used as $\Y_j$.
Let $\rho_n(r,J)$ be the relative density of the digraph based on $\X_n$ and $\Y$ which yields $J$ Delaunay triangles, and let $w_j = A(T_j) / A(C_H(\Y))$ for $j=1,\ldots,J$, where $A(C_H(\Y))=\sum_{j=1}^{J}A(T_j)$ with $A(\cdot)$ being the area functional. Then we obtain the
following as a corollary to Theorem 2.

{\bf Corollary 1}
\label{cor:asy-norm-rho-NYr}
The asymptotic null distribution for $\rho_n(r,J)$ conditional on $\mathcal W=\{w_1,\ldots,w_J\}$
for $r \in [1,\infty)$
is given by $\N(\mu(r,J),\nu(r,J)/n)$ provided that $\nu(r,J)>0$ with
\begin{equation}
\mu(r,J):=\mu(r) \,\sum_{j=1}^{J}w_j^2\; \; \text{ and } \;\; \nu(r,J):= \nu(r) \,\sum_{j=1}^{J}w_j^3 +4\,\mu(r)^2\left[\sum_{j=1}^{J}w_j^3-\left(\sum_{j=1}^{J}w_j^2 \right)^2\right],
\end{equation}
where $\mu(r)$ and $\nu(r)$ are given in Equations (\ref{eq:Asymean}) and (\ref{eq:Asyvar}),
respectively.

{\bfseries Proof:} See Appendix 5 for the proof $\blacksquare$

By an appropriate application of Jensen's inequality, we see that $\sum_{j=1}^{J}w_j^3 \ge \left(\sum_{j=1}^{J}w_j^2 \right)^2.$
Therefore, the covariance  $\nu(r,J)=0$ iff both $\nu(r)=0$ and $\sum_{j=1}^{J}w_j^3=\left(\sum_{j=1}^{J}w_j^2 \right)^2$ hold, so asymptotic normality may hold even when $\nu(r)=0$.

Similarly, for the segregation (association) alternatives with $4\,\,\epsilon^2/3 \times 100 \%$ of the triangles around the vertices of each triangle is forbidden (allowed), we obtain the above asymptotic distribution of $\rho_n(r,J)$ with $\mu(r,J)$ being replaced by $\mu(r,J,\epsilon)$, $\nu(r,J)$ by $\nu(r,J,\epsilon)$, $\mu(r)$ by $\mu(r,\epsilon)$, and $\nu(r)$ by $\nu(r,\epsilon)$.

Thus in the case of $J>1$, we have a (conditional) test of
$H_0: X_i \stackrel{iid}{\sim} \U(C_H(\Y))$
which once again
rejects against segregation for large values of $\rho_n(r,J)$ and
rejects against association for small values of $\rho_n(r,J)$.

The segregation (with $\delta=1/16$ i.e. $\epsilon=\sqrt{3}/8$), null, and association (with $\delta=1/4$ i.e. $\epsilon=\sqrt{3}/12$) realizations (from left to right) are depicted in Figure \ref{fig:deldata1} with $n=100$ and in Figure \ref{fig:deldata} with $n=1000$.  For both $n=100$ and $n=1000$, for the null realization,
the $p$-value is greater than 0.1 for all $r$ values and both alternatives.
For the segregation realization,
with $n=100$ we obtain
$p < 0.001$ for $r \le 3$,
$p = 0.025$ for $r=5$, and
$p > 0.1$ for $r \ge 10$
and
with $n=1000$ we obtain
$p < 0.0031$ for $ 1< r \le 5$
and
$p > 0.24$ for $r=1$ and $r\ge 10$.
For the association realization,
with $n=100$,
we obtain
$p < 0.05$ for $r=1.5,\,2$,
and
$p > 0.06$ for for other values of $r$
and
with $n=1000$,
we obtain
$p < 0.0135$ for $1<r \le 3$,
$p=.14$ for $r=1$, and
$p > 0.25$ for for $r \ge 5$. Note that this is only for one realization of $\X_n$.

\begin{figure}[]
\centering
\psfrag{power}{ \Huge{\bfseries{power}}}
\psfrag{r}{\Huge{$r$}}
\rotatebox{-90}{ \resizebox{1.8 in}{!}{ \includegraphics{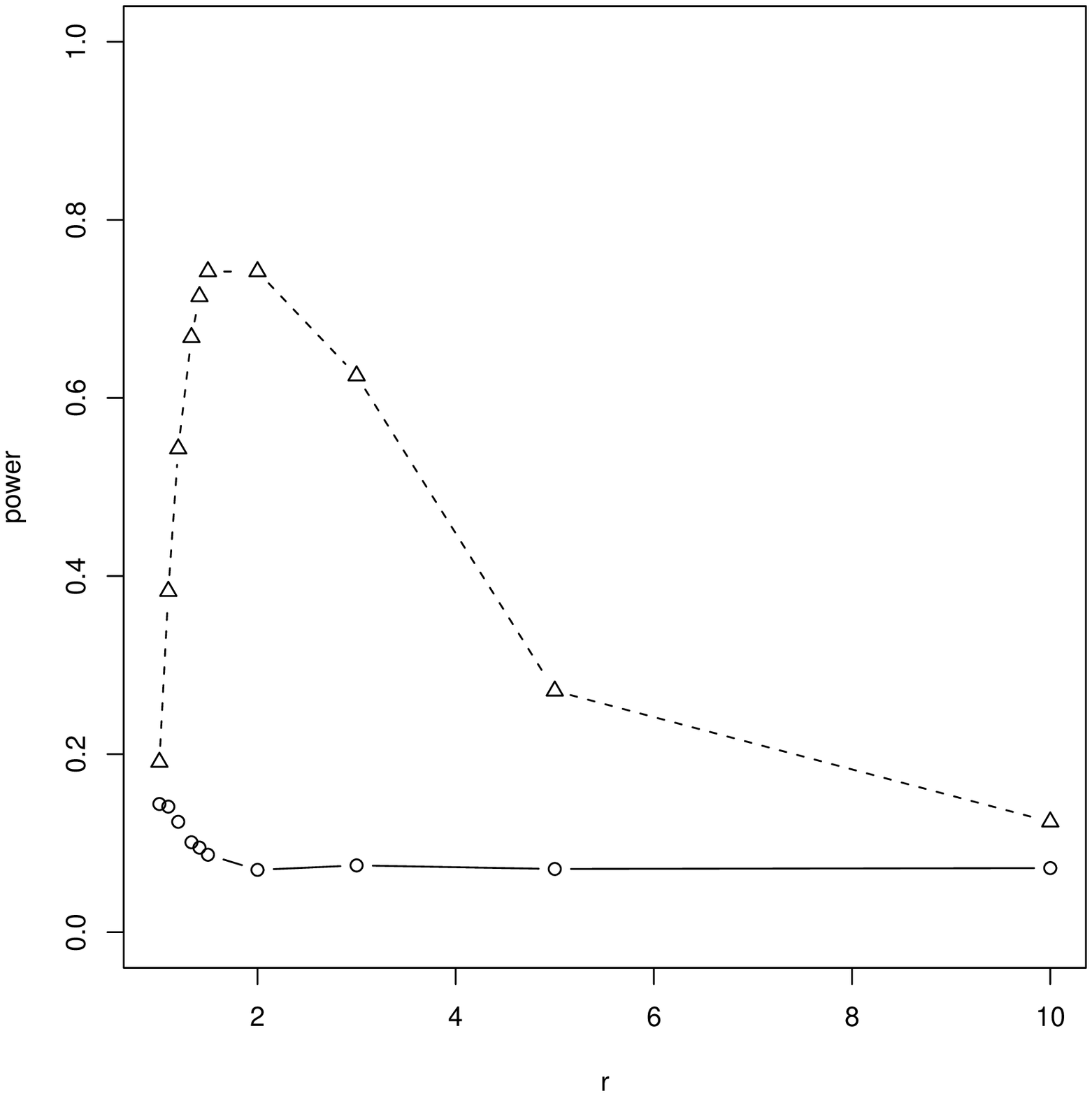}}}
\rotatebox{-90}{ \resizebox{1.8 in}{!}{ \includegraphics{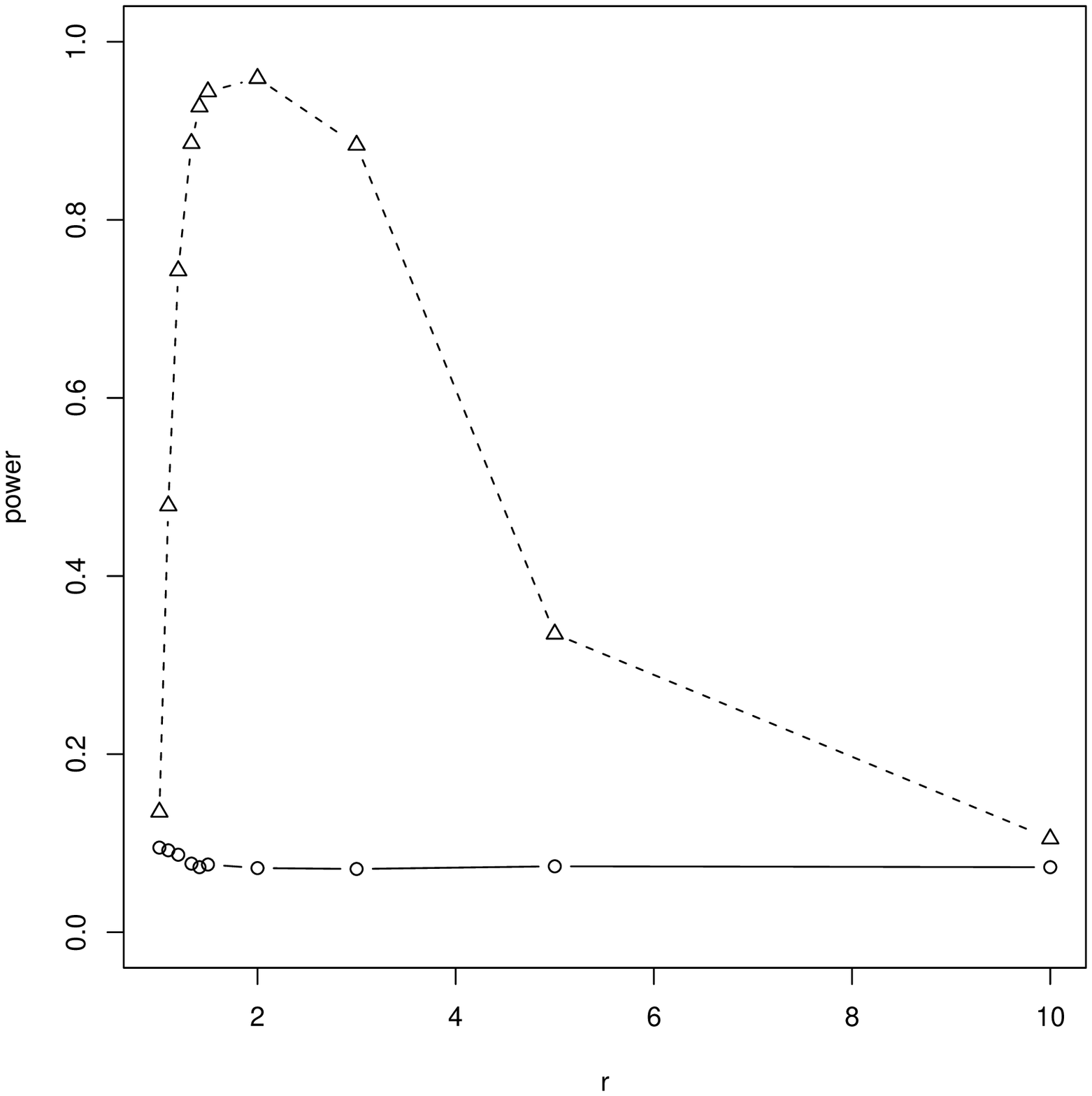}}}
\rotatebox{-90}{ \resizebox{1.8 in}{!}{ \includegraphics{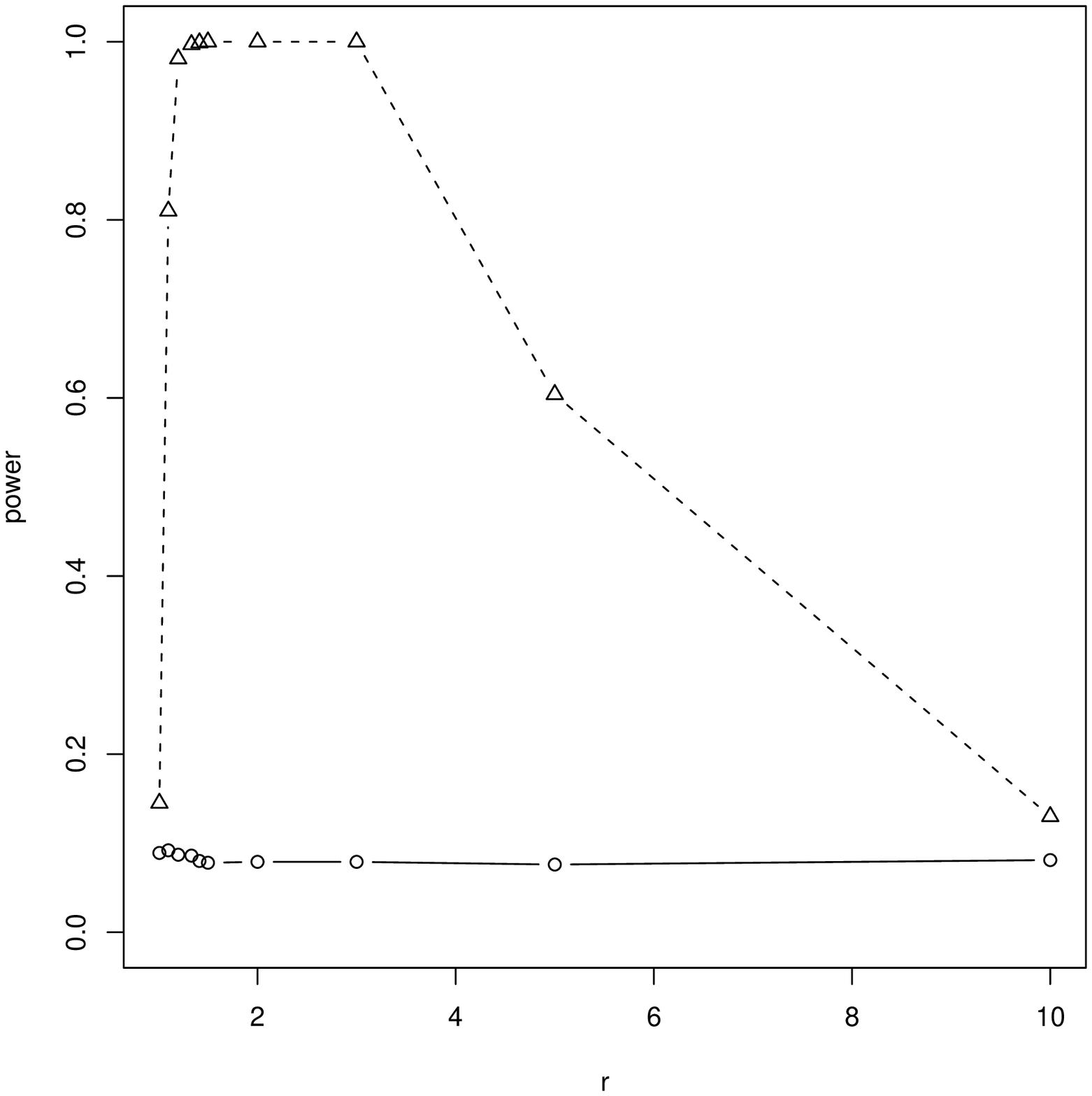}}}
\caption{
 \label{fig:MTSegSim}
Monte Carlo power using the asymptotic critical value against $H^S_{\sqrt{3}/8}$,
as a function of $r$, for $n=100$ (left), $n=200$ (middle), and $n=500$ (right) conditional on the realization of $\Y$ in Figure \ref{fig:deldata1}.  The circles represent the empirical significance levels while triangles represent the empirical power values.}
\end{figure}

\begin{figure}[]
\centering
\psfrag{power}{ \Huge{\bfseries{power}}}
\psfrag{r}{\Huge{$r$}}
\rotatebox{-90}{ \resizebox{1.8 in}{!}{ \includegraphics{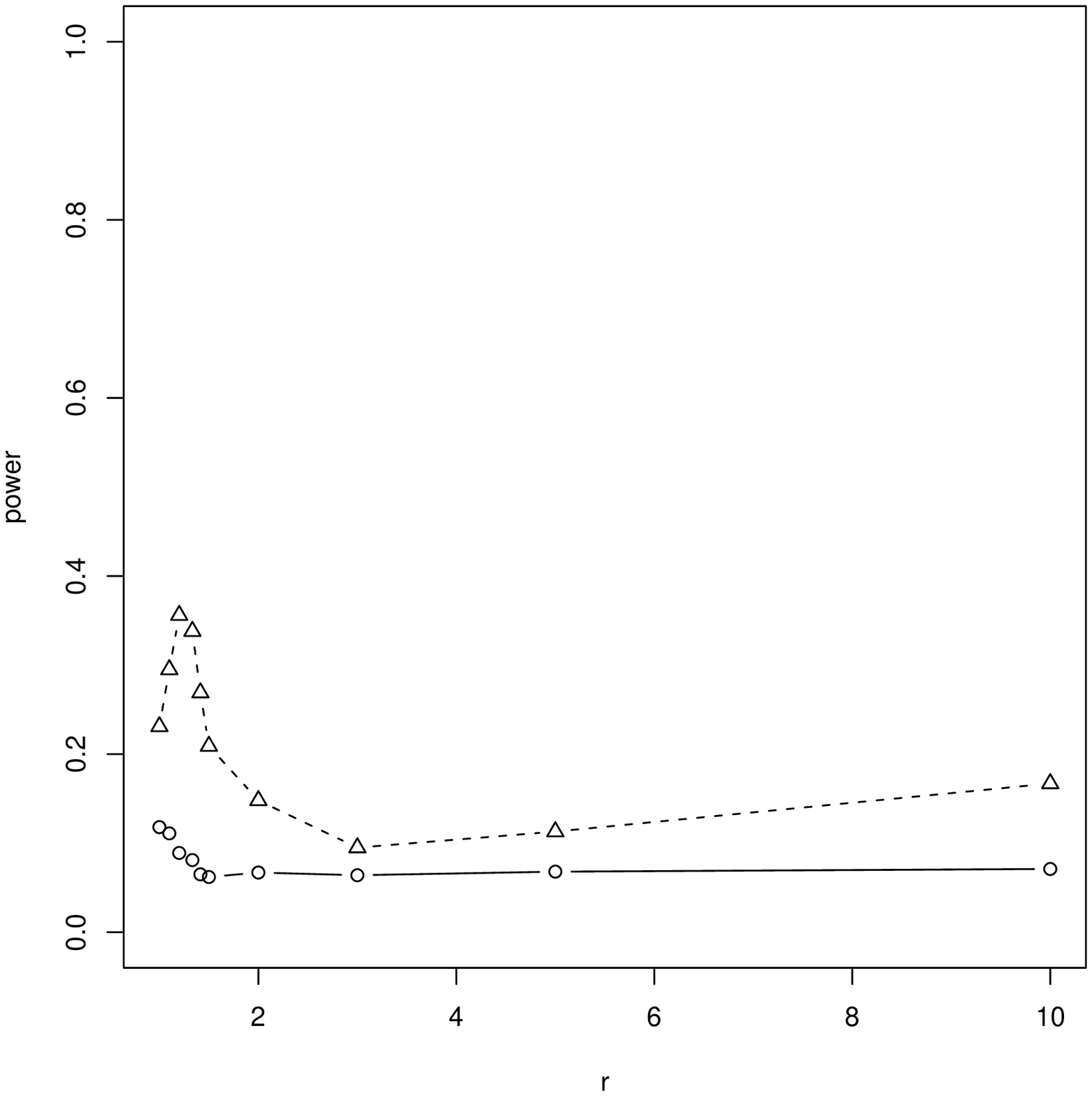}}}
\rotatebox{-90}{ \resizebox{1.8 in}{!}{ \includegraphics{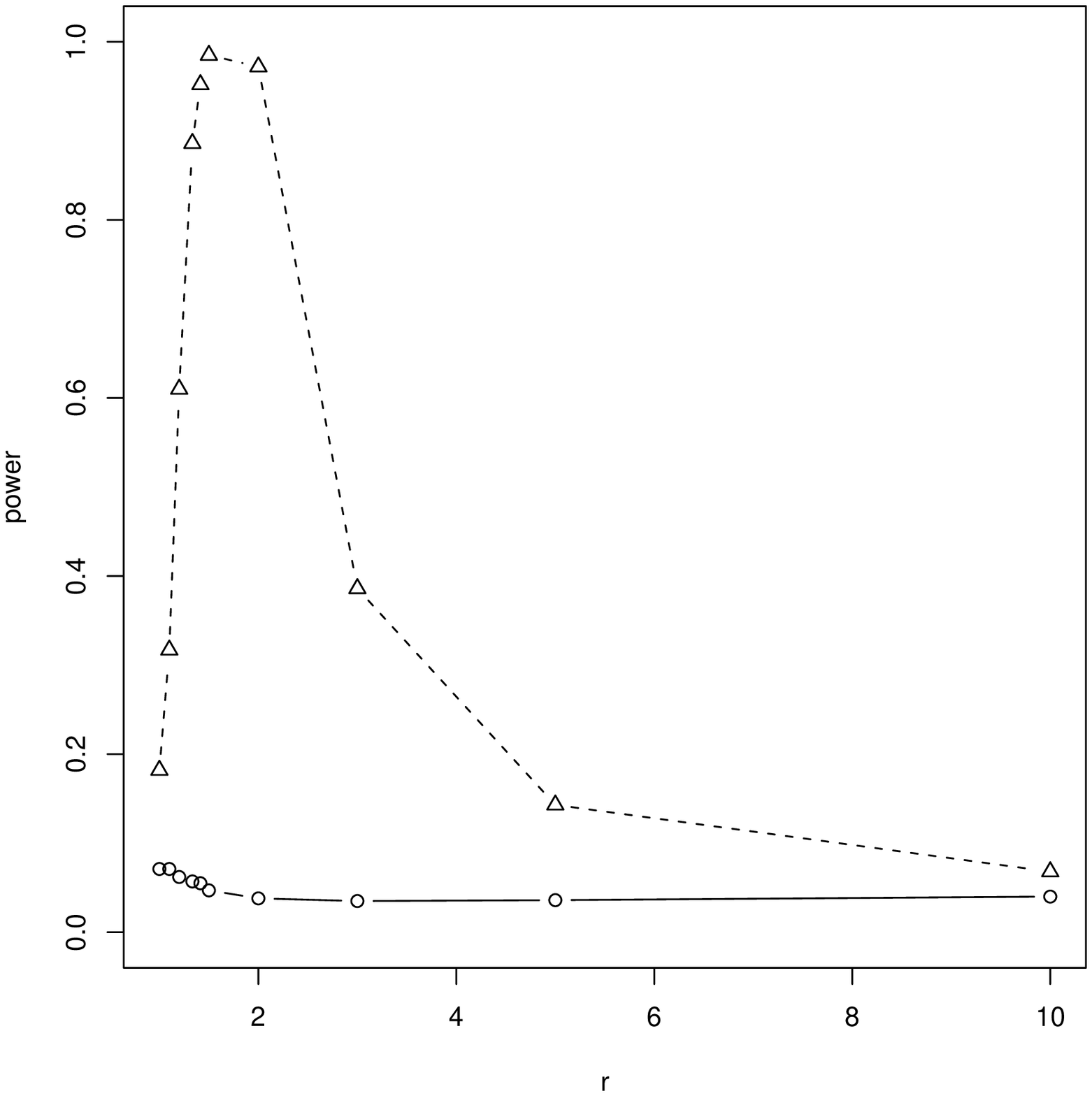}}}
\rotatebox{-90}{ \resizebox{1.8 in}{!}{ \includegraphics{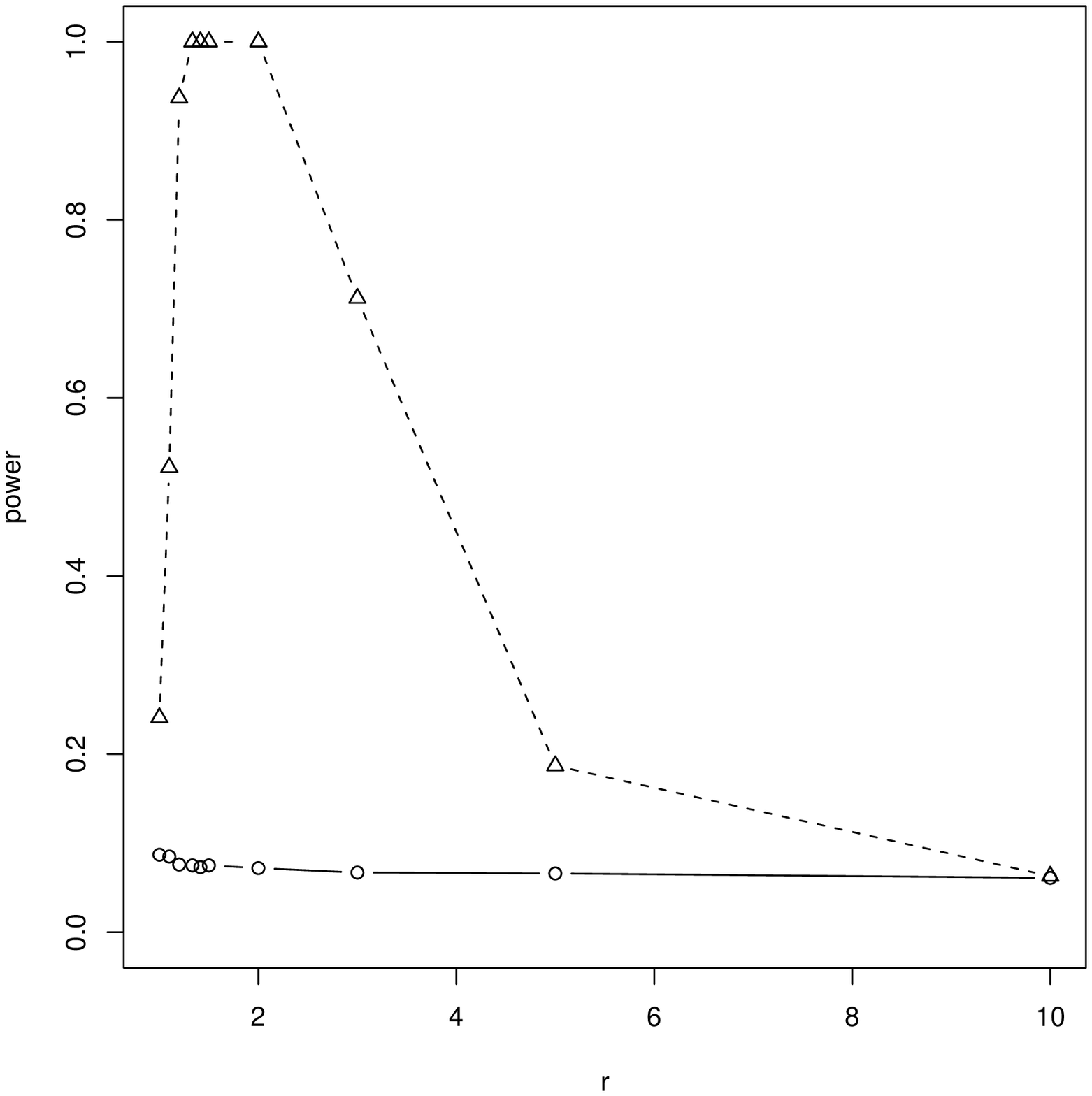}}}
\caption{
 \label{fig:MTAggSim}
Monte Carlo power using the asymptotic critical value against $H^A_{\sqrt{3}/12}$
as a function of $r$, for $n=100$ (left), $n=200$ (middle), and $n=500$ (right) conditional on the realization of $\Y$ in Figure \ref{fig:deldata}.  The circles represent the empirical significance levels while triangles represent the empirical power values.}
\end{figure}

We implement the above described Monte Carlo experiment $1000$ times with $n=100$, $n=200$, and $n=500$ and find the empirical significance levels $\widehat{\alpha}_S(n,J)$ and $\widehat{\alpha}_A(n,J)$ and the empirical powers $\widehat{\beta}^S_{n}(r,\sqrt{3}/8,J)$ and $\widehat{\beta}^A_{n}(r,\sqrt{3}/12,J)$.  These empirical estimates are presented in Table \ref{tab:MT-asy-emp-val} and plotted in Figures \ref{fig:MTSegSim} and \ref{fig:MTAggSim}.  Notice that the empirical significance levels are all larger than .05 for both alternatives, so this test is liberal in rejecting $H_0$ against both alternatives for the given realization of $\Y$ and $n$ values.  The smallest empirical significance levels and highest empirical power estimates occur at moderate $r$ values ($r=3/2,\,2,\,3$) against segregation and at smaller $r$ values ($r=\sqrt{2},\,3/2$) against association. Based on this analysis, for the given realization of $\Y$, we suggest the use of moderate $r$ values for segregation and slightly smaller for association.  Notice also that as $n$ increases, the empirical power estimates gets larger for both alternatives.

\begin{table}[t]
\centering
\begin{tabular}{|c|c|c|c|c|c|c|c|c|c|c|}
\hline
$r$  & 1 & 11/10 &6/5 & 4/3 & $\sqrt{2}$ &3/2 & 2 & 3 & 5 & 10\\
\hline
\multicolumn{11}{|c|}{$n=100,\,N=1000$} \\
\hline
$\widehat{\alpha}_S(n,J)$ & .144 & .141 & .124 & .101 & .095 & .087 & .070 & .075 & .071 & .072 \\
\hline
$\widehat{\beta}^S_{n}(r,\sqrt{3}/8,J)$ & .191 & .383 & .543 & .668 & .714 & .742 & .742 & .625 & .271 & .124 \\
\hline
\hline
$\widehat{\alpha}_A(n,J)$ & .118 & .111 & .089 & .081 & .065 & .062 & .067 & .064 & .068 & .071  \\
\hline
$\widehat{\beta}^A_{n}(r,\sqrt{3}/12,J)$ & .231 & .295 & .356 & .338 & .269 & .209 & .148 & .095 & .113 & .167 \\
\hline
\multicolumn{11}{|c|}{$n=200,\,N=1000$} \\
\hline
$\widehat{\alpha}_S(n,J)$ & .095 & .092 & .087 & .077 & .073 & .076 & .072 & .071 & .074 & .073 \\
\hline
$\widehat{\beta}^S_{n}(r,\sqrt{3}/8,J)$ & .135 & .479 & .743 & .886 & .927 & .944 & .959 & .884 & .335 & .105 \\
\hline
\hline
$\widehat{\alpha}_A(n,J)$ & .071 & .071 & .062 & .057 & .055 & .047 & .038 & .035 & .036 & .040  \\
\hline
$\widehat{\beta}^A_{n}(r,\sqrt{3}/12,J)$ & .182 & .317 & .610 & .886 & .952 & .985 & .972 & .386 & .143 & .068 \\
\hline
\multicolumn{11}{|c|}{$n=500,\,N=1000$} \\
\hline
$\widehat{\alpha}_S(n,J)$ & .089 & .092 & .087 & .086 & .080 & .078 & .079 & .079 & .076 & .081 \\
\hline
$\widehat{\beta}^S_{n}(r,\sqrt{3}/8,J)$ & .145 & .810 & .981 & .997 & .999 & 1.000 & 1.000 & 1.000 & .604 & .130 \\
\hline
\hline
$\widehat{\alpha}_A(n,J)$ & .087 & .085 & .076 & .075 & .073 & .075 & .072 & .067 & .066 & .061  \\
\hline
$\widehat{\beta}^A_{n}(r,\sqrt{3}/12,J)$ & .241 & .522 & .937 & 1.000 & 1.000 & 1.000 & 1.000 & .712 & .187 & .063 \\
\hline
\end{tabular}
\caption{
\label{tab:MT-asy-emp-val}
The empirical significance level and empirical power values under $H^S_{\sqrt{3}/8}$ and $H^A_{\sqrt{3}/12}$, $N=1000$, $n=100$, and $J=13$, at $\alpha=.05$ for the realization of $\Y$ in Figure \ref{fig:deldata}.}
\end{table}

{\bf Remark}
\label{rem:conditional-test}
The conditional test presented here is appropriate
when $w_j \in \mathcal W$ are fixed, not random.
An unconditional version requires the joint distribution
of the number and relative size of Delaunay triangles
when $\Y$ is, for instance, a Poisson point pattern.
Alas, this joint distribution is not available (see \cite{okabe:2000}). $\square$

\subsubsection{Related Test Statistics in Multiple Triangle Case}
\label{sec:rel-test-stat-MT}
For $J>1$, we have derived the asymptotic distribution of $\rho_n(r,J)=\frac{|\A|}{(n\,(n-1))}$. Let $\A_j$ be the number of arcs and $\rho_{n_j}(r)$ be the relative density for triangle $T_j$ and $n_j:=|\X_n \cap T_j|$, for $j=1,\ldots,J$.  So
$$\sum_{j=1}^{J}\frac{n_j\,(n_j-1)}{n\,(n-1)} \rho_{n_j}(r)= \rho_n(r,J),$$
since
$$\sum_{j=1}^{J}\frac{n_j\,(n_j-1)}{n\,(n-1)} \rho_{n_j}(r)=\frac{\sum_{j=1}^{J} |\A_j|}{n\,(n-1)}=\frac{|\A|}{n\,(n-1)}=\rho_n(r,J).$$

Let $\widehat{U}_n:=\sum_{j=1}^{J}w_j^2\cdot\rho_{n_j}(r)$ where $w_j=A(T_j)/A(C_H(\Y))$.  Since $\rho_{n_j}(r)$ are asymptotically independent, $\sqrt{n}\left( \widehat{U}_n-\mu(r,J) \right)$ and $\sqrt{n}\left( \rho_n(r,J)-\mu(r,J) \right)$ both converge in distribution to  $\N(0,\nu(r,J))$.

In the denominator of $\rho_n(r,J)$, we use $n(n-1)$ as the maximum number of arcs possible.  However, by definition, we can at most have a digraph with $J$ complete symmetric components of order $n_j$, for $j=1,\ldots,J$.  Then the maximum number possible is $n_t:=\sum_{j=1}^{J}n_j\,(n_j-1)$.  So the (adjusted) relative density is $\rho^{adj}_{n,J}(r):=\frac{|\A|}{n_t}$ and $\rho^{adj}_n(r)=\frac{\sum_{j=1}^{J} |\A_j|}{n_t}=\sum_{j=1}^{J}\frac{n_j\,(n_j-1)}{n_t}\,\rho_{n_j}(r)$. Since $\frac{n_j\,(n_j-1)}{n_t} \ge 0$ for each $j$, and $\sum_{j=1}^{J}\frac{n_j\,(n_j-1)}{n_t}=1$, $\rho^{adj}_{n,J}(r)$ is a mixture of $\rho_{n_j}(r)$. Then $\E\left[ \rho^{adj}_{n,J}(r) \right]=\mu(r,J)$ and the asymptotic variance of $\rho^{adj}_{n,J}(r)$ is
$$\frac{1}{n} \left[\nu(r) \left(\sum_{j=1}^{J}w_j^3/\left(\sum_{j=1}^{J}w_j^2 \right)^2 \right)+4\,\mu(r)^2\left(\sum_{j=1}^{J}w_j^3/\left( \sum_{j=1}^{J}w_j^2 \right)^2-1\right) \right]. \square $$

\subsubsection{Asymptotic Efficacy Analysis for $J>1$}
The PAE, HLAE, and asymptotic power function analysis are given for $J=1$. For $J>1$, the analysis will depend on both the number of triangles as well as the relative sizes of the triangles.  So the optimal $r$ values with respect to these efficacy criteria for $J=1$ do not necessarily hold for $J>1$, so the analysis need to be updated, given the values of $J$ and $\mathcal W$.

Under segregation alternative $H^S_{\epsilon}$, the PAE is given by
\begin{equation}
\label{eqn:PAE-NYR-J>1}
\PAE_J^S(r) = \frac{\bigl( \mu_S^{\prime\prime}(r,J,\epsilon=0)\bigr)^2}{\nu(r,J)}=
   \frac{\left( \mu_S^{\prime\prime}(r,\epsilon=0)\,\sum_{j=1}^{J}w_j^2 \right)^2}{\nu(r) \,\sum_{j=1}^{J}w_j^3 +4\,\mu(r)^2\left(\sum_{j=1}^{J}w_j^3-\left(\sum_{j=1}^{J}w_j^2 \right)^2\right)}.
\end{equation}
Under association alternative $H^A_{\epsilon}$ the PAE is similar.
\begin{figure}[]
\centering
\psfrag{r}{\scriptsize{$r$}}
\epsfig{figure=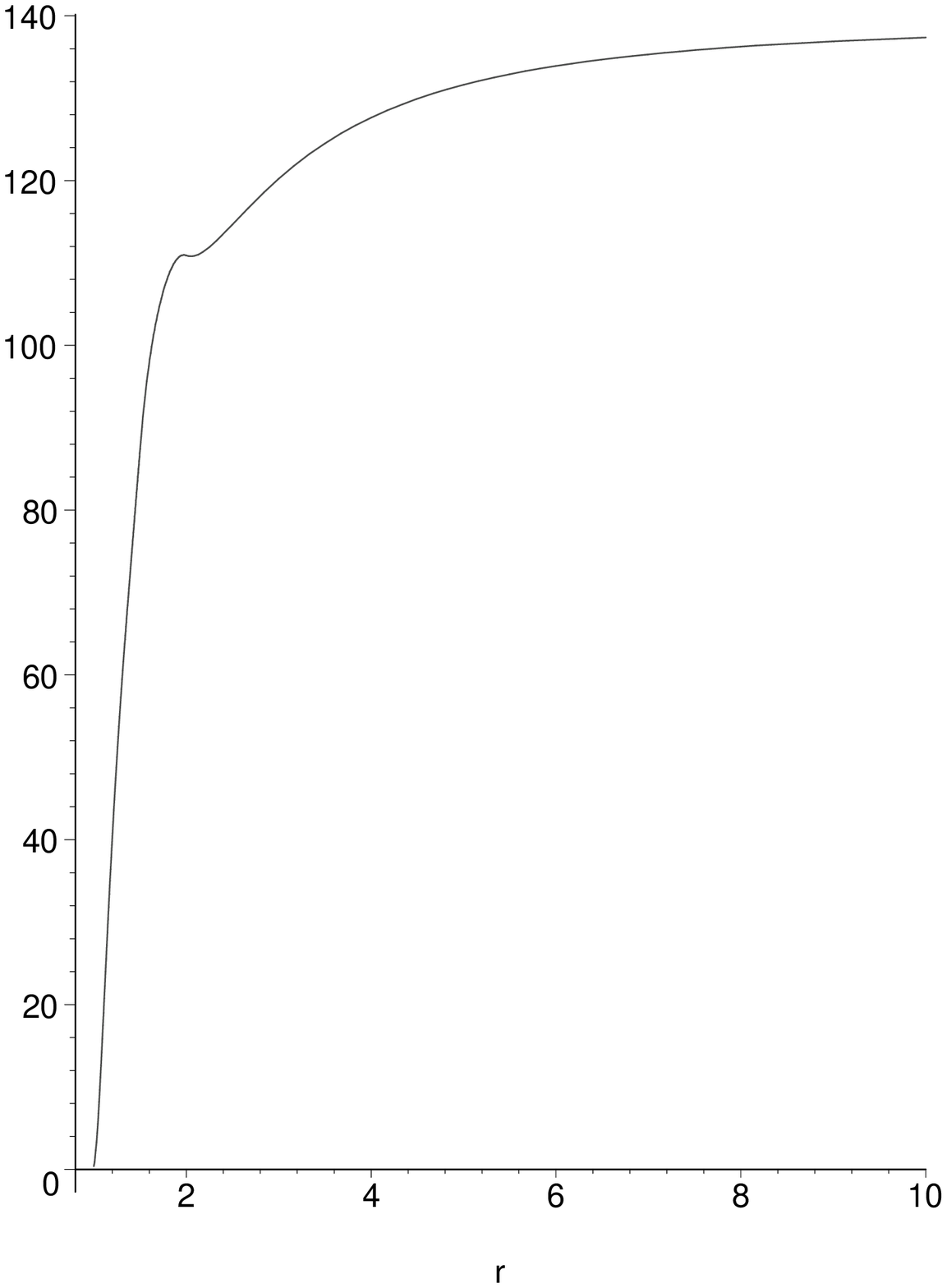, height=180pt, width=180pt}
\epsfig{figure=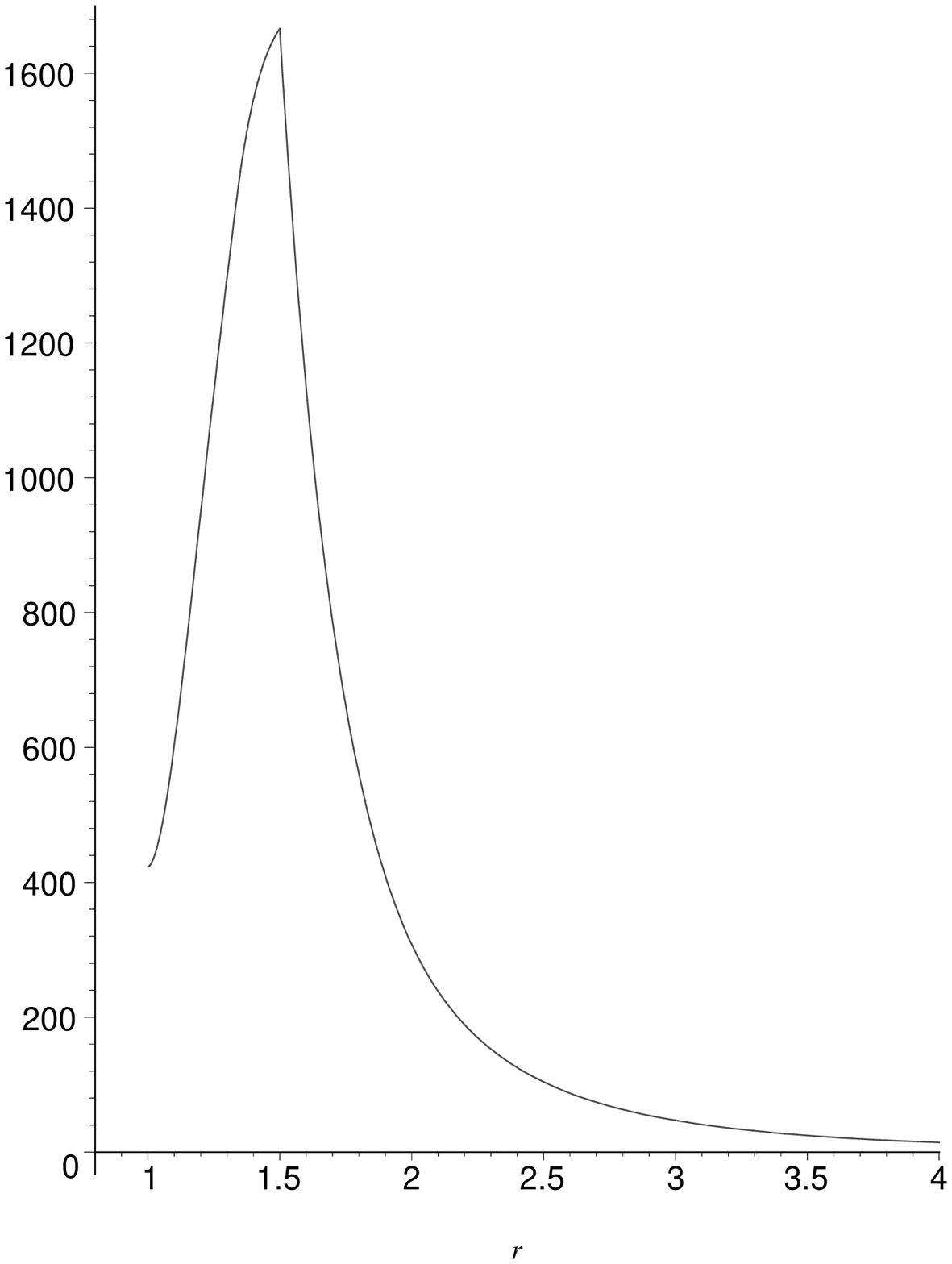, height=180pt, width=180pt}
\caption{\label{fig:MT-PAECurves}
Pitman asymptotic efficacy against segregation (left) and association (right)
as a function of $r$ with the realization of $\Y$ in Figure \ref{fig:deldata1}.
Notice that vertical axes are differently scaled.}
\end{figure}

In Figure \ref{fig:MT-PAECurves}, we present the PAE as a function of $r$ for both segregation and association conditional on the realization of $\Y$ in Figure \ref{fig:deldata}. Notice that, unlike $J=1$ case, $\PAE_J^S(r)$ is bounded. Some values of interest are
$\PAE_J^S(\rho_n(1)) = .3884$, $\lim_{r \rightarrow \infty} \PAE_J^S(r) = \frac{8\,\sum_{j=1}^{J}w_j^2}{256\,\left(\sum_{j=1}^{J}w_j^3-\left(\sum_{j=1}^{J}w_j^2 \right)^2\right)} \approx 139.34$, $\argsup_{r \in [1,2]} \PAE_J^S(r) \approx 1.974$. As for association, $\PAE_J^A(r=1) = 422.9551$, $ \lim_{r \rightarrow \infty} \PAE_J^A(r) = 0$,
$\argsup_{r \ge 1} \PAE_J^A(r) \approx 1.5$ with $\PAE_J^A(r=1.5) \approx 1855.9672$.
Based on the asymptotic efficacy analysis, we suggest,
for large $n$ and small $\epsilon$,
choosing moderate $r$ for testing against segregation and association.

Under segregation, the HLAE is given by
\begin{equation}
\label{eqn:HLAE-NYr-J>1}
\HLAE_J^S(r,\epsilon):=\frac{\left( \mu_S(r,J,\epsilon)-\mu(r,J)\right)^2}{\nu_S(r,J,\epsilon)}=\frac{\left(\mu_S(r,\epsilon)\,\left(\sum_{j=1}^{J}w_j^2 \right)-\mu(r)\,\left(\sum_{j=1}^{J}w_j^2 \right)\right)^2}{\nu_S(r,\epsilon) \,\sum_{j=1}^{J}w_j^3 +4\,\mu_S(r,\epsilon)^2\left(\sum_{j=1}^{J}w_j^3-\left(\sum_{j=1}^{J}w_j^2 \right)^2\right)}.
\end{equation}
 Notice that $\HLAE_J^S(r,\epsilon=0)=0$ and $\lim_{\rightarrow \infty}\HLAE_J^S(r,\epsilon)=0$.

\begin{figure}[]
\centering
\psfrag{r}{\scriptsize{$r$}}
\epsfig{figure=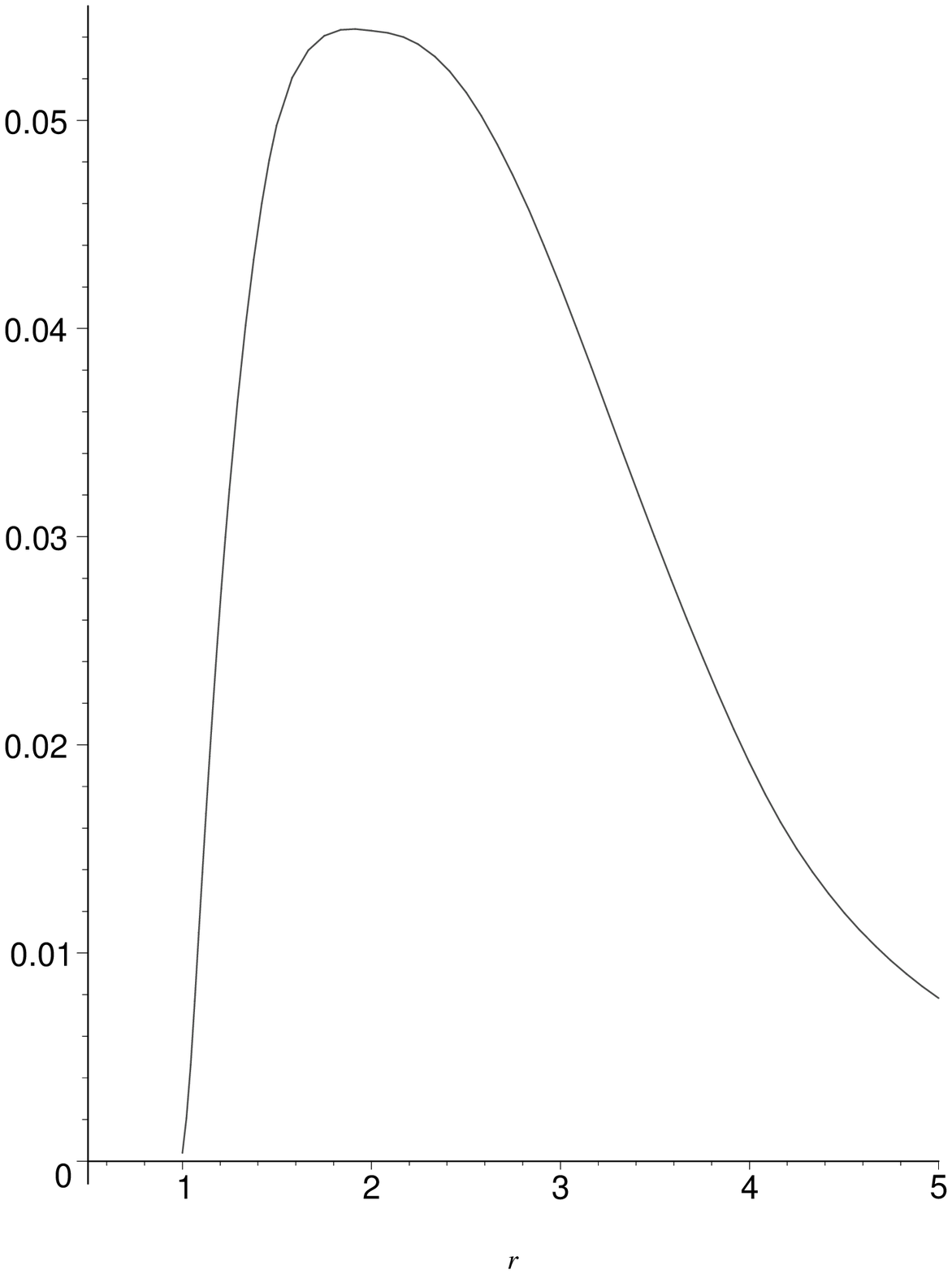, height=125pt, width=125pt}
\epsfig{figure=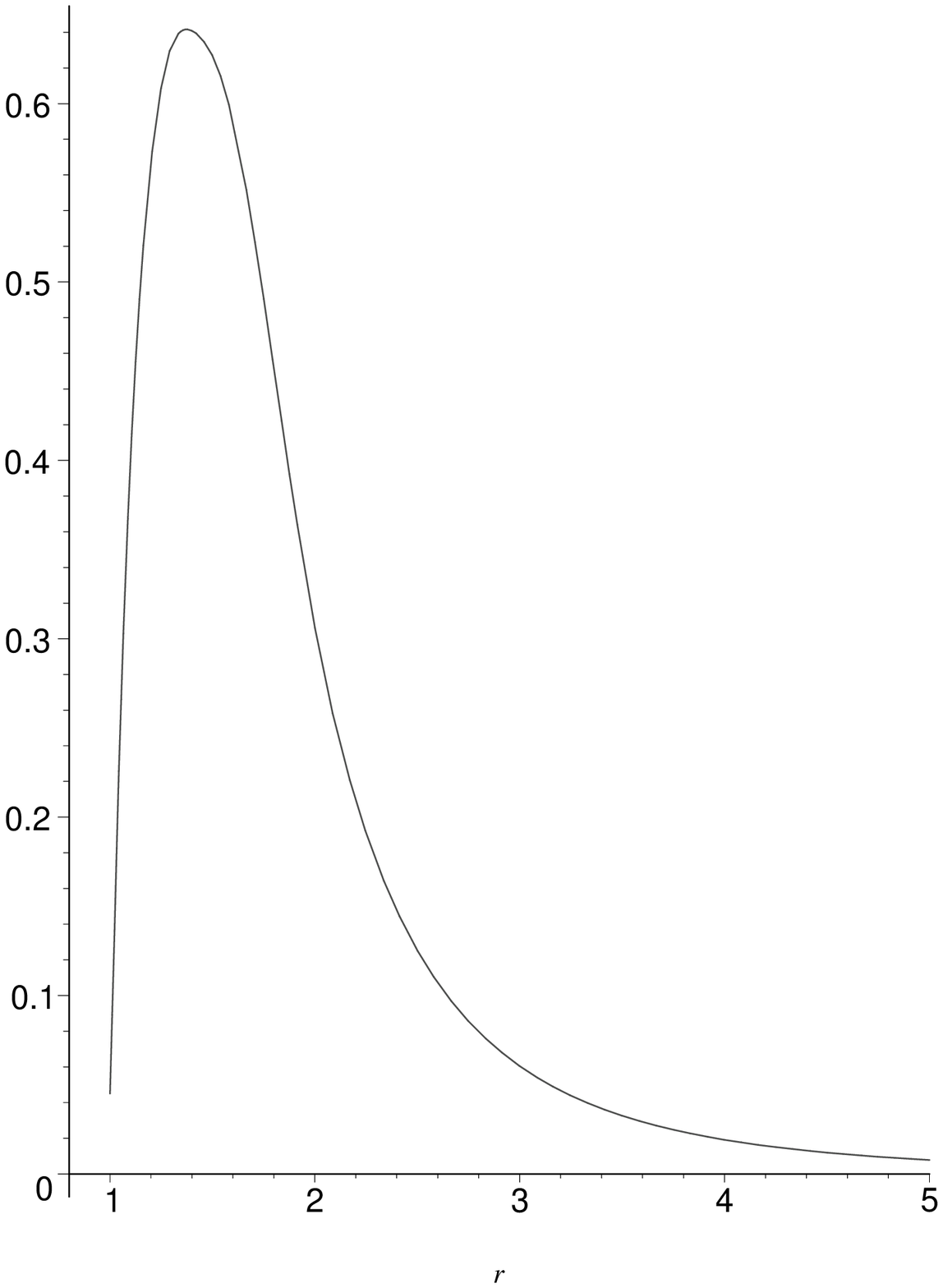, height=125pt, width=125pt}
\epsfig{figure=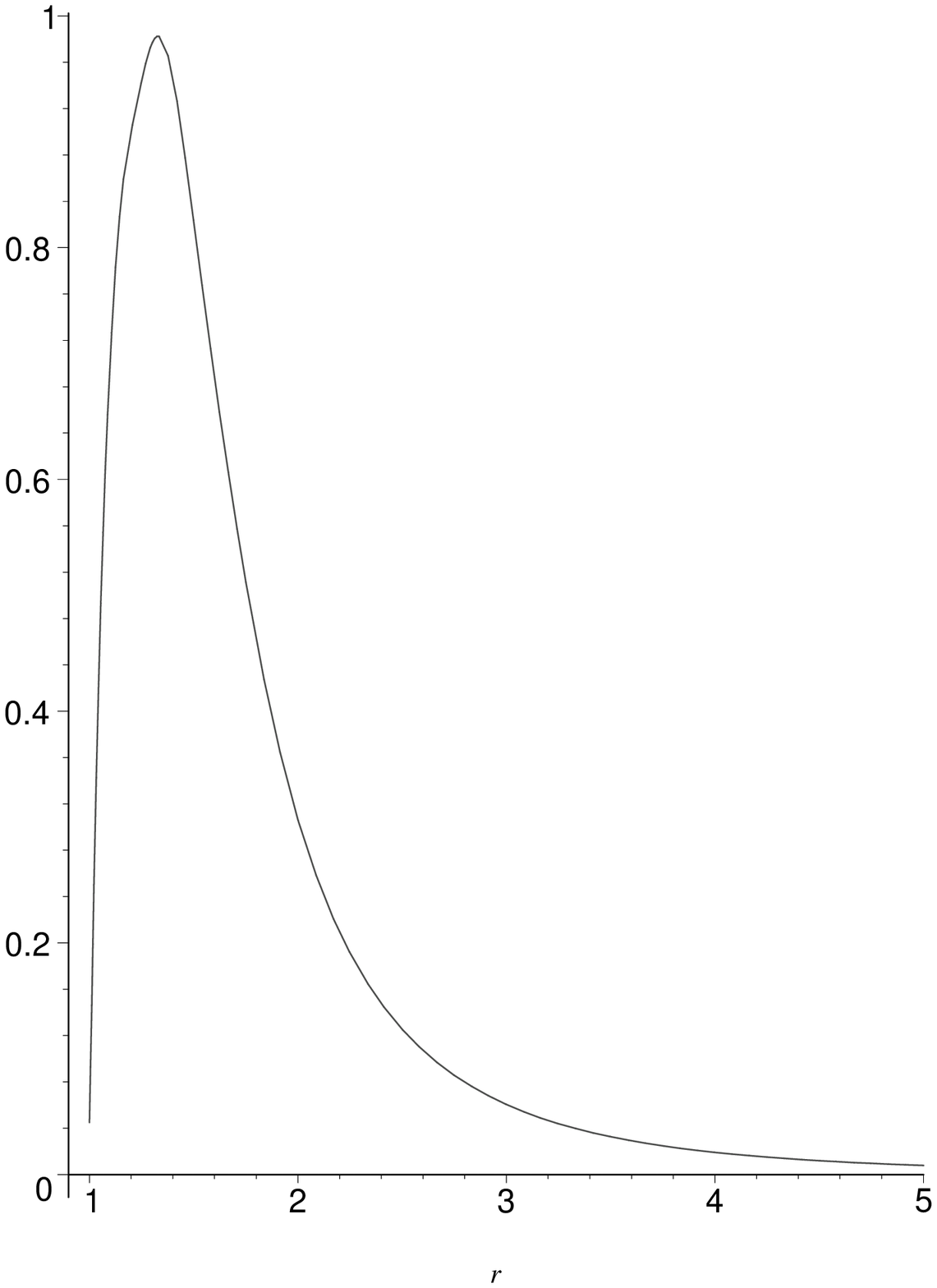, height=125pt, width=125pt}
\caption{
\label{fig:MT-HLAE-Seg}
Hodges-Lehmann asymptotic efficacy
against segregation alternative $H^S_{\epsilon}$
as a function of $r$
for $\epsilon= \sqrt{3}/8, \sqrt{3}/4, 2\,\sqrt{3}/7$ (left to right) conditional on the realization of $\Y$ in Figure \ref{fig:deldata1}.}
\end{figure}

We calculate HLAE of $\rho_n(r,J)$ under $H^S_{\epsilon}$ for $\epsilon=\sqrt{3}/8$, $\epsilon=\sqrt{3}/4$, and $\epsilon=2\,\sqrt{3}/7$. In Figure \ref{fig:MT-HLAE-Seg} we present $\HLAE_J^S(r,\epsilon)$ for these $\epsilon$ values conditional on the realization of $\Y$ in Figure \ref{fig:deldata1}.

Note that with $\epsilon=\sqrt{3}/8$, $\HLAE_J^S \left( r=1,\sqrt{3}/8 \right) \approx.0004$ and $\argsup_{r \in [1,\infty]} \HLAE_J^S\left( r,\sqrt{3}/8 \right) \approx 1.8928$ with the supremum $\approx .0544$.  With $\epsilon=\sqrt{3}/4$, $\HLAE_J^S\left( r=1,\sqrt{3}/4 \right) \approx .0450$ and $\argsup_{r \in [1,\infty]} \HLAE_J^S\left( r,\sqrt{3}/4 \right) \approx 1.3746$ with the supremum $\approx .6416$. With $\epsilon=2\,\sqrt{3}/7$, $\HLAE_J^S\left( r=1,2\,\sqrt{3}/7 \right) \approx .045$ and $\argsup_{r \in [1,\infty]} \\
\HLAE_J^S\left( r,2\,\sqrt{3}/7 \right) \approx 1.3288$ with the supremum $ \approx .9844$.  Furthermore, we observe that $\HLAE_J^S\left( r,2\,\sqrt{3}/7 \right)>\HLAE_J^S \left( r,\sqrt{3}/4 \right)>\HLAE_J^S\left( r,\sqrt{3}/8 \right)$ at each $r$.  Based on the HLAE analysis for the given $\Y$ we suggest moderate $r$ values for moderate segregation, and small $r$ values for severe segregation.

The explicit form of $\HLAE_J^A(r,\epsilon)$ is similar which implies $\HLAE_J^A(r,\epsilon=0)=0$ and $\lim_{r\rightarrow \infty}\HLAE_J^A(r,\epsilon)=0$.

\begin{figure}[]
\centering
\psfrag{r}{\scriptsize{$r$}}
\epsfig{figure=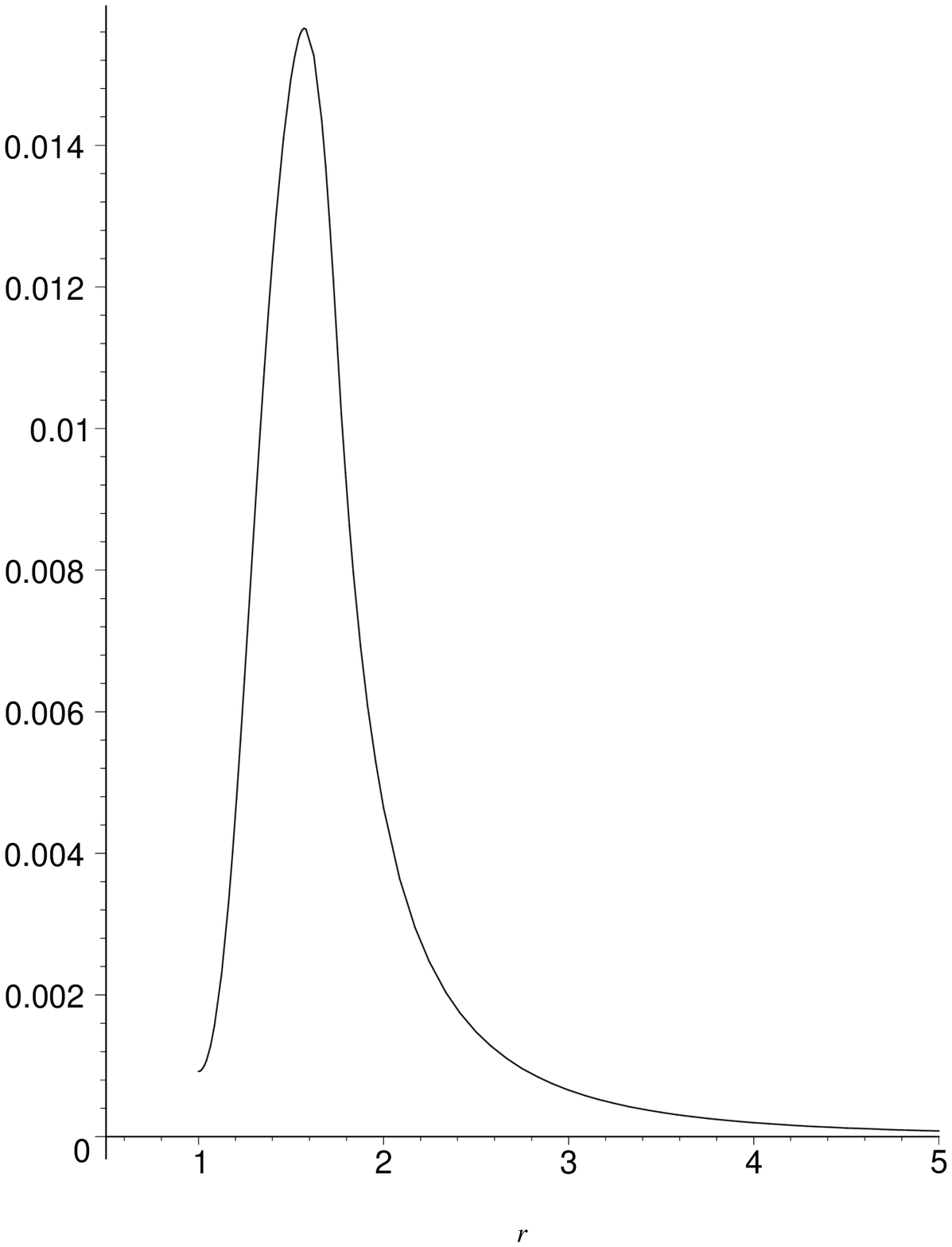, height=125pt, width=125pt}
\epsfig{figure=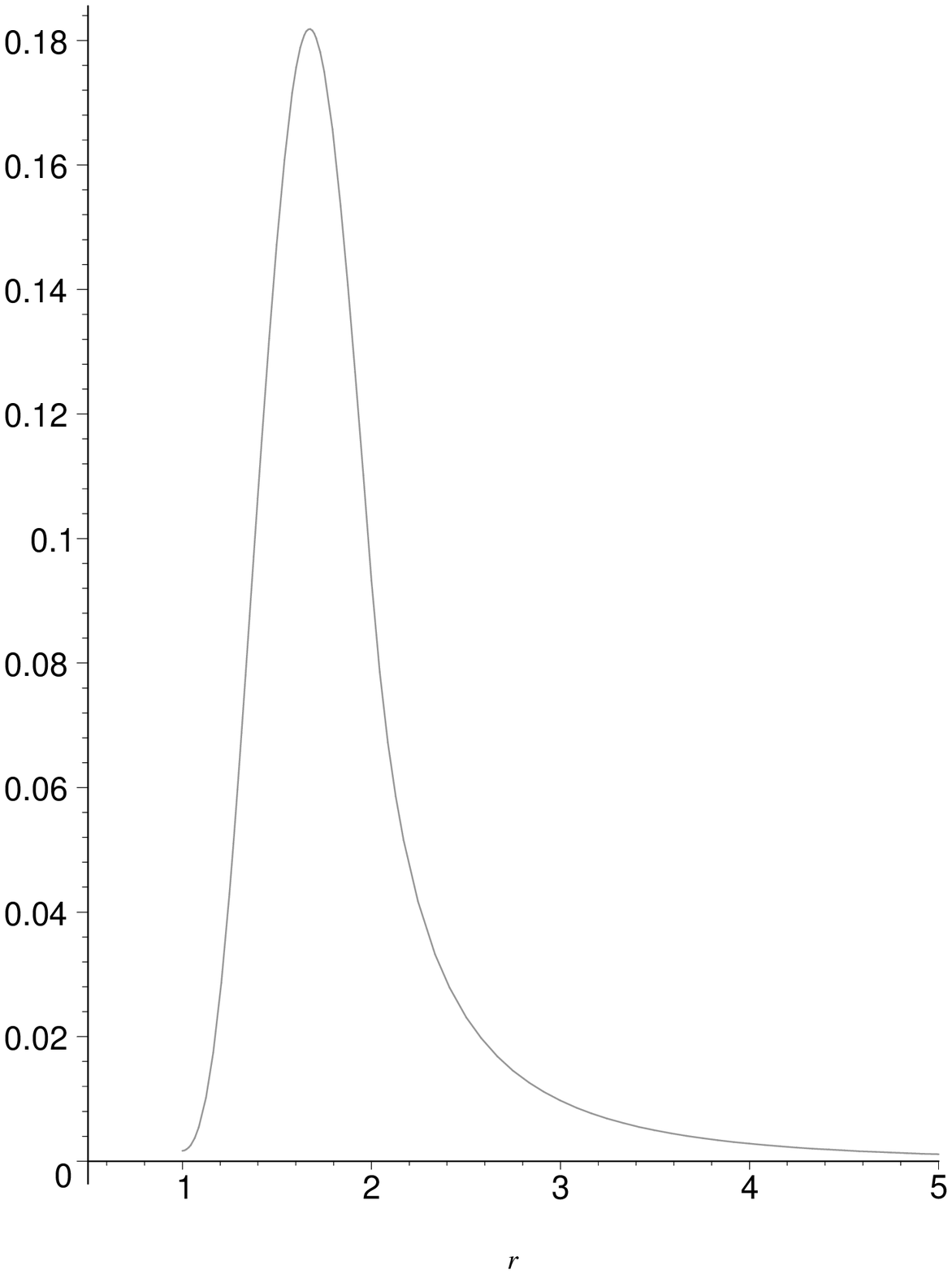, height=125pt, width=125pt}
\epsfig{figure=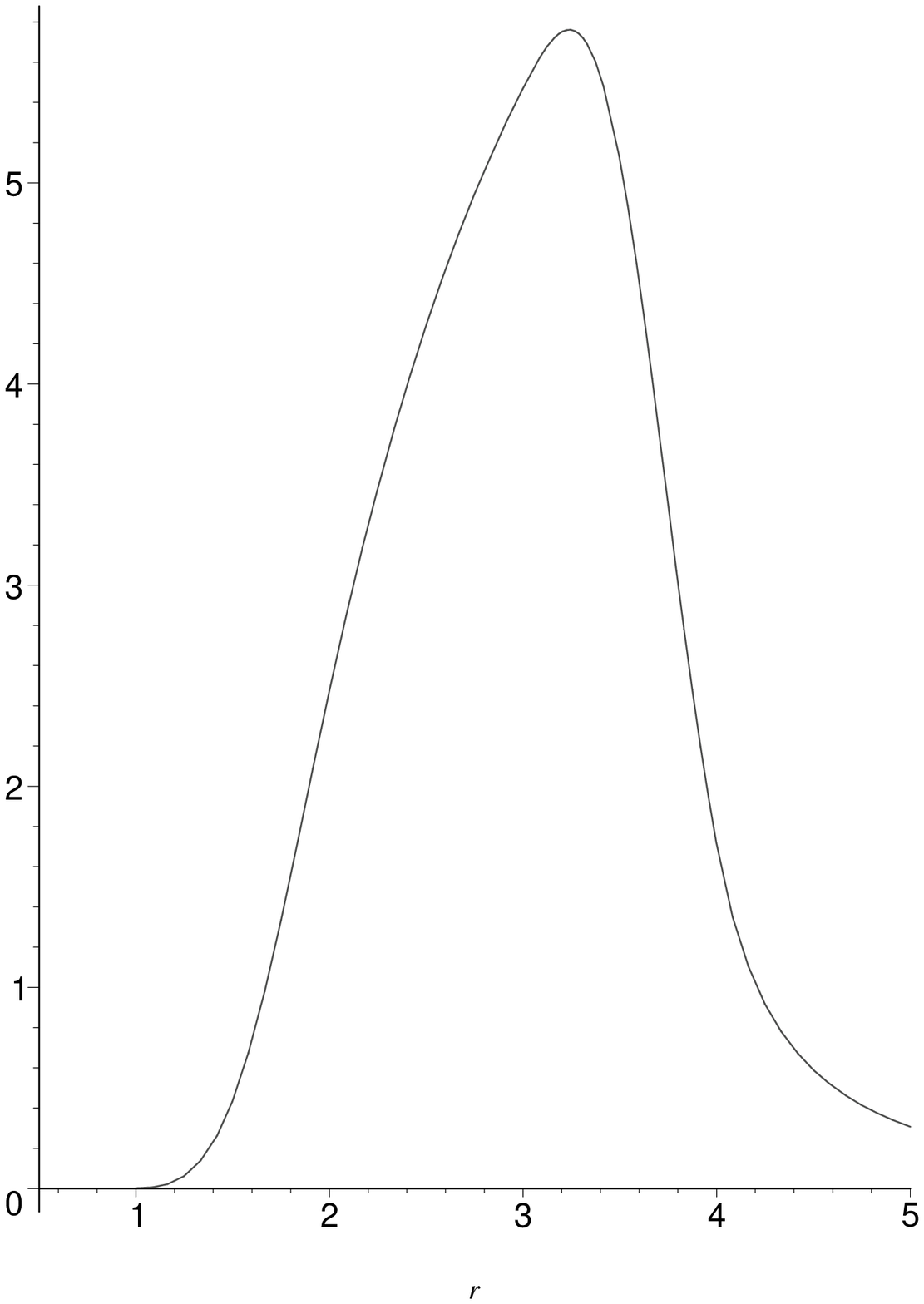, height=125pt, width=125pt}
\caption{
\label{fig:MT-HLAE-Agg}
Hodges-Lehmann asymptotic efficacy
against association alternative $H^A_{\epsilon}$
as a function of $r$
for $\epsilon= \sqrt{3}/21, \sqrt{3}/12, 5\,\sqrt{3}/24$ (left to right) conditional on the realization of $\Y$ in Figure \ref{fig:deldata1}.}
\end{figure}

We calculate HLAE of $\rho_n(r,J)$ under $H^A_{\epsilon}$ for $\epsilon=\sqrt{3}/21$, $\epsilon=\sqrt{3}/12$, and $\epsilon=5\,\sqrt{3}/24$. In Figure \ref{fig:MT-HLAE-Agg} we present $\HLAE_J^S(r,\epsilon)$ for these $\epsilon$ values conditional on the realization of $\Y$ in Figure \ref{fig:deldata}

Note that with $\epsilon=\sqrt{3}/21$, $\HLAE_J^A\left( r=1,\sqrt{3}/21 \right) \approx.0009$ and $\argsup_{r \in [1,\infty]} \HLAE_J^A\left( r,\sqrt{3}/21 \right) \approx 1.5734$ with the supremum $\approx .0157$.  With $\epsilon=\sqrt{3}/12$, $\HLAE_J^A\left( r=1,\sqrt{3}/12 \right) \approx .0168$ and $\argsup_{r \in [1,\infty]}\\
 \HLAE_J^A\left( r,\sqrt{3}/12 \right) \approx 1.6732$ with the supremum $\approx .1818$. With $\epsilon=5\,\sqrt{3}/24$, $\HLAE_J^A\left( r=1,5\,\sqrt{3}/24 \right) \approx .0017$ and $\argsup_{r \in [1,\infty]} \HLAE_J^A\left( r,5\,\sqrt{3}/24 \right) \approx 3.2396$ with the supremum $\approx 5.7616$.  Furthermore, we observe that $\HLAE_J^A\left( r,5\,\sqrt{3}/24 \right)>\HLAE_J^A\left( r,\sqrt{3}/12 \right)>\HLAE_J^A\left( r,\sqrt{3}/21 \right)$ at each $r$.  Based on the HLAE analysis for the given $\Y$ we suggest moderate $r$ values for moderate association and larger $r$ values for severe association.

\section{Discussion}
\label{sec:discussion}
The extension to $\R^d$ for $d > 2$ is straightforward.
See \cite{ceyhan:TR-dom-num-NPE-spatial} for more detail.
Moreover, the geometry invariance,
asymptotic normality of the $U$-statistic and consistency of the tests hold for $d>2$.

The first proximity map similar to the $r$-factor proximity map $\NY^r$ in literature
is the spherical proximity map $N_S(x):=B(x,r(x))$, (see the references for CCCD in the Introduction).
A slight variation of $N_S$ is the arc-slice proximity map $N_{AS}(x):=B(x,r(x)) \cap T(x)$ where
$T(x)$ is the Delaunay cell that contains $x$ (see \cite{ceyhan:CS-JSM-2003}).
Furthermore, Ceyhan and Priebe introduced the central similarity proximity map
$N_{CS}$ in \cite{ceyhan:CS-JSM-2003} and $\NY^r$ in \cite{ceyhan:TR-dom-num-NPE-spatial}.
The $r$-factor proximity map, when compared to the others,
has the advantages that the asymptotic distribution of the
domination number $\gamma_n(\NY^r)$ is tractable (see \cite{ceyhan:TR-dom-num-NPE-spatial}),
the exact minimum dominating sets can be found in polynomial time.
Moreover $\NY^r$ and $N_{CS}$ are geometry invariant for uniform data over triangles.
Additionally, the mean and variance of $\rho_n$ is not analytically tractable
for $N_S$ and $N_{AS}$. While $\NY^r(x)$, $N_{CS}(x)$, and $N_{AS}(x)$ are well
defined only for $x \in C_H(\Y)$, the convex hull of $\Y$, $N_S(x)$ is well defined for all $x \in \R^d$.
$N_S$ and $N_{AS}$ require no effort to extend to higher dimensions.

There are many tests available for segregation and association in literature.
See \cite{dixon:1994} for a survey on these tests and relevant references.
The most prevalent of these tests are Pielou's $\chi^2$ test of independence
and Ripley's test based on $K(t)$ and $L(t)$ functions.
However, the test we introduce here is not comparable to either of them,
since it is a conditional test --- conditional on a realization of $J=|\Y|$ and $\mathcal W$ and
we require the number of triangles $J$ is fixed and relatively small compared to $n=|\X_n|$.
The null hypothesis for testing spatial patterns has two major forms:
\begin{itemize}
\item[(i)] assuming random labeling of locations, i.e. spatial randomness does not necessarily hold, as in Pielou's test which only tests for the association between classes,
\item[(ii)] assuming not only random labeling but also complete spatial randomness, that is, each class is distributed randomly throughout the area of interest, as in Ripley's test.
\end{itemize}
Our conditional test is closer to the latter in this regard.

The test based on the mean domination number in \cite{ceyhan:TR-dom-num-NPE-spatial} is not a conditional test,
but requires both $n$ and number of Delaunay triangles $J$ to be large.  The comparison for a large but fixed $J$ is possible.
Furthermore, under segregation alternatives, the Pitman asymptotic efficacy is not applicable to the mean domination number case, however, for large $n$ and $J$ we suggest the use of it over arc density since for each $\epsilon>0$, Hodges-Lehmann asymptotic efficacy is unbounded for the mean domination number case, while it is bounded for arc density case with $J>1$. As for the association alternative, HLAE suggests moderate $r$ values which has finite Hodges-Lehmann asymptotic efficacy. So again, for large $J$ and $n$ mean domination number is preferable.  The basic advantage of $\rho_n(r)$ is that, it does not require $J$ to be large, so for small $J$ it is preferable.

\section*{Acknowledgments}
This work was partially sponsored by the Defense Advanced Research Projects Agency as administered
by the Air Force Office of Scientific Research under contract DOD F49620-99-1-0213.


\section*{Appendix 1: Derivation of $\mu(r)$ and $\nu(r)$}
In the standard equilateral triangle, let $\y_1=(0,0)$, $\y_2=(1,0)$, $\y_3=\bigl( 1/2,\sqrt{3}/2 \bigr)$, $M_C$ be the center of mass, $M_j$ be the midpoints of the edges $e_j$ for $j=1,2,3$. Then $M_C=\bigl(1/2,\sqrt{3}/6\bigr)$, $M_1=\bigl(3/4,\sqrt{3}/4 \bigr)$, $M_2=\bigl(1/4,\sqrt{3}/4\bigr)$, $M_3=(1/2,0)$.

Recall that $\E[\rho_n(r)]=\frac{1}{n\,(n-1)}\sum \sum_{i < j } \,\E[h_{ij}]=\frac{1}{2}\E[h_{12}]=\mu(r)=P\bigl(X_j \in \NY^r(X_i)\bigr)$.

Let $\X_n$ be a random sample of size $n$ from $\U(T(\Y))$. For $x_1=(u,v)$, $\ell_r(x_1)=r\,v+r\,\sqrt{3}\,u-\sqrt{3}\,x.$  Next, let $N_1:=\ell_r(x_1)\cap e_3$ and $N_2:=\ell_r(x_1)\cap e_2$. Then for $z_1 \in T_s:=T(\y_1,M_3,M_C)$, $N_{\Y}^r(z_1)=T(\y_1,N_1,N_2)$ provided that $\ell_r(x_1)$ is not outside of $T(\Y)$, where
 $$N_1=\bigl(r\,\bigl(y_1+\sqrt{3}\,x_1\bigr)\sqrt{3}/3,0\bigr)\text{ and } N_2=\bigl(r\,\bigl(y_1+\sqrt{3}\,x_1\bigr)\sqrt{3}/6,\bigl(y_1+\sqrt{3}\,x_1\bigr) r/2\bigr).$$

\subsection*{Derivation of $\mu(r)$ in Theorem 2}
Now we find $\mu(r)$ for $r \in [1,\infty)$.
Observe that, by symmetry,
 $$\mu(r)=P\bigl( X_2 \in \NY^r(X_1) \bigr)=6\,P\bigl( X_2 \in \NY^r(X_1), X_1 \in T_s \bigr).$$
 Let $\ell_s(r,x)$ be the line such that $r\,d(\y_1,\ell_s(r,x))=d(\y_1,e_1)$, so $\ell_s(r,x)=\sqrt{3}\,(1/r-x)$.  Then if $x_1 \in T_s$ is above $\ell_s(r,x)$ then $\NY^r(x_1)=T(\Y)$, otherwise, $\NY^r(x_1)\subsetneq T(\Y)$.

For $r \in [1,3/2)$, $\ell_s(r,x)\cap T_s=\emptyset$, so $\NY^r(x)\subsetneq T(\Y)$ for all $x \in T_s$.  Then
$$
\mu(r)=6\,P\bigl( X_2 \in \NY^r(X_1), X_1 \in T_s \bigr)=6\,\int_0^{1/2}\int_0^{x/\sqrt{3}} \frac{A(\NY^r(x_1))}{A(T(\Y))^2}dydx = 6\,\left(\frac{37}{1296}\,r^2 \right)=\frac{37}{216}\,r^2.
$$
where $A(\NY^r(x_1))=\frac{\sqrt{3}}{12}\,r^2\left( y+\sqrt{3}\,x \right)^2 $ and $A(T(\Y))=\sqrt{3}/4$.

\begin{figure} [ht]
    \centering
   \scalebox{.4}{\input{ls_lam_cases.pstex_t}}
    \caption{The cases for relative position of $\ell_s(r,x)$ with various $r$ values.}
    \label{fig:ls-lam-cases}
\end{figure}
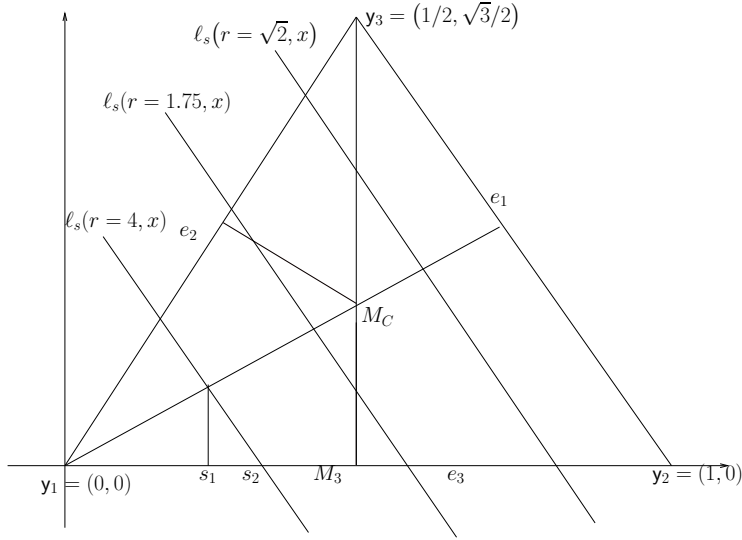

For $r \in [3/2,2)$, $\ell_s(r,x)$ crosses through $\overline{M_3M}_C$.  Let the $x$ coordinate of $\ell_s(r,x)\cap \overline{\y_1M}_C$ be $s_1$, then $s_1=\frac{3}{4\,r}$.  See Figure \ref{fig:ls-lam-cases}.

 Then
\begin{align*}
&P\bigl( X_2 \in \NY^r(X_1), X_1 \in T_s \bigr)= \int_0^{s_1}\int_0^{x/\sqrt{3}} \frac{A(\NY^r(x_1))}{A(T(\Y))^2}dydx+\int_{s_1}^{1/2}\int_0^{\ell_s(r,x)} \frac{A(\NY^r(x_1))}{A(T(\Y))^2}dydx+\\
 &\int_{s_1}^{1/2}\int_{\ell_s(r,x)}^{x/\sqrt{3}} \frac{1}{A(T(\Y))}dydx = -\frac{-36+r^4+64\,r-32\,r^2}{48\,r^2}.
\end{align*}
Hence for $r \in [3/2,2)$, $\mu(r)=-\frac{1}{8}\,r^2-8\,r^{-1}+\frac{9}{2}\,r^{-2}+4$.

For $r \in [2,\infty)$, $\ell_s(r,x)$ crosses through $\overline{\y_1M}_3$.  Let the $x$ coordinate of $\ell_s(r,x)\cap \overline{\y_1M}_3$ be $s_2$, then $s_2=1/r$. See Figure \ref{fig:ls-lam-cases}.

Then
\begin{eqnarray*}
\lefteqn{P\bigl( X_2 \in \NY^r(X_1), X_1 \in T_s \bigr)= \int_0^{s_1}\int_0^{x/\sqrt{3}} \frac{A(\NY^r(x_1))}{A(T(\Y))^2}dydx+\int_{s_1}^{s_2}\int_0^{\ell_s(r,x)} \frac{A(\NY^r(x_1))}{A(T(\Y))^2}dydx}\\
& &+\int_{s_1}^{s_2}\int_{\ell_s(r,x)}^{x/\sqrt{3}} \frac{1}{A(T(\Y))}dydx+\int_{s_2}^{1/2}\int_{0}^{x/\sqrt{3}} \frac{1}{A(T(\Y))}dydx=\frac{-3+2\,r^2}{12\,r^2}.
\end{eqnarray*}
Hence for $r \in [2,\infty)$, $\mu(r)=1-\frac{3}{2}\,r^{-2}$.

For $r=\infty $, $\mu(r)=1$ follows trivially.

\subsection*{Derivation of $\nu(r)$ in Theorem 2}
To find $\Cov[h_{12},h_{13}] $, we introduce a related concept.

{\bf Definition:}
Let $(\Omega,\mathcal{M})$ be a measurable space and
consider the proximity map $N:\Omega \times \wp(\Omega) \rightarrow \wp(\Omega)$,
where $\wp(\cdot)$ represents the power set functional.
For $B \subset \Omega$, the {\em $\G_1$-region}, $\G_1(\cdot)=\G_1(\cdot,N):\Omega \rightarrow \wp(\Omega)$ associates the region $\G_1(B):=\{z \in \Omega: B \subseteq  N(z)\}$ with each set $B \subset \Omega$. For $x \in \Omega$, we denote $\G_1(\{x\})$ as $\G_1(x)$.  Note that $\G_1$-region depends on proximity region $N(\cdot)$.

Furthermore, let $\G_1(\cdot,\NY^r)$ be the $\G_1$-region associated with $\NY^r(\cdot)$, let $A_{ij}$ be the event that $\{X_iX_j \in \A\} =\{X_i \in \NY^r(X_j)\}$, then $h_{ij}=I(A_{ij})+I(A_{ji})$. Let
$$
P^r_{2N}:=P(\{X_2,X_3\} \subset \NY^r(X_1)),\;\;
P^r_M:=P(X_2 \in \NY^r(X_1), X_3 \in \G_1(X_1,\NY^r),\;\;
P^r_{2G}:=P(\{X_2,X_3\} \subset \G_1(X_1,\NY^r)).
$$
 Then
$\Cov[h_{12},h_{13}]=\E[h_{12}\,h_{13}]-\E[h_{12}]\E[h_{13}]$ where
\begin{eqnarray*}
\E[h_{12}\,h_{13}] & = & \E[(\I(A_{12})+\I(A_{21}))\,(\I(A_{13})+\I(A_{31})] \\
                  & = & P(A_{12} \cap A_{13})+P(A_{12} \cap A_{31})+P(A_{21} \cap A_{13})+P(A_{21}\cap A_{31}). \\
                  &=&P(\{X_2,X_3\} \subset \NY^r(X_1))+2\,P(X_2 \in \NY^r(X_1), X_3 \in \G_1(X_1,\NY^r))+P(\{X_2,X_3\} \subset \G_1(X_1,\NY^r))\\
 &=&P^r_{2N}+2\,P^r_M+P^r_{2G}.
\end{eqnarray*}
So $\nu(r)=\Cov[h_{12},h_{13}]  = \left(P^r_{2N}+2\,P^r_M+P^r_{2G}\right)-[2\,\mu(r)]^2.$

Furthermore, for any $x_1=(u,v) \in T(\Y)$, $\G_1(x_1,\NY^r)$ is a convex or nonconvex polygon.  Let $\xi_j(r,x)$ be the line  between $x_1$ and the vertex $\y_j$ parallel to the edge $e_j$ such that $r\,d(\y_j,\xi_j(r,x))=d(\y_j,\ell_r(x_1)) \text{ for } j=1,2,3.$
Then   $\G_1(x_1,\NY^r)\cap R(\y_j)$ is bounded by $\xi_j(r,x)$ and the median lines.

For $x_1=(u,v)$, $\xi_1(r,x)=-\sqrt{3}\,x+(v+\sqrt{3}\,u)/r,\; \xi_2(r,x)=(v+\sqrt{3}r\,(x-1)+\sqrt{3}(1-u))/r \text{ and } \xi_3(r,x)=(\sqrt{3}(r-1)+2\,v)/(2\,r).$
 To find the covariance, we need to find the possible types of $\G_1\left( x_1,\NY^r \right)$ and $\NY^r(x_1)$ for $r \in [1,\infty)$.

We partition $[1,\infty)$ with respect to the types of $\NY^r(x_1)$ and $\G_1\left( x_1,\NY^r \right)$ and obtain $[1,4/3)$, $[4/3,3/2)$, $[3/2,2)$, $[2,\infty)$.

For $r\in [1,4/3)$, there are six cases regarding $\G_1\left( x_1,\NY^r \right)$ and one case for $\NY^r(x_1)$. See Figure \ref{fig:G1-NYr-Cases-1} for the prototypes of these six cases of $\G_1\left( x_1,\NY^r \right)$. Each case $j$, corresponds to the region $R_j$ in Figure \ref{regions-for-N_nu2}, where
$\ell_{am}\,(x) =x/\sqrt{3},\;\; q_1(x)=(2\,r+3\,x-3)/\sqrt{3},\;\; q_2(x)=\sqrt{3}\,(1/2-r/3),\\
q_3(x)=\sqrt{3}\,(x-1+r/2),\;\;  q_4(x)=\sqrt{3}\,(1/2-r/4),\;\;   q_{12}(x)=\sqrt{3}\,(r/2-x)$ and $s_1=1-2\,r/3, \; s_2=3/2-r,\; s_3=1-r/2, \; s_4=3/2-5\,r/6, \; s_5=3\,r/8.$
\begin{figure} [ht]
   \centering
   \scalebox{.27}{\input{G1ofxCase1.pstex_t}} 
   \scalebox{.27}{\input{G1ofxCase2.pstex_t}}
   \scalebox{.27}{\input{G1ofxCase3.pstex_t}}
   \scalebox{.27}{\input{G1ofxCase4.pstex_t}}
   \scalebox{.27}{\input{G1ofxCase5.pstex_t}}
   \scalebox{.27}{\input{G1ofxCase6.pstex_t}}
   \caption{The prototypes of the six cases of $\G_1\left( x_1,\NY^r \right)$ for $x_1 \in T_s$ for $r \in [1,4/3)$.}
\label{fig:G1-NYr-Cases-1}
\end{figure}
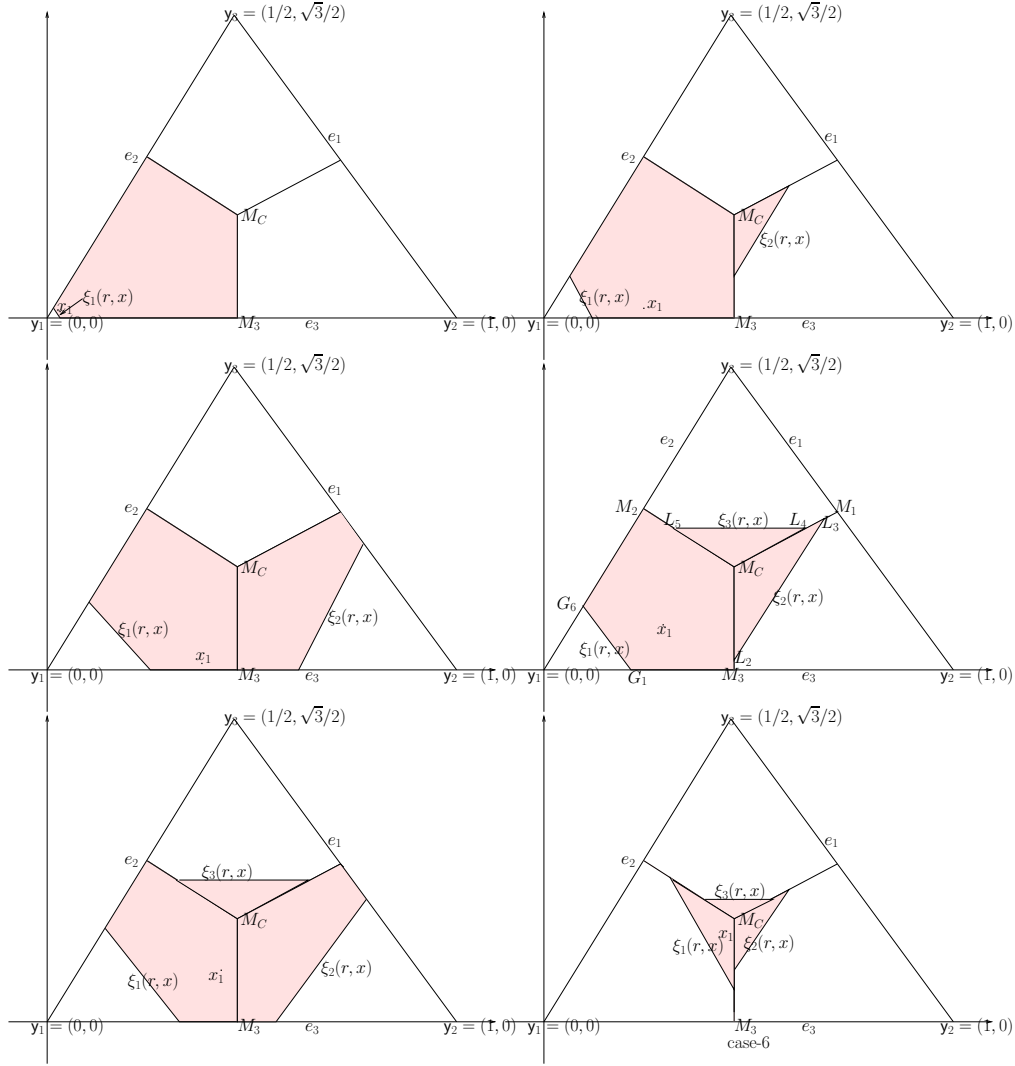
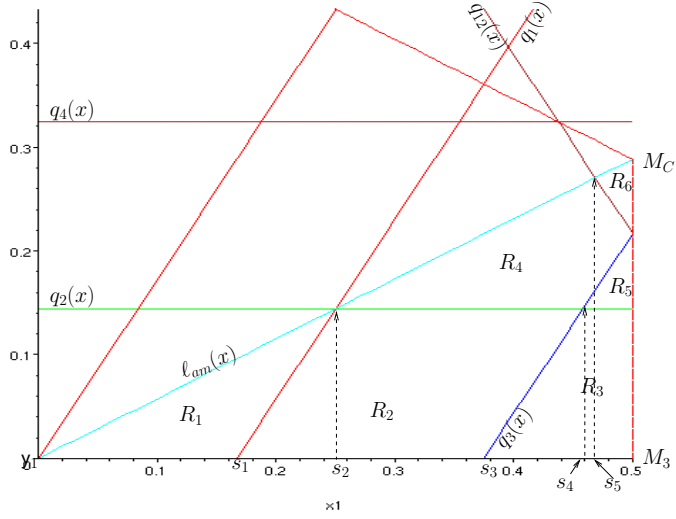
\begin{figure} []
    \centering
   \scalebox{.4}{\input{N_nu2GamRegions.pstex_t}}
    \caption{The regions corresponding to the prototypes of the six cases for $r \in [1,4/3)$ with $r=1.25$.}
    \label{regions-for-N_nu2}
\end{figure}
The explicit forms of $R_j$, $j=1,\ldots,6$ are as follows:
\begin{align*}
R_1&=\bigl\{(x,y)\in  [0,s_1]\times [0,\ell_{am}\,(x)] \cup [s_1, s_2]\times [q_1(x), \ell_{am}\,(x)]\bigr\},\\
R_2&=\bigl\{(x,y)\in [s_1,s_2] \times [0,q_1(x)]\cup [s_2,s_3] \times [0,q_2(x)] \cup [s_3,s_4] \times [q_3(x),q_2(x)]\bigr\},\\
R_3&=\bigl\{(x,y)\in [s_3,s_4]\times [0,q_3(x)] \cup [s_4,1/2] \times [0,q_2(x)]\bigr\},\\
R_4&=\bigl\{(x,y)\in [s_1,s_2] \times [0,q_1(x)] \cup  [s_4,s_5] \times [q_3(x),\ell_{am}\,(x)]\cup  [s_5,1/2] \times [q_3(x),q_{12}\,(x)]\bigr\},\\
R_5&=\bigl\{(x,y)\in [ s_4,1/2] \times [q_2(x),q_3(x)] \bigr\},\;\;R_6=\bigl\{(x,y)\in [s_5,1/2] \times [q_{12}\,(x),\ell_{am}\,(x)]\bigr\}.
\end{align*}

By symmetry,  $P^r_{2N}=6\,P\bigl( \{X_2,X_3\} \subset \NY^r(X_1),\; X_1 \in T_s \bigr).$

For $r \in [1,4/3)$,
$$ P\bigl( \{X_2,X_3\} \subset \NY^r(X_1),\; X_1 \in T_s \bigr) = \int_0^{1/2}\int_0^{\ell_{am}\,(x)} \frac{A(\NY^r(x_1))^2}{A(T(\Y))^3}dydx=\frac{781}{116640}\,r^4, $$
where $A(\NY^r(x_1))=\frac{\sqrt{3}}{12}\,r^2\left( y+\sqrt{3}\,x \right)^2 $.  Hence for $r \in [1,4/3)$, $P^r_{2N}=\frac{781\,r^4}{19440}.$  Note that the same results also hold for $r \in [4/3,3/2)$.

Next, by symmetry, $P^r_{2G}=6\,P\bigl( \{X_2,X_3\} \subset \G_1\left( X_1,\NY^r \right),\; X_1 \in T_s \bigr),$
and
 $$P\bigl( \{X_2,X_3\} \subset \G_1\left( X_1,\NY^r \right),\; X_1 \in T_s \bigr)=\sum_{j=1}^6 P\bigl( \{X_2,X_3\} \subset \G_1\left( X_1,\NY^r \right),\; X_1 \in R_j \bigr).$$
For $x_1 \in R_1$,
\begin{align*}
&P\bigl( \{X_2,X_3\} \subset \G_1\left( X_1,\NY^r \right),\; X_1 \in R_1 \bigr)=\int_0^{s_1}\int_0^{\ell_{am}\,(x)} \frac{A\left( \G_1\left( x_1,\NY^r \right) \right)^2}{A(T(\Y))^3}dydx\\
&+\int_{s_1}^{s_2}\int_{q_1(x)}^{\ell_{am}\,(x)} \frac{A\left( \G_1\left( x_1,\NY^r \right) \right)^2}{A(T(\Y))^3}dydx = \frac{(211\,r^4-1716\,r^3+5751\,r^2-6696\,r+2511)(2\,r-3)^2}{10935\,r^4},
\end{align*}
where $A\left( \G_1\left( x_1,\NY^r \right) \right)=-\frac{\sqrt{3}\,\left( r^2-\left(\sqrt{3}\,x+y \right)^2\right)}{12\,r^2}$.

For $x_1 \in R_2$,
\begin{eqnarray*}
\lefteqn{P\bigl( \{X_2,X_3\} \subset \G_1\left( X_1,\NY^r \right),\; X_1 \in R_2 \bigr) = \int_{s_1}^{s_2}\int_0^{\ell_{am}\,(x)} \frac{A\left( \G_1\left( x_1,\NY^r \right) \right)^2}{A(T(\Y))^3}dydx}\\
& &+\int_{s_2}^{s_3}\int_{0}^{q_2(x)} \frac{A\left( \G_1\left( x_1,\NY^r \right) \right)^2}{A(T(\Y))^3}dydx+\int_{s_3}^{s_4}\int_{q_3(x)}^{q_2(x)} \frac{A\left( \G_1\left( x_1,\NY^r \right) \right)^2}{A(T(\Y))^3}dydx\\
& = &-\frac{(2\,r-3)(440\,r^4-4091\,r^3+13476\,r^2-16506\,r+6696)}{9720\,r^3}.
\end{eqnarray*}
where $A\left( \G_1\left( x_1,\NY^r \right) \right)=\frac{\sqrt{3}\,\left(-4\,\sqrt{3}\,r\,y-12\,r+12\,r\,x+5\,r^2+2\,y^2+6\,\sqrt{3}\,y-8\,x\,\sqrt{3}\,y+9-18\,x+6\,x^2\right)}{12\,r^2}$.

For $x_1 \in R_3$,
\begin{eqnarray*}
\lefteqn{P\left( \{X_2,X_3\} \subset \G_1\left( X_1,\NY^r \right),\; X_1 \in R_3 \right)}\\
& = &\int_{s_3}^{s_4}\int_0^{q_3(x)} \frac{A\left( \G_1\left( x_1,\NY^r \right) \right)^2}{A(T(\Y))^3}dydx+\int_{s_4}^{1/2}\int_{0}^{q_2(x)} \frac{A\left( \G_1\left( x_1,\NY^r \right) \right)^2}{A(T(\Y))^3}dydx\\
& = &-\frac{(2\,r-3)(21056\,r^5-7845\,r^4+231300\,r^3-943650\,r^2+1127520\,r-428652)}{262440\,r^4}.
\end{eqnarray*}
where $A\left( \G_1\left( x_1,\NY^r \right) \right)=-\frac{\sqrt{3}\,\left( 2\,y^2+2\,\sqrt{3}\,y+3-6\,x+6\,x^2-2\,r^2 \right)}{12\,r^2}$.

For $x_1 \in R_4$,
\begin{eqnarray*}
\lefteqn{P(\{X_2,X_3\} \subset \G_1\left( X_1,\NY^r \right),\; X_1 \in R_4)=\int_{s_2}^{s_4}\int_{q_2(x)}^{\ell_{am}\,(x)} \frac{A\left( \G_1\left( x_1,\NY^r \right) \right)^2}{A(T(\Y))^3}dydx}\\
& &+\int_{s_4}^{s_5}\int_{q_3(x)}^{\ell_{am}\,(x)} \frac{A\left( \G_1\left( x_1,\NY^r \right) \right)^2}{A(T(\Y))^3}dydx+\int_{s_5}^{1/2}\int_{q_3(x)}^{q_{12}\,(x)} \frac{A\left( \G_1\left( x_1,\NY^r \right) \right)^2}{A(T(\Y))^3}dydx\\
& = &-{\frac{12873091}{699840}}\,r^2+{\frac{81239}{648}}\,r+{\frac{14714}{27}}\,r^{-1}-{\frac{4238}{9}}\,r^{-2}+{\frac{656}{3}}\,r^{-3}-{\frac{128}{3}}\,r^{-4}-{\frac{77123}{216}}.
\end{eqnarray*}
where $A\left( \G_1\left( x_1,\NY^r \right) \right)=\frac{\sqrt{3}\,\left( 9\,r^2+18-24\,r+4\,\sqrt{3}\,r
\,y-18\,x+6\,x^2+14\,y^2+12\,r\,x-8\,x\,\sqrt{3}\,y-6\,\sqrt{3}\,y \right)}{12\,r^2}$.

For $x_1 \in R_5$,
\begin{eqnarray*}
\lefteqn{P\left( \{X_2,X_3\} \subset \G_1\left( X_1,\NY^r \right),\; X_1 \in R_5 \right)=\int_{s_4}^{1/2}\int_{q_2(x)}^{q_3(x)} \frac{A\left( \G_1\left( x_1,\NY^r \right) \right)^2}{A(T(\Y))^3}dydx}\\
& = &\frac{(89305\,r^4-364080\,r^3+598320\,r^2-468288\,r+145152)(-6+5\,r)^2}{262440\,r^4}.
\end{eqnarray*}
where $A\left( \G_1\left( x_1,\NY^r \right) \right)=\frac{\sqrt{3}\,\left( 9\,r^2+18-24\,r+4\,\sqrt{3}\,r\,y-18\,x+6\,x^2+14\,y^2+12\,r\,x-8\,x\,\sqrt{3}\,y-6\,\sqrt{3}\,y\right)}{12\,r^2}$.

For $x_1 \in R_6$,
\begin{align*}
&P(\{X_2,X_3\} \subset \G_1\left( X_1,\NY^r \right),\; X_1 \in R_6)=\int_{s_5}^{1/2}\int_{q_{12}{x}}^{\ell_{am}\,(x)} \frac{A\left( \G_1\left( x_1,\NY^r \right) \right)^2}{A(T(\Y))^3}dydx\\
& = \frac{(1081\,r^4-4672\,r^3+7624\,r^2-5568\,r+1536)(-4+3\,r)^2}{960\,r^4}.
\end{align*}
where $A\left( \G_1\left( x_1,\NY^r \right) \right)=-\frac{\sqrt{3}\,\left( \sqrt{3}\,y-2\,r^2-3+3\,x+4\,r-3\,x^2-3\,y^2 \right)}{2\,r^2}$.

So
\begin{eqnarray*}
P(\{X_2,X_3\} \subset \G_1\left( X_1,\NY^r \right))&=&6\,\left( \frac{25687}{349920}\,r^2-\frac{133}{972}\,r+\frac{14}{81}\,r^{-1}-\frac{1}{9}\,r^{-2}+\frac{1}{90}\,r^{-4}-\frac{1}{324} \right)\\
&=&\frac {25687\,r^6-47880\,r^5-1080\,r^4+60480\,r^3-38880\,r^2+3888}{58320\,r^4}.
\end{eqnarray*}

Furthermore, by symmetry, $P^r_M=6\,P(X_2 \in \NY^r(X_1),\;X_3\in \G_1\left( X_1,\NY^r \right),\; X_1 \in T_s),$ and
$$P(X_2 \in \NY^r(X_1),\;X_3\in \G_1\left( X_1,\NY^r \right),\; X_1 \in T_s)\\
=\sum_{j=1}^6 P(X_2 \in \NY^r(X_1),\;X_3\in \G_1\left( X_1,\NY^r \right),\; X_1 \in R_j).
$$
For $x_1 \in R_1$,
\begin{eqnarray*}
\lefteqn{P(X_2 \in \NY^r(X_1),\;X_3\in \G_1\left( X_1,\NY^r \right),\; X_1 \in R_1)=\int_0^{s_1}\int_0^{\ell_{am}\,(x)} \frac{A(\NY^r(x_1))\,A\left( \G_1\left( x_1,\NY^r \right) \right)}{A(T(\Y))^3}dydx}\\
& &+\int_{s_1}^{s_2}\int_{q_1(x)}^{\ell_{am}\,(x)} \frac{A(\NY^r(x_1))\,A\left( \G_1\left( x_1,\NY^r \right) \right)}{A(T(\Y))^3}dydx=-\frac{1}{21870}\,(143\,r^2-744\,r+558)(2\,r-3)^4.
\end{eqnarray*}
For $x_1 \in R_2$,
{\small
\begin{eqnarray*}
\lefteqn{P(X_2 \in \NY^r(X_1),\;X_3\in \G_1\left( X_1,\NY^r \right),\; X_1 \in R_2)=\int_{s_1}^{s_2}\int_0^{\ell_{am}\,(x)} \frac{A(\NY^r(x_1))\,A\left( \G_1\left( x_1,\NY^r \right) \right)}{A(T(\Y))^3}dydx}\\
& &+\int_{s_2}^{s_3}\int_{0}^{q_2(x)} \frac{A(\NY^r(x_1))\,A\left( \G_1\left( x_1,\NY^r \right) \right)}{A(T(\Y))^3}dydx+\int_{s_3}^{s_4}\int_{q_3(x)}^{q_2(x)} \frac{A(\NY^r(x_1))\,A\left( \G_1\left( x_1,\NY^r \right) \right)}{A(T(\Y))^3}dydx\\
& = &\frac{1}{349920}\,r\,(2\,r-3)(23014\,r^4-187311\,r^3+517896\,r^2-594216\,r+241056).\end{eqnarray*}
}
For $x_1 \in R_3$,
{\small
\begin{eqnarray*}
\lefteqn{P(X_2 \in \NY^r(X_1),\;X_3\in \G_1\left( X_1,\NY^r \right),\; X_1 \in R_3)}\\
& = &\int_{s_3}^{s_4}\int_0^{q_3(x)} \frac{A(\NY^r(x_1))\,A\left( \G_1\left( x_1,\NY^r \right) \right)}{A(T(\Y))^3}dydx+\int_{s_4}^{1/2}\int_{0}^{q_2(x)} \frac{A(\NY^r(x_1))\,A\left( \G_1\left( x_1,\NY^r \right) \right)}{A(T(\Y))^3}dydx\\
& = &{\frac{1}{1049760}}\,(2\,r-3)(874\,r^5-297327\,r^4+1858392\,r^3-4298832\,r^2+4280202\,r-1546209).
\end{eqnarray*}
}
For $x_1 \in R_4$,
\begin{eqnarray*}
\lefteqn{P(X_2 \in \NY^r(X_1),\;X_3\in \G_1\left( X_1,\NY^r \right),\; X_1 \in R_4)=\int_{s_2}^{s_4}\int_{q_2(x)}^{\ell_{am}\,(x)} \frac{A(\NY^r(x_1))\,A\left( \G_1\left( x_1,\NY^r \right) \right)}{A(T(\Y))^3}dydx}\\
& &+\int_{s_4}^{s_5}\int_{q_3(x)}^{\ell_{am}\,(x)} \frac{A(\NY^r(x_1))\,A\left( \G_1\left( x_1,\NY^r \right) \right)}{A(T(\Y))^3}dydx+\int_{s_5}^{1/2}\int_{q_3(x)}^{q_{12}\,(x)} \frac{A(\NY^r(x_1))\,A\left( \G_1\left( x_1,\NY^r \right) \right)}{A(T(\Y))^3}dydx\\
& = &-\frac{1}{466560}\,r\,(1762560\,r-497664-2661120\,r^2+201395\,r^5-1017720\,r^4+2212560\,r^3).
\end{eqnarray*}

For $x_1 \in R_5$,
\begin{align*}
&P(X_2 \in \NY^r(X_1),\;X_3\in \G_1\left( X_1,\NY^r \right),\; X_1 \in R_5)=\int_{s_4}^{1/2}\int_{q_2(x)}^{q_3(x)} \frac{A(\NY^r(x_1))\,A\left( \G_1\left( x_1,\NY^r \right) \right)}{A(T(\Y))^3}dydx\\
& = \frac{1}{262440}\,(1570\,r^4-1380\,r^3-11205\,r^2+29700\,r-19116)(-6+5\,r)^2.
\end{align*}

For $x_1 \in R_6$,
\begin{align*}
&P(X_2 \in \NY^r(X_1),\;X_3\in \G_1\left( X_1,\NY^r \right),\; X_1 \in R_6)=\int_{s_5}^{1/2}\int_{q_{12}{x}}^{\ell_{am}\,(x)} \frac{A(\NY^r(x_1))\,A\left( \G_1\left( x_1,\NY^r \right) \right)}{A(T(\Y))^3}dydx\\
& = \frac{1}{51840}\,(1485\,r^4-2064\,r^3+16\,r^2-128\,r+768)(-4+3\,r)^2.
\end{align*}

Thus
\begin{eqnarray*}
&P(X_2 \in \NY^r(X_1),\;X_3\in \G_1\left( X_1,\NY^r \right))=6\,\left( {\frac{3007}{699840}}\,r^6-{\frac{8}{405}}\,r^5+{\frac{5}{648}}\,r^4+\frac{1}{9}\,r^3-{\frac{133}{648}}\,r^2+{\frac{56}{405}}\,r-{\frac{143}{4320}} \right)\\
&={\frac{3007}{116640}}\,r^6-{\frac{16}{135}}\,r^5+{\frac{5}{108}}\,r^4+\frac{2}{3}\,r^3-{\frac{133}{108}}\,r^2+{\frac{112}{135}}\,r-{\frac{143}{720}}.
\end{eqnarray*}

Hence
{\small
\begin{multline*}
\E[h_{12}\,h_{13}]=\Bigl[3007\,r^{10}-13824\,r^9+7743\,r^8+77760\,r^7-117953\,r^6+48888\,r^5-24246\,r^4+60480\,r^3\\
-38880\,r^2+3888\Bigr]/\Bigl[58320\,r^4\Bigr].
\end{multline*}
}
Thus
\begin{multline*}
\nu(r)=\Bigl[3007\,r^{10}-13824\,r^9+898\,r^8+77760\,r^7-117953\,r^6+48888\,r^5-24246\,r^4+60480\,r^3-38880\,r^2\\
+3888\Bigr]/\Bigl[58320\,r^4\Bigr].
\end{multline*}

\begin{figure} []
   \centering
   \scalebox{.3}{\input{G1ofxCase7.pstex_t}}
   \caption{The prototype of the new case for $\G_1\left( x_1,\NY^r \right)$ for $x_1 \in T_s$ for $r \in [4/3,3/2)$.}
\label{fig:G1-NYr-Cases-2}
\end{figure}
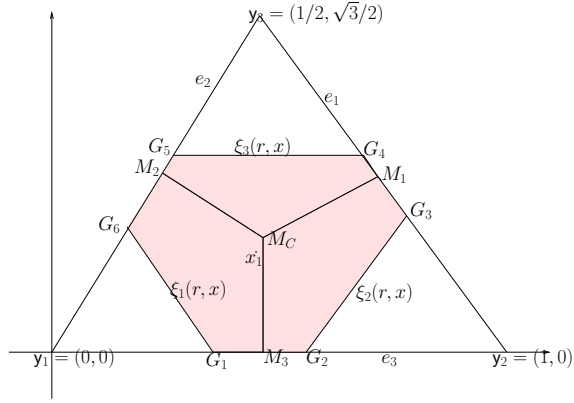

\begin{figure} [ht]
    \centering
   \scalebox{.4}{\input{N_nu6GamRegions.pstex_t}}
    \caption{The regions corresponding to the six cases for $r \in [4/3,3/2)$}
    \label{regions for N_nu6}
\end{figure}
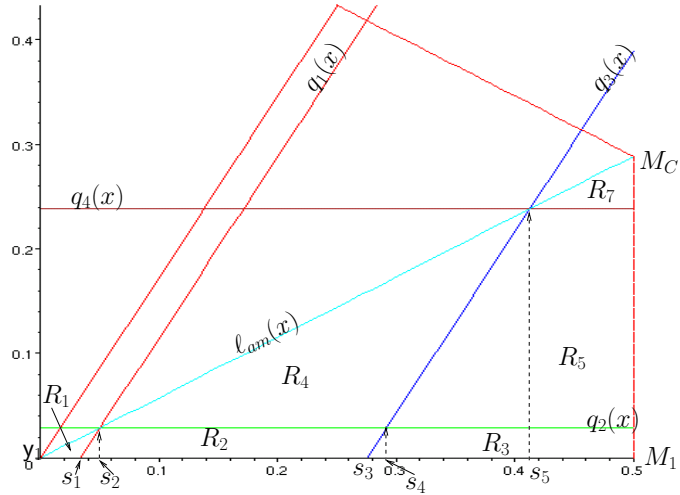

For $r \in [4/3,3/2)$, there are six cases regarding $\G_1\left( x_1,\NY^r \right)$ and one case for $\NY^r(x_1)$. Prototypes of the five of the cases for $\G_1\left( x_1,\NY^r \right)$ are as in case$-j$ for $j=1,\ldots,5$ in Figure \ref{fig:G1-NYr-Cases-1} and the new case, case-7, is depicted in Figure \ref{fig:G1-NYr-Cases-2}.
Each case $j$ corresponds to the region $R_j$ in Figure \ref{regions for N_nu6} where  $s_1=1-2\,r/3, \; s_2=3/2-r,\; s_3=1-r/2, \; s_4=3/2-5\,r/6, \; s_5=3/2-3\,r/4.$
The explicit forms of $R_j$, $j=1,2,3$ are same as before, for $j=4,5,7$ are given below:
{\small
\begin{align*}
R_4&=\bigl\{(x,y)\in [s_2,s_4] \times [q_2(x),\ell_{am}\,(x)] \cup  [s_4,s_6] \times [q_3(x),\ell_{am}\,(x)]\bigr\}\\
R_5&=\bigl\{(x,y)\in [ s_4,s_6] \times [q_2(x),q_3(x)] \cup [s_6,1/2] \times [q_2(x),q_4(x)]\bigr\}\\
R_7&=\bigl\{(x,y)\in [s_6,1/2] \times [q_4(x),\ell_{am}\,(x)]\bigr\}
\end{align*}
}

where $\ell_{am}\,(x)=x/\sqrt{3}$, $q_1(x)=(2\,r-3)/\sqrt{3}+\sqrt{3}\,x$, $q_2(x)=\sqrt{3}\,(1/2-r/3)$, $q_3(x)=\sqrt{3}\,(x-1+r/2)$, and $q_4(x)=\sqrt{3}\,(1/2-r/4)$.

Then $P^r_{2N}=\frac{781\,r^4}{19440}.$ We use the same limits of integration in $\mu(r)$ calculations with the integrand $A(\NY^r(x_1))^2/A(T(\Y))^3$.

Next, by symmetry, $P^r_{2G}=6\,P\bigl( \{X_2,X_3\} \subset \G_1\left( X_1,\NY^r \right),\; X_1 \in T_s \bigr),$ and, let $S_I:=\{1,2,3,4,5,7\}$, then
 $$P\bigl( \{X_2,X_3\} \subset \G_1\left( X_1,\NY^r \right),\; X_1 \in T_s \bigr)=\sum_{j=S_I} P\left( \{X_2,X_3\} \subset \G_1\left( X_1,\NY^r \right),\; X_1 \in R_j \right).$$
For $x_1 \in R_j$, $j=1,2,3 $ we get the same result as before.

For $x_1 \in R_4$,
\begin{eqnarray*}
\lefteqn{P(\{X_2,X_3\} \subset \G_1\left( X_1,\NY^r \right),\; X_1 \in R_4) = \int_{s_2}^{s_4}\int_{q_2(x)}^{\ell_{am}\,(x)} \frac{A\left( \G_1\left( x_1,\NY^r \right) \right)^2}{A(T(\Y))^3}dydx}\\
& &+\int_{s_4}^{s_6}\int_{q_3(x)}^{\ell_{am}\,(x)} \frac{A\left( \G_1\left( x_1,\NY^r \right) \right)^2}{A(T(\Y))^3}dydx =\frac {9637\,r^4-89640\,r^3+288360\,r^2-362880\,r+155520}{349920\,r^2}.
\end{eqnarray*}
where $A\left( \G_1\left( x_1,\NY^r \right) \right)=\frac{\sqrt{3}\,\left( 9\,r^2+18-24\,r+4\,\sqrt{3}\,r\,y-18\,x+6\,x^2+14\,y^2+12\,r\,x-8\,x\,\sqrt{3}\,y-6\,\sqrt{3}\,y \right)}{12\,r^2}$.

For $x_1 \in R_5$,
\begin{eqnarray*}
\lefteqn{P\left( \{X_2,X_3\} \subset \G_1\left( X_1,\NY^r \right),\; X_1 \in R_5 \right)}\\
& = &\int_{s_4}^{s_6}\int_{q_2(x)}^{q_3(x)} \frac{A\left( \G_1\left( x_1,\NY^r \right) \right)^2}{A(T(\Y))^3}dydx+\int_{s_6}^{1/2}\int_{q_2(x)}^{q_4(x)} \frac{A\left( \G_1\left( x_1,\NY^r \right) \right)^2}{A(T(\Y))^3}dydx\\
& = &\frac {87251\,r^5+13219200\,r-11214720\,r^2-5225472+3377160\,r^3-261288\,r^4}{2099520\,r^3}.
\end{eqnarray*}
where $A\left( \G_1\left( x_1,\NY^r \right) \right)$ is same as before.

For $x_1 \in R_7$,
\begin{eqnarray*}
&P\left( \{X_2,X_3\} \subset \G_1\left( X_1,\NY^r \right),\; X_1 \in R_7 \right)=\int_{s_6}^{1/2}\int_{q_4(x)}^{\ell_{am}\,(x)} \frac{A\left( \G_1\left( x_1,\NY^r \right) \right)^2}{A(T(\Y))^3}dydx\\
& =\frac{(57\,r^4+96\,r^3-72\,r^2-576\,r+512)(-4+3\,r)^2}{2880\,r^4}.
\end{eqnarray*}
where $A\left( \G_1\left( x_1,\NY^r \right) \right)=-\frac{\sqrt{3}\,\left( 6\,y^2-2\,\sqrt{3}\,y+6-6\,x+6\,x^2-3\,r^2 \right)}{12\,r^2}$.

So,
\begin{eqnarray*}
P^r_{2G}&=&6\,\Biggl(\frac {-47880\,r^5-38880\,r^2+25687\,r^6-1080\,r^4+60480\,r^3+3888}{349920\,r^4}\Biggr)\\
&=&\frac {-47880\,r^5-38880\,r^2+25687\,r^6-1080\,r^4+60480\,r^3+3888}{58320\,r^4}.
\end{eqnarray*}

Furthermore,
$$
P^r_M=\sum_{j\in S_I} P(X_2 \in \NY^r(X_1),\;X_3\in \G_1\left( X_1,\NY^r \right),\; X_1 \in R_j).
$$

For $x_1 \in R_j$, $j=1,2,3 $ we get the same result as before.

For $x_1 \in R_4$,
\begin{eqnarray*}
\lefteqn{P\left( X_2 \in \NY^r(X_1),\;X_3\in \G_1\left( X_1,\NY^r \right),\; X_1 \in R_4 \right)}\\
& = &\int_{s_2}^{s_4}\int_{q_2(x)}^{\ell_{am}\,(x)} \frac{A(\NY^r(x_1))\,A\left( \G_1\left( x_1,\NY^r \right) \right)}{A(T(\Y))^3}dydx+\int_{s_4}^{s_6}\int_{q_3(x)}^{\ell_{am}\,(x)} \frac{A(\NY^r(x_1))\,A\left( \G_1\left( x_1,\NY^r \right) \right)}{A(T(\Y))^3}dydx\\
& = &-{\frac{1}{466560}}\,r^2(207360+404640\,r^2-483840\,r-142920\,r^3+17687\,r^4).
\end{eqnarray*}

For $x_1 \in R_5$,
\begin{eqnarray*}
\lefteqn{P\left( X_2 \in \NY^r(X_1),\;X_3\in \G_1\left( X_1,\NY^r \right),\; X_1 \in R_5 \right)}\\
& = &\int_{s_4}^{s_6}\int_{q_2(x)}^{q_3(x)} \frac{A(\NY^r(x_1))\,A\left( \G_1\left( x_1,\NY^r \right) \right)}{A(T(\Y))^3}dydx+\int_{s_6}^{1/2}\int_{q_2(x)}^{q_4(x)} \frac{A(\NY^r(x_1))\,A\left( \G_1\left( x_1,\NY^r \right) \right)}{A(T(\Y))^3}dydx\\
& = &-\frac{r\,(399064320\,r-150792192+171990000\,r^3-391461120\,r^2-31140648\,r^4+1230359\,r^5}{67184640}.
\end{eqnarray*}

For $x_1 \in R_7$,
\begin{eqnarray*}
&P\left( X_2 \in \NY^r(X_1),\;X_3\in \G_1\left( X_1,\NY^r \right),\; X_1 \in R_7 \right)=\int_{s_6}^{1/2}\int_{q_4(x)}^{\ell_{am}\,(x)} \frac{A(\NY^r(x_1))\,A\left( \G_1\left( x_1,\NY^r \right) \right)}{A(T(\Y))^3}dydx\\
& = \frac{1}{829440}\,(2727\,r^4-3648\,r^3-52736\,r^2+166656\,r-121600)(-4+3\,r)^2.
\end{eqnarray*}

Then,
\begin{eqnarray*}
&P^r_M=6\,\Biggl({\frac{5467}{2799360}}\,r^6-{\frac{35}{2592}}\,r^5+{\frac{37}{1296}}\,r^4-{\frac{13}{648}}\,r^2+{\frac{83}{12960}}\Biggr)={\frac{5467}{466560}}\,r^6-{\frac{35}{432}}\,r^5+{\frac{37}{216}}\,r^4-{\frac{13}{108}}\,r^2+\frac{83}{216}.
\end{eqnarray*}

So,
\begin{multline*}
\E[h_{12}\,h_{13}]=\Bigl[5467\,r^{10}-37800\,r^9+89292\,r^8+46588\,r^6-191520\,r^5+13608\,r^4+241920\,r^3-155520\,r^2\\
+15552\Bigr]/\Bigl[233280\,r^4\Bigr].
\end{multline*}

Thus, for $r \in [4/3,3/2)$
\begin{multline*}
\nu(r)=\Bigl[5467\,r^{10}-37800\,r^9+61912\,r^8+46588\,r^6-191520\,r^5+13608\,r^4+241920\,r^3-155520\,r^2\\
+15552\Bigr]/\Bigl[233280\,r^4\Bigr].
\end{multline*}

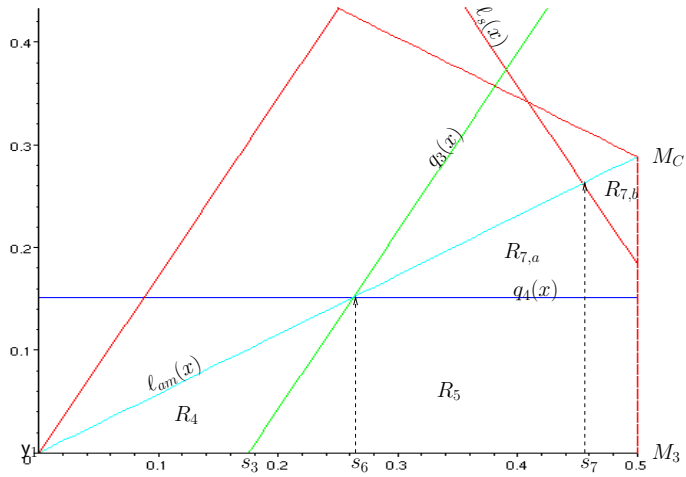
\begin{figure} []
    \centering
   \scalebox{.4}{\input{N_nu8GamRegions.pstex_t}}
    \caption{The regions corresponding to the three cases for $r \in [3/2,2)$ with $r=1.65 $}
    \label{regions for N_nu8}
\end{figure}

For $r \in [3/2,2)$, there are three cases regarding $\G_1\left( x_1,\NY^r \right)$ and two cases for $\NY^r(x_1)$.  The prototypes of these three cases as in cases 4,5, and 7 of Figures \ref{fig:G1-NYr-Cases-1} and \ref{fig:G1-NYr-Cases-2}.
Each case $j$, corresponds to the region $R_j$ in Figure \ref{regions for N_nu8} where $q_j(x)$ are same as before for $j=3,4$, and $s_j$, $j=3,4,6$ are same as before and $s_7=3/(4\,r)$. Observe that for $x_1 \in R_4 \cup R_5 \cup R_{7a}$, $\NY^r(x_1)=T_r(x_1)\subsetneq T(\Y)$, and for $x_1 \in R_{7b}$, $\NY^r(x_1)=T(\Y)$. So there are four regions to consider to calculate the covariance.

Then, for $x_1=(x,y)\in R_j$, $\G_1\left( x_1,\NY^r \right)$ are same as before for $j=4,5,7$. The explicit forms of $R_j$, $j=4, 5, 7a, 7b$ are given below (the explicit form of $R_7$ is same as before):
\begin{align*}
R_4&=\bigl\{(x,y)\in [0,s_3] \times [0,\ell_{am}\,(x)] \cup  [s_3,s_6] \times [q_3(x),\ell_{am}\,(x)]\bigr\}\\
R_5&=\bigl\{(x,y)\in [s_3,s_6] \times [0,q_3(x)] \cup [s_6,1/2] \times [0,q_4(x)]\bigr\}\\
R_{7,a}&=\bigl\{(x,y)\in [s_6,s_7] \times [q_4(x),\ell_{am}\,(x)] \cup [s_7,1/2] \times [q_4(x),\ell_s(r,x)]\bigr\}\\
R_{7,b}&=\bigl\{(x,y)\in [s_7,1/2] \times [\ell_s(r,x),\ell_{am}\,(x)]\bigr\}
\end{align*}

Now,
{\small
\begin{align*}
&P(\{X_2,X_3\} \subset \NY^r(X_1), X_1 \in T_s)= \int_0^{1/2}\int_0^{\ell_{am}\,(x)} \frac{A(\NY^r(x_1))^2}{A(T(\Y))^3}dydx=\int_0^{s_7}\int_0^{\ell_{am}\,(x)} \frac{A(\NY^r(x_1))^2}{A(T(\Y))^3}dydx\\
&  +\int_{s_7}^{1/2}\int_0^{\ell_s(x)} \frac{A(\NY^r(x_1))^2}{A(T(\Y))^3}dydx +\int_{s_7}^{1/2}\int_{\ell_s(x)}^{\ell_{am}\,(x)} \frac{1}{A(T(\Y))}dydx =-\frac{-480+r^6+768\,r-320\,r^2}{480\,r^2}.
\end{align*}
}
Hence  $P^r_{2N}=-\frac{-480+r^6+768\,r-320\,r^2}{80\,r^2}.$

Next, by symmetry, $P^r_{2G}=6\,P\bigl( \{X_2,X_3\} \subset \G_1\left( X_1,\NY^r \right),\; X_1 \in T_s \bigr),$ and, let $S_I:=\{4,5,7\}$, then
 $$P\bigl( \{X_2,X_3\} \subset \G_1\left( X_1,\NY^r \right),\; X_1 \in T_s \bigr)=\sum_{j=S_I} P(\{X_2,X_3\} \subset \G_1\left( X_1,\NY^r \right),\; X_1 \in R_j).$$
For $x_1 \in R_4$,
\begin{eqnarray*}
\lefteqn{P\left( \{X_2,X_3\} \subset \G_1\left( X_1,\NY^r \right),\; X_1 \in R_4 \right)=\int_{0}^{s_3}\int_{0}^{\ell_{am}\,(x)} \frac{A\left( \G_1\left( x_1,\NY^r \right) \right)^2}{A(T(\Y))^3}dydx}\\
& &+\int_{s_3}^{s_6}\int_{q_3(x)}^{\ell_{am}\,(x)} \frac{A\left( \G_1\left( x_1,\NY^r \right) \right)^2}{A(T(\Y))^3}dydx = \frac{(237\,r^4-956\,r^3+1728\,r^2-1584\,r+592)(-2+r)^2}{480\,r^4}.
\end{eqnarray*}
For $x_1 \in R_5$,
\begin{align*}
&P\left( \{X_2,X_3\} \subset \G_1\left( X_1,\NY^r \right),\; X_1 \in R_5 \right) = \int_{s_3}^{s_6}\int_{0}^{q_3(x)} \frac{A\left( \G_1\left( x_1,\NY^r \right) \right)^2}{A(T(\Y))^3}dydx\\
& +\int_{s_6}^{1/2}\int_{0}^{q_4(x)} \frac{A\left( \G_1\left( x_1,\NY^r \right) \right)^2}{A(T(\Y))^3}dydx = \frac{(r-2)(1909\,r^5-6142\,r^4+10036\,r^3-14808\,r^2+15024\,r-6048)}{2880\,r^4}.
\end{align*}

For $x_1 \in R_7$ the result is same as before. So
\begin{eqnarray*}
P^r_{2G}&=&6\,\left(\frac{7320\,r^4-984\,r^5+13\,r^6-20480\,r^3+27840\,r^2-18816\,r+5152}{1440\,r^4}\right)\\
&=&\frac{7320\,r^4-984\,r^5+13\,r^6-20480\,r^3+27840\,r^2-18816\,r+5152}{240\,r^4}.
\end{eqnarray*}

Furthermore,
 $$P\left( X_2 \in \NY^r(X_1),\;X_3\in \G_1\left( X_1,\NY^r \right),\; X_1 \in T_s \right)=\int_{0}^{1/2}\int_{0}^{\ell_{am}\,(x)} \frac{A(\NY^r(x_1))\,A\left( \G_1\left( x_1,\NY^r \right) \right)}{A(T(\Y))^3}dydx.
$$
For $x_1 \in R_4$,
\begin{eqnarray*}
\lefteqn{P(X_2 \in \NY^r(X_1),\;X_3\in \G_1\left( X_1,\NY^r \right),\; X_1 \in R_4)=\int_{0}^{s_3}\int_{0}^{\ell_{am}\,(x)} \frac{A(\NY^r(x_1))\,A\left( \G_1\left( x_1,\NY^r \right) \right)}{A(T(\Y))^3}dydx}\\
& & +\int_{s_3}^{s_6}\int_{q_3(x)}^{\ell_{am}\,(x)} \frac{A(\NY^r(x_1))\,A\left( \G_1\left( x_1,\NY^r \right) \right)}{A(T(\Y))^3}dydx = \frac{1}{1920}\,(99\,r^2-16\,r-84)(-2+r)^4.
\end{eqnarray*}
For $x_1 \in R_5$,
\begin{eqnarray*}
\lefteqn{P\left( X_2 \in \NY^r(X_1),\;X_3\in \G_1\left( X_1,\NY^r \right),\; X_1 \in R_5\right)}\\
& = &\int_{s_3}^{s_6}\int_{0}^{q_3(x)} \frac{A(\NY^r(x_1))\,A\left( \G_1\left( x_1,\NY^r \right) \right)}{A(T(\Y))^3}dydx+\int_{s_6}^{1/2}\int_{0}^{q_4(x)} \frac{A(\NY^r(x_1))\,A\left( \G_1\left( x_1,\NY^r \right) \right)}{A(T(\Y))^3}dydx\\
& = &-{\frac{1}{92160}}\,(-2+r)(7535\,r^5-35210\,r^4+9500\,r^3+181560\,r^2-308880\,r+147168).
\end{eqnarray*}
For $x_1 \in R_{7a}$,
\begin{eqnarray*}
\lefteqn{P\left( X_2 \in \NY^r(X_1),\;X_3\in \G_1\left( X_1,\NY^r \right),\; X_1 \in R_{7a} \right)=\int_{s_6}^{s_7}\int_{q_4(x)}^{\ell_{am}\,(x)} \frac{A(\NY^r(x_1))\,A\left( \G_1\left( x_1,\NY^r \right) \right)}{A(T(\Y))^3}dydx}\\
& &+\int_{s_7}^{1/2}\int_{q_4(x)}^{\ell_s(x)} \frac{A(\NY^r(x_1))\,A\left( \G_1\left( x_1,\NY^r \right) \right)}{A(T(\Y))^3}dydx = {\frac{303}{10240}}\,r^6-{\frac{91}{768}}\,r^5-{\frac{53}{128}}\,r^4+{\frac{235}{72}}\,r^3-{\frac{173}{24}}\,r^2\\
& & +{\frac{101}{15}}\,r+\frac{2}{3}\,r^{-1}-\frac{3}{4}\,r^{-2}-{\frac{16}{9}}\,r^{-3}+4\,r^{-4}-{\frac{18}{5}}\,r^{-5}+\frac{3}{2}\,r^{-6}-{\frac{34}{15}}.
\end{eqnarray*}

For $x_1 \in R_{7b}$,
\begin{align*}
&P\left( X_2 \in \NY^r(X_1),\;X_3\in \G_1\left( X_1,\NY^r \right),\; X_1 \in R_{7b} \right) =\\
& \int_{s_7}^{1/2}\int_{\ell_s(x)}^{\ell_{am}\,(x)} \frac{A(\NY^r(x_1))\,A\left( \G_1\left( x_1,\NY^r \right) \right)}{A(T(\Y))^3}dydx = \frac{(2\,r^4-4\,r^2+4\,r-3)(2\,r-3)^2}{12\,r^6}.
\end{align*}
Hence
\begin{multline*}
P^r_{2G}=6\,\Bigl(-\Bigl[7\,r^{12}-72\,r^{11}+240\,r^{10}-1440\,r^8+3456\,r^7-10296\,r^6+15360\,r^5+6720\,r^4-40960\,r^3\\
+46080\,r^2-27648\,r+8640\Bigr]/\Bigl[11520\,r^6\Bigr]\Bigr)=-\Bigl[7\,r^{12}-72\,r^{11}+240\,r^{10}-1440\,r^8+3456\,r^7\\
-10296\,r^6+15360\,r^5+6720\,r^4-40960\,r^3+46080\,r^2-27648\,r+8640\Bigr]/\Bigl[1920\,r^6\Bigr].
\end{multline*}
Thus
\begin{multline*}
\E[h_{12}\,h_{13}]=-\Bigl[7\,r^{12}-72\,r^{11}+252\,r^{10}-1492\,r^8+7392\,r^7-43416\,r^6+106496\,r^5-110400\,r^4+\\
34304\,r^3+25472\,r^2-27648\,r+8640\Bigr]/\Bigl[960\,r^6\Bigr].
\end{multline*}
Therefore, for $r \in [3/2,2)$
\begin{multline*}
\nu(r)=-\Bigl[7\,r^{12}-72\,r^{11}+312\,r^{10}-5332\,r^8+15072\,r^7+13704\,r^6-139264\,r^5+273600\,r^4-242176\,r^3\\
+103232\,r^2-27648\,r+8640\Bigr]/\Bigl[960\,r^6\Bigr].
\end{multline*}

For $r \in [2,\infty)$, there is only one case regarding $\G_1\left( x_1,\NY^r \right)$, namely $R_7$, and two cases regarding $\NY^r(x_1)$. Furthermore, $s_7$, is same as before and $s_8=1/r$.  Observe that for $x_1 \in R_{7a}$, $\NY^r(x_1)=T_r(x_1)\subsetneq T(\Y)$, and for $x_1 \in R_{7b}$, $\NY^r(x_1)=T(\Y)$. So there are two regions to consider to calculate the covariance.

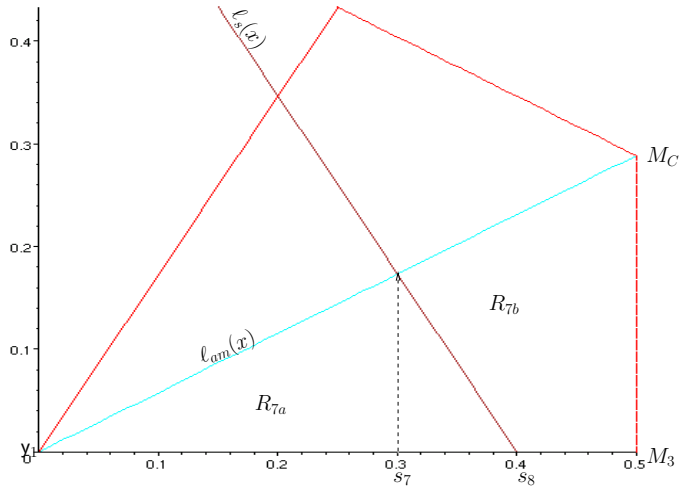
\begin{figure} [ht]
    \centering
   \scalebox{.4}{\input{N_nu10GamRegions.pstex_t}}
    \caption{The regions corresponding to the two cases for $\NY^r(x_1)$ for $r \in [2,\infty)$ with $r=2.5$}
    \label{regions for N_nu10}
\end{figure}
For $x_1=(x,y)\in R_7$, $\G_1\left( x_1,\NY^r \right)$ is same as before. The explicit form of $R_7$, is same as $T_s$.  For $R_{7a}$ and $R_{7b}$, see below:
\begin{align*}
R_{7,a}&=\bigl\{(x,y)\in [0,s_7] \times [0,\ell_{am}\,(x)] \cup [s_7,s_8] \times [0,\ell_s(r,x)]\bigr\}\\
R_{7,b}&=\bigl\{(x,y)\in [s_7,s_8] \times [\ell_s(r,x),\ell_{am}\,(x)]\cup [s_8,1/2] \times [0,\ell_{am}\,(r,x)]\bigr\}
\end{align*}
Now,
\begin{align*}
&P(\{X_2,X_3\} \subset \NY^r(X_1), X_1 \in T_s)= \int_0^{s_7}\int_0^{\ell_{am}\,(x)} \frac{A(\NY^r(x_1))^2}{A(T(\Y))^3}dydx+\int_{s_7}^{s_8}\int_0^{\ell_s(x)} \frac{A(\NY^r(x_1))^2}{A(T(\Y))^3}dydx\\
& +\int_{s_7}^{s_8}\int_{\ell_s(x)}^{\ell_{am}\,(x)} \frac{1}{A(T(\Y))}dydx+\int_{s_8}^{1/2}\int_{0}^{\ell_{am}\,(x)} \frac{1}{A(T(\Y))}dydx = -\frac{1}{3}\,r^{-2}+\frac{1}{6}.
\end{align*}
Hence $P^r_{2N}=1-2\,r^{-2}.$
Next,
\begin{align*}
&P\bigl( \{X_2,X_3\} \subset \G_1\left( X_1,\NY^r \right),\; X_1 \in T_s \bigr)=P(\{X_2,X_3\} \subset \G_1\left( X_1,\NY^r \right),\; X_1 \in R_7)\\
& = \int_{0}^{1/2}\int_{0}^{\ell_{am}\,(x)} \frac{A\left( \G_1\left( x_1,\NY^r \right) \right)^2}{A(T(\Y))^3}dydx=\frac{34-45\,r^2+15\,r^4}{90\,r^4}.
\end{align*}
So $P^r_{2G}=\frac{34-45\,r^2+15\,r^4}{15\,r^4}.$

Furthermore,
\begin{eqnarray*}
&P(X_2 \in \NY^r(X_1),\;X_3\in \G_1\left( X_1,\NY^r \right), X_1 \in T_s)= \int_0^{1/2}\int_0^{\ell_{am}\,(x)} \frac{A(\NY^r(x_1))\,A\left( \G_1\left( x_1,\NY^r \right) \right)}{A(T(\Y))^3}dydx\\
& = \int_0^{s_7}\int_0^{\ell_{am}\,(x)} \frac{A(\NY^r(x_1))\,A\left( \G_1\left( x_1,\NY^r \right) \right)}{A(T(\Y))^3}dydx+\int_{s_7}^{s_8}\int_0^{\ell_s(x)} \frac{A(\NY^r(x_1))\,A\left( \G_1\left( x_1,\NY^r \right) \right)}{A(T(\Y))^3}dydx\\
& +\int_{s_7}^{s_8}\int_{\ell_s(x)}^{\ell_{am}\,(x)} \frac{A\left( \G_1\left( x_1,\NY^r \right) \right)}{A(T(\Y))^2}dydx+\int_{s_8}^{1/2}\int_{0}^{\ell_{am}\,(x)} \frac{A\left( \G_1\left( x_1,\NY^r \right) \right)}{A(T(\Y))^2}dydx=\frac{25-48\,r+90\,r^2-90\,r^4+30\,r^6}{180\,r^6}.
\end{eqnarray*}
So $P^r_M=\frac{25-48\,r+90\,r^2-90\,r^4+30\,r^6}{30\,r^6}. $

Hence, $\E[h_{12}\,h_{13}]=\frac {60\,r^6-165\,r^4+124\,r^2-48\,r+25}{15\,r^6}.$ Thus, for $r \in [2,\infty)$,
$$\nu(r)=\frac {15\,r^4-11\,r^2-48\,r+25}{15\,r^6}.$$

For $r=\infty $, it is trivial to see that $\nu(r)=0$.

\section*{Appendix 2: The Mean $\mu(r,\epsilon)$ Under Segregation and Association Alternatives}
Derivation of $\mu(r,\epsilon)$ involves detailed geometric calculations and partitioning of the space of $(r,\epsilon,x_1)$ for $r \in [1,\infty)$, $\epsilon \in \bigl[ 0,\sqrt{3}/3 \bigr)$, and $x_1 \in T_s=T_s$ See Appendix 3 for the derivation of $\mu(r,\epsilon)$ at a demonstrative interval.

\subsection*{$\mu_S(r,\epsilon)$ Under Segregation Alternatives}
Under segregation, we compute $\mu_S(r,\epsilon)$ explicitly.
For $\epsilon \in \bigl[ 0,\sqrt{3}/8 \bigr)$, $\mu_S(r,\epsilon)=\sum_{j=1}^7 \varpi_{1,j}(r,\epsilon)\,\I(r \in \mI_j)$
where
{\small
\begin{align*}
\varpi_{1,1}(r,\epsilon)&=-\frac{576\,r^2\epsilon^4-1152\,\epsilon^4-37\,r^2+288\,\epsilon^2}{216\,(2\,\epsilon+1)^2(2\,\epsilon-1)^2},\\
\varpi_{1,2}(r,\epsilon)&=-\Bigl[ 576\,r^4\,\epsilon^4-1152\,r^2\epsilon^4+91\,r^4+512\,\sqrt{3}\,r^3\epsilon+2592\,r^2\epsilon^2+1536\,\sqrt{3}\,r\,\epsilon^3+1152\,\epsilon^4\,\\
&-768\,r^3-2304\,\,\sqrt{3}\,r^2\epsilon-6912\,r\,\epsilon^2-2304\,\,\sqrt{3}\,\epsilon^3+1728\,r^2+3456\,\sqrt{3}\,r\,\epsilon+5184\,\,\epsilon^2\\
&-1728\,r-1728\,\sqrt{3}\,\epsilon+648\Bigl]/\Bigl[ 216\,r^2(2\,\epsilon+1)^2(2\,\epsilon-1)^2\Bigl],
\end{align*}
}
{\small
\begin{align*}
\varpi_{1,3}(r,\epsilon)&=-\Bigl[ 192\,r^4\,\epsilon^4-384\,\,r^2\epsilon^4+9\,r^4+864\,\,r^2\epsilon^2+512\,\sqrt{3}\,r\,\epsilon^3+384\,\,\epsilon^4-2304\,\,r\,\epsilon^2-768\,\sqrt{3}\,\epsilon^3\\
&-288\,r^2+1728\,\epsilon^2+576\,r-324\,\Bigl]/\Bigl[ 72\,r^2(2\,\epsilon+1)^2(2\,\epsilon-1)^2\Bigl],\\
\varpi_{1,4}(r,\epsilon)&=-\Bigl[ 192\,r^4\,\epsilon^4-384\,\,r^2\epsilon^4-9\,r^4-96\,\sqrt{3}\,r^3\epsilon+288\,r^2\epsilon^2-128\,\epsilon^4+144\,r^3+576\,\sqrt{3}\,r^2\epsilon+256\\
&\sqrt{3}\,\epsilon^3-720\,r^2-1152\,\sqrt{3}\,r\,\epsilon-576\,\epsilon^2+1152\,r+768\,\sqrt{3}\,\epsilon-612\Bigl]/\Bigl[ 72\,r^2(2\,\epsilon+1)^2(2\,\epsilon-1)^2\Bigl],\\
\varpi_{1,5}(r,\epsilon)&=-\frac{48\,r^4{\epsilon}^4-96\,r^2{\epsilon}^4+72\,r^2{\epsilon}^2-32\,{\epsilon}^4+64\,\,\sqrt{3}{\epsilon}^3-18\,r^2-144\,{\epsilon}^2+27}{18\,r^2(2\,\epsilon+1)^2(2\,\epsilon-1)^2},\\
\varpi_{1,6}(r,\epsilon)&=\frac{48\,r^4{\epsilon}^4+256\,r^3{\epsilon}^4-128\,\sqrt{3}r^3{\epsilon}^3+288\,r^2{\epsilon}^4-192\,\sqrt{3}r^2{\epsilon}^3+72\,r^2{\epsilon}^2+18\,r^2+48\,\sqrt{3}\,\epsilon-45}{18\,(2\,\epsilon+1)^2(2\,\epsilon-1)^2r^2},\\
\varpi_{1,7}(r,\epsilon)&=1,
\end{align*}
}
with the corresponding intervals $\mI_1=\Bigl[ 1,3/2-\sqrt{3}\,\epsilon \Bigr)$, $\mI_2=\Bigl[ 3/2-\sqrt{3}\,\epsilon,3/2\Bigr)$, $\mI_3=\Bigl[ 3/2,2-4\,\,\epsilon/\sqrt{3}\Bigr)$, $\mI_4=\Bigl[ 2-4\,\,\epsilon/\sqrt{3},2 \Bigr)$, $\mI_5=\Bigl[ 2,\sqrt{3}/(2\,\epsilon)-1 \Bigr)$, $\mI_6=\Bigl[ \sqrt{3}/(2\,\epsilon)-1,\sqrt{3}/(2\,\epsilon) \Bigr)$, and $\mI_7=\Bigl[ \sqrt{3}/(2\,\epsilon),\infty \Bigr)$.

For $\epsilon \in \Bigl[ \sqrt{3}/8,\sqrt{3}/6 \Bigr)$,
$\mu_S(r,\epsilon)=\sum_{j=1}^7 \varpi_{2,j}(r,\epsilon)\,\I(r \in \mI_j)$ where $\varpi_{2,j}(r,\epsilon)=\varpi_{1,j}(r,\epsilon)$ for $j=1,2,4,5,6$, and for $j=3,7$,
\begin{align*}
\varpi_{2,3}(r,\epsilon)&=-\Bigl[576\,r^4{\epsilon}^4-1152\,r^2{\epsilon}^4+37\,r^4+224\,\,\sqrt{3}r^3\,\epsilon+864\,\,r^2{\epsilon}^2-384\,\,{\epsilon}^4-336\,r^3-576\,\sqrt{3}r^2\,\epsilon\\
&+768\,\sqrt{3}{\epsilon}^3+432\,r^2-1728\,{\epsilon}^2+576\,\sqrt{3}\,\epsilon-216\Bigr]/\Bigl[216\,r^2(2\,\epsilon+1)^2(2\,\epsilon-1)^2\Bigr],\\
\varpi_{2,7}(r,\epsilon)&=1,
\end{align*}
with the corresponding intervals $\mI_1=\Bigl[ 1,3/2-\sqrt{3}\,\epsilon \Bigr)$, $\mI_2=\Bigl[ 3/2-\sqrt{3}\,\epsilon,2-4\,\,\epsilon/\sqrt{3} \Bigr)$, $\mI_3=\Bigl[ 2-4\,\,\epsilon/\sqrt{3},3/2 \Bigr)$, $\mI_4=[3/2,2)$, $\mI_5=\Bigl[2,\sqrt{3}/(2\,\epsilon)-1\Bigr)$, $\mI_6=\Bigl[ \sqrt{3}/(2\,\epsilon)-1,\sqrt{3}/(2\,\epsilon) \Bigr)$, and $\mI_5=\Bigl[\sqrt{3}/(2\,\epsilon),\infty \Bigr)$.

For $\epsilon \in \Bigl[ \sqrt{3}/6,\sqrt{3}/4 \Bigr)$,
{\small
$\mu_S(r,\epsilon)=\sum_{j=1}^6 \varpi_{3,j}(r,\epsilon)\,\I(r \in \mI_j)$ where $\varpi_{3,1}(r,\epsilon)=\varpi_{1,2}(r,\epsilon)$ and
\begin{align*}
\varpi_{3,2}(r,\epsilon)&=-\Bigl[576\,r^4\,\epsilon^4-1152\,r^2\epsilon^4+37\,r^4+224\,\,\sqrt{3}\,r^3\epsilon+864\,\,r^2\epsilon^2-384\,\,\epsilon^4-336\,r^3-576\,\sqrt{3}\,r^2\epsilon\\
&+768\,\sqrt{3}\,\epsilon^3+432\,r^2-1728\,\epsilon^2+576\,\sqrt{3}\,\epsilon-216\Bigr]/\Bigl[216\,r^2(2\,\epsilon+1)^2(2\,\epsilon-1)^2\Bigr],\\
\varpi_{3,3}(r,\epsilon)&=\Bigl[576\,r^2\epsilon^4+3072\,r\,\epsilon^4-1536\,\sqrt{3}\,r\,\epsilon^3+3456\,\epsilon^4-2304\,\,\sqrt{3}\,\epsilon^3-37\,r^2-224\,\,\sqrt{3}\,r\,\epsilon\\
&+864\,\,\epsilon^2+336\,r+576\,\sqrt{3}\,\epsilon-432\Bigr]/\Bigl[216\,(2\,\epsilon+1)^2(2\,\epsilon-1)^2\Bigr],
\end{align*}

\begin{align*}
\varpi_{3,4}(r,\epsilon)&=\Bigl[192\,r^4\,\epsilon^4+1024\,\,r^3\epsilon^4-512\,\sqrt{3}\,r^3\epsilon^3+1152\,r^2\epsilon^4-768\,\sqrt{3}\,r^2\epsilon^3+9\,r^4+96\,\sqrt{3}\,r^3\epsilon+288\,r^2\epsilon^2\\
&-144\,r^3-576\,\sqrt{3}\,r^2\epsilon+720\,r^2+1152\,\sqrt{3}\,r\,\epsilon-1152\,r-576\,\sqrt{3}\,\epsilon+540\Bigr]/\Bigl[72\,r^2(2\,\epsilon+1)^2(2\,\epsilon-1)^2\Bigr],\\
\varpi_{3,5}(r,\epsilon)&=\frac{48\,r^4\,\epsilon^4+256\,r^3\epsilon^4-128\,\sqrt{3}\,r^3\epsilon^3+288\,r^2\epsilon^4-192\,\sqrt{3}\,r^2\epsilon^3+72\,r^2\epsilon^2+18\,r^2+48\,\sqrt{3}\,\epsilon-45}{18\,r^2(2\,\epsilon+1)^2(2\,\epsilon-1)^2},\\
\varpi_{3,6}(r,\epsilon)&=1,
\end{align*}
}
with the corresponding intervals $\mI_1=\Bigl[1,2-4\,\,\epsilon/\sqrt{3}\Bigr)$, $\mI_2=\Bigl[ 2-4\,\,\epsilon/\sqrt{3},\sqrt{3}/(2\,\epsilon)-1 \Bigr)$, $\mI_3=\Bigl[\sqrt{3}/(2\,\epsilon)-1,3/2\Bigr)$, $\mI_4=[3/2,2)$, $\mI_5=\Bigl[ 2,\sqrt{3}/(2\,\epsilon) \Bigr)$, and $\mI_5=\Bigl[ \sqrt{3}/(2\,\epsilon),\infty \Bigr)$.

For $\epsilon \in \Bigl[ \sqrt{3}/4,\sqrt{3}/3 \Bigr)$,
$\mu_S(r,\epsilon)=\sum_{j=1}^3 \varpi_{4,j}(r,\epsilon)\,\I(r \in \mI_j)$ where
\begin{align*}
\varpi_{4,1}(r,\epsilon)&=-\frac{9\,r^2\epsilon^2+2\,\sqrt{3}\,r^2\epsilon+48\,r\,\epsilon^2+r^2-16\,\sqrt{3}\,r\,\epsilon-90\,\epsilon^2-12\,r+36\,\sqrt{3}\,\epsilon}{18\,\left( 3\,\epsilon-\sqrt{3} \right)^2},\\
\varpi_{4,2}(r,\epsilon)&=-\Bigl[ 9\,r^4\,\epsilon^4-4\,\,\sqrt{3}\,r^4\,\epsilon^3+48\,r^3\epsilon^4-48\,\sqrt{3}\,r^3\epsilon^3-90\,r^2\epsilon^4+36\,r^3\epsilon^2+96\,\sqrt{3}\,r^2\epsilon^3-126\,r^2\epsilon^2\\
&-32\,\sqrt{3}\,r\,\epsilon^3-48\,\epsilon^4+36\,\sqrt{3}\,r^2\epsilon+144\,r\,\epsilon^2+96\,\sqrt{3}\,\epsilon^3-18\,r^2-72\,\sqrt{3}\,r\,\epsilon-216\,\epsilon^2+36\,r\\
&+72\,\sqrt{3}\,\epsilon-27 \Bigr]/ \Bigl[ 2\,\left( 3\,\epsilon-\sqrt{3} \right)^4r^2 \Bigr],\\
\varpi_{4,3}(r,\epsilon)&=1,
\end{align*}
with the corresponding intervals $\mI_1=\Bigl[ 1,3-2\,\epsilon/\sqrt{3} \Bigr)$, $\mI_2=\Bigl[ 3-2\,\epsilon/\sqrt{3},\sqrt{3}/\epsilon-2 \Bigr)$, and $\mI_3=\Bigl[ \sqrt{3}/\epsilon-2,\infty \Bigr)$.

\subsection*{$\mu_A(r,\epsilon)$ Under Association Alternatives}
Under association, we compute $\mu_A(r,\epsilon)$ explicitly.
For $\epsilon \in \left[ 0,\left( 7\,\sqrt{3}-3\,\sqrt{15} \right)/12 \approx .042 \right)$,
$\mu_A(r,\epsilon)=\sum_{j=1}^6 \varpi_{1,j}(r,\epsilon)\,\I(r \in \mI_j)$ where
{\small
\begin{align*}
\varpi_{1,1}(r,\epsilon)&=-\Bigl[3456\,{\epsilon}^4r^4+9216\,{\epsilon}^4r^3-3072\,\sqrt{3}{\epsilon}^3r^4-17280\,{\epsilon}^4r^2-3072\,\sqrt{3}{\epsilon}^3r^3+2304\,\,{\epsilon}^2r^4\\
&+4608\,\sqrt{3}{\epsilon}^3r^2-2304\,\,{\epsilon}^2r^3+6336\,{\epsilon}^4+6144\,\sqrt{3}{\epsilon}^3\,r+6912\,{\epsilon}^2r^2+512\,\sqrt{3}\,\epsilon\,r^3\\
&-101\,r^4-6144\,\sqrt{3}{\epsilon}^3-11520\,{\epsilon}^2\,r-1536\,\sqrt{3}\,\epsilon\,r^2+256\,r^3+5760\,{\epsilon}^2+1536\,\sqrt{3}\,\epsilon\,r\\ &-384\,\,r^2-512\,\sqrt{3}\,\epsilon+256\,r-64\,\Bigr]/\Bigl[24\,\,\left( 6\,\epsilon+\sqrt{3} \right)^2\left( 6\,\epsilon-\sqrt{3} \right)^2r^2\Bigr],\\
\varpi_{1,2}(r,\epsilon)&=-\Bigl[1728\,{\epsilon}^4r^4-1536\,\sqrt{3}{\epsilon}^3r^4-31104\,\,{\epsilon}^4r^2+1152\,{\epsilon}^2r^4+15552\,{\epsilon}^4+10368\,{\epsilon}^2r^2-37\,r^4\\
&-20736\,{\epsilon}^2\,r+10368\,{\epsilon}^2\Bigr]/\Bigl[24\,\,\left( 6\,\epsilon+\sqrt{3} \right)^2\left( 6\,\epsilon-\sqrt{3} \right)^2r^2\Bigr],
\end{align*}
\begin{align*}
\varpi_{1,3}(r,\epsilon)&=\Bigl[-2592\,{\epsilon}^4r^4-2304\,\,\sqrt{3}{\epsilon}^3r^4-46656\,{\epsilon}^4r^2+1728\,{\epsilon}^2r^4+10656\,{\epsilon}^4-9216\,\sqrt{3}{\epsilon}^3\,r\\
&+9072\,{\epsilon}^2r^2-432\,\sqrt{3}\,\epsilon\,r^3-15\,r^4+12288\,\sqrt{3}{\epsilon}^3-13824\,\,{\epsilon}^2\,r+1728\,\sqrt{3}\,\epsilon\,r^2-216\,r^3\\
&+4032\,{\epsilon}^2-2304\,\,\sqrt{3}\,\epsilon\,r+432\,r^2+1024\,\,\sqrt{3}\,\epsilon-384\,\,r+128\Bigr]/\Bigl[36\,\left( 6\,\epsilon+\sqrt{3} \right)^2\left( 6\,\epsilon-\sqrt{3} \right)^2r^2\Bigr],\\
\varpi_{1,4}(r,\epsilon)&=-\frac{1728\,{\epsilon}^4r^4-1536\,\sqrt{3}{\epsilon}^3r^4-31104\,\,{\epsilon}^4r^2+1152\,{\epsilon}^2r^4-5184\,\,{\epsilon}^4+2592\,{\epsilon}^2r^2-37\,r^4-3456\,{\epsilon}^2}{24\,\,\left( 6\,\epsilon+\sqrt{3} \right)^2\left( 6\,\epsilon-\sqrt{3} \right)^2r^2},\\
\varpi_{1,5}(r,\epsilon)&={\frac{9}{8}}\,{\frac{1152\,{\epsilon}^4r^2+192\,{\epsilon}^4-192\,{\epsilon}^2r^2-r^4+128\,{\epsilon}^2+32\,r^2-64\,\,r+36}{\left( 6\,\epsilon+\sqrt{3} \right)^2\left( 6\,\epsilon-\sqrt{3} \right)^2r^2}},\\
\varpi_{1,6}(r,\epsilon)&=-{\frac{9}{8}}\,{\frac{(r+6)(r-2)^3}{\left( 6\,\epsilon+\sqrt{3} \right)^2\left( 6\,\epsilon-\sqrt{3} \right)^2r^2}},
\end{align*}
with the corresponding intervals $\mI_1=\left[1,\frac{1+2\,\sqrt{3}\,\epsilon}{1-\sqrt{3}\,\epsilon}\right)$, $\mI_2=\left[\frac{1+2\,\sqrt{3}\,\epsilon}{1-\sqrt{3}\,\epsilon},\frac{4\,\,\left( 1-\sqrt{3}\,\epsilon \right)}{3}\right)$, $\mI_3=\left[\frac{4\,\,\left( 1-\sqrt{3}\,\epsilon \right)}{3},\frac{4\,\,\left(1+2\,\sqrt{3}\,\epsilon \right)}{3}\right)$, $\mI_4=\left[\frac{4\,\,\left( 1+2\,\sqrt{3}\,\epsilon \right)}{3},\frac{3}{2\,\left( 1-\sqrt{3}\,\epsilon \right)}\right)$, $\mI_5=\left[\frac{3}{2\,\left( 1-\sqrt{3}\,\epsilon \right)},2\right)$ and $\mI_6=[2,\infty)$.
}

For $\epsilon \in \Bigl[ \left( 7\,\sqrt{3}-3\,\sqrt{15} \right)/12,\sqrt{3}/12 \Bigr)$,
$\mu_A(r,\epsilon)=\sum_{j=1}^6 \varpi_{2,j}(r,\epsilon)\,\I(r \in \mI_j)$ where $\varpi_{2,j}(r,\epsilon)=\varpi_{1,j}(r,\epsilon)$ for $j=1,3,4,5,6$ and
\begin{align*}
\varpi_{2,2}(r,\epsilon)&=\Bigl[-3456\,\epsilon^2r^4+111\,r^4-5184\,\,\epsilon^4r^4+4608\,\sqrt{3}\,\epsilon^3r^4-336\,\sqrt{3}\,\epsilon\,r^3-168\,r^3-13824\,\,\epsilon^4r^3\\
&+4608\,\sqrt{3}\,\epsilon^3r^3+3456\,\epsilon^2r^3+144\,r^2-6912\,\sqrt{3}\,\epsilon^3r^2-3888\,\epsilon^2r^2+576\,\sqrt{3}\,\epsilon\,r^2+25920\,\epsilon^4r^2\\
&+3168\,\epsilon^4+2880\,\epsilon^2-256\,\sqrt{3}\,\epsilon-32-3072\,\sqrt{3}\,\epsilon^3\Bigr]/\Bigl[36\,\left( \sqrt{3}+6\,\epsilon \right)^2 \left( -6\,\epsilon+\sqrt{3} \right)^2r^2 \Bigr]
\end{align*}
with the corresponding intervals $\mI_1=\Bigl[1,\frac{4\,\,\left( 1-\sqrt{3}\,\epsilon \right)}{3} \Bigr)$, $\mI_2=\Bigl[\frac{4\,\,\left( 1-\sqrt{3}\,\epsilon \right)}{3},\frac{1+2\,\sqrt{3}\,\epsilon}{1-\sqrt{3}\,\epsilon}\Bigr)$, $\mI_3=\Bigl[\frac{1+2\,\sqrt{3}\,\epsilon}{1-\sqrt{3}\,\epsilon},\frac{4\,\,(1+2\,\sqrt{3}\,\epsilon}{3} \Bigr)$, $\mI_4=\Bigl[\frac{4\,\,(1+2\,\sqrt{3}\,\epsilon}{3},\frac{3}{2\,\left( 1-\sqrt{3}\,\epsilon \right)}\Bigr)$, $\mI_5=\Bigl[\frac{3}{2\,\left( 1-\sqrt{3}\,\epsilon \right)},2\Bigr)$ and $\mI_6=[2,\infty)$.

For $\epsilon \in \Bigl[ \sqrt{3}/12,\sqrt{3}/3 \Bigr)$,
$\mu_A(r,\epsilon)=\sum_{j=1}^3 \varpi_{3,j}(r,\epsilon)\,\I(r \in \mI_j)$ where
\begin{align*}
\varpi_{3,1}(r,\epsilon)&=\frac{2\,r^2-1}{6\,r^2},\\
\varpi_{3,2}(r,\epsilon)&=\Bigl[432\,{\epsilon}^4r^4+1152\,{\epsilon}^4r^3-576\,\sqrt{3}{\epsilon}^3r^4+1296\,{\epsilon}^4r^2-960\,\sqrt{3}{\epsilon}^3r^3+864\,\,{\epsilon}^2r^4-864\,\,\sqrt{3}{\epsilon}^3r^2\\
&+576\,{\epsilon}^2r^3-192\,\sqrt{3}\,\epsilon\,r^4-360\,{\epsilon}^4+648\,{\epsilon}^2r^2+64\,\,\sqrt{3}\,\epsilon\,r^3+48\,r^4+192\,\sqrt{3}{\epsilon}^3-144\,\sqrt{3}\,\epsilon\,r^2\\
&-64\,\,r^3-504\,\,{\epsilon}^2+72\,r^2+88\,\sqrt{3}\,\epsilon-25\Bigr]/\Bigl[16\,\left( 3\,\epsilon-\sqrt{3} \right)^4r^2\Bigr],\\
\varpi_{3,3}(r,\epsilon)&=-\frac{-54\,\,{\epsilon}^2r^2+36\,\sqrt{3}\,\epsilon\,r^2+15\,{\epsilon}^2-18\,r^2+2\,\sqrt{3}\,\epsilon+20}{6\,\left( -3\,\epsilon+\sqrt{3} \right)^2r^2},
\end{align*}
with the corresponding intervals $\mI_1=\Bigl[1,\frac{1+2\,\sqrt{3}\,\epsilon}{2\,\left( 1-\sqrt{3}\,\epsilon \right)}\Bigr)$, $\mI_3=\Bigl[\frac{1+2\,\sqrt{3}\,\epsilon}{2\,\left( 1-\sqrt{3}\,\epsilon \right)},\frac{3}{2\,\left( 1-\sqrt{3}\,\epsilon \right)}\Bigr)$, $\mI_5=\Bigl[\frac{3}{2\,\left( 1-\sqrt{3}\,\epsilon \right)},\infty\Bigr)$.

\section*{Appendix 3: Derivation of $\mu(r,\epsilon)$}
We demonstrate the derivation of $\mu_S(r,\epsilon)$ for segregation with $\epsilon \in \Bigl[0,\sqrt{3}/8\Bigr)$ and among the intervals of $r$ that do not vanish as $\epsilon \rightarrow 0$.  So the resultant expressions can be used in PAE analysis.

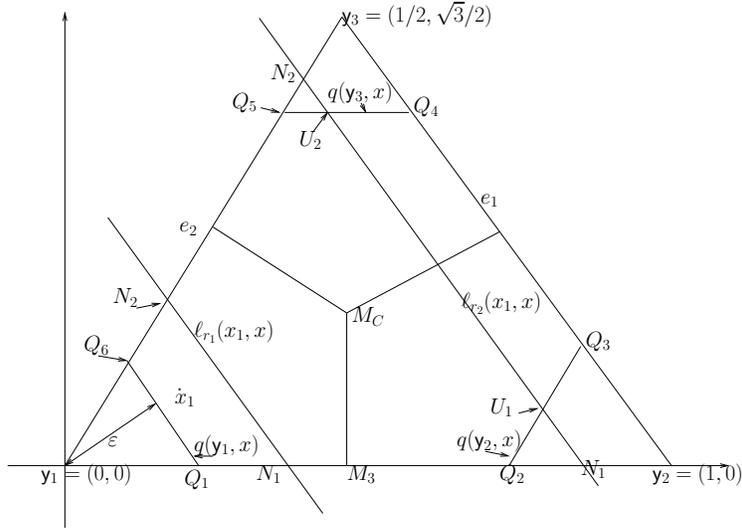
\begin{figure} [ht]
   \centering
   \scalebox{.4}{\input{ArcProbSeg1.pstex_t}}
   \caption{The support under $H^S_{\epsilon}$ for $\epsilon \in (0,\sqrt{3}/4)$ and the two types $\ell_r(x_1,x)$  for $r_1<r_2$.}
\label{fig:arc-prob-seg1}
\end{figure}

First, observe that, by symmetry,
 $$\mu_S(r,\epsilon)=P\bigl(X_2 \in \NY^r(X_1,\epsilon)\bigr)=6\,P\bigl(X_2 \in \NY^r(X_1,\epsilon), X_1 \in T_s \setminus T(\y_1,\epsilon)\bigr).$$
 Let $q(\y_j,x)$ be the line parallel to $e_j$  and crossing $T(\Y)$ such that $d(\y_j,q(\y_j,x))=\epsilon$ for $j=1,2,3$. Furthermore, let $T_{\epsilon}:=T(\Y)\setminus \cup_{j=1}^3 T(\y_j,\epsilon)$. Then $q(\y_1,x)=2\,\epsilon-\sqrt{3}\,x $, $q(\y_2,x)=\sqrt{3}\,x-\sqrt{3}+2\,\epsilon$, and $q(\y_3,x)=\sqrt{3}/2-\epsilon$.  Now, let
\begin{align*}
Q_1&=q(\y_1,x)\cap \overline{\y_1\y_2}=\left(2\,\epsilon/\sqrt{3},0 \right), &
Q_2&=q(\y_2,x)\cap \overline{\y_1\y_2}=\left(1-2\,\epsilon/\sqrt{3},0\right),\\
Q_3&=q(\y_2,x)\cap \overline{\y_2\y_3}=\left(1-\epsilon/\sqrt{3},\epsilon\right),&
Q_4&=q(\y_3,x)\cap \overline{\y_2\y_3}=\left(1/2+\epsilon/\sqrt{3},\sqrt{3}/2-\epsilon\right), \\
Q_5&=q(\y_3,x)\cap \overline{\y_1\y_3}=\left(1/2-\epsilon/\sqrt{3},\sqrt{3}/2-\epsilon\right),&
Q_6&=q(\y_1,x)\cap \overline{\y_1\y_3}=\left(\epsilon/\sqrt{3},\epsilon\right).
\end{align*}
See Figure \ref{fig:arc-prob-seg1}. Then $T(\y_1,\epsilon)=T(\y_1,Q_1,Q_6)$, $T(\y_2,\epsilon)=T(Q_2,\y_2,Q_3)$, and $T(\y_3,\epsilon)=T(Q_4,Q_5,\y_3)$, and for $\epsilon \in [0,\sqrt{3}/4)$, $T_{\epsilon}$ is the hexagon with vertices, $Q_j,\;j=1,\ldots,6$.
Now, let $q_2(x)$ be the line such that $r \, d(\y_1,q_2(x))=d(\y_1,\ell(Q_2))=d(\y_1,\overline{\y_2\y_3})-\epsilon$, $q_3(x)$ be the line such that $r \, d(\y_1,q_3(x))=d(\y_1,\overline{\y_2\y_3})$.  Then $q_2(x)=-\sqrt{3}\,x+\left( \sqrt{3}-2\,\epsilon \right)/r$ and $q_3(x)$ is the same as $\ell_s(x)$ before.  Let the $x$ coordinate of $q(\y_1,x)\cap \ell_{am}\,(x)$ be $s_1$, $q_2(x)\cap \ell_{am}\,(x)$ be $s_3$, and $\ell_s(x)\cap \ell_{am}\,(x)$ be $s_5$ and $Q_1=(s_2,0)$, $q_2(x)\cap \overline{\y_1\y_2}=(s_4,0)$, and $\ell_s(x)\cap \overline{\y_1\y_2}=(s_6,0)$. So $s_1:=\sqrt{3}\,\epsilon/2$, $s_2=2\,\epsilon/\sqrt{3}$, $s_3=\left( 3-2\,\epsilon\,\sqrt{3} \right)/(4\,r)$, $s_4=\left( 3-2\,\epsilon\,\sqrt{3} \right)/(3\,r)$, $s_5=3/(4\,r)$, and $s_6=1/r$.

\begin{figure} [ht]
   \centering
   \scalebox{.4}{\input{PartSegrgt2.pstex_t}}
   \caption{The partition of $T_s$ for different types of $\NY^r(\cdot,\epsilon)$ under $H^S_{\epsilon}$ with $r \in [2,\sqrt{3}/(2\,\epsilon))$.}
\label{fig:part-seg-r-gt-2}
\end{figure}
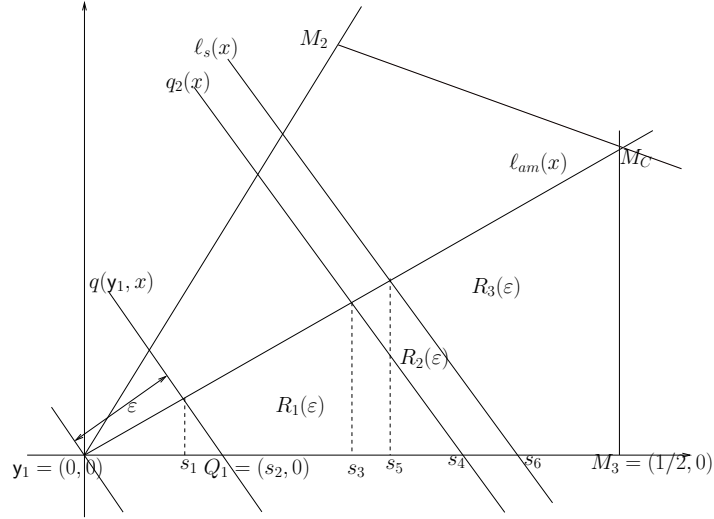

See Figure \ref{fig:part-seg-r-gt-2} for an $r \in [2,\sqrt{3}/(2\,\epsilon))$.
Furthermore, for $x_1=(x,y)\in R_{CM}\,(\y_1)$, let
\begin{align*}
U_1&:=q(\y_2,x)\cap \ell_r(x_1,x)=\left(\left(\sqrt{3}\,y/6+x/2\right)\, r+1/2-\epsilon/\sqrt{3},\left(y/2+\sqrt{3}\,x/2\right)\, r-\sqrt{3}/2+\epsilon\right) \text{, and}\\
U_2&:=q(\y_3,x)\cap \ell_r(x_1,x)=\left(\left(y/\sqrt{3}+x\right)\, r+\epsilon/\sqrt{3}-1/2,\sqrt{3}/2-\epsilon\right).
\end{align*}

Let $\msP(a_1,a_2,\ldots, a_n)$ denote the polygon with vertices $a_1,a_2,\ldots, a_n $. If $x_1$ is below $q_2(x)$, then $\NY^r(x_1,\epsilon)=A(\NY^r(x_1)) \setminus T(\y_1,\epsilon)=\msP(Q_1,N_1,N_2,Q_6)$, if $x_1$ is between $q_2(x)$ and $\ell_s(x)$, then $\NY^r(x_1,\epsilon)=\msP(Q_1,Q_2,U_1,\\
U_2,Q_5,Q_6)$, and if $x_1$ is above $\ell_s(x)$, then $\NY^r(x_1,\epsilon)=T_{\epsilon}$.

For $r \in \Bigl[1,3/2-\sqrt{3}\,\epsilon\Bigr)$, since $\epsilon$ small enough that $q_2(x)\cap T_s=\emptyset$, then $N(x,\epsilon)\subsetneq T_{\epsilon}$  for all $x \in T_s \setminus T(\y_1,\epsilon)$.  Then
\begin{eqnarray*}
\lefteqn{P(X_2 \in \NY^r(X_1,\epsilon), X_1 \in T_s \setminus T(\y_1,\epsilon))=\int_{s_1}^{s_2}\int_{q(\y_1,x)}^{\ell_{am}\,(x)} \frac{A(\NY^r(x_1,\epsilon))}{A(T_{\epsilon})^2}dydx}\\
& & + \int_{s_2}^{1/2}\int_{0}^{\ell_{am}\,(x)} \frac{A(\NY^r(x_1,\epsilon))}{A(T_{\epsilon})^2}dydx=-\frac{(576\,r^2-1152){\epsilon}^4+288\,{\epsilon}^2-37\,r^2}{1296\,(2\,\epsilon-1)^2(2\,\epsilon+1)^2}.
\end{eqnarray*}
where $A(\NY^r(x_1,\epsilon))=A(\msP(Q_1,N_1,N_2,Q_6))=\left( \frac{\sqrt{3}}{12}\,y^2+\frac{1}{2}\,x\,y+\frac{\sqrt{3}}{4}x^2 \right)\, r^2-\frac{\sqrt{3}}{3}\,{\epsilon}^2$ and $A(T_{\epsilon})=\sqrt{3}/4-\sqrt{3}\,{\epsilon}^2$ and $\ell_{am}\,(x)=x/\sqrt{3} $ is the equation of the line segment $\overline{\y_1M}_C$. Hence for $r \in [1,3/2)$, $\mu(r)=-\frac{(576\,r^2-1152){\epsilon}^4+288\,{\epsilon}^2-37\,r^2}{216\,(2\,\epsilon-1)^2(2\,\epsilon+1)^2}$.

For $r \in \Bigl[ 3/2,2-4\,\epsilon/\sqrt{3} \Bigr)$, $\ell_s(x)$ crosses through $\overline{M_3M}_C$. Since $\epsilon$ small enough so that $q_2(x)$ does the same. So $x_1$ below $q_2(x)$ is equivalent to $x_1 \in R_1(\epsilon)$ where
$$
R_1(\epsilon)=\bigl\{ (x,y)\in  [s_3,s_5]\times [q_2(x),\ell_{am}\,(x)] \cup [s_2, s_3]\times [q_2(x), \ell_{am}\,(x)] \cup [s_3, 1/2]\times [0,q_2(x)] \bigr\}.
$$
Then
\begin{eqnarray*}
\lefteqn{P(X_2 \in \NY^r(X_1,\epsilon), X_1 \in R_1(\epsilon)) = \int_0^{1/2}\int_0^{\ell_{am}\,(x)} \frac{A(\NY^r(x_1,\epsilon))}{A(T_{\epsilon})^2}dydx}\\
& = &\int_{s_1}^{s_2}\int_{q(\y_1,x)}^{\ell_{am}\,(x)} \frac{A(\msP(Q_1,N_1,N_2,Q_6))}{A(T_{\epsilon})^2}dydx+\int_{s_2}^{s_3}\int_{0}^{\ell_{am}\,(x)} \frac{A(\msP(Q_1,N_1,N_2,Q_6))}{A(T_{\epsilon})^2}dydx\\
&  & +\int_{s_3}^{1/2}\int_{0}^{q_2(x)} \frac{A(\msP(Q_1,N_1,N_2,Q_6))}{A(T_{\epsilon})^2}dydx =\Bigl[(384\,r^2+576-192\,r^4){\epsilon}^4+512\,{\epsilon}^3\sqrt{3}\,r+\\
& &(288\,r^2-1728){\epsilon}^2+\Bigl( -576\,\sqrt{3}\,r+864\,\sqrt{3}\Bigr)\,\epsilon-9\,r^4-324+288\,r\Bigr]/\Bigl[432\,((2\,\epsilon-1)^2(2\,\epsilon+1)^2r^2)\Bigr]
\end{eqnarray*}
where $A(\msP(Q_1,N_1,N_2,Q_6))$ is same as before.

Next, $x_1$ between $q_2(x)$ and $\ell_s(x)$ is equivalent to $x_1 \in R_2(\epsilon)$ where
$$
R_2(\epsilon)=\bigl\{ (x,y)\in  [s_{9},s_2]\times [q_1(x),\ell_{am}\,(x)] \cup [s_5, 1/2]\times [q_2(x),\ell_s(x)] \bigr\}.
$$
Then
\begin{eqnarray*}
&P(X_2 \in \NY^r(X_1,\epsilon), X_1 \in R_2(\epsilon)) = \int_{s_3}^{s_5}\int_{q_2(x)}^{\ell_{am}\,(x)} \frac{A(\msP(Q_1,Q_2,U_1,U_2,Q_5,Q_6))}{A(T_{\epsilon})^2}dydx\\
& +\int_{s_5}^{1/2}\int_{0}^{q_2(x)} \frac{A(\msP(Q_1,Q_2,U_1,U_2,Q_5,Q_6))}{A(T_{\epsilon})^2}dydx= -\frac{2\,\sqrt{3}\,\epsilon\,(10\,{\epsilon}^3\sqrt{3}+(32\,r-24){\epsilon}^2+\left( -27\,\sqrt{3}+12\,\sqrt{3}\,r \right)\,\epsilon-18\,r+27)}{27\,(4\,{\epsilon}^2-1)^2r^2}
\end{eqnarray*}
where
\begin{multline*}
A(\msP(Q_1,Q_2,U_1,U_2,Q_5,Q_6))=-\sqrt{3}\,{\epsilon}^2+\left( 2-2\,r\,y/\sqrt{3}-2\,r\,x \right)\,\epsilon+r\,y+\sqrt{3}\,r\,x-\sqrt{3}/2-\sqrt{3}\,r^2\,y^2/12\\
-r^2\,x\,y/2-\sqrt{3}\,r^2\,x^2/4.
\end{multline*}

Furthermore, $x_1$ above $\ell_s(x)$ is equivalent to $x_1 \in R_3(\epsilon)$ where
$
R_3(\epsilon)=\bigl\{ (x,y)\in  [s_5, 1/2] \times [\ell_{am}\,(x),\ell_s(x)] \bigr\}.
$

Then $P(X_2 \in \NY^r(X_1,\epsilon), X_1 \in R_3(\epsilon)) = \int_{s_5}^{1/2}\int_{\ell_s(x)}^{\ell_{am}\,(x)} \frac{1}{A(T_{\epsilon})}dydx = -\frac{(2\,r-3)^2}{6\,(2\,\epsilon-1)(2\,\epsilon+1)\, r^2}.$

Hence for $r \in [3/2,2)$,
\begin{eqnarray*}
\mu_S(r,\epsilon)&=&
6\,\Bigl(-\Bigl[(-384\,r^2+384+192\,r^4){\epsilon}^4+(-768\,\sqrt{3}+512\,\sqrt{3}\,r){\epsilon}^3+(1728-2304\,r+864\,r^2){\epsilon}^2\\
& &-288\,r^2-324+9\,r^4+576\,r)\Bigr]/\Bigl[432\,((2\,\epsilon-1)^2(2\,\epsilon+1)^2\,r^2)\Bigr]\Bigr)\\
&=&-\Bigl[(-384\,r^2+384+192\,r^4){\epsilon}^4+(-768\,\sqrt{3}+512\,\sqrt{3}\,r){\epsilon}^3+(1728-2304\,r+864\,r^2){\epsilon}^2\\
& &-288\,r^2-324+9\,r^4+576\,r\Bigr]/\Bigl[72\, ((2\,\epsilon+1)^2(2\,\epsilon-1)^2r^2)\Bigr].
\end{eqnarray*}

For $r \in [2,\infty)$, $\ell_s(x)$ crosses through $\overline{\y_1M}_3$, so the same types of $\NY^r(x_1,\epsilon)$ occur as above. The explicit forms of $R_j(\epsilon)$, $j=1,2,3 $ change and are given below:
\begin{align*}
R_1(\epsilon)&=\bigl\{(x,y)\in  [s_1,s_2]\times [q(\y_1,x),\ell_{am}\,(x)] \cup [s_2, s_3] \times [0, \ell_{am}\,(x)] \cup [s_3, s_4] \times [0, q_2(x)(x)]\bigr\}\\
R_2(\epsilon)&=\bigl\{(x,y)\in [s_3,s_5] \times [q_2(x),\ell_{am}\,(x)]\cup [s_5,s_4] \times [q_2(x),,\ell_s(x)] \cup [s_4,s_6] \times [0,\ell_s(x)]\bigr\}\\
R_3(\epsilon)&=\bigl\{(x,y)\in [s_5,s_6]\times [\ell_s(x),\ell_{am}\,(x)] \cup [s_6,1/2] \times [0,\ell_{am}\,(x)]\bigr\}.
\end{align*}
Then
\begin{eqnarray*}
\lefteqn{P\left( X_2 \in \NY^r(X_1,\epsilon), X_1 \in R_1(\epsilon)\right) =\int_{s_1}^{s_2}\int_{q(\y_1,x)}^{\ell_{am}\,(x)} \frac{A(\msP(Q_1,N_1,N_2,Q_6))}{A(T_{\epsilon})^2}dydx}\\
& &+\int_{s_2}^{s_3}\int_{0}^{\ell_{am}\,(x)} \frac{A(\msP(Q_1,N_1,N_2,Q_6))}{A(T_{\epsilon})^2}dydx +\int_{s_3}^{s_4}\int_{0}^{q_2(x)} \frac{A(\msP(Q_1,N_1,N_2,Q_6))}{A(T_{\epsilon})^2}dydx\\
& = & -\frac{-9+(-32\,r^2+16\,r^4+16){\epsilon}^4-48\,{\epsilon}^2+24\,\epsilon\,\sqrt{3}}{36\,(4\,{\epsilon}^2-1)^2r^2}.
\end{eqnarray*}
where $A(\msP(Q_1,N_1,N_2,Q_6))$ is same as before.

Furthermore,
\begin{eqnarray*}
\lefteqn{P\left( X_2 \in \NY^r(X_1,\epsilon), X_1 \in R_2(\epsilon) \right) = \int_{s_3}^{s_5}\int_{q_2(x)}^{\ell_{am}\,(x)} \frac{A(\msP(Q_1,Q_2,U_1,U_2,Q_5,Q_6))}{A(T_{\epsilon})^2}dydx}\\
& & +\int_{s_5}^{s_4}\int_{q_2(x)}^{\ell_s(x)} \frac{A(\msP(Q_1,Q_2,U_1,U_2,Q_5,Q_6))}{A(T_{\epsilon})^2}dydx +\int_{s_4}^{s_6}\int_{0}^{\ell_s(x)} \frac{A(\msP(Q_1,Q_2,U_1,U_2,Q_5,Q_6))}{A(T_{\epsilon})^2}dydx\\
& = &\frac{2\,\sqrt{3}\,\epsilon\,(-27\,\epsilon\,\sqrt{3}-24\,{\epsilon}^2+10\,{\epsilon}^3\sqrt{3}+27)}{81\,(4\,{\epsilon}^2-1)^2r^2}.
\end{eqnarray*}
where $A(\msP(Q_1,Q_2,U_1,U_2,Q_5,Q_6))$ is same as before.

Next,
\begin{eqnarray*}
&P\left( X_2 \in \NY^r(X_1,\epsilon), X_1 \in R_3(\epsilon) \right) = \int_{s_5}^{s_6}\int_{\ell_s(x)}^{\ell_{am}\,(x)} \frac{1}{A(T_{\epsilon})}dydx+\int_{s_6}^{1/2}\int_{0}^{\ell_{am}\,(x)} \frac{1}{A(T_{\epsilon})}dydx=\frac{3-r^2}{6\,(4\,{\epsilon}^2-1)\, r^2}.
\end{eqnarray*}
Hence for $r \in [2, \sqrt{3}/(2\,\epsilon)-1)$,
\begin{eqnarray*}
\mu_S(r,\epsilon)&=&6\,\left(-\frac{(48\,r^4-32-96\,r^2){\epsilon}^4+64\,{\epsilon}^3\sqrt{3}+(72\,r^2-144){\epsilon}^2+27-18\,r^2}{108\,(4\,{\epsilon}^2-1)^2r^2}\right)\\
&=&-\frac{(48\,r^4-32-96\,r^2){\epsilon}^4+64\,{\epsilon}^3\sqrt{3}+(72\,r^2-144){\epsilon}^2+27-18\,r^2}{18\,(4\,{\epsilon}^2-1)^2r^2}.
\end{eqnarray*}
For $r=\infty $, it is trivial to see that $\mu(r)=1$. In fact, for fixed $\epsilon>0$, $\mu(r)=1$ for $r \ge \sqrt{3}/(2\,\epsilon)$.

\section*{Appendix 4:Derivation of $\mu_S(r,\epsilon)$ and $\nu_S(r,\epsilon)$ for Segregation with $\epsilon=\sqrt{3}/8$}
For the segregation alternative with $\epsilon=\sqrt{3}/8$, $\mu_S\left( r,\epsilon=\sqrt{3}/8 \right)=1$ for $r \ge 4 $, so we find $\mu_S\left( r,\epsilon=\sqrt{3}/8\right)$ for $r \in [1,4)$.  In particular, for the mean we partition $[1,4)$ into five intervals, $[1,9/8),\,[9/8,3/2),\,[3/2,2),\,[2,3),\,[3,4)$, and for the covariance into twelve intervals, $[1,12/11),\,[12/11,9/8),\,[9/8,\sqrt{6}/11),\,[\sqrt{6}/11,21/16),\,[21/16,4/3)$,\\
$[4/3,3/2),\,[3/2,\sqrt{3})$, $[\sqrt{3},7/4),\,[7/4,2),\,[2,3),\,[3,7/2),\,[7/2,4)$. We pick the sample intervals $[3/2,2)$ and $[7/4,2)$ to demonstrate the calculations of the mean and the variance, respectively.  Then observe that, by symmetry,
 $$\mu_S(r,\epsilon)=P(X_2 \in \NY^r(X_1,\epsilon))=6\,P(X_2 \in \NY^r(X_1,\epsilon), X_1 \in R_{\epsilon}\,(\y_1)).$$
 Then $q(\y_1,x)=\sqrt{3}\,(1/4-x)$, $q(\y_2,x)=\sqrt{3}\,(4\,x-3)/4$, and $q(\y_3,x)=3\,\sqrt{3}/8$. See Section \ref{sec:derivation-par-eps} for the definition of $q(\y_j,x)$.  Hence,
$
Q_1=q(\y_1,x)\cap \overline{\y_1\y_2}=(1/4,0),\;\;\;
Q_2=q(\y_2,x)\cap \overline{\y_1\y_2}=(3/4,0), \;\;\;
Q_3=q(\y_2,x)\cap \overline{\y_2\y_3}=(7/8,3/8),\;\;\;
Q_4=q(\y_3,x)\cap \overline{\y_2\y_3}=\left(5/8,3\,\sqrt{3}/8 \right),\;\;\;
Q_5=q(\y_3,x)\cap \overline{\y_1\y_3}=\left(3/8,3\,\sqrt{3}/8 \right),\;\;\;
Q_6=q(\y_1,x)\cap \overline{\y_1\y_3}=(1/8,3/8).
$
Then $T(\y_1,\epsilon)=T(\y_1,Q_1,Q_6)$, $T(\y_2,\epsilon)=T(Q_2,y_2,Q_3)$, and $T(\y_3,\epsilon)=T(Q_4,Q_5,\y_3)$, and for $\epsilon =\sqrt{3}/8$, $T_{\epsilon}$ is the hexagon with vertices, $Q_j,\;j=1,\ldots,6$.

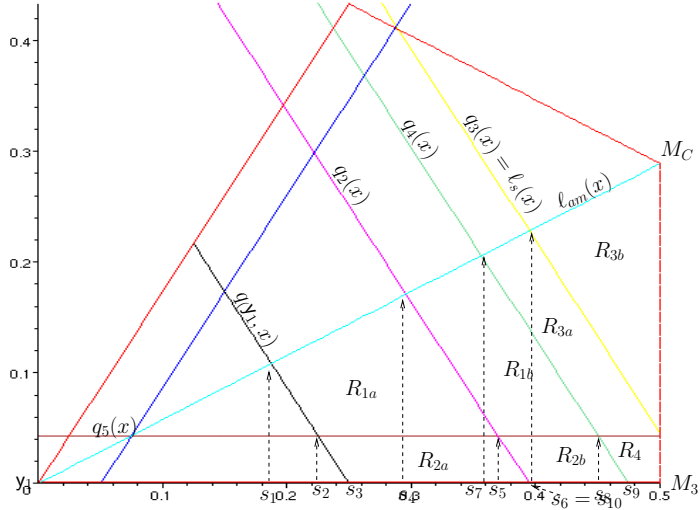
\begin{figure} [ht]
    \centering
   \scalebox{.4}{\input{HL_seg_eps1.pstex_t}}
    \caption{The regions corresponding to the seven cases for $r \in [1,4)$ with $r=1.9 $}
    \label{regions for HL_seg_eps1}
\end{figure}

Furthermore, $q_2(x)=-\sqrt{3}\,(x-3/(4\,r))$ and $q_3(x)=\sqrt{3}\,(1/r-x)$ (i.e. the same as $\ell_s(x)$ before). See Section \ref{sec:derivation-par-eps} for the definition of $q_j(x)$. Let the $x$ coordinate of $q(\y_1,x)\cap \ell_{am}\,(x)$ be $s_1$, $q_2(x)\cap \ell_{am}\,(x)$ be $s_4$, and $\ell_s(x)\cap \ell_{am}\,(x)$ be $s_{10}$ and $Q_1=(s_3,0)$, $q_2(x)\cap \overline{\y_1\y_2}=(s_6,0)$, and $\ell_s(x)\cap \overline{\y_1\y_2}=(s_{12},0)$. So $s_1=3/16$ and $s_3=1/4$, $s_4:=9/(16\,r)$, $s_6=3/(4\,r)=s_{10}$, and $s_{12}=1/r$.  See Figure \ref{regions for HL_seg_eps1}.

For $x_1=(x,y)\in R_{CM}\,(\y_1)$, let
\begin{align*}
U_1&:=q(\y_2,x)\cap \ell_r(x_1,x)=\left(\sqrt{3}\left( 4\,r\,y+4\,\sqrt{3}\,r\,x+3\,\sqrt{3}\right)/24,r\,y/2+\sqrt{3}\,r\,x/2-3\,\sqrt{3}/8\right)\text{, and}\\
U_2&:=q(\y_3,x)\cap \ell_r(x_1,x)=\left(\sqrt{3}\,\left( 8\,r\,y+8\,\sqrt{3}\,r\,x-3\,\sqrt{3}\right)/24,3\,\sqrt{3}/8 \right).
\end{align*}

If $x_1$ is below $q_2(x)$, then $\NY^r\left( x_1,\epsilon \right)=\NY^r(x_1) \setminus T(\y_1,\sqrt{3}/8)=\msP(Q_1,N_1,N_2,Q_6)$, if $x_1$ is between $q_2(x)$ and $\ell_s(x)$, then $\NY^r\left( x_1,\epsilon \right)=\msP(Q_1,Q_2,U_1,U_2,Q_5,Q_6)$, and if $x_1$ is above $\ell_s(x)$, then $\NY^r\left( x_1,\epsilon \right)=T_{\epsilon}$.

For $r \in [3/2,2)$, $q_3(x)$ crosses through $\overline{M_3M}_C$ and $q_2(x)$ crosses $\overline{\y_1M}_3$. So $x_1$ below $q_2(x)$ is equivalent to $x_1 \in R_1\left( \sqrt{3}/8 \right)$ where
$$
R_1\left( \sqrt{3}/8 \right)=\bigl\{ (x,y)\in  [s_{1},s_3]\times [q(\y_1,x),\ell_{am}\,(x)] \cup [s_{3}, s_{4}]\times [0, \ell_{am}\,(x)] \cup [s_{4},s_6]\times [0,q_2(x)] \bigr\}.
$$
Then
\begin{align*}
&P\left( X_2 \in \NY^r\left( X_1,\epsilon \right), X_1 \in R_1\left( \sqrt{3}/8 \right) \right) = \int_{T_s\setminus T(\y_1,\sqrt{3}/8)} \frac{A(\NY^r\left( x_1,\epsilon \right))}{A(T_{\epsilon})^2}dydx=\\
&\int_{s_1}^{s_{3}}\int_{q(\y_1,x)}^{\ell_{am}\,(x)} \frac{A(\msP(Q_1,N_1,N_2,Q_6))}{A(T_{\epsilon})^2}dydx +\int_{s_{3}}^{s_{4}}\int_{0}^{\ell_{am}\,(x)} \frac{A(\msP(Q_1,N_1,N_2,Q_6))}{A(T_{\epsilon})^2}dydx \\
&+\int_{s_{4}}^{s_6}\int_{0}^{q_2(x)} \frac{A(\msP(Q_1,N_1,N_2,Q_6))}{A(T_{\epsilon})^2}dydx = -\frac{r^4-2\,r^2-63}{676\,r^2}.
\end{align*}

where $A(\msP(Q_1,N_1,N_2,Q_6))=\frac{\sqrt{3}}{576}\,\left( 12\,r\,x+4\,\sqrt{3}\,r\,y+3 \right) \left( 12\,r\,x+4\,\sqrt{3}\,r\,y-3 \right)$ is same as before.

Next, $x_1$ between $q_2(x)$ and $\ell_s(x)$ is equivalent to $x_1 \in R_2\left( \sqrt{3}/8 \right)$ where
$$
R_2\left( \sqrt{3}/8 \right)=\bigl\{ (x,y)\in  [s_{4},s_6]\times [q_2(x),\ell_{am}\,(x)] \cup [s_{6},1/2]\times [0, \ell_s(x)] \bigr\}.
$$
 Then
\begin{eqnarray*}
\lefteqn{P\left( X_2 \in \NY^r\left( X_1,\epsilon \right), X_1 \in R_2\left( \sqrt{3}/8 \right) \right)= \int_{s_{4}}^{s_{6}}\int_{q_2(x)}^{\ell_{am}\,(x)} \frac{A(\msP(Q_1,Q_2,U_1,U_2,Q_5,Q_6))}{A(T_{\epsilon})^2}dydx}\\
& &+\int_{s_6}^{1/2}\int_{0}^{\ell_s(x)} \frac{A(\msP(Q_1,Q_2,U_1,U_2,Q_5,Q_6))}{A(T_{\epsilon})^2}dydx = \frac{64\,r^4-768\,r^3+1824\,r^2-128\,r-1467}{2028\,r^2}.
\end{eqnarray*}
where $A(\msP(Q_1,Q_2,U_1,U_2,Q_5,Q_6))=\left( -\frac{\sqrt{3}}{12}y^2-\frac{1}{2}\,x\,y-\frac{\sqrt{3}}{4}x^2\right)\, r^2+\left( \frac{3}{4}\,y+\frac{3\,\sqrt{3}}{4}\,x \right)\, r-{\frac{19}{64}}\,\sqrt{3}$.

Furthermore, $x_1$ above $q_3(x)$ is $x_1 \in R_3\left( \sqrt{3}/8 \right)$ where $R_3\left( \sqrt{3}/8 \right)=\bigl\{ (x,y)\in  [s_{10},1/2]\times [\ell_s(x),\ell_{am}\,(x)]\bigr\}.$

Then
$$
P\left( X_2 \in \NY^r\left( X_1,\epsilon \right), X_1 \in R_3\left( \sqrt{3}/8 \right) \right) = \int_{s_{10}}^{1/2}\int_{q_3(x)}^{\ell_{am}\,(x)} \frac{1}{A(T_{\epsilon})}dydx = \frac{8\,(4\,r^2-12\,r+9)}{39\,r^2}.
$$
Hence for $r \in [3/2,2)$,
{\small
\begin{eqnarray*}
\mu_S\left( r,\sqrt{3}/8 \right)=6\,\left(\frac{61\,r^4-768\,r^3+3494\,r^2-5120\,r+2466}{2028\,r^2}\right)=\frac{61\,r^4-768\,r^3+3494\,r^2-5120\,r+2466}{338\,r^2}.
\end{eqnarray*}
}

For $r \ge 4 $, it is trivial to see that $\mu(r)=1$.

To find the covariance, we need to find the possible types of $\G_1\left( x_1,\NY^r,\epsilon \right)$ and $\NY^r\left( x_1,\epsilon \right)$ for $r \in [1,4)$.  The intersection points of $\xi_j(r,x)$ with $\partial(T(\Y))$ and $\partial(R(\y_j))$ for $j=1,2,3$, i.e. $G_1-G_6$ and $L_1-L_6$ are same as before. Recall also $M_C,\,M_1,M_2,M_3$ and $\y_1,\y_2,\y_3$.  Then $\G_1\left( x_1,\NY^r,\epsilon \right)$ is a polygon whose vertices are a subset of the above points.

There are six cases regarding $\G_1\left( x_1,\NY^r,\epsilon \right)$ and one case for $\NY^r\left( x_1,\epsilon \right)$.  Each case $j$, corresponds to the region $R_j\left( \sqrt{3}/8 \right)$ in Figure \ref{regions for HL_seg_eps1} where $q(\y_1,x)$, $q_2(x)$, $q_3(x)$, $s_j$ for $j=1,3,4,6,10$ are same as before and
 $q_4(x)=-\sqrt{3}\,(4\,x-r)/4,\;q_5(x)=\sqrt{3}\,(2-r)/4$ and $s_2=(r-1)/4, \; s_5=(r^2-2\,r+3)/(4\,r),\; s_7=3\,r/16, \; s_8=(r-1)/2$. (see Figure \ref{regions for HL_seg_eps1}).

Then, for\\
$
x_1=(x,y)\in R_1\left( \sqrt{3}/8 \right),  \G_1\left( x_1,\NY^r \right)=\msP(Q_1,G_2,G_3,G_4,G_5,Q_6)\\
x_1\in R_2\left( \sqrt{3}/8 \right),  \G_1\left( x_1,\NY^r \right)=\msP(Q_1,G_2,G_3,M_2,L_4,L_5,M_3,Q_6)\\
x_1\in R_3\left( \sqrt{3}/8 \right),  \G_1\left( x_1,\NY^r \right)=\msP(G_1,G_2,G_3,G_4,G_5,G_6)\\
x_1\in R_4\left( \sqrt{3}/8 \right),  \G_1\left( x_1,\NY^r \right)=\msP(G_1,G_2,G_3,G_4.G_5,G_6)
$

The explicit forms of $R_j\left( \sqrt{3}/8 \right)$, $j=1,\ldots,4$ are as follows:
\begin{align*}
R_1\left( \sqrt{3}/8 \right)&=\bigl\{(x,y)\in  [s_1,s_2]\times [q(\y_1,x),\ell_{am}\,(x)] \cup [s_2, s_7]\times [q_5(x), \ell_{am}\,(x)] \cup [s_7, s_8]\times [q_5(x), q_4(x)]\bigr\}\\
R_2\left( \sqrt{3}/8 \right)&=\bigl\{(x,y)\in [s_2,s_3] \times [q(\y_1,x),q_5(x)]\cup [s_3,s_8] \times [0,q_5(x)] \cup [s_8,s_9] \times [0,q_4(x)]\bigr\}\\
R_3\left( \sqrt{3}/8 \right)&=\bigl\{(x,y)\in [s_7,s_8]\times [q_3(x),\ell_{am}\,(x)] \cup [s_8,1/2] \times [q_5(x),\ell_{am}\,(x)]\bigr\}\\
R_4\left( \sqrt{3}/8 \right)&=\bigl\{(x,y)\in [s_8,s_9] \times [q_4(x),q_5(x)] \cup [s_9,1/2] \times [0,q_5(x)]\bigr\}.
\end{align*}
Let $P^r_{2N}(\epsilon):=P\left( \{X_2,X_3\} \subset \NY^r(X_1,\epsilon) \right)$, $P^r_{2G}(\epsilon):=P\left( \{X_2,X_3\} \subset \G_1^r(X_1,\epsilon) \right)$, and $P^r_M(\epsilon):=P\bigl( X_2 \in \NY^r(X_1,\epsilon),X_3 \in \G_1^r(X_1,\epsilon) \bigr)$.

Now, by symmetry, $P^r_{2N}\left( \sqrt{3}/8 \right) =6\,P\left( \{X_2,X_3\} \subset \NY^r\left( X_1,\epsilon \right),\; X_1 \in T(\y_1,\sqrt{3}/8) \right).$

For $r \in [7/4,2)$,
\begin{eqnarray*}
\lefteqn{P\left( \{X_2,X_3\} \subset \NY^r\left( X_1,\epsilon \right),\; X_1 \in T(\y_1,\sqrt{3}/8)\right) = \int_0^{1/2}\int_0^{\ell_{am}\,(x)} \frac{A(\NY^r\left( x_1,\epsilon \right))^2}{A(T_{\epsilon})^3}dydx}\\
& = &-\frac{261\,r^6-4608\,r^5+29105\,r^4-72960\,r^3+32575\,r^2+78848\,r-67620}{65910\,r^2}.
\end{eqnarray*}
where $\NY^r\left( x_1,\epsilon \right)=\msP(Q_1,N_1,N_2,Q_6)$ for $x_1 \in R_{1a}\left( \sqrt{3}/8 \right) \cup R_{2a}\left( \sqrt{3}/8 \right)$, $\NY^r\left( x_1,\epsilon \right)=\msP(Q_1,Q_2,U_1,U_2,Q_5,Q_6)$ for $x_1 \in R_{1b}\left( \sqrt{3}/8 \right) \cup R_{2b}\left( \sqrt{3}/8 \right) \cup R_{3a}\left( \sqrt{3}/8 \right) \cup R_{4}$, and $\NY^r\left( x_1,\epsilon \right)=T_{\epsilon}$ for $x_1 \in R_{3b}\left( \sqrt{3}/8 \right)$, all of whose areas are given above. Hence for $r \in [7/4,2)$,
 $$P(\{X_2,X_3\} \subset \NY^r\left( X_1,\epsilon \right))=-\frac{261\,r^6-4608\,r^5+29105\,r^4-72960\,r^3+32575\,r^2+78848\,r-67620}{10985\,r^2}.$$

Next, by symmetry, $P^r_{2G}\left( \sqrt{3}/8 \right)=6\,P\left( \{X_2,X_3\} \subset \G_1^r\left( X_1,\epsilon \right),\; X_1 \in T_s\setminus T(\y_1,\sqrt{3}/8) \right),$
and
 $$P\left( \{X_2,X_3\} \subset \G_1^r\left( X_1,\epsilon \right),\; X_1 \in T_s\setminus T(\y_1,\sqrt{3}/8) \right)=\sum_{j=1}^4 P\left( \{X_2,X_3\} \subset \G_1^r\left( X_1,\epsilon \right),\; X_1 \in R_j\left( \sqrt{3}/8 \right) \right).$$

For $x_1 \in R_1\left( \sqrt{3}/8 \right)$,
\begin{eqnarray*}
\lefteqn{P\left( \{X_2,X_3\} \subset \G_1^r\left( X_1,\epsilon \right),\; X_1 \in R_1\left( \sqrt{3}/8 \right) \right) = \int_{s_1}^{s_2}\int_{q(\y_1,x)}^{\ell_{am}\,(x)} \frac{A\left( \G_1\left( x_1,\NY^r,\epsilon \right) \right)^2}{A(T_{\epsilon})^3}dydx}\\
& &+\int_{s_2}^{s_7}\int_{q_5(x)}^{\ell_{am}\,(x)} \frac{A\left( \G_1\left( x_1,\NY^r,\epsilon \right) \right)^2}{A(T_{\epsilon})^3}dydx +\int_{s_7}^{s_8}\int_{q_5(x)}^{q_4(x)} \frac{A\left( \G_1\left( x_1,\NY^r,\epsilon \right) \right)^2}{A(T_{\epsilon})^3}dydx\\
& = &\frac{18894\,r^6-12248\,r^5-131375\,r^4+45360\,r^3+584030\,r^2-841816\,r+337155}{65910\,r^4}.
\end{eqnarray*}
where $A\left( \G_1\left( x_1,\NY^r,\epsilon \right) \right)=\frac{\sqrt{3}\,\left( 45\,r^2+32\,y\,\sqrt{3}\,x+96\,x-48\,x^2-80\,y^2+32\,\sqrt{3}\,y-96 \right)}{192\,r^2}$.

For $x_1 \in R_2\left( \sqrt{3}/8 \right)$,
\begin{eqnarray*}
\lefteqn{P\left( \{X_2,X_3\} \subset \G_1^r\left( X_1,\epsilon \right),\; X_1 \in R_2\left( \sqrt{3}/8 \right) \right) = \int_{s_2}^{s_3}\int_{q(\y_1,x)}^{q_5(x)} \frac{A\left( \G_1\left( x_1,\NY^r,\epsilon \right) \right)^2}{A(T_{\epsilon})^3}dydx}\\
& &+\int_{s_3}^{s_8}\int_{0}^{q_5(x)} \frac{A\left( \G_1\left( x_1,\NY^r,\epsilon \right) \right)^2}{A(T_{\epsilon})^3}dydx +\int_{s_8}^{s_9}\int_{0}^{q_4(x)} \frac{A\left( \G_1\left( x_1,\NY^r,\epsilon \right) \right)^2}{A(T_{\epsilon})^3}dydx\\
& = &-\frac{4\,(4865\,r^6-24063\,r^5+49460\,r^4-73210\,r^3+98045\,r^2-84107\,r+29010)}{32955\,r^4}.
\end{eqnarray*}
where $A\left( \G_1\left( x_1,\NY^r,\epsilon \right) \right)=\frac{\sqrt{3}\,\left( 96+32\,y\,\sqrt{3}\,x+128\,\sqrt{3}\,r\,y-224\,\sqrt{3}\,y-192\,r+96\,x+93\,r^2-48\,x^2+176\,y^2 \right)}{192\,r^2}$.

For $x_1 \in R_3\left( \sqrt{3}/8 \right)$,
\begin{eqnarray*}
\lefteqn{P\left( \{X_2,X_3\} \subset \G_1^r\left( X_1,\epsilon \right),\; X_1 \in R_3\left( \sqrt{3}/8 \right) \right)}\\
& = &\int_{s_7}^{s_8}\int_{q_4(x)}^{\ell_{am}\,(x)} \frac{A\left( \G_1\left( x_1,\NY^r,\epsilon \right) \right)^2}{A(T_{\epsilon})^3}dydx+\int_{s_8}^{1/2}\int_{q_5(x)}^{\ell_{am}\,(x)} \frac{A\left( \G_1\left( x_1,\NY^r,\epsilon \right) \right)^2}{A(T_{\epsilon})^3}dydx\\
& = &-\frac{46293\,r^6-100944\,r^5-254880\,r^4+506880\,r^3+829440\,r^2-2064384\,r+1048576}{98865\,r^4}.
\end{eqnarray*}
where $A\left( \G_1\left( x_1,\NY^r,\epsilon \right) \right)=\frac{\sqrt{3}\,\left( 3\,r^2+6\,x-6\,x^2-6\,y^2+2\,\sqrt{3}\,y-6 \right)}{12\,r^2}$.

For $x_1 \in R_4\left( \sqrt{3}/8 \right)$,
\begin{eqnarray*}
\lefteqn{P\left( \{X_2,X_3\} \subset \G_1^r\left( X_1,\epsilon \right),\; X_1 \in R_4\left( \sqrt{3}/8 \right) \right)}\\
& = &\int_{s_2}^{s_4}\int_{q_2(x)}^{\ell_{am}\,(x)} \frac{A\left( \G_1\left( x_1,\NY^r,\epsilon \right) \right)^2}{A(T_{\epsilon})^3}dydx+\int_{s_4}^{s_5}\int_{q_3(x)}^{\ell_{am}\,(x)} \frac{A\left( \G_1\left( x_1,\NY^r,\epsilon \right) \right)^2}{A(T_{\epsilon})^3}dydx\\
& = &\frac{8\,(3577\,r^6-20548\,r^5+45620\,r^4-61760\,r^3+79040\,r^2-77824\,r+32256)}{32955\,r^4}.
\end{eqnarray*}
where $A\left( \G_1\left( x_1,\NY^r,\epsilon \right) \right)=\frac{\sqrt{3}\,(-6\,r+3\,x-3\,x^2+5\,y^2+3\,r^2+3+4\,\sqrt{3}\,r\,y-7\,\sqrt{3}\,y)}{6\,r^2}$.

So
\begin{align*}
&P^r_{2G}\left( \sqrt{3}/8 \right)=6\,\left(\frac{19032\,r^6-243648\,r^5+1118355\,r^4-2085120\,r^3+1534050\,r^2-113664\,r-233639}{197730\,r^4}\right)\\
&=\frac{19032\,r^6-243648\,r^5+1118355\,r^4-2085120\,r^3+1534050\,r^2-113664\,r-233639}{32955\,r^4}.
\end{align*}
Furthermore, by symmetry,
 $$P^r_M\left( \sqrt{3}/8 \right)=6\,P\left( X_2 \in \NY^r\left( X_1,\epsilon \right),\;X_3\in \G_1^r\left( X_1,\epsilon \right),\; X_1 \in T_s\setminus T(\y_1,\sqrt{3}/8) \right),$$
and
\begin{multline*}
P\left( X_2 \in \NY^r\left( X_1,\epsilon \right),\;X_3\in \G_1^r\left( X_1,\epsilon \right),\; X_1 \in T_s\setminus T(\y_1,\sqrt{3}/8) \right)\\
=\sum_{j} P\left( X_2 \in \NY^r(X_1),\;X_3\in \G_1\left( X_1,\NY^r \right),\; X_1 \in R_j\left( \sqrt{3}/8 \right) \right)
\end{multline*}
where $j \in \{1a,\,1b,\,2a,\,2b,3a,\,3b,4\}$. The explicit forms of these regions are
{\small
\begin{align*}
R_{1a}\left( \sqrt{3}/8 \right)&=\bigl\{(x,y)\in  [s_1,s_2]\times [q(\y_1,x),\ell_{am}\,(x)] \cup [s_2, s_4]\times [q_5(x), \ell_{am}\,(x)] \cup [s_4, s_5]\times [q_5(x), q_2(x)]\bigr\},\\
R_{1b}\left( \sqrt{3}/8 \right)&=\bigl\{(x,y)\in [s_4,s_7] \times [q_2(x),\ell_{am}\,(x)]\cup [s_7,s_5] \times [q_2(x),q_4(x)] \cup [s_5,s_8] \times [q_5(x),q_4(x)]\bigr\},\\
R_{2a}\left( \sqrt{3}/8 \right)&=\bigl\{(x,y)\in [s_2,s_3]\times [q(\y_1,x),q_5(x)] \cup [s_3,s_5] \times [0,q_5(x)] \cup [s_5,s_6] \times [0,q_2(x)] \bigr\},\\
R_{2b}\left( \sqrt{3}/8 \right)&=\bigl\{(x,y)\in [s_5,s_6] \times [q_2(x),q_5(x)] \cup [s_6,s_8] \times [0,q_5(x)] \cup [s_8,s_9] \times [0,q_4(x)]\bigr\},\\
R_{3a}\left( \sqrt{3}/8 \right)&=\bigl\{(x,y)\in [s_7,s_{10}] \times [q_4(x),\ell_{am}\,(x)] \cup [s_{10},s_8] \times [q_4(x),q_3(x)] \cup [s_8,1/2] \times [q_5(x),q_4(x)]\bigr\},\\
R_{3b}\left( \sqrt{3}/8 \right)&=\bigl\{(x,y)\in [s_{10},1/2] \times [q_3(x),\ell_{am}\,(x)]\bigr\}.
\end{align*}
}
 $R_4\left( \sqrt{3}/8 \right)$ is the same as before.

For $x_1 \in R_{1a}\left( \sqrt{3}/8 \right)$,
{\small
\begin{align*}
&P\left( X_2 \in \NY^r\left( X_1,\epsilon \right),\;X_3\in \G_1^r\left( X_1,\epsilon \right),\; X_1 \in R_{1a}\left( \sqrt{3}/8 \right)\right) = \int_{s_1}^{s_2}\int_{q(\y_1,x)}^{\ell_{am}\,(x)} \frac{A(\NY^r\left( x_1,\epsilon \right))\,A\left( \G_1\left( x_1,\NY^r,\epsilon \right) \right)}{A(T_{\epsilon})^3}dydx\\
&+\int_{s_2}^{s_4}\int_{q_5(x)}^{\ell_{am}\,(x)} \frac{A(\NY^r\left( x_1,\epsilon \right))\,A\left( \G_1\left( x_1,\NY^r,\epsilon \right) \right)}{A(T_{\epsilon})^3}dydx+\int_{s_4}^{s_5}\int_{q_5(x)}^{q_2(x)} \frac{A(\NY^r\left( x_1,\epsilon \right))\,A\left( \G_1\left( x_1,\NY^r,\epsilon \right) \right)}{A(T_{\epsilon})^3}dydx\\
& =-\Bigl[2960\,r^9-3045\,r^8-21504\,r^7-28554\,r^6+101040\,r^5+205785\,r^4-550080\,r^3+391392\,r^2\\
& -114048\,r+12150\Bigr]/\Bigl[395460\,r^6\Bigr].
\end{align*}
}
For $x_1 \in R_{1b}\left( \sqrt{3}/8 \right)$,
{\small
\begin{align*}
&P\left(X_2 \in \NY^r\left( X_1,\epsilon \right),\;X_3\in \G_1^r\left( X_1,\epsilon \right),\; X_1 \in R_{1b}\left( \sqrt{3}/8 \right)\right) = \int_{s_4}^{s_7}\int_{q_2(x)}^{\ell_{am}\,(x)} \frac{A(\NY^r\left( x_1,\epsilon \right))\,A\left( \G_1\left( x_1,\NY^r,\epsilon \right) \right)}{A(T_{\epsilon})^3}dydx\\
& +\int_{s_7}^{s_5}\int_{q_2(x)}^{q_4(x)} \frac{A(\NY^r\left( x_1,\epsilon \right))\,A\left( \G_1\left( x_1,\NY^r,\epsilon \right) \right)}{A(T_{\epsilon})^3}dydx+\int_{s_5}^{s_8}\int_{q_5(x)}^{q_4(x)} \frac{A(\NY^r\left( x_1,\epsilon \right))\,A\left( \G_1\left( x_1,\NY^r,\epsilon \right) \right)}{A(T_{\epsilon})^3}dydx\\
& = -\Bigl[(r^2-3)(3111\,r^{10}+96\,r^9-67659\,r^8+26528\,r^7+341223\,r^6-352896\,r^5-177291\,r^4+229632\,r^3\\
&+23922\,r^2-22464\,r+3078)\Bigr]/\Bigl[395460\, r^6\Bigr].
\end{align*}
}
For $x_1 \in R_{2a}\left( \sqrt{3}/8 \right)$,
{\small
\begin{align*}
&P\left( X_2 \in \NY^r\left( X_1,\epsilon \right),\;X_3\in \G_1^r\left( X_1,\epsilon \right),\; X_1 \in R_{2a}\left( \sqrt{3}/8 \right) \right) = \int_{s_2}^{s_3}\int_{q(\y_1,x)}^{q_5(x)} \frac{A(\NY^r\left( x_1,\epsilon \right))\,A\left( \G_1\left( x_1,\NY^r,\epsilon \right) \right)}{A(T_{\epsilon})^3}dydx\\
& +\int_{s_3}^{s_5}\int_{0}^{q_5(x)} \frac{A(\NY^r\left( x_1,\epsilon \right))\,A\left( \G_1\left( x_1,\NY^r,\epsilon \right) \right)}{A(T_{\epsilon})^3}dydx+\int_{s_5}^{s_6}\int_{0}^{q_2(x)} \frac{A(\NY^r\left( x_1,\epsilon \right))\,A\left( \G_1\left( x_1,\NY^r,\epsilon \right) \right)}{A(T_{\epsilon})^3}dydx\\
& = \frac{2\,(r-2)(r-3)(530\,r^6+1225\,r^5+1697\,r^4-3399\,r^3-3027\,r^2+5274\,r-1188)}{98865\,r^5}.
\end{align*}
}
For $x_1 \in R_{2b}\left( \sqrt{3}/8 \right)$,
{\small
\begin{align*}
&P\left( X_2 \in \NY^r\left( X_1,\epsilon \right),\;X_3\in \G_1^r\left( X_1,\epsilon \right),\; X_1 \in R_{2b}\left( \sqrt{3}/8 \right) \right) = \int_{s_5}^{s_6}\int_{q_2(x)}^{q_5(x)} \frac{A(\NY^r\left( x_1,\epsilon \right))\,A\left( \G_1\left( x_1,\NY^r,\epsilon \right) \right)}{A(T_{\epsilon})^3}dydx\\
& +\int_{s_6}^{s_8}\int_{0}^{q_5(x)} \frac{A(\NY^r\left( x_1,\epsilon \right))\,A\left( \G_1\left( x_1,\NY^r,\epsilon \right) \right)}{A(T_{\epsilon})^3}dydx+\int_{s_8}^{s_9}\int_{0}^{q_4(x)} \frac{A(\NY^r\left( x_1,\epsilon \right))\,A\left( \G_1\left( x_1,\NY^r,\epsilon \right) \right)}{A(T_{\epsilon})^3}dydx\\
& = \frac{2\,(r-2)(r^2-3)(467\,r^8-50\,r^7-7789\,r^6+1290\,r^5+16083\,r^4-4110\,r^3-7311\,r^2-810\,r+702)}{98865\,r^5}.
\end{align*}
}
For $x_1 \in R_{3a}\left( \sqrt{3}/8 \right)$,
{\small
\begin{align*}
&P\left( X_2 \in \NY^r\left( X_1,\epsilon \right),\;X_3\in \G_1^r\left( X_1,\epsilon \right),\; X_1 \in R_{3a}\left( \sqrt{3}/8 \right) \right) = \int_{s_7}^{s_{10}}\int_{q_4(x)}^{\ell_{am}\,(x)} \frac{A(\NY^r\left( x_1,\epsilon \right))\,A\left( \G_1\left( x_1,\NY^r,\epsilon \right) \right)}{A(T_{\epsilon})^3}dydx\\
& +\int_{s_{10}}^{s_8}\int_{q_4(x)}^{q_3(x)} \frac{A(\NY^r\left( x_1,\epsilon \right))\,A\left( \G_1\left( x_1,\NY^r,\epsilon \right) \right)}{A(T_{\epsilon})^3}dydx+\int_{s_8}^{1/2}\int_{q_5(x)}^{q_3(x)} \frac{A(\NY^r\left( x_1,\epsilon \right))\,A\left( \G_1\left( x_1,\NY^r,\epsilon \right) \right)}{A(T_{\epsilon})^3}dydx\\
& = \Bigl[3934\,r^{12}-11040\,r^{11}-2352\,r^{10}-283680\,r^9+1239855\,r^8-751008\,r^7-3225344\,r^6\\
& +6125568\,r^5-3847680\,r^4+81920\,r^3+1843200\,r^2-2045952\,r+905472\Bigr]/\Bigl[395460\,r^6\Bigr].
\end{align*}
}
For $x_1 \in R_{3b}\left( \sqrt{3}/8 \right)$,
\begin{align*}
&P\left( X_2 \in \NY^r\left( X_1,\epsilon \right),\;X_3\in \G_1^r\left( X_1,\epsilon \right),\; X_1 \in R_{3b}\left( \sqrt{3}/8 \right) \right)=\\
& \int_{s_{10}}^{1/2}\int_{q_3(x)}^{\ell_{am}\,(x)} \frac{A(\NY^r\left( x_1,\epsilon \right))\,A\left( \G_1\left( x_1,\NY^r,\epsilon \right) \right)}{A(T_{\epsilon})^3}dydx = \frac{64\,(2\,r^4-4\,r^2+4\,r-3)(2\,r-3)^2}{507\,r^6}.
\end{align*}

For $x_1 \in R_4\left( \sqrt{3}/8 \right)$,
{\small
\begin{eqnarray*}
\lefteqn{P\left( X_2 \in \NY^r\left( X_1,\epsilon \right),\;X_3\in \G_1^r\left( X_1,\epsilon \right),\; X_1 \in R_4\left( \sqrt{3}/8 \right) \right)}\\
& = &\int_{s_8}^{s_9}\int_{q_4(x)}^{q_5(x)} \frac{A(\NY^r\left( x_1,\epsilon \right))\,A\left( \G_1\left( x_1,\NY^r,\epsilon \right) \right)}{A(T_{\epsilon})^3}dydx+\int_{s_9}^{1/2}\int_{0}^{q_5(x)} \frac{A(\NY^r\left( x_1,\epsilon \right))\,A\left( \G_1\left( x_1,\NY^r,\epsilon \right) \right)}{A(T_{\epsilon})^3}dydx\\
& = &-\frac{8\,(48\,r^6-22\,r^5-73\,r^4-3114\,r^3+4391\,r^2+1574\,r-2850)(r-2)^2}{32955\,r^2}.
\end{eqnarray*}
}
Then
\begin{eqnarray*}
P^r_M(\epsilon)&=&6\,\Bigl(-\Bigl[49\,r^{12}-1536\,r^{11}+17952\,r^{10}-129280\,r^9+609420\,r^8-1728768\,r^7+2757670\,r^6\\
& &-3013632\,r^5+3418140\,r^4-3829760\,r^3+3026880\,r^2-1571616\,r+445284\Bigr]/\Bigl[395460\,r^6\Bigr]\Bigr)\\
&=&-\Bigl[49\,r^{12}-1536\,r^{11}+17952\,r^{10}-129280\,r^9+609420\,r^8-1728768\,r^7+2757670\,r^6\\
& &-3013632\,r^5+3418140\,r^4-3829760\,r^3+3026880\,r^2-1571616\,r+445284\Bigr]/\Bigl[65910\,r^6\Bigr].
\end{eqnarray*}
So
\begin{multline*}
\E^S_{\sqrt{3}/8}[h_{12}\,h_{13}]=-\Bigl[49\,r^{12}-1536\,r^{11}+18735\,r^{10}-143104\,r^9+677703\,r^8-1704000\,r^7+1737040\,r^6\\
-691968\,r^5+1681230\,r^4-3716096\,r^3+3260519\,r^2-1571616\,r+445284\Bigr]/\Bigl[3295\,r^6\Bigr].
\end{multline*}
Hence
\begin{multline*}
\nu_S\left( r,\sqrt{3}/8 \right)=\\
-\Bigl[637\,r^{12}-19968\,r^{11}+299370\,r^{10}-3265792\,r^9+24051519\,r^8-112023360\,r^7+328179640\,r^6\\
-602490624\,r^5+673558110\,r^4-427086848\,r^3+133604087\,r^2-20431008\,r+5788692\Bigr]/\Bigl[428415\,r^6\Bigr].
\end{multline*}

Derivation of $\mu_S(r,\epsilon)$ and $\nu_S(r,\epsilon)$ for segregation with $\epsilon=\sqrt{3}/4$ and with $\epsilon=2\,\sqrt{3}/7$ are similar.

\section*{Appendix 5: Proof of Corollary 1}
In the multiple triangle case,
\begin{multline*}
\mu(r,J)=\E[\rho_n(r,J)]=\frac{1}{n\,(n-1)}\sum\hspace*{-0.1 in}\sum_{i < j \hspace*{0.25 in}}   \hspace*{-0.1 in} \,\E[h_{ij}]=\frac{1}{2}\E[h_{12}]=\E[\I\left( A_{12} \right)] =P\bigl( X_2 \in \NY^r(X_1) \bigr).
\end{multline*}
By definition of $N_{\Y}^r(\cdot)$, $P\bigl( X_2 \in \NY^r(X_1) \bigr)=0$ if $X_1$ and $X_2$ are in different triangles. So by the law of total probability
\begin{eqnarray*}
\mu(r,J)&:=&P\bigl( X_2 \in \NY^r(X_1) \bigr)= \sum_{j=1}^{J}P\bigl( X_2 \in \NY^r(X_1)\,|\,\{X_1,X_2\} \subset T_j \bigr)\,P\bigl( \{X_1,X_2\} \subset T_j \bigr)\\
&=& \sum_{j=1}^{J}\mu(r)\,P\bigl( \{X_1,X_2\} \subset T_j \bigr) \;\;\;\text{ (since $P\bigl( X_2 \in \NY^r(X_1)\,|\,\{X_1,X_2\} \subset T_j \bigr)=\mu(r)$)}\\
&=& \mu(r) \, \sum_{j=1}^{J}\Bigl[A(T_j) / A(C_H(\Y))\Bigr]^2\;\;\; \text{ (since $P\bigl( \{X_1,X_2\} \subset T_j \bigr)=\bigl( A(T_j) / A(C_H(\Y)) \bigr)^2$)}
\end{eqnarray*}
Then $\mu(r,J)=\mu(r)\cdot \left(\sum_{j=1}^{J}w_j^2\right)$ where $\mu(r)$ is given in Equation (\ref{eq:Asymean}).

Furthermore, the asymptotic variance is
\begin{eqnarray*}
\nu(r,J)&=&\E[h_{12}\,h_{13}]-\E[h_{12}]\,\E[h_{13}]\\
& = & P\bigl( \{X_2,X_3\} \subset \NY^r(X_1) \bigr)+2\,P\bigl( X_2 \in \NY^r(X_1), X_3 \in \G_1(X_1,\NY^r) \bigr)\\
& &+P\bigl(\{X_2,X_3\} \subset \G_1(X_1,\NY^r)\bigr)-4\,\,(\mu(r,J))^2.
\end{eqnarray*}

Let $P^r_{2N}:=P\bigl( \{X_2,X_3\} \subset \NY^r(X_1) \bigr)$, $P^r_{2G}:=P\bigl(\{X_2,X_3\} \subset \G_1(X_1,\NY^r)\bigr)$, and $P^r_{M}:=P\bigl( X_2 \in \NY^r(X_1), X_3 \in \G_1(X_1,\NY^r) \bigr)$. Then for $J>1$, we have
\begin{eqnarray*}
P\bigl( \{X_2,X_3\} \subset \NY^r(X_1) \bigr)&=&\sum_{j=1}^{J}P\bigl( \{X_2,X_3\} \subset \NY^r(X_1)\,|\, \{X_1,X_2,X_3\} \subset T_j \bigr)\, P\bigl( \{X_1,X_2,X_3\} \subset T_j \bigr)\\
& = &\sum_{j=1}^{J}P^r_{2N}\, \bigl( A(T_j) / A(C_H(\Y)) \bigr)^3 =P^r_{2N}\, \left(\sum_{j=1}^{J}w_j^3 \right).
\end{eqnarray*}
Similarly, $P\bigl( X_2 \in \NY^r(X_1), X_3 \in \G_1(X_1,\NY^r) \bigr)=P^r_{M}\,\left(\sum_{j=1}^{J}w_j^3 \right) \text{  and  }P\bigl( \{X_2,X_3\} \subset \G_1(X_1,\NY^r) \bigr)=P^r_{2G}\,\Bigl(\sum_{j=1}^{J}w_j^3 \Bigr)$, hence,
{\small
$$\nu(r,J)=\bigl( P^r_{2N}+2\,P^r_M+P^r_{2G} \bigr)\,\left(\sum_{j=1}^{J}w_j^3 \right)-4\,\,\mu(r,J)^2=\nu(r)\,\left(\sum_{j=1}^{J}w_j^3 \right)+4\,\,\mu(r)^2\,\left(\sum_{j=1}^{J}w_j^3-\left( \sum_{j=1}^{J}w_j^2 \right)^2\right),$$
}
so conditional on $\mathcal W$, if $\nu(r,J)>0$ then $\sqrt{n}\,\bigl( \rho_n(r,J)-\mu(r,J)\bigr) \stackrel {\mathcal L}{\longrightarrow} \N(0,\nu(r,J))$.
$\blacksquare$

\end{document}

%% file: Nofnu.pstex_t
\begin{picture}(0,0)%
\includegraphics{Nofnu.pstex}%
\end{picture}%
\setlength{\unitlength}{3947sp}%
\begingroup\makeatletter\ifx\SetFigFont\undefined%
\gdef\SetFigFont#1#2#3#4#5{%
  \reset@font\fontsize{#1}{#2pt}%
  \fontfamily{#3}\fontseries{#4}\fontshape{#5}%
  \selectfont}%
\fi\endgroup%
\begin{picture}(10684,8043)(589,-7348)
\put(6975,-1389){\rotatebox{310.0}{\makebox(0,0)[lb]{\smash{{\SetFigFont{22}{14.4}{\rmdefault}{\bfdefault}{\updefault}{\color[rgb]{0,0,0}$e(x)$}%
}}}}}
\put(1405,-2508){\rotatebox{320.0}{\makebox(0,0)[lb]{\smash{{\SetFigFont{22}{14.4}{\rmdefault}{\mddefault}{\updefault}{\color[rgb]{0,0,0}$\ell(x)$}%
}}}}}
\put(4876,539){\makebox(0,0)[lb]{\smash{{\SetFigFont{22}{14.4}{\rmdefault}{\mddefault}{\updefault}{\color[rgb]{0,0,0}$\y_3$}%
}}}}
\put(11251,-6136){\makebox(0,0)[lb]{\smash{{\SetFigFont{22}{14.4}{\rmdefault}{\mddefault}{\updefault}{\color[rgb]{0,0,0}$\y_2$}%
}}}}
\put(2326,-5461){\makebox(0,0)[lb]{\smash{{\SetFigFont{22}{14.4}{\rmdefault}{\mddefault}{\updefault}{\color[rgb]{0,0,0}$R(\y_1)$}%
}}}}
\put(4276,-2161){\makebox(0,0)[lb]{\smash{{\SetFigFont{22}{14.4}{\rmdefault}{\mddefault}{\updefault}{\color[rgb]{0,0,0}$R(\y_3)$}%
}}}}
\put(7351,-5236){\makebox(0,0)[lb]{\smash{{\SetFigFont{22}{14.4}{\rmdefault}{\mddefault}{\updefault}{\color[rgb]{0,0,0}$R(\y_2)$}%
}}}}
\put(2551,164){\rotatebox{315.0}{\makebox(0,0)[lb]{\smash{{\SetFigFont{22}{14.4}{\rmdefault}{\mddefault}{\updefault}{\color[rgb]{0,0,0}$\ell_2(x)$}%
}}}}}
\put(960,-5299){\rotatebox{45.0}{\makebox(0,0)[lb]{\smash{{\SetFigFont{22}{14.4}{\rmdefault}{\mddefault}{\updefault}{\color[rgb]{0,0,0}$d(v(x),\ell(x))$}%
}}}}}
\put(1051,-6511){\makebox(0,0)[lb]{\smash{{\SetFigFont{22}{14.4}{\rmdefault}{\mddefault}{\updefault}{\color[rgb]{0,0,0}$\y_1=v(x)$}%
}}}}
\put(3301,-4936){\makebox(0,0)[lb]{\smash{{\SetFigFont{22}{14.4}{\rmdefault}{\bfdefault}{\updefault}{\color[rgb]{0,0,0}$x$}%
}}}}
\put(3451,-6736){\rotatebox{45.0}{\makebox(0,0)[lb]{\smash{{\SetFigFont{22}{14.4}{\rmdefault}{\mddefault}{\updefault}{\color[rgb]{0,0,0}$d(v(x),\ell_2(x))=2\,d(v(x),\ell(x))$}%
}}}}}
\end{picture}%

%% file: ArcProbEps.pstex_t
\begin{picture}(0,0)%
\includegraphics{ArcProbEps.pstex}%
\end{picture}%
\setlength{\unitlength}{3947sp}%
\begingroup\makeatletter\ifx\SetFigFont\undefined%
\gdef\SetFigFont#1#2#3#4#5{%
  \reset@font\fontsize{#1}{#2pt}%
  \fontfamily{#3}\fontseries{#4}\fontshape{#5}%
  \selectfont}%
\fi\endgroup%
\begin{picture}(10204,10004)(-2,-10763)
\put(5101,-10036){\makebox(0,0)[lb]{\smash{{\SetFigFont{25}{16.8}{\rmdefault}{\mddefault}{\updefault}{\color[rgb]{0,0,0}$r$}%
}}}}
\put(8101,-3061){\makebox(0,0)[lb]{\smash{{\SetFigFont{25}{16.8}{\rmdefault}{\mddefault}{\updefault}{\color[rgb]{0,0,0}$\mu(\NPE^r)$}%
}}}}
\put(6901,-5236){\makebox(0,0)[lb]{\smash{{\SetFigFont{25}{16.8}{\rmdefault}{\mddefault}{\updefault}{\color[rgb]{0,0,0}$\mu_S\left( \NPE^r,2\,\sqrt{3}/7 \right)$}%
}}}}
\put(6901,-4411){\makebox(0,0)[lb]{\smash{{\SetFigFont{25}{16.8}{\rmdefault}{\mddefault}{\updefault}{\color[rgb]{0,0,0}$\mu_S\left( \NPE^r,\sqrt{3}/4 \right)$}%
}}}}
\put(6826,-3661){\makebox(0,0)[lb]{\smash{{\SetFigFont{25}{16.8}{\rmdefault}{\mddefault}{\updefault}{\color[rgb]{0,0,0}$\mu_S\left( \NPE^r,\sqrt{3}/8 \right)$}%
}}}}
\end{picture}%

%% file: AsyVarEps.pstex_t
\begin{picture}(0,0)%
\includegraphics{AsyVarEps.pstex}%
\end{picture}%
\setlength{\unitlength}{3947sp}%
\begingroup\makeatletter\ifx\SetFigFont\undefined%
\gdef\SetFigFont#1#2#3#4#5{%
  \reset@font\fontsize{#1}{#2pt}%
  \fontfamily{#3}\fontseries{#4}\fontshape{#5}%
  \selectfont}%
\fi\endgroup%
\begin{picture}(10204,10004)(-2,-10763)
\put(5026,-10111){\makebox(0,0)[lb]{\smash{{\SetFigFont{25}{16.8}{\rmdefault}{\mddefault}{\updefault}{\color[rgb]{0,0,0}$r$}%
}}}}
\put(7801,-5086){\makebox(0,0)[lb]{\smash{{\SetFigFont{25}{16.8}{\rmdefault}{\mddefault}{\updefault}{\color[rgb]{0,0,0}$\nu(\NPE^r)$}%
}}}}
\put(1126,-6811){\makebox(0,0)[lb]{\smash{{\SetFigFont{25}{16.8}{\rmdefault}{\mddefault}{\updefault}{\color[rgb]{0,0,0}$\nu_S\left( \NPE^r,2\,\sqrt{3}/7 \right)$}%
}}}}
\put(4126,-8761){\makebox(0,0)[lb]{\smash{{\SetFigFont{25}{16.8}{\rmdefault}{\mddefault}{\updefault}{\color[rgb]{0,0,0}$\nu_S\left( \NPE^r,\sqrt{3}/4 \right)$}%
}}}}
\put(5551,-7111){\makebox(0,0)[lb]{\smash{{\SetFigFont{25}{16.8}{\rmdefault}{\mddefault}{\updefault}{\color[rgb]{0,0,0}$\nu_S\left( \NPE^r,\sqrt{3}/8 \right)$}%
}}}}
\end{picture}%

%% file: hlae_agg3.pstex_t
\begin{picture}(0,0)%
\includegraphics{hlae_agg3.pstex}%
\end{picture}%
\setlength{\unitlength}{3947sp}%
\begingroup\makeatletter\ifx\SetFigFont\undefined%
\gdef\SetFigFont#1#2#3#4#5{%
  \reset@font\fontsize{#1}{#2pt}%
  \fontfamily{#3}\fontseries{#4}\fontshape{#5}%
  \selectfont}%
\fi\endgroup%
\begin{picture}(10204,10004)(-2,-10763)
\put(5054,-10349){\makebox(0,0)[lb]{\smash{{\SetFigFont{24}{16.8}{\rmdefault}{\mddefault}{\updefault}{\color[rgb]{0,0,0}$r$}%
}}}}
\end{picture}%

%% file: hlae_agg2.pstex_t
\begin{picture}(0,0)%
\includegraphics{hlae_agg2.pstex}%
\end{picture}%
\setlength{\unitlength}{3947sp}%
\begingroup\makeatletter\ifx\SetFigFont\undefined%
\gdef\SetFigFont#1#2#3#4#5{%
  \reset@font\fontsize{#1}{#2pt}%
  \fontfamily{#3}\fontseries{#4}\fontshape{#5}%
  \selectfont}%
\fi\endgroup%
\begin{picture}(10204,10004)(-2,-10763)
\put(5054,-10349){\makebox(0,0)[lb]{\smash{{\SetFigFont{24}{16.8}{\rmdefault}{\mddefault}{\itdefault}{\color[rgb]{0,0,0}$r$}%
}}}}
\end{picture}%

%% file: hlae_agg1.pstex_t
\begin{picture}(0,0)%
\includegraphics{hlae_agg1.pstex}%
\end{picture}%
\setlength{\unitlength}{3947sp}%
\begingroup\makeatletter\ifx\SetFigFont\undefined%
\gdef\SetFigFont#1#2#3#4#5{%
  \reset@font\fontsize{#1}{#2pt}%
  \fontfamily{#3}\fontseries{#4}\fontshape{#5}%
  \selectfont}%
\fi\endgroup%
\begin{picture}(10204,10004)(-2,-10763)
\put(5054,-10349){\makebox(0,0)[lb]{\smash{{\SetFigFont{24}{16.8}{\rmdefault}{\mddefault}{\updefault}{\color[rgb]{0,0,0}$r$}%
}}}}
\end{picture}%

%% file: ArcProbEpsAgg.pstex_t
\begin{picture}(0,0)%
\includegraphics{ArcProbEpsAgg.pstex}%
\end{picture}%
\setlength{\unitlength}{3947sp}%
\begingroup\makeatletter\ifx\SetFigFont\undefined%
\gdef\SetFigFont#1#2#3#4#5{%
  \reset@font\fontsize{#1}{#2pt}%
  \fontfamily{#3}\fontseries{#4}\fontshape{#5}%
  \selectfont}%
\fi\endgroup%
\begin{picture}(10204,10004)(-2,-10763)
\put(5101,-10111){\makebox(0,0)[lb]{\smash{{\SetFigFont{25}{16.8}{\rmdefault}{\mddefault}{\updefault}{\color[rgb]{0,0,0}$r$}%
}}}}
\put(1726,-3136){\makebox(0,0)[lb]{\smash{{\SetFigFont{25}{16.8}{\rmdefault}{\mddefault}{\updefault}{\color[rgb]{0,0,0}$\mu(\NPE^r)$}%
}}}}
\put(7126,-5686){\makebox(0,0)[lb]{\smash{{\SetFigFont{25}{16.8}{\rmdefault}{\mddefault}{\updefault}{\color[rgb]{0,0,0}$\mu_A\left( \NPE^r,\sqrt{3}/21 \right)$}%
}}}}
\put(7426,-4861){\makebox(0,0)[lb]{\smash{{\SetFigFont{25}{16.8}{\rmdefault}{\mddefault}{\updefault}{\color[rgb]{0,0,0}$\mu_A\left( \NPE^r,\sqrt{3}/12 \right)$}%
}}}}
\put(7651,-4036){\makebox(0,0)[lb]{\smash{{\SetFigFont{25}{16.8}{\rmdefault}{\mddefault}{\updefault}{\color[rgb]{0,0,0}$\mu_A\left( \NPE^r,5\,\sqrt{3}/24 \right)$}%
}}}}
\end{picture}%

%% file: AsyVarEpsAgg.pstex_t
\begin{picture}(0,0)%
\includegraphics{AsyVarEpsAgg.pstex}%
\end{picture}%
\setlength{\unitlength}{3947sp}%
\begingroup\makeatletter\ifx\SetFigFont\undefined%
\gdef\SetFigFont#1#2#3#4#5{%
  \reset@font\fontsize{#1}{#2pt}%
  \fontfamily{#3}\fontseries{#4}\fontshape{#5}%
  \selectfont}%
\fi\endgroup%
\begin{picture}(10204,10004)(-2,-10763)
\put(4651,-10111){\makebox(0,0)[lb]{\smash{{\SetFigFont{25}{16.8}{\rmdefault}{\mddefault}{\updefault}{\color[rgb]{0,0,0}$r$}%
}}}}
\put(7876,-8911){\makebox(0,0)[lb]{\smash{{\SetFigFont{25}{16.8}{\rmdefault}{\mddefault}{\updefault}{\color[rgb]{0,0,0}$\nu(\NPE^r)$}%
}}}}
\put(6676,-4936){\makebox(0,0)[lb]{\smash{{\SetFigFont{25}{16.8}{\rmdefault}{\mddefault}{\updefault}{\color[rgb]{0,0,0}$\nu_A\left( \NPE^r,5\,\sqrt{3}/24 \right)$}%
}}}}
\put(7051,-6361){\makebox(0,0)[lb]{\smash{{\SetFigFont{25}{16.8}{\rmdefault}{\mddefault}{\updefault}{\color[rgb]{0,0,0}$\nu_A\left( \NPE^r,\sqrt{3}/12 \right)$}%
}}}}
\put(7426,-7336){\makebox(0,0)[lb]{\smash{{\SetFigFont{25}{16.8}{\rmdefault}{\mddefault}{\updefault}{\color[rgb]{0,0,0}$\nu_A\left( \NPE^r,\sqrt{3}/21 \right)$}%
}}}}
\end{picture}%

%% file: ls_lam_cases.pstex_t
\begin{picture}(0,0)%
\includegraphics{ls_lam_cases.pstex}%
\end{picture}%
\setlength{\unitlength}{3947sp}%
\begingroup\makeatletter\ifx\SetFigFont\undefined%
\gdef\SetFigFont#1#2#3#4#5{%
  \reset@font\fontsize{#1}{#2pt}%
  \fontfamily{#3}\fontseries{#4}\fontshape{#5}%
  \selectfont}%
\fi\endgroup%
\begin{picture}(11349,8280)(289,-7498)
\put(10426,-6586){\makebox(0,0)[lb]{\smash{{\SetFigFont{22}{16.8}{\rmdefault}{\mddefault}{\updefault}{\color[rgb]{0,0,0}$\y_2=(1,0)$}%
}}}}
\put(7201,-6586){\makebox(0,0)[lb]{\smash{{\SetFigFont{22}{16.8}{\rmdefault}{\mddefault}{\updefault}{\color[rgb]{0,0,0}$e_3$}%
}}}}
\put(5101,-6586){\makebox(0,0)[lb]{\smash{{\SetFigFont{22}{16.8}{\rmdefault}{\mddefault}{\updefault}{\color[rgb]{0,0,0}$M_3$}%
}}}}
\put(3301,-6586){\makebox(0,0)[lb]{\smash{{\SetFigFont{22}{14.4}{\rmdefault}{\mddefault}{\updefault}{\color[rgb]{0,0,0}$s_1$}%
}}}}
\put(1201,-2611){\makebox(0,0)[lb]{\smash{{\SetFigFont{22}{16.8}{\rmdefault}{\mddefault}{\updefault}{\color[rgb]{0,0,0}$\ell_s(r=4,x)$}%
}}}}
\put(1801,-736){\makebox(0,0)[lb]{\smash{{\SetFigFont{22}{16.8}{\rmdefault}{\mddefault}{\updefault}{\color[rgb]{0,0,0}$\ell_s(r=1.75,x)$}%
}}}}
\put(3226,314){\makebox(0,0)[lb]{\smash{{\SetFigFont{22}{16.8}{\rmdefault}{\mddefault}{\updefault}{\color[rgb]{0,0,0}$\ell_s\bigl(r=\sqrt{2},x\bigr)$}%
}}}}
\put(5926,614){\makebox(0,0)[lb]{\smash{{\SetFigFont{22}{16.8}{\rmdefault}{\mddefault}{\updefault}{\color[rgb]{0,0,0}$\y_3=\bigl(1/2,\sqrt{3}/2\bigr)$}%
}}}}
\put(7876,-2236){\makebox(0,0)[lb]{\smash{{\SetFigFont{22}{16.8}{\rmdefault}{\mddefault}{\updefault}{\color[rgb]{0,0,0}$e_1$}%
}}}}
\put(3976,-6586){\makebox(0,0)[lb]{\smash{{\SetFigFont{22}{14.4}{\rmdefault}{\mddefault}{\updefault}{\color[rgb]{0,0,0}$s_2$}%
}}}}
\put(826,-6736){\makebox(0,0)[lb]{\smash{{\SetFigFont{22}{16.8}{\rmdefault}{\mddefault}{\updefault}{\color[rgb]{0,0,0}$\y_1=(0,0)$}%
}}}}
\put(5851,-4111){\makebox(0,0)[lb]{\smash{{\SetFigFont{22}{16.8}{\rmdefault}{\mddefault}{\updefault}{\color[rgb]{0,0,0}$M_{C}$}%
}}}}
\put(3001,-2761){\makebox(0,0)[lb]{\smash{{\SetFigFont{22}{16.8}{\rmdefault}{\mddefault}{\updefault}{\color[rgb]{0,0,0}$e_2$}%
}}}}
\end{picture}%

%% file: G1ofxCase1.pstex_t
\begin{picture}(0,0)%
\includegraphics{G1ofxCase1.pstex}%
\end{picture}%
\setlength{\unitlength}{3947sp}%
\begingroup\makeatletter\ifx\SetFigFont\undefined%
\gdef\SetFigFont#1#2#3#4#5{%
  \reset@font\fontsize{#1}{#2pt}%
  \fontfamily{#3}\fontseries{#4}\fontshape{#5}%
  \selectfont}%
\fi\endgroup%
\begin{picture}(11349,8130)(289,-7348)
\put(10426,-6586){\makebox(0,0)[lb]{\smash{{\SetFigFont{25}{16.8}{\rmdefault}{\mddefault}{\updefault}{\color[rgb]{0,0,0}$\y_2=(1,0)$}%
}}}}
\put(826,-6586){\makebox(0,0)[lb]{\smash{{\SetFigFont{25}{16.8}{\rmdefault}{\mddefault}{\updefault}{\color[rgb]{0,0,0}$\y_1=(0,0)$}%
}}}}
\put(5326,614){\makebox(0,0)[lb]{\smash{{\SetFigFont{25}{16.8}{\rmdefault}{\mddefault}{\updefault}{\color[rgb]{0,0,0}$\y_3=(1/2,\sqrt{3}/2)$}%
}}}}
\put(7726,-2236){\makebox(0,0)[lb]{\smash{{\SetFigFont{25}{16.8}{\rmdefault}{\mddefault}{\updefault}{\color[rgb]{0,0,0}$e_1$}%
}}}}
\put(3001,-2686){\makebox(0,0)[lb]{\smash{{\SetFigFont{25}{16.8}{\rmdefault}{\mddefault}{\updefault}{\color[rgb]{0,0,0}$e_2$}%
}}}}
\put(7201,-6586){\makebox(0,0)[lb]{\smash{{\SetFigFont{25}{16.8}{\rmdefault}{\mddefault}{\updefault}{\color[rgb]{0,0,0}$e_3$}%
}}}}
\put(5626,-6586){\makebox(0,0)[lb]{\smash{{\SetFigFont{25}{16.8}{\rmdefault}{\mddefault}{\updefault}{\color[rgb]{0,0,0}$M_3$}%
}}}}
\put(2026,-5986){\makebox(0,0)[lb]{\smash{{\SetFigFont{25}{16.8}{\rmdefault}{\mddefault}{\updefault}{\color[rgb]{0,0,0}$\xi_1(r,x)$}%
}}}}
\put(5701,-4111){\makebox(0,0)[lb]{\smash{{\SetFigFont{25}{16.8}{\rmdefault}{\mddefault}{\updefault}{\color[rgb]{0,0,0}$M_{C}$}%
}}}}
\put(1426,-6211){\makebox(0,0)[lb]{\smash{{\SetFigFont{25}{16.8}{\rmdefault}{\mddefault}{\updefault}{\color[rgb]{0,0,0}$x_1$}%
}}}}
\end{picture}%

%% file: G1ofxCase2.pstex_t
\begin{picture}(0,0)%
\includegraphics{G1ofxCase2.pstex}%
\end{picture}%
\setlength{\unitlength}{3947sp}%
\begingroup\makeatletter\ifx\SetFigFont\undefined%
\gdef\SetFigFont#1#2#3#4#5{%
  \reset@font\fontsize{#1}{#2pt}%
  \fontfamily{#3}\fontseries{#4}\fontshape{#5}%
  \selectfont}%
\fi\endgroup%
\begin{picture}(11349,8130)(289,-7348)
\put(10426,-6586){\makebox(0,0)[lb]{\smash{{\SetFigFont{25}{16.8}{\rmdefault}{\mddefault}{\updefault}{\color[rgb]{0,0,0}$\y_2=(1,0)$}%
}}}}
\put(826,-6586){\makebox(0,0)[lb]{\smash{{\SetFigFont{25}{16.8}{\rmdefault}{\mddefault}{\updefault}{\color[rgb]{0,0,0}$\y_1=(0,0)$}%
}}}}
\put(5326,614){\makebox(0,0)[lb]{\smash{{\SetFigFont{25}{16.8}{\rmdefault}{\mddefault}{\updefault}{\color[rgb]{0,0,0}$\y_3=(1/2,\sqrt{3}/2)$}%
}}}}
\put(7726,-2236){\makebox(0,0)[lb]{\smash{{\SetFigFont{25}{16.8}{\rmdefault}{\mddefault}{\updefault}{\color[rgb]{0,0,0}$e_1$}%
}}}}
\put(3001,-2686){\makebox(0,0)[lb]{\smash{{\SetFigFont{25}{16.8}{\rmdefault}{\mddefault}{\updefault}{\color[rgb]{0,0,0}$e_2$}%
}}}}
\put(7201,-6586){\makebox(0,0)[lb]{\smash{{\SetFigFont{25}{16.8}{\rmdefault}{\mddefault}{\updefault}{\color[rgb]{0,0,0}$e_3$}%
}}}}
\put(5626,-6586){\makebox(0,0)[lb]{\smash{{\SetFigFont{25}{16.8}{\rmdefault}{\mddefault}{\updefault}{\color[rgb]{0,0,0}$M_3$}%
}}}}
\put(2026,-5986){\makebox(0,0)[lb]{\smash{{\SetFigFont{25}{16.8}{\rmdefault}{\mddefault}{\updefault}{\color[rgb]{0,0,0}$\xi_1(r,x)$}%
}}}}
\put(5701,-4111){\makebox(0,0)[lb]{\smash{{\SetFigFont{25}{16.8}{\rmdefault}{\mddefault}{\updefault}{\color[rgb]{0,0,0}$M_{C}$}%
}}}}
\put(3601,-6136){\makebox(0,0)[lb]{\smash{{\SetFigFont{25}{16.8}{\rmdefault}{\mddefault}{\updefault}{\color[rgb]{0,0,0}$x_1$}%
}}}}
\put(6226,-4636){\makebox(0,0)[lb]{\smash{{\SetFigFont{25}{16.8}{\rmdefault}{\mddefault}{\updefault}{\color[rgb]{0,0,0}$\xi_2(r,x)$}%
}}}}
\end{picture}%

%% file: G1ofxCase3.pstex_t
\begin{picture}(0,0)%
\includegraphics{G1ofxCase3.pstex}%
\end{picture}%
\setlength{\unitlength}{3947sp}%
\begingroup\makeatletter\ifx\SetFigFont\undefined%
\gdef\SetFigFont#1#2#3#4#5{%
  \reset@font\fontsize{#1}{#2pt}%
  \fontfamily{#3}\fontseries{#4}\fontshape{#5}%
  \selectfont}%
\fi\endgroup%
\begin{picture}(11349,8130)(289,-7348)
\put(10426,-6586){\makebox(0,0)[lb]{\smash{{\SetFigFont{25}{16.8}{\rmdefault}{\mddefault}{\updefault}{\color[rgb]{0,0,0}$\y_2=(1,0)$}%
}}}}
\put(826,-6586){\makebox(0,0)[lb]{\smash{{\SetFigFont{25}{16.8}{\rmdefault}{\mddefault}{\updefault}{\color[rgb]{0,0,0}$\y_1=(0,0)$}%
}}}}
\put(5326,614){\makebox(0,0)[lb]{\smash{{\SetFigFont{25}{16.8}{\rmdefault}{\mddefault}{\updefault}{\color[rgb]{0,0,0}$\y_3=(1/2,\sqrt{3}/2)$}%
}}}}
\put(7726,-2236){\makebox(0,0)[lb]{\smash{{\SetFigFont{25}{16.8}{\rmdefault}{\mddefault}{\updefault}{\color[rgb]{0,0,0}$e_1$}%
}}}}
\put(3001,-2686){\makebox(0,0)[lb]{\smash{{\SetFigFont{25}{16.8}{\rmdefault}{\mddefault}{\updefault}{\color[rgb]{0,0,0}$e_2$}%
}}}}
\put(7201,-6586){\makebox(0,0)[lb]{\smash{{\SetFigFont{25}{16.8}{\rmdefault}{\mddefault}{\updefault}{\color[rgb]{0,0,0}$e_3$}%
}}}}
\put(5626,-6586){\makebox(0,0)[lb]{\smash{{\SetFigFont{25}{16.8}{\rmdefault}{\mddefault}{\updefault}{\color[rgb]{0,0,0}$M_3$}%
}}}}
\put(5701,-4111){\makebox(0,0)[lb]{\smash{{\SetFigFont{25}{16.8}{\rmdefault}{\mddefault}{\updefault}{\color[rgb]{0,0,0}$M_{C}$}%
}}}}
\put(7726,-5236){\makebox(0,0)[lb]{\smash{{\SetFigFont{25}{16.8}{\rmdefault}{\mddefault}{\updefault}{\color[rgb]{0,0,0}$\xi_2(r,x)$}%
}}}}
\put(4651,-6136){\makebox(0,0)[lb]{\smash{{\SetFigFont{25}{16.8}{\rmdefault}{\mddefault}{\updefault}{\color[rgb]{0,0,0}$x_1$}%
}}}}
\put(2851,-5536){\makebox(0,0)[lb]{\smash{{\SetFigFont{25}{16.8}{\rmdefault}{\mddefault}{\updefault}{\color[rgb]{0,0,0}$\xi_1(r,x)$}%
}}}}
\end{picture}%

%% file: G1ofxCase4.pstex_t
\begin{picture}(0,0)%
\includegraphics{G1ofxCase4.pstex}%
\end{picture}%
\setlength{\unitlength}{3947sp}%
\begingroup\makeatletter\ifx\SetFigFont\undefined%
\gdef\SetFigFont#1#2#3#4#5{%
  \reset@font\fontsize{#1}{#2pt}%
  \fontfamily{#3}\fontseries{#4}\fontshape{#5}%
  \selectfont}%
\fi\endgroup%
\begin{picture}(11349,8130)(289,-7348)
\put(10426,-6586){\makebox(0,0)[lb]{\smash{{\SetFigFont{25}{16.8}{\rmdefault}{\mddefault}{\updefault}{\color[rgb]{0,0,0}$\y_2=(1,0)$}%
}}}}
\put(826,-6586){\makebox(0,0)[lb]{\smash{{\SetFigFont{25}{16.8}{\rmdefault}{\mddefault}{\updefault}{\color[rgb]{0,0,0}$\y_1=(0,0)$}%
}}}}
\put(5326,614){\makebox(0,0)[lb]{\smash{{\SetFigFont{25}{16.8}{\rmdefault}{\mddefault}{\updefault}{\color[rgb]{0,0,0}$\y_3=(1/2,\sqrt{3}/2)$}%
}}}}
\put(7201,-6586){\makebox(0,0)[lb]{\smash{{\SetFigFont{25}{16.8}{\rmdefault}{\mddefault}{\updefault}{\color[rgb]{0,0,0}$e_3$}%
}}}}
\put(2026,-5986){\makebox(0,0)[lb]{\smash{{\SetFigFont{25}{16.8}{\rmdefault}{\mddefault}{\updefault}{\color[rgb]{0,0,0}$\xi_1(r,x)$}%
}}}}
\put(5701,-4111){\makebox(0,0)[lb]{\smash{{\SetFigFont{25}{16.8}{\rmdefault}{\mddefault}{\updefault}{\color[rgb]{0,0,0}$M_{C}$}%
}}}}
\put(3826,-5536){\makebox(0,0)[lb]{\smash{{\SetFigFont{25}{16.8}{\rmdefault}{\mddefault}{\updefault}{\color[rgb]{0,0,0}$x_1$}%
}}}}
\put(3901,-1111){\makebox(0,0)[lb]{\smash{{\SetFigFont{25}{16.8}{\rmdefault}{\mddefault}{\updefault}{\color[rgb]{0,0,0}$e_2$}%
}}}}
\put(6901,-1111){\makebox(0,0)[lb]{\smash{{\SetFigFont{25}{16.8}{\rmdefault}{\mddefault}{\updefault}{\color[rgb]{0,0,0}$e_1$}%
}}}}
\put(5326,-6586){\makebox(0,0)[lb]{\smash{{\SetFigFont{25}{16.8}{\rmdefault}{\mddefault}{\updefault}{\color[rgb]{0,0,0}$M_3$}%
}}}}
\put(3151,-6661){\makebox(0,0)[lb]{\smash{{\SetFigFont{25}{16.8}{\rmdefault}{\mddefault}{\updefault}{\color[rgb]{0,0,0}$G_1$}%
}}}}
\put(1501,-4936){\makebox(0,0)[lb]{\smash{{\SetFigFont{25}{16.8}{\rmdefault}{\mddefault}{\updefault}{\color[rgb]{0,0,0}$G_6$}%
}}}}
\put(2851,-2686){\makebox(0,0)[lb]{\smash{{\SetFigFont{25}{16.8}{\rmdefault}{\mddefault}{\updefault}{\color[rgb]{0,0,0}$M_2$}%
}}}}
\put(3976,-2986){\makebox(0,0)[lb]{\smash{{\SetFigFont{25}{16.8}{\rmdefault}{\mddefault}{\updefault}{\color[rgb]{0,0,0}$L_5$}%
}}}}
\put(5251,-2986){\makebox(0,0)[lb]{\smash{{\SetFigFont{25}{16.8}{\rmdefault}{\mddefault}{\updefault}{\color[rgb]{0,0,0}$\xi_3(r,x)$}%
}}}}
\put(6901,-2986){\makebox(0,0)[lb]{\smash{{\SetFigFont{25}{16.8}{\rmdefault}{\mddefault}{\updefault}{\color[rgb]{0,0,0}$L_4$}%
}}}}
\put(7651,-3061){\makebox(0,0)[lb]{\smash{{\SetFigFont{25}{16.8}{\rmdefault}{\mddefault}{\updefault}{\color[rgb]{0,0,0}$L_3$}%
}}}}
\put(5626,-6211){\makebox(0,0)[lb]{\smash{{\SetFigFont{25}{16.8}{\rmdefault}{\mddefault}{\updefault}{\color[rgb]{0,0,0}$L_2$}%
}}}}
\put(6526,-4786){\makebox(0,0)[lb]{\smash{{\SetFigFont{25}{16.8}{\rmdefault}{\mddefault}{\updefault}{\color[rgb]{0,0,0}$\xi_2(r,x)$}%
}}}}
\put(7951,-2686){\makebox(0,0)[lb]{\smash{{\SetFigFont{25}{16.8}{\rmdefault}{\mddefault}{\updefault}{\color[rgb]{0,0,0}$M_1$}%
}}}}
\end{picture}%

%% file: G1ofxCase5.pstex_t
\begin{picture}(0,0)%
\includegraphics{G1ofxCase5.pstex}%
\end{picture}%
\setlength{\unitlength}{3947sp}%
\begingroup\makeatletter\ifx\SetFigFont\undefined%
\gdef\SetFigFont#1#2#3#4#5{%
  \reset@font\fontsize{#1}{#2pt}%
  \fontfamily{#3}\fontseries{#4}\fontshape{#5}%
  \selectfont}%
\fi\endgroup%
\begin{picture}(11349,8130)(289,-7348)
\put(10426,-6586){\makebox(0,0)[lb]{\smash{{\SetFigFont{25}{16.8}{\rmdefault}{\mddefault}{\updefault}{\color[rgb]{0,0,0}$\y_2=(1,0)$}%
}}}}
\put(826,-6586){\makebox(0,0)[lb]{\smash{{\SetFigFont{25}{16.8}{\rmdefault}{\mddefault}{\updefault}{\color[rgb]{0,0,0}$\y_1=(0,0)$}%
}}}}
\put(5326,614){\makebox(0,0)[lb]{\smash{{\SetFigFont{25}{16.8}{\rmdefault}{\mddefault}{\updefault}{\color[rgb]{0,0,0}$\y_3=(1/2,\sqrt{3}/2)$}%
}}}}
\put(7726,-2236){\makebox(0,0)[lb]{\smash{{\SetFigFont{25}{16.8}{\rmdefault}{\mddefault}{\updefault}{\color[rgb]{0,0,0}$e_1$}%
}}}}
\put(3001,-2686){\makebox(0,0)[lb]{\smash{{\SetFigFont{25}{16.8}{\rmdefault}{\mddefault}{\updefault}{\color[rgb]{0,0,0}$e_2$}%
}}}}
\put(7201,-6586){\makebox(0,0)[lb]{\smash{{\SetFigFont{25}{16.8}{\rmdefault}{\mddefault}{\updefault}{\color[rgb]{0,0,0}$e_3$}%
}}}}
\put(5626,-6586){\makebox(0,0)[lb]{\smash{{\SetFigFont{25}{16.8}{\rmdefault}{\mddefault}{\updefault}{\color[rgb]{0,0,0}$M_3$}%
}}}}
\put(5701,-4111){\makebox(0,0)[lb]{\smash{{\SetFigFont{25}{16.8}{\rmdefault}{\mddefault}{\updefault}{\color[rgb]{0,0,0}$M_{C}$}%
}}}}
\put(4801,-2986){\makebox(0,0)[lb]{\smash{{\SetFigFont{25}{16.8}{\rmdefault}{\mddefault}{\updefault}{\color[rgb]{0,0,0}$\xi_3(r,x)$}%
}}}}
\put(4951,-5386){\makebox(0,0)[lb]{\smash{{\SetFigFont{25}{16.8}{\rmdefault}{\mddefault}{\updefault}{\color[rgb]{0,0,0}$x_1$}%
}}}}
\put(7576,-5236){\makebox(0,0)[lb]{\smash{{\SetFigFont{25}{16.8}{\rmdefault}{\mddefault}{\updefault}{\color[rgb]{0,0,0}$\xi_2(r,x)$}%
}}}}
\put(3076,-5536){\makebox(0,0)[lb]{\smash{{\SetFigFont{25}{16.8}{\rmdefault}{\mddefault}{\updefault}{\color[rgb]{0,0,0}$\xi_1(r,x)$}%
}}}}
\end{picture}%

%% file: G1ofxCase6.pstex_t
\begin{picture}(0,0)%
\includegraphics{G1ofxCase6.pstex}%
\end{picture}%
\setlength{\unitlength}{3947sp}%
\begingroup\makeatletter\ifx\SetFigFont\undefined%
\gdef\SetFigFont#1#2#3#4#5{%
  \reset@font\fontsize{#1}{#2pt}%
  \fontfamily{#3}\fontseries{#4}\fontshape{#5}%
  \selectfont}%
\fi\endgroup%
\begin{picture}(11349,8130)(289,-7348)
\put(10426,-6586){\makebox(0,0)[lb]{\smash{{\SetFigFont{25}{16.8}{\rmdefault}{\mddefault}{\updefault}{\color[rgb]{0,0,0}$\y_2=(1,0)$}%
}}}}
\put(826,-6586){\makebox(0,0)[lb]{\smash{{\SetFigFont{25}{16.8}{\rmdefault}{\mddefault}{\updefault}{\color[rgb]{0,0,0}$\y_1=(0,0)$}%
}}}}
\put(5326,614){\makebox(0,0)[lb]{\smash{{\SetFigFont{25}{16.8}{\rmdefault}{\mddefault}{\updefault}{\color[rgb]{0,0,0}$\y_3=(1/2,\sqrt{3}/2)$}%
}}}}
\put(7726,-2236){\makebox(0,0)[lb]{\smash{{\SetFigFont{25}{16.8}{\rmdefault}{\mddefault}{\updefault}{\color[rgb]{0,0,0}$e_1$}%
}}}}
\put(3001,-2686){\makebox(0,0)[lb]{\smash{{\SetFigFont{25}{16.8}{\rmdefault}{\mddefault}{\updefault}{\color[rgb]{0,0,0}$e_2$}%
}}}}
\put(7201,-6586){\makebox(0,0)[lb]{\smash{{\SetFigFont{25}{16.8}{\rmdefault}{\mddefault}{\updefault}{\color[rgb]{0,0,0}$e_3$}%
}}}}
\put(5626,-6586){\makebox(0,0)[lb]{\smash{{\SetFigFont{25}{16.8}{\rmdefault}{\mddefault}{\updefault}{\color[rgb]{0,0,0}$M_3$}%
}}}}
\put(5701,-4111){\makebox(0,0)[lb]{\smash{{\SetFigFont{25}{16.8}{\rmdefault}{\mddefault}{\updefault}{\color[rgb]{0,0,0}$M_{C}$}%
}}}}
\put(5251,-4411){\makebox(0,0)[lb]{\smash{{\SetFigFont{25}{16.8}{\rmdefault}{\mddefault}{\updefault}{\color[rgb]{0,0,0}$x_1$}%
}}}}
\put(5851,-4636){\makebox(0,0)[lb]{\smash{{\SetFigFont{25}{16.8}{\rmdefault}{\mddefault}{\updefault}{\color[rgb]{0,0,0}$\xi_2(r,x)$}%
}}}}
\put(5176,-3436){\makebox(0,0)[lb]{\smash{{\SetFigFont{25}{16.8}{\rmdefault}{\mddefault}{\updefault}{\color[rgb]{0,0,0}$\xi_3(r,x)$}%
}}}}
\put(4201,-4711){\makebox(0,0)[lb]{\smash{{\SetFigFont{25}{16.8}{\rmdefault}{\mddefault}{\updefault}{\color[rgb]{0,0,0}$\xi_1(r,x)$}%
}}}}
\put(5476,-6961){\makebox(0,0)[lb]{\smash{{\SetFigFont{25}{14.4}{\rmdefault}{\mddefault}{\updefault}{\color[rgb]{0,0,0}case-6}%
}}}}
\end{picture}%

%% file: N_nu2GamRegions.pstex_t
\begin{picture}(0,0)%
\includegraphics{N_nu2GamRegions.pstex}%
\end{picture}%
\setlength{\unitlength}{3947sp}%
\begingroup\makeatletter\ifx\SetFigFont\undefined%
\gdef\SetFigFont#1#2#3#4#5{%
  \reset@font\fontsize{#1}{#2pt}%
  \fontfamily{#3}\fontseries{#4}\fontshape{#5}%
  \selectfont}%
\fi\endgroup%
\begin{picture}(10923,8028)(301,-7186)
\put(3226,-5761){\makebox(0,0)[lb]{\smash{{\SetFigFont{22}{18}{\rmdefault}{\bfdefault}{\updefault}{\color[rgb]{0,0,0}$R_1$}%
}}}}
\put(6226,-5686){\makebox(0,0)[lb]{\smash{{\SetFigFont{22}{18}{\rmdefault}{\bfdefault}{\updefault}{\color[rgb]{0,0,0}$R_2$}%
}}}}
\put(9526,-5311){\makebox(0,0)[lb]{\smash{{\SetFigFont{22}{18}{\rmdefault}{\bfdefault}{\updefault}{\color[rgb]{0,0,0}$R_3$}%
}}}}
\put(8251,-3361){\makebox(0,0)[lb]{\smash{{\SetFigFont{22}{18}{\rmdefault}{\bfdefault}{\updefault}{\color[rgb]{0,0,0}$R_4$}%
}}}}
\put(9976,-3736){\makebox(0,0)[lb]{\smash{{\SetFigFont{22}{18}{\rmdefault}{\bfdefault}{\updefault}{\color[rgb]{0,0,0}$R_5$}%
}}}}
\put(9976,-2086){\makebox(0,0)[lb]{\smash{{\SetFigFont{22}{18}{\rmdefault}{\bfdefault}{\updefault}{\color[rgb]{0,0,0}$R_6$}%
}}}}
\put(7801,764){\rotatebox{300.0}{\makebox(0,0)[lb]{\smash{{\SetFigFont{22}{18}{\rmdefault}{\bfdefault}{\updefault}{\color[rgb]{0,0,0}$q_{12}(x)$}%
}}}}}
\put(8701,164){\rotatebox{60.0}{\makebox(0,0)[lb]{\smash{{\SetFigFont{22}{18}{\rmdefault}{\bfdefault}{\updefault}{\color[rgb]{0,0,0}$q_1(x)$}%
}}}}}
\put(1201,-3886){\makebox(0,0)[lb]{\smash{{\SetFigFont{22}{18}{\rmdefault}{\bfdefault}{\updefault}{\color[rgb]{0,0,0}$q_2(x)$}%
}}}}
\put(1201,-961){\makebox(0,0)[lb]{\smash{{\SetFigFont{22}{18}{\rmdefault}{\bfdefault}{\updefault}{\color[rgb]{0,0,0}$q_4(x)$}%
}}}}
\put(751,-6436){\makebox(0,0)[lb]{\smash{{\SetFigFont{22}{18}{\rmdefault}{\bfdefault}{\updefault}{\color[rgb]{0,0,0}$\y_1$}%
}}}}
\put(10501,-6436){\makebox(0,0)[lb]{\smash{{\SetFigFont{22}{18}{\rmdefault}{\bfdefault}{\updefault}{\color[rgb]{0,0,0}$M_3$}%
}}}}
\put(10501,-1786){\makebox(0,0)[lb]{\smash{{\SetFigFont{22}{18}{\rmdefault}{\bfdefault}{\updefault}{\color[rgb]{0,0,0}$M_{C}$}%
}}}}
\put(8328,-6188){\rotatebox{50.0}{\makebox(0,0)[lb]{\smash{{\SetFigFont{22}{18}{\rmdefault}{\bfdefault}{\updefault}{\color[rgb]{0,0,0}$q_3(x)$}%
}}}}}
\put(5626,-6586){\makebox(0,0)[lb]{\smash{{\SetFigFont{22}{18}{\rmdefault}{\bfdefault}{\updefault}{\color[rgb]{0,0,0}$s_2$}%
}}}}
\put(4051,-6511){\makebox(0,0)[lb]{\smash{{\SetFigFont{22}{18}{\rmdefault}{\bfdefault}{\updefault}{\color[rgb]{0,0,0}$s_1$}%
}}}}
\put(7951,-6586){\makebox(0,0)[lb]{\smash{{\SetFigFont{22}{18}{\rmdefault}{\bfdefault}{\updefault}{\color[rgb]{0,0,0}$s_3$}%
}}}}
\put(9901,-6811){\makebox(0,0)[lb]{\smash{{\SetFigFont{22}{18}{\rmdefault}{\bfdefault}{\updefault}{\color[rgb]{0,0,0}$s_5$}%
}}}}
\put(9151,-6811){\makebox(0,0)[lb]{\smash{{\SetFigFont{22}{18}{\rmdefault}{\bfdefault}{\updefault}{\color[rgb]{0,0,0}$s_4$}%
}}}}
\put(3301,-5086){\rotatebox{20.0}{\makebox(0,0)[lb]{\smash{{\SetFigFont{22}{18}{\rmdefault}{\bfdefault}{\updefault}{\color[rgb]{0,0,0}$\ell_{am}(x)$}%
}}}}}
\end{picture}%

%% file: G1ofxCase7.pstex_t
\begin{picture}(0,0)%
\includegraphics{G1ofxCase7.pstex}%
\end{picture}%
\setlength{\unitlength}{3947sp}%
\begingroup\makeatletter\ifx\SetFigFont\undefined%
\gdef\SetFigFont#1#2#3#4#5{%
  \reset@font\fontsize{#1}{#2pt}%
  \fontfamily{#3}\fontseries{#4}\fontshape{#5}%
  \selectfont}%
\fi\endgroup%
\begin{picture}(11349,8130)(289,-7348)
\put(10426,-6586){\makebox(0,0)[lb]{\smash{{\SetFigFont{25}{16.8}{\rmdefault}{\mddefault}{\updefault}{\color[rgb]{0,0,0}$\y_2=(1,0)$}%
}}}}
\put(826,-6586){\makebox(0,0)[lb]{\smash{{\SetFigFont{25}{16.8}{\rmdefault}{\mddefault}{\updefault}{\color[rgb]{0,0,0}$\y_1=(0,0)$}%
}}}}
\put(5326,614){\makebox(0,0)[lb]{\smash{{\SetFigFont{25}{16.8}{\rmdefault}{\mddefault}{\updefault}{\color[rgb]{0,0,0}$\y_3=(1/2,\sqrt{3}/2)$}%
}}}}
\put(5701,-4111){\makebox(0,0)[lb]{\smash{{\SetFigFont{25}{16.8}{\rmdefault}{\mddefault}{\updefault}{\color[rgb]{0,0,0}$M_{C}$}%
}}}}
\put(7576,-5236){\makebox(0,0)[lb]{\smash{{\SetFigFont{25}{16.8}{\rmdefault}{\mddefault}{\updefault}{\color[rgb]{0,0,0}$\xi_2(r,x)$}%
}}}}
\put(5251,-4486){\makebox(0,0)[lb]{\smash{{\SetFigFont{25}{16.8}{\rmdefault}{\mddefault}{\updefault}{\color[rgb]{0,0,0}$x_1$}%
}}}}
\put(3676,-5161){\makebox(0,0)[lb]{\smash{{\SetFigFont{25}{16.8}{\rmdefault}{\mddefault}{\updefault}{\color[rgb]{0,0,0}$\xi_1(r,x)$}%
}}}}
\put(5026,-2161){\makebox(0,0)[lb]{\smash{{\SetFigFont{25}{16.8}{\rmdefault}{\mddefault}{\updefault}{\color[rgb]{0,0,0}$\xi_3(r,x)$}%
}}}}
\put(5626,-6586){\makebox(0,0)[lb]{\smash{{\SetFigFont{25}{16.8}{\rmdefault}{\mddefault}{\updefault}{\color[rgb]{0,0,0}$M_3$}%
}}}}
\put(4426,-6661){\makebox(0,0)[lb]{\smash{{\SetFigFont{25}{16.8}{\rmdefault}{\mddefault}{\updefault}{\color[rgb]{0,0,0}$G_1$}%
}}}}
\put(6526,-6586){\makebox(0,0)[lb]{\smash{{\SetFigFont{25}{16.8}{\rmdefault}{\mddefault}{\updefault}{\color[rgb]{0,0,0}$G_2$}%
}}}}
\put(8101,-6586){\makebox(0,0)[lb]{\smash{{\SetFigFont{25}{16.8}{\rmdefault}{\mddefault}{\updefault}{\color[rgb]{0,0,0}$e_3$}%
}}}}
\put(6901,-1111){\makebox(0,0)[lb]{\smash{{\SetFigFont{25}{16.8}{\rmdefault}{\mddefault}{\updefault}{\color[rgb]{0,0,0}$e_1$}%
}}}}
\put(4201,-736){\makebox(0,0)[lb]{\smash{{\SetFigFont{25}{16.8}{\rmdefault}{\mddefault}{\updefault}{\color[rgb]{0,0,0}$e_2$}%
}}}}
\put(8101,-2761){\makebox(0,0)[lb]{\smash{{\SetFigFont{25}{16.8}{\rmdefault}{\mddefault}{\updefault}{\color[rgb]{0,0,0}$M_1$}%
}}}}
\put(2926,-2611){\makebox(0,0)[lb]{\smash{{\SetFigFont{25}{16.8}{\rmdefault}{\mddefault}{\updefault}{\color[rgb]{0,0,0}$M_2$}%
}}}}
\put(8701,-3511){\makebox(0,0)[lb]{\smash{{\SetFigFont{25}{16.8}{\rmdefault}{\mddefault}{\updefault}{\color[rgb]{0,0,0}$G_3$}%
}}}}
\put(7726,-2236){\makebox(0,0)[lb]{\smash{{\SetFigFont{25}{16.8}{\rmdefault}{\mddefault}{\updefault}{\color[rgb]{0,0,0}$G_4$}%
}}}}
\put(3226,-2161){\makebox(0,0)[lb]{\smash{{\SetFigFont{25}{16.8}{\rmdefault}{\mddefault}{\updefault}{\color[rgb]{0,0,0}$G_5$}%
}}}}
\put(2176,-3811){\makebox(0,0)[lb]{\smash{{\SetFigFont{25}{16.8}{\rmdefault}{\mddefault}{\updefault}{\color[rgb]{0,0,0}$G_6$}%
}}}}
\end{picture}%

%% file: N_nu6GamRegions.pstex_t
\begin{picture}(0,0)%
\includegraphics{N_nu6GamRegions.pstex}%
\end{picture}%
\setlength{\unitlength}{3947sp}%
\begingroup\makeatletter\ifx\SetFigFont\undefined%
\gdef\SetFigFont#1#2#3#4#5{%
  \reset@font\fontsize{#1}{#2pt}%
  \fontfamily{#3}\fontseries{#4}\fontshape{#5}%
  \selectfont}%
\fi\endgroup%
\begin{picture}(10738,8112)(301,-7348)
\put(751,-6361){\makebox(0,0)[lb]{\smash{{\SetFigFont{25}{14.4}{\rmdefault}{\mddefault}{\updefault}{\color[rgb]{0,0,0}$\y_1$}%
}}}}
\put(10501,-6511){\makebox(0,0)[lb]{\smash{{\SetFigFont{25}{14.4}{\rmdefault}{\mddefault}{\updefault}{\color[rgb]{0,0,0}$M_1$}%
}}}}
\put(10426,-1861){\makebox(0,0)[lb]{\smash{{\SetFigFont{25}{14.4}{\rmdefault}{\mddefault}{\updefault}{\color[rgb]{0,0,0}$M_C$}%
}}}}
\put(9601,-5911){\makebox(0,0)[lb]{\smash{{\SetFigFont{25}{14.4}{\rmdefault}{\mddefault}{\updefault}{\color[rgb]{0,0,0}$q_2(x)$}%
}}}}
\put(4126,-4786){\rotatebox{25.0}{\makebox(0,0)[lb]{\smash{{\SetFigFont{25}{14.4}{\rmdefault}{\mddefault}{\updefault}{\color[rgb]{0,0,0}$\ell_{am}(x)$}%
}}}}}
\put(9826,-661){\rotatebox{50.0}{\makebox(0,0)[lb]{\smash{{\SetFigFont{25}{14.4}{\rmdefault}{\mddefault}{\updefault}{\color[rgb]{0,0,0}$q_3(x)$}%
}}}}}
\put(5330,-661){\rotatebox{55.0}{\makebox(0,0)[lb]{\smash{{\SetFigFont{25}{14.4}{\rmdefault}{\mddefault}{\updefault}{\color[rgb]{0,0,0}$q_1(x)$}%
}}}}}
\put(1501,-2415){\makebox(0,0)[lb]{\smash{{\SetFigFont{25}{14.4}{\rmdefault}{\mddefault}{\updefault}{\color[rgb]{0,0,0}$q_4(x)$}%
}}}}
\put(3526,-6211){\makebox(0,0)[lb]{\smash{{\SetFigFont{25}{14.4}{\rmdefault}{\bfdefault}{\updefault}{\color[rgb]{0,0,0}$R_2$}%
}}}}
\put(7951,-6286){\makebox(0,0)[lb]{\smash{{\SetFigFont{25}{14.4}{\rmdefault}{\bfdefault}{\updefault}{\color[rgb]{0,0,0}$R_3$}%
}}}}
\put(4801,-5236){\makebox(0,0)[lb]{\smash{{\SetFigFont{25}{14.4}{\rmdefault}{\bfdefault}{\updefault}{\color[rgb]{0,0,0}$R_4$}%
}}}}
\put(9151,-4936){\makebox(0,0)[lb]{\smash{{\SetFigFont{25}{14.4}{\rmdefault}{\bfdefault}{\updefault}{\color[rgb]{0,0,0}$R_5$}%
}}}}
\put(6751,-6886){\makebox(0,0)[lb]{\smash{{\SetFigFont{25}{14.4}{\rmdefault}{\mddefault}{\updefault}{\color[rgb]{0,0,0}$s_4$}%
}}}}
\put(1051,-5536){\makebox(0,0)[lb]{\smash{{\SetFigFont{25}{14.4}{\rmdefault}{\bfdefault}{\updefault}{\color[rgb]{0,0,0}$R_1$}%
}}}}
\put(1351,-6811){\makebox(0,0)[lb]{\smash{{\SetFigFont{25}{14.4}{\rmdefault}{\mddefault}{\updefault}{\color[rgb]{0,0,0}$s_1$}%
}}}}
\put(1951,-6811){\makebox(0,0)[lb]{\smash{{\SetFigFont{25}{14.4}{\rmdefault}{\mddefault}{\updefault}{\color[rgb]{0,0,0}$s_2$}%
}}}}
\put(5926,-6661){\makebox(0,0)[lb]{\smash{{\SetFigFont{25}{14.4}{\rmdefault}{\mddefault}{\updefault}{\color[rgb]{0,0,0}$s_3$}%
}}}}
\put(8701,-6736){\makebox(0,0)[lb]{\smash{{\SetFigFont{25}{14.4}{\rmdefault}{\mddefault}{\updefault}{\color[rgb]{0,0,0}$s_5$}%
}}}}
\put(9601,-2311){\makebox(0,0)[lb]{\smash{{\SetFigFont{25}{14.4}{\rmdefault}{\bfdefault}{\updefault}{\color[rgb]{0,0,0}$R_7$}%
}}}}
\end{picture}%

%% file: N_nu8GamRegions.pstex_t
\begin{picture}(0,0)%
\includegraphics{N_nu8GamRegions.pstex}%
\end{picture}%
\setlength{\unitlength}{3947sp}%
\begingroup\makeatletter\ifx\SetFigFont\undefined%
\gdef\SetFigFont#1#2#3#4#5{%
  \reset@font\fontsize{#1}{#2pt}%
  \fontfamily{#3}\fontseries{#4}\fontshape{#5}%
  \selectfont}%
\fi\endgroup%
\begin{picture}(10919,8041)(301,-7348)
\put(751,-6361){\makebox(0,0)[lb]{\smash{{\SetFigFont{22}{18}{\rmdefault}{\mddefault}{\updefault}{\color[rgb]{0,0,0}$\y_1$}%
}}}}
\put(3151,-5836){\makebox(0,0)[lb]{\smash{{\SetFigFont{22}{18}{\rmdefault}{\bfdefault}{\updefault}{\color[rgb]{0,0,0}$R_4$}%
}}}}
\put(4201,-6586){\makebox(0,0)[lb]{\smash{{\SetFigFont{22}{18}{\rmdefault}{\mddefault}{\updefault}{\color[rgb]{0,0,0}$s_3$}%
}}}}
\put(2776,-5386){\rotatebox{25.0}{\makebox(0,0)[lb]{\smash{{\SetFigFont{22}{18}{\rmdefault}{\mddefault}{\updefault}{\color[rgb]{0,0,0}$\ell_{am}(x)$}%
}}}}}
\put(5926,-6586){\makebox(0,0)[lb]{\smash{{\SetFigFont{22}{18}{\rmdefault}{\mddefault}{\updefault}{\color[rgb]{0,0,0}$s_6$}%
}}}}
\put(7276,-5461){\makebox(0,0)[lb]{\smash{{\SetFigFont{22}{18}{\rmdefault}{\bfdefault}{\updefault}{\color[rgb]{0,0,0}$R_5$}%
}}}}
\put(9526,-6586){\makebox(0,0)[lb]{\smash{{\SetFigFont{22}{18}{\rmdefault}{\mddefault}{\updefault}{\color[rgb]{0,0,0}$s_7$}%
}}}}
\put(8476,-3886){\makebox(0,0)[lb]{\smash{{\SetFigFont{22}{18}{\rmdefault}{\mddefault}{\updefault}{\color[rgb]{0,0,0}$q_4(x)$}%
}}}}
\put(8326,-3286){\makebox(0,0)[lb]{\smash{{\SetFigFont{22}{18}{\rmdefault}{\bfdefault}{\updefault}{\color[rgb]{0,0,0}$R_{7,a}$}%
}}}}
\put(7236,-1873){\rotatebox{50.0}{\makebox(0,0)[lb]{\smash{{\SetFigFont{22}{18}{\rmdefault}{\mddefault}{\updefault}{\color[rgb]{0,0,0}$q_3(x)$}%
}}}}}
\put(9901,-2311){\makebox(0,0)[lb]{\smash{{\SetFigFont{22}{18}{\rmdefault}{\bfdefault}{\updefault}{\color[rgb]{0,0,0}$R_{7,b}$}%
}}}}
\put(10651,-1786){\makebox(0,0)[lb]{\smash{{\SetFigFont{22}{18}{\rmdefault}{\mddefault}{\updefault}{\color[rgb]{0,0,0}$M_C$}%
}}}}
\put(7876,614){\rotatebox{300.0}{\makebox(0,0)[lb]{\smash{{\SetFigFont{22}{18}{\rmdefault}{\mddefault}{\updefault}{\color[rgb]{0,0,0}$\ell_s(x)$}%
}}}}}
\put(10651,-6436){\makebox(0,0)[lb]{\smash{{\SetFigFont{22}{18}{\rmdefault}{\mddefault}{\updefault}{\color[rgb]{0,0,0}$M_3$}%
}}}}
\end{picture}%

%% file: N_nu10GamRegions.pstex_t
\begin{picture}(0,0)%
\includegraphics{N_nu10GamRegions.pstex}%
\end{picture}%
\setlength{\unitlength}{3947sp}%
\begingroup\makeatletter\ifx\SetFigFont\undefined%
\gdef\SetFigFont#1#2#3#4#5{%
  \reset@font\fontsize{#1}{#2pt}%
  \fontfamily{#3}\fontseries{#4}\fontshape{#5}%
  \selectfont}%
\fi\endgroup%
\begin{picture}(10844,8037)(301,-7348)
\put(751,-6361){\makebox(0,0)[lb]{\smash{{\SetFigFont{22}{18}{\rmdefault}{\bfdefault}{\updefault}{\color[rgb]{0,0,0}$\y_1$}%
}}}}
\put(10576,-6511){\makebox(0,0)[lb]{\smash{{\SetFigFont{22}{18}{\rmdefault}{\bfdefault}{\updefault}{\color[rgb]{0,0,0}$M_3$}%
}}}}
\put(10576,-1786){\makebox(0,0)[lb]{\smash{{\SetFigFont{22}{18}{\rmdefault}{\bfdefault}{\updefault}{\color[rgb]{0,0,0}$M_C$}%
}}}}
\put(6601,-6811){\makebox(0,0)[lb]{\smash{{\SetFigFont{22}{18}{\rmdefault}{\mddefault}{\updefault}{\color[rgb]{0,0,0}$s_7$}%
}}}}
\put(8551,-6811){\makebox(0,0)[lb]{\smash{{\SetFigFont{22}{18}{\rmdefault}{\mddefault}{\updefault}{\color[rgb]{0,0,0}$s_8$}%
}}}}
\put(4426,-5686){\makebox(0,0)[lb]{\smash{{\SetFigFont{22}{18}{\rmdefault}{\bfdefault}{\updefault}{\color[rgb]{0,0,0}$R_{7a}$}%
}}}}
\put(8101,-4111){\makebox(0,0)[lb]{\smash{{\SetFigFont{22}{18}{\rmdefault}{\bfdefault}{\updefault}{\color[rgb]{0,0,0}$R_{7b}$}%
}}}}
\put(3601,-4936){\rotatebox{20.0}{\makebox(0,0)[lb]{\smash{{\SetFigFont{22}{18}{\rmdefault}{\mddefault}{\updefault}{\color[rgb]{0,0,0}$\ell_{am}(x)$}%
}}}}}
\put(4051,539){\rotatebox{305.0}{\makebox(0,0)[lb]{\smash{{\SetFigFont{22}{18}{\rmdefault}{\mddefault}{\updefault}{\color[rgb]{0,0,0}$\ell_s(x)$}%
}}}}}
\end{picture}%

%% file: ArcProbSeg1.pstex_t
\begin{picture}(0,0)%
\includegraphics{ArcProbSeg1.pstex}%
\end{picture}%
\setlength{\unitlength}{3947sp}%
\begingroup\makeatletter\ifx\SetFigFont\undefined%
\gdef\SetFigFont#1#2#3#4#5{%
  \reset@font\fontsize{#1}{#2pt}%
  \fontfamily{#3}\fontseries{#4}\fontshape{#5}%
  \selectfont}%
\fi\endgroup%
\begin{picture}(11349,8130)(289,-7348)
\put(10426,-6586){\makebox(0,0)[lb]{\smash{{\SetFigFont{22}{16.8}{\rmdefault}{\mddefault}{\updefault}{\color[rgb]{0,0,0}$\y_2=(1,0)$}%
}}}}
\put(826,-6586){\makebox(0,0)[lb]{\smash{{\SetFigFont{22}{16.8}{\rmdefault}{\mddefault}{\updefault}{\color[rgb]{0,0,0}$\y_1=(0,0)$}%
}}}}
\put(7726,-2236){\makebox(0,0)[lb]{\smash{{\SetFigFont{22}{16.8}{\rmdefault}{\mddefault}{\updefault}{\color[rgb]{0,0,0}$e_1$}%
}}}}
\put(3001,-2686){\makebox(0,0)[lb]{\smash{{\SetFigFont{22}{16.8}{\rmdefault}{\mddefault}{\updefault}{\color[rgb]{0,0,0}$e_2$}%
}}}}
\put(5626,-6586){\makebox(0,0)[lb]{\smash{{\SetFigFont{22}{16.8}{\rmdefault}{\mddefault}{\updefault}{\color[rgb]{0,0,0}$M_3$}%
}}}}
\put(5701,-4111){\makebox(0,0)[lb]{\smash{{\SetFigFont{22}{16.8}{\rmdefault}{\mddefault}{\updefault}{\color[rgb]{0,0,0}$M_{C}$}%
}}}}
\put(3076,-6661){\makebox(0,0)[lb]{\smash{{\SetFigFont{22}{16.8}{\rmdefault}{\mddefault}{\updefault}{\color[rgb]{0,0,0}$Q_1$}%
}}}}
\put(8026,-6586){\makebox(0,0)[lb]{\smash{{\SetFigFont{22}{16.8}{\rmdefault}{\mddefault}{\updefault}{\color[rgb]{0,0,0}$Q_2$}%
}}}}
\put(9376,-4486){\makebox(0,0)[lb]{\smash{{\SetFigFont{22}{16.8}{\rmdefault}{\mddefault}{\updefault}{\color[rgb]{0,0,0}$Q_3$}%
}}}}
\put(6676,-811){\makebox(0,0)[lb]{\smash{{\SetFigFont{22}{16.8}{\rmdefault}{\mddefault}{\updefault}{\color[rgb]{0,0,0}$Q_4$}%
}}}}
\put(3226,-4336){\makebox(0,0)[lb]{\smash{{\SetFigFont{22}{16.8}{\rmdefault}{\mddefault}{\updefault}{\color[rgb]{0,0,0}$\ell_{r_1}(x_1,x)$}%
}}}}
\put(7426,-3886){\makebox(0,0)[lb]{\smash{{\SetFigFont{22}{16.8}{\rmdefault}{\mddefault}{\updefault}{\color[rgb]{0,0,0}$\ell_{r_2}(x_1,x)$}%
}}}}
\put(2926,-5386){\makebox(0,0)[lb]{\smash{{\SetFigFont{22}{16.8}{\rmdefault}{\mddefault}{\updefault}{\color[rgb]{0,0,0}$x_1$}%
}}}}
\put(1876,-6061){\makebox(0,0)[lb]{\smash{{\SetFigFont{22}{16.8}{\rmdefault}{\mddefault}{\updefault}{\color[rgb]{0,0,0}$\varepsilon$}%
}}}}
\put(3226,-6136){\makebox(0,0)[lb]{\smash{{\SetFigFont{22}{16.8}{\rmdefault}{\mddefault}{\updefault}{\color[rgb]{0,0,0}$q(\y_1,x)$}%
}}}}
\put(7351,-6061){\makebox(0,0)[lb]{\smash{{\SetFigFont{22}{16.8}{\rmdefault}{\mddefault}{\updefault}{\color[rgb]{0,0,0}$q(\y_2,x)$}%
}}}}
\put(5326,-586){\makebox(0,0)[lb]{\smash{{\SetFigFont{22}{16.8}{\rmdefault}{\mddefault}{\updefault}{\color[rgb]{0,0,0}$q(\y_3,x)$}%
}}}}
\put(4426,-286){\makebox(0,0)[lb]{\smash{{\SetFigFont{22}{16.8}{\rmdefault}{\mddefault}{\updefault}{\color[rgb]{0,0,0}$N_2$}%
}}}}
\put(4201,-6586){\makebox(0,0)[lb]{\smash{{\SetFigFont{22}{16.8}{\rmdefault}{\mddefault}{\updefault}{\color[rgb]{0,0,0}$N_1$}%
}}}}
\put(9301,-6511){\makebox(0,0)[lb]{\smash{{\SetFigFont{22}{16.8}{\rmdefault}{\mddefault}{\updefault}{\color[rgb]{0,0,0}$N_1$}%
}}}}
\put(5551,614){\makebox(0,0)[lb]{\smash{{\SetFigFont{22}{16.8}{\rmdefault}{\mddefault}{\updefault}{\color[rgb]{0,0,0}$\y_3=(1/2,\sqrt{3}/2)$}%
}}}}
\put(3826,-736){\makebox(0,0)[lb]{\smash{{\SetFigFont{22}{16.8}{\rmdefault}{\mddefault}{\updefault}{\color[rgb]{0,0,0}$Q_5$}%
}}}}
\put(4876,-1336){\makebox(0,0)[lb]{\smash{{\SetFigFont{22}{16.8}{\rmdefault}{\mddefault}{\updefault}{\color[rgb]{0,0,0}$U_2$}%
}}}}
\put(7876,-5536){\makebox(0,0)[lb]{\smash{{\SetFigFont{22}{16.8}{\rmdefault}{\mddefault}{\updefault}{\color[rgb]{0,0,0}$U_1$}%
}}}}
\put(1501,-4561){\makebox(0,0)[lb]{\smash{{\SetFigFont{22}{16.8}{\rmdefault}{\mddefault}{\updefault}{\color[rgb]{0,0,0}$Q_6$}%
}}}}
\put(1951,-3811){\makebox(0,0)[lb]{\smash{{\SetFigFont{22}{16.8}{\rmdefault}{\mddefault}{\updefault}{\color[rgb]{0,0,0}$N_2$}%
}}}}
\end{picture}%

%% file: PartSegrgt2.pstex_t
\begin{picture}(0,0)%
\includegraphics{PartSegrgt2.pstex}%
\end{picture}%
\setlength{\unitlength}{3947sp}%
\begingroup\makeatletter\ifx\SetFigFont\undefined%
\gdef\SetFigFont#1#2#3#4#5{%
  \reset@font\fontsize{#1}{#2pt}%
  \fontfamily{#3}\fontseries{#4}\fontshape{#5}%
  \selectfont}%
\fi\endgroup%
\begin{picture}(10662,8124)(76,-7348)
\put(3076,-6661){\makebox(0,0)[lb]{\smash{{\SetFigFont{22}{16.8}{\rmdefault}{\mddefault}{\updefault}{\color[rgb]{0,0,0}$Q_1=(s_2,0)$}%
}}}}
\put(2701,-6586){\makebox(0,0)[lb]{\smash{{\SetFigFont{22}{16.8}{\rmdefault}{\mddefault}{\updefault}{\color[rgb]{0,0,0}$s_1$}%
}}}}
\put(5326,-6661){\makebox(0,0)[lb]{\smash{{\SetFigFont{22}{16.8}{\rmdefault}{\mddefault}{\updefault}{\color[rgb]{0,0,0}$s_3$}%
}}}}
\put(5926,-6586){\makebox(0,0)[lb]{\smash{{\SetFigFont{22}{16.8}{\rmdefault}{\mddefault}{\updefault}{\color[rgb]{0,0,0}$s_5$}%
}}}}
\put(7876,-1861){\makebox(0,0)[lb]{\smash{{\SetFigFont{22}{16.8}{\rmdefault}{\mddefault}{\updefault}{\color[rgb]{0,0,0}$\ell_{am}(x)$}%
}}}}
\put(9601,-1786){\makebox(0,0)[lb]{\smash{{\SetFigFont{22}{16.8}{\rmdefault}{\mddefault}{\updefault}{\color[rgb]{0,0,0}$M_{C}$}%
}}}}
\put(9151,-6586){\makebox(0,0)[lb]{\smash{{\SetFigFont{22}{16.8}{\rmdefault}{\mddefault}{\updefault}{\color[rgb]{0,0,0}$M_3=(1/2,0)$}%
}}}}
\put(4576, 89){\makebox(0,0)[lb]{\smash{{\SetFigFont{22}{16.8}{\rmdefault}{\mddefault}{\updefault}{\color[rgb]{0,0,0}$M_2$}%
}}}}
\put(1876,-5686){\makebox(0,0)[lb]{\smash{{\SetFigFont{22}{16.8}{\rmdefault}{\mddefault}{\updefault}{\color[rgb]{0,0,0}$\varepsilon$}%
}}}}
\put(1276,-3736){\makebox(0,0)[lb]{\smash{{\SetFigFont{22}{16.8}{\rmdefault}{\mddefault}{\updefault}{\color[rgb]{0,0,0}$q(\y_1,x)$}%
}}}}
\put(2926,-61){\makebox(0,0)[lb]{\smash{{\SetFigFont{22}{16.8}{\rmdefault}{\mddefault}{\updefault}{\color[rgb]{0,0,0}$\ell_s(x)$}%
}}}}
\put(2476,-586){\makebox(0,0)[lb]{\smash{{\SetFigFont{22}{16.8}{\rmdefault}{\mddefault}{\updefault}{\color[rgb]{0,0,0}$q_2(x)$}%
}}}}
\put(8101,-6511){\makebox(0,0)[lb]{\smash{{\SetFigFont{22}{16.8}{\rmdefault}{\mddefault}{\updefault}{\color[rgb]{0,0,0}$s_6$}%
}}}}
\put(6901,-6511){\makebox(0,0)[lb]{\smash{{\SetFigFont{22}{16.8}{\rmdefault}{\mddefault}{\updefault}{\color[rgb]{0,0,0}$s_4$}%
}}}}
\put( 76,-6661){\makebox(0,0)[lb]{\smash{{\SetFigFont{22}{16.8}{\rmdefault}{\mddefault}{\updefault}{\color[rgb]{0,0,0}$\y_1=(0,0)$}%
}}}}
\put(4201,-5686){\makebox(0,0)[lb]{\smash{{\SetFigFont{22}{16.8}{\rmdefault}{\mddefault}{\updefault}{\color[rgb]{0,0,0}$R_1(\varepsilon)$}%
}}}}
\put(6151,-4936){\makebox(0,0)[lb]{\smash{{\SetFigFont{22}{16.8}{\rmdefault}{\mddefault}{\updefault}{\color[rgb]{0,0,0}$R_2(\varepsilon)$}%
}}}}
\put(7276,-3811){\makebox(0,0)[lb]{\smash{{\SetFigFont{22}{16.8}{\rmdefault}{\mddefault}{\updefault}{\color[rgb]{0,0,0}$R_3(\varepsilon)$}%
}}}}
\end{picture}%

%% file: HL_seg_eps1.pstex_t
\begin{picture}(0,0)%
\includegraphics{HL_seg_eps1.pstex}%
\end{picture}%
\setlength{\unitlength}{3947sp}%
\begingroup\makeatletter\ifx\SetFigFont\undefined%
\gdef\SetFigFont#1#2#3#4#5{%
  \reset@font\fontsize{#1}{#2pt}%
  \fontfamily{#3}\fontseries{#4}\fontshape{#5}%
  \selectfont}%
\fi\endgroup%
\begin{picture}(10963,8100)(301,-7336)
\put(451,-6961){\makebox(0,0)[lb]{\smash{{\SetFigFont{22}{14.4}{\rmdefault}{\bfdefault}{\updefault}{\color[rgb]{0,0,0}$\y_1$}%
}}}}
\put(10726,-7036){\makebox(0,0)[lb]{\smash{{\SetFigFont{22}{14.4}{\rmdefault}{\bfdefault}{\updefault}{\color[rgb]{0,0,0}$M_3$}%
}}}}
\put(10576,-1786){\makebox(0,0)[lb]{\smash{{\SetFigFont{22}{14.4}{\rmdefault}{\bfdefault}{\updefault}{\color[rgb]{0,0,0}$M_C$}%
}}}}
\put(9001,-2611){\rotatebox{25.0}{\makebox(0,0)[lb]{\smash{{\SetFigFont{22}{14.4}{\rmdefault}{\mddefault}{\updefault}{\color[rgb]{0,0,0}$\ell_{am}(x)$}%
}}}}}
\put(1651,-6136){\makebox(0,0)[lb]{\smash{{\SetFigFont{22}{14.4}{\rmdefault}{\mddefault}{\updefault}{\color[rgb]{0,0,0}$q_5(x)$}%
}}}}
\put(4276,-7186){\makebox(0,0)[lb]{\smash{{\SetFigFont{22}{14.4}{\rmdefault}{\mddefault}{\updefault}{\color[rgb]{0,0,0}$s_1$}%
}}}}
\put(5626,-7111){\makebox(0,0)[lb]{\smash{{\SetFigFont{22}{14.4}{\rmdefault}{\mddefault}{\updefault}{\color[rgb]{0,0,0}$s_3$}%
}}}}
\put(5101,-7111){\makebox(0,0)[lb]{\smash{{\SetFigFont{22}{14.4}{\rmdefault}{\mddefault}{\updefault}{\color[rgb]{0,0,0}$s_2$}%
}}}}
\put(6451,-7186){\makebox(0,0)[lb]{\smash{{\SetFigFont{22}{14.4}{\rmdefault}{\mddefault}{\updefault}{\color[rgb]{0,0,0}$s_4$}%
}}}}
\put(9526,-7186){\makebox(0,0)[lb]{\smash{{\SetFigFont{22}{14.4}{\rmdefault}{\mddefault}{\updefault}{\color[rgb]{0,0,0}$s_8$}%
}}}}
\put(8851,-7261){\makebox(0,0)[lb]{\smash{{\SetFigFont{22}{14.4}{\rmdefault}{\mddefault}{\updefault}{\color[rgb]{0,0,0}$s_6=s_{10}$}%
}}}}
\put(7501,-7111){\makebox(0,0)[lb]{\smash{{\SetFigFont{22}{14.4}{\rmdefault}{\mddefault}{\updefault}{\color[rgb]{0,0,0}$s_7$}%
}}}}
\put(7876,-7111){\makebox(0,0)[lb]{\smash{{\SetFigFont{22}{14.4}{\rmdefault}{\mddefault}{\updefault}{\color[rgb]{0,0,0}$s_5$}%
}}}}
\put(9976,-7111){\makebox(0,0)[lb]{\smash{{\SetFigFont{22}{14.4}{\rmdefault}{\mddefault}{\updefault}{\color[rgb]{0,0,0}$s_9$}%
}}}}
\put(5626,-5536){\makebox(0,0)[lb]{\smash{{\SetFigFont{22}{14.4}{\rmdefault}{\bfdefault}{\updefault}{\color[rgb]{0,0,0}$R_{1a}$}%
}}}}
\put(8101,-5236){\makebox(0,0)[lb]{\smash{{\SetFigFont{22}{14.4}{\rmdefault}{\bfdefault}{\updefault}{\color[rgb]{0,0,0}$R_{1b}$}%
}}}}
\put(6751,-6661){\makebox(0,0)[lb]{\smash{{\SetFigFont{22}{14.4}{\rmdefault}{\bfdefault}{\updefault}{\color[rgb]{0,0,0}$R_{2a}$}%
}}}}
\put(8926,-6586){\makebox(0,0)[lb]{\smash{{\SetFigFont{22}{14.4}{\rmdefault}{\bfdefault}{\updefault}{\color[rgb]{0,0,0}$R_{2b}$}%
}}}}
\put(9526,-3361){\makebox(0,0)[lb]{\smash{{\SetFigFont{22}{14.4}{\rmdefault}{\bfdefault}{\updefault}{\color[rgb]{0,0,0}$R_{3b}$}%
}}}}
\put(8701,-4561){\makebox(0,0)[lb]{\smash{{\SetFigFont{22}{14.4}{\rmdefault}{\bfdefault}{\updefault}{\color[rgb]{0,0,0}$R_{3a}$}%
}}}}
\put(9901,-6511){\makebox(0,0)[lb]{\smash{{\SetFigFont{22}{14.4}{\rmdefault}{\bfdefault}{\updefault}{\color[rgb]{0,0,0}$R_4$}%
}}}}
\put(3901,-3961){\rotatebox{300.0}{\makebox(0,0)[lb]{\smash{{\SetFigFont{22}{14.4}{\rmdefault}{\mddefault}{\updefault}{\color[rgb]{0,0,0}$q_(\y_1,x)$}%
}}}}}
\put(5446,-1993){\rotatebox{305.0}{\makebox(0,0)[lb]{\smash{{\SetFigFont{22}{14.4}{\rmdefault}{\mddefault}{\updefault}{\color[rgb]{0,0,0}$q_2(x)$}%
}}}}}
\put(6526,-1261){\rotatebox{305.0}{\makebox(0,0)[lb]{\smash{{\SetFigFont{22}{14.4}{\rmdefault}{\mddefault}{\updefault}{\color[rgb]{0,0,0}$q_4(x)$}%
}}}}}
\put(7501,-1186){\rotatebox{305.0}{\makebox(0,0)[lb]{\smash{{\SetFigFont{22}{14.4}{\rmdefault}{\mddefault}{\updefault}{\color[rgb]{0,0,0}$q_3(x)=\ell_s(x)$}%
}}}}}
\end{picture}%